\documentclass[a4paper, 12pt]{article}
\usepackage[utf8]{inputenc}
\usepackage[english,russian]{babel}
\usepackage{amsmath}
\usepackage{amsfonts}
\usepackage{amsbib}
\DeclareMathAlphabet{\mathpzc}{OT1}{pzc}{m}{it}
\renewcommand{\Im}{\operatorname{Im}}

\newcommand{\codim}{\operatorname{codim}}
\newcommand{\Rank}{\operatorname{rank}}
\newcommand{\sgn}{\operatorname{sign}}
\newcommand{\inter}{\operatorname{int}}
\newcommand{\e}{\mathfrak e}
\renewcommand{\l}{\mathfrak l}
\renewcommand{\k}{\mathfrak k}
\newcommand{\m}{\mathfrak m}
\newcommand{\n}{\mathfrak n}
\newcommand{\M}{\mathfrak M}
\newcommand{\Iota}{\mathfrak J}
\newcommand{\card}{\operatorname{card}}
\newcommand{\diam}{\operatorname{diam}}
\newcommand{\supp}{\operatorname{supp}}
\newcommand{\supvrai}{\operatornamewithlimits{sup\,vrai}}
\newcommand{\N}{\mathbb N}
\newcommand{\Z}{\mathbb Z}
\newcommand{\R}{\mathbb R}
\newcommand{\Nu}{\mathcal N}
\newcommand{\D}{\mathcal D}
\newcommand{\J}{\mathcal J}
\newcommand{\I}{\mathcal I}

\newcommand{\mes}{\operatorname{mes}}

\allowdisplaybreaks

\begin{document}

\author{ С. Н. Кудрявцев }
\title{Продолжение функций из изотропных пространств Никольского -- Бесова
и их приближение вместе с производными}
\date{}
\maketitle
\begin{abstract}
В статье рассмотрены изотропные пространства Никольского и Бесова с
нормами, в определении которых вместо модуля непрерывности известного порядка
частных производных функций фиксированного порядка используется
"$L_p$-усреднённый" модуль непрерывности функций соответствующего
порядка. Для таких пространств функций, заданных в ограниченных областях
$(1,\ldots,1)$-типа (в широком смысле), построены непрерывные линейные
отображения их в обычные изотропные пространства Никольского и Бесова в
$ \R^d, $ являющиеся операторами продолжения функций, что влечёт
совпадение тех и других пространств в упомянутых областях. Установлено, что всякая ограниченная область в $ \R^d $ с
липшицевой границей является областью
$(1,\ldots,1)$-типа (в широком смысле). В работе также найдена слабая
асимптотика аппроксимационных характеристик,
относящихся к задаче восстановления функций вместе с их производными
по значениям функций в заданном числе точек, задаче С.Б. Стечкина для
оператора дифференцирования, задаче описания асимптотики поперечников для
изотропных классов Никольского и Бесова в этих областях.
\end{abstract}

Ключевые слова: изотропные пространства Никольского -- Бесова,
продолжение функций, эквивалентные нормы, восстановление функций,
приближение оператора, поперечник
\bigskip

\centerline{Введение}
\bigskip

В работе рассмотрен ряд задач теории функциональных пространств и теории
приближений для изотропных пространств Никольского и Бесова функций, заданных
в областях из определённого класса. Эти задачи объединяет общая техника вывода
метрических соотношений, используемых для получения верхних оценок величин,
являющихся предметом изучения в этих задачах. Остановимся подробнее на
содержании работы.

При $ d \in \N, \alpha \in \R_+, 1 \le p < \infty, 1 \le \theta < \infty $ для
области $ D \subset \R^d $ вводятся в рассмотрение
пространства $ (B_{p,\theta}^\alpha)^\prime(D) ((H_p^\alpha)^\prime(D)) $
с нормами
$$
\| f \|_{(B_{p,\theta}^\alpha)^\prime(D)} = \max\biggl(\| f \|_{L_p(D)},
\left(\int_0^\infty t^{-1 -\theta \alpha}
(\Omega^{\prime l}(f, t)_{L_p(D)})^{\theta} \,dt \right)^{1/\theta}\biggr),
$$
$$
\| f \|_{(H_p^\alpha)^\prime(D)} = \max(\| f \|_{L_p(D)},
\sup_{t \in \R_+} t^{-\alpha} \Omega^{\prime l}(f,t)_{L_p(D)}),
$$
где
\begin{multline*}
\Omega^{\prime l}(f,t)_{L_p(D)} =
\biggl((2 t)^{-d} \int_{ \{\xi \in \R^d: |\xi_j| \le t, \ j =1,\ldots,d\}}
\| \Delta_{\xi}^{l} f\|_{L_p(D_{l \xi})}^p d\xi\biggr)^{1 /p},\\
t \in \R_+, l = \min \{m \in \N: \alpha < m \}.
\end{multline*}

Для ограниченных областей $D $ так называемого $ \e $-типа (см. п. 1.2.)
построены непрерывные линейные отображения пространств $ (B_{p,\theta}^\alpha)^\prime(D)
((H_p^\alpha)^\prime(D)) $ в пространства Бесова $ B_{p,\theta}^\alpha(\R^d) $
(Никольского $ H_p^\alpha(\R^d)), $ являющиеся операторами продолжения
функций, что влечёт совпадение соответствующих пространств (см. п. 2.1.).
Публикации, в которых рассматривается такая задача, автору не известны. Из
близких работ по этой тематике приведём [1] -- [4] (см. также имеющуюся
там литературу), в которых изучается вопрос о продолжении за пределы
области определения гладких функций из изотропных пространств с
сохранением класса. Отметим, что средства построения операторов
продолжения функций и схемы доказательства непрерывности таких
операторов, применяемые ниже, отличаются от тех, что использовались в упомянутых
работах.

В случае ограниченной области $ D \ \e $-типа в статье при соответствующих
условиях на параметры, определяющие рассматриваемые объекты,
установлена слабая асимптотика, т.е. найден порядок величины
наилучшей точности восстановления в пространстве $ W_q^{\mathpzc m}(D) $
по значениям в $ n $ точках функций $ f $ из классов Никольского
$ (\mathcal H_p^\alpha)^\prime(D) $ и Бесова
$ (\mathcal B_{p,\theta}^\alpha)^\prime(D), $ кроме того,
получена слабая асимптотика поведения в зависимости от $ \rho $
наилучшей точности приближения в $ L_q(D) $ оператора частного
дифференцирования $ \D^\lambda $ на классах
$ (\mathcal H_p^\alpha)^\prime(D) $ и
$ (\mathcal B_{p,\theta}^\alpha)^\prime(D) $ операторами, действующими из
$ L_s(D) $ в $ L_q(D), $ норма которых не превосходит $ \rho, $
и, наконец, описана слабая асимптотика колмогоровского, гельфандовского,
линейного, александровского и энтропийного $n$-поперечников в
пространствах Соболева $ W_q^{\mathpzc m}(D) $ единичного шара пространств
$ (H_p^\alpha)^\prime(D) $ и $ (B_{p,\theta}^\alpha)^\prime(D). $ В этой части
настоящая работа продолжает исследования, проводившиеся автором в
[5] -- [7] для упомянутых задач в
отношении изотропных классов Никольского и Бесова функций, заданных
на кубе $ I^d. $ Проблематика, к которой относятся рассматриваемые задачи,
занимает значительное место в теории приближений, и обзор результатов по ней
в рамках этой статьи невозможен.

Отметим, что общие подходы при решении упомянутых задач сохранились такими же,
как в [5] -- [7]. Однако для их реализации потребовались другие
средства приближения и существенные изменения схем вывода вспомогательных
соотношений, используемых для получения верхних оценок изучаемых величин,
в сравнении с теми, что применялись в [5] -- [7]. Нижние оценки
рассматриваемых ниже величин устанавливаются, опираясь на нижние
оценки соответствующих величин, полученные в [5] -- [7]. Как и
следовало ожидать, слабая асимптотика изучаемых величин в задачах о
восстановлении функций, приближении оператора дифференцирования и нахождении
порядков поперечников совпадает со слабой асимптотикой соответствующих им
величин, установленной в [5] -- [7].

Работа состоит из введения и пяти параграфов, в первом из которых
рассматриваются некоторые средства для решения объявленных задач, во втором --
задача о продолжении, в третьем -- задача восстановления, в четвёртом --
задача приближения оператора дифференцирования, в пятом -- поперечники.
\bigskip

\centerline{\S 1. Предварительные сведения и вспомогательные утверждения}
\bigskip

1.1. В этом пункте вводятся обозначения, относящиеся к
функциональным пространствам и классам, рассматриваемым в настоящей работе,
а также приводятся некоторые факты, необходимые в дальнейшем.

Для $ d \in \N $ через $ \Z_+^d $ обозначим множество
$$
\Z_+^d = \{\lambda = (\lambda_1, \ldots, \lambda_d) \in \Z^d:
\lambda_j \ge0, j =1, \ldots, d\}.
$$

При $ d \in \N $ для $ l \in \Z_+ $ через $ \Z_{+ l}^d $ обозначим множество
$$
\{\lambda \in \Z_+^d: |\lambda| \le l\},
$$
где $ |\lambda| = \sum_{j=1}^d \lambda_j, $
а через $ \Z_{+ l}^{\prime\,d} $ --- множество
$$
\{\lambda \in \Z_+^d: |\lambda| = l\}.
$$

Обозначим также при $ d \in \N $ для $ l \in \Z_+^d $ через
$ \Z_+^d(l) $ множество
$$
\Z_+^d(l) = \{ \lambda \in \Z_+^d: \lambda_j \le l_j, j =1, \ldots,d\}.
$$

Для $  m,n \in \Z: m \le n, $ и $ d \in \N $ обозначим через
$$
\Nu_{m,n}^d = \{ \nu \in \Z^d: m \le \nu_j \le n, j =1,\ldots,d\}.
$$

Напомним, что при $ n \in \N $ и $ 1 \le p  \le \infty $ через $ l_p^n $
обозначается пространство $ \R^n $ с фиксированной в нём нормой
$$
\|x\|_{l_p^n} = \begin{cases} (\sum_{j =1}^n |x_j|^p)^{1/p} \text{ при } p < \infty; \\
\max_{j =1}^n |x_j| \text{ при } p = \infty, \end{cases} \ x \in \R^n.
$$

Отметим здесь неоднократно используемое в дальнейшем неравенство
\begin{equation*} \tag{1.1.1}
| \sum_{j =1}^n x_j|^a \le \sum_{j =1}^n |x_j|^a, \ x_j \in \R, \ j =1,\ldots,n, \
n \in \N, \ 0 \le a \le 1.
\end{equation*}

Для $ d \in \N $ в пространстве $\R^d$ фиксируем норму
$$
\| x\| = \|x\|_{l_\infty^d},\ x \in \R^d.
$$

Для $ d \in \N, l \in \Z_+ $ через $ \mathcal P^{l,d} $ будем
обозначать пространство вещественных  полиномов, состоящее из всех
функций $ f: \R^d \mapsto \R $ вида
$$
f(x) = \sum_{\lambda \in \Z_{+ l}^d} a_\lambda x^\lambda,
$$
где $ a_\lambda \in \R,\  x^\lambda = x_1^{\lambda_1} \cdot \ldots \cdot
x_d^{\lambda_d},\
\lambda \in \Z_{+ l}^d,$
а для области $ D \subset \R^d $ через
$ \mathcal P^{l,d}(D) $ обозначим пространство функций $ f, $ определённых на
$ D, $ для каждой из которых существует  $ g \in \mathcal P^{l,d} $  такой,
что $ f = g \mid_D. $

Для топологического пространства $ T $ через $ C(T) $ будем
обозначать пространство непрерывных вещественных функций на $ T. $
Для множества $ A \subset T $ через $ \overline A $ обозначается замыкание
множества $ A, $ а через $ \inter A $  -- его внутренность.

Для измеримого по Лебегу множества $ D \subset \R^d $ при $ 1 \le p \le \infty $
через $ L_p(D),$ как обычно, обозначается
пространство всех вещественных измеримых на $ D $ функций $ f,$
для которых определена норма
$$
\|f\|_{L_p(D)} = \begin{cases} (\int_D |f(x)|^p dx)^{1/p}, \ 1 \le p < \infty;\\
\supvrai_{x \in D}|f(x)|, \ p = \infty. \end{cases}
$$

При $ d \in \N $ для $ \lambda \in \Z_+^d $ через $ \D^\lambda $
будем обозначать оператор дифференцирования $ \D^\lambda =
\frac{\D^{|\lambda|}} {\D x_1^{\lambda_1} \ldots \D x_d^{\lambda_d}}. $

При обозначении известных пространств дифференцируемых функций будем
ориентироваться на [8].
Для области $ D \subset \R^d $ при $ 1 \le p \le \infty, l \in \Z_+ $ через
$ W_p^l(D) $ обозначается пространство всех функций $ f \in L_p(D), $
для которых для каждого $ \lambda \in \Z_{+ l}^{\prime \,d} $
обобщённая производная $ \D^\lambda f $ принадлежит $ L_p(D) $ с нормой
$$
\|f\|_{W_p^l(D)} = \max(\|f\|_{L_p(D)}, \max_{\lambda \in \Z_{+ l}^{\prime \,d}}
\|\D^\lambda f\|_{L_p(D)}).
$$

При $ d \in \N $ для $ x,y \in \R^d $ будем обозначать
$$
(x,y) = \sum_{j =1}^d x_j y_j.
$$

Для $ x,y \in \R^d $ положим $ xy = x \cdot y = (x_1 y_1, \ldots,x_d y_d), $
а для $ x \in \R^d $ и $ A \subset \R^d $ определим
$$
x A = x \cdot A = \{xy: y \in A\}.
$$

Напомним, что $ t_+ = \frac{1} {2} (t +|t|), \ t \in \R. $

Обозначим через $ \R_+^d $ множество $ x \in \R^d, $ для которых
$ x_j >0 $ при $ j =1,\ldots,d.$

Будем обозначать через $ \chi_A $ характеристическую функцию
множества $ A \subset \R^d. $

При $ d \in \N $ определим множества
$$
I^d = \{x \in \R^d: 0 < x_j < 1, \ j =1,\ldots,d\},
$$
$$
\overline I^d = \{x \in \R^d: 0 \le x_j \le 1, \ j =1,\ldots,d\},
$$
$$
B_0^d = \{ x \in \R^d: \|x\| < 1\}, \
B^d = \{x \in \R^d: \|x\| \le 1\}.
$$

Через $ \e $ будем обозначать вектор в $ \R^d, $ задаваемый
равенством $ \e = (1,\ldots,1). $

При $ d \in \N $ для $ j =1,\ldots,d$ через $ e_j $ будем обозначать вектор
$ e_j = (0,\ldots,0,1_j,0,\ldots,0).$

Далее, напомним, что для области $ D \subset \R^d $ и вектора $ h \in \R^d $
через $ D_h $ обозначается множество
$$
D_h = \{x \in D: x +th \in D \  \forall \ t \in \ \overline I\}.
$$

Для функции $ f, $ заданной в области $ D \subset \R^d, $ и
вектора $ h \in \R^d $ определим в $ D_h $ её разность $ \Delta_h f $ с шагом
$ h, $ полагая
$$
(\Delta_h f)(x) = f(x +h) -f(x), \ x \in D_h,
$$
а для $ l \in \N: l \ge 2, $ в $ D_{lh} $ определим $l$-ую разность
$ \Delta_h^l f $ функции $ f $ с шагом $ h $ равенством
$$
(\Delta_h^l f)(x) = (\Delta_h (\Delta_h^{l -1} f))(x), \ x \in D_{lh},
$$
положим также $ \Delta_h^0 f = f. $

Как известно, справедливо равенство
$$
(\Delta_h^l f)(\cdot) = \sum_{k=0}^l C_l^k (-1)^{l -k} f(\cdot +kh), \
C_l^k = \frac{l!} {k! (l -k)!}.
$$

Отметим, что для $ l \in \Z_+,\ d \in \N$ для $ f \in \mathcal P^{l,d} $ и
$ h \in \R^d$ справедливо включение $ \Delta_h f \in \mathcal P^{l -1,d},$
где $ \mathcal P^{-1,d} $ --- множество, состоящее из единственного элемента ---
функции тождественно равной нулю
в $ \R^d. $ Отсюда вытекает, что при  $ l,d \in \N $ для любой области
$ D \subset \R^d $ для $ f \in \mathcal P^{l-1,d}(D) $ и $ h \in \R^d $ верно равенство
$$
(\Delta_h^l f)(x) =0 \ \text{ при } \ x \in D_{lh}.
$$

Определим пространства и классы функций, изучаемые в настоящей работе (ср. с
[8]). Но прежде введём некоторые обозначения.

Пусть $ D$ --- область в $ \R^d$ и $ 1 \le p < \infty.$
Для $ f \in L_p(D) $ при $ l \in \Z_+ $ обозначим модуль
непрерывности в $ L_p(D) $ порядка $ l $ через
$$
\Omega^l(f,t)_{L_p(D)} =
\supvrai_{h \in t B^d} \| \Delta_{h}^l f\|_{L_p(D_{l h })}, \ t \in \R_+,
$$

а также введём в рассмотрение "усреднённый" модуль непрерывности
в $ L_p(D) $ порядка $ l, $ полагая
\begin{multline*}
\Omega^{\prime l}(f,t)_{L_p(D)} =
\biggl((2 t)^{-d} \int_{ t B^d} \| \Delta_{\xi}^l f\|_{L_p(D_{l \xi})}^p d\xi\biggr)^{1 /p} = \\
\biggl((2 t)^{-d} \int_{ t B^d} \int_{D_{l \xi}} |\Delta_{\xi}^l f(x)|^p dx
d\xi\biggr)^{1 /p}, \ t \in \R_+.
\end{multline*}

Из приведенных определений видно, что
\begin{multline*} \tag{1.1.2}
\Omega^{\prime l} (f, t)_{L_p(D)} \le
\Omega^l (f, t)_{L_p(D)}, \ t \in \R_+,
f \in L_p(D), 1 \le p < \infty, \\
l \in \Z_+, D \text{ -- произвольная область в } \R^d.
\end{multline*}

Пусть теперь $ d \in \N, \alpha \in \R_+, 1 \le p < \infty $ и $ D $ --
область в $ \R^d. $ Тогда зададим число $ l = l(\alpha) \in \N, $ полагая
$ l(\alpha) = \min \{m \in \N: \alpha < m \}, $ а также выберем число $ \l \in \Z_+ $
такое, что $ \l < \alpha, $ и через
$ (H_p^{\alpha})^{\l}(D) $ обозначим пространство всех функций
$ f \in W_p^{\l}(D), $ для которых при $ \lambda \in \Z_{+ \l}^{\prime \,d} $
конечна величина
$$
\sup_{t \in \R_+} t^{-(\alpha -\l)} \cdot \Omega^{l -\l}(\D^{\lambda} f, t)_{L_p(D)} < \infty,
$$
а через $ (\mathcal H_p^{\alpha})^{\l}(D) $ -- множество функций $ f \in
W_p^{\l}(D), $ для которых при $ \lambda \in \Z_{+ \l}^{\prime \,d} $ выполняется неравенство
$$
\sup_{t \in \R_+} t^{-(\alpha -\l)} \Omega^{l -\l}(\D^{\lambda} f,
t)_{L_p(D)} \le 1.
$$

В пространстве $ (H_p^\alpha)^{\l}(D) $ вводится норма
$$
\| f \|_{(H_p^\alpha)^{\l}(D)} = \max(\| f \|_{W_p^{\l}(D)},
\max_{ \lambda \in \Z_{+ \l}^{\prime \,d}}
\sup_{t \in \R_+} t^{-(\alpha -\l)} \Omega^{l -\l}(\D^{\lambda} f, t)_{L_p(D)}).
$$

При тех же условиях на $ \alpha, p, D $ обозначим через
$ (H_p^\alpha)^\prime(D) $ пространство всех функций
$ f \in L_p(D), $ для которых выполняется условие
$$
\sup_{t \in \R_+} t^{-\alpha} \cdot \Omega^{\prime l}(f, t)_{L_p(D)} < \infty,
$$
а через $ (\mathcal H_p^\alpha)^\prime(D) $ -- множество функций
$ f \in L_p(D), $ обладающих тем свойством, что соблюдается неравенство
$$
\sup_{t \in \R_+} t^{-\alpha} \cdot \Omega^{\prime l}(f, t)_{L_p(D)} \le 1,
\text{ где } l = l(\alpha).
$$

В пространстве $ (H_p^\alpha)^\prime(D) $ задаётся норма
$$
\| f \|_{(H_p^\alpha)^\prime(D)} = \max(\| f \|_{L_p(D)}, \sup_{t \in \R_+}
t^{-\alpha} \Omega^{\prime l}(f, t)_{L_p(D)}).
$$

Для области $ D \subset \R^d $ при $ \alpha \in \R_+,\ 1 \le p < \infty,
\theta \in \R: 1 \le \theta < \infty, $
полагая, как и выше, $ l = l(\alpha) \in \N, $ и выбирая $ \l \in \Z_+:
\l < \alpha, $ через $ (B_{p,\theta}^\alpha)^{\l}(D) $ обозначим пространство
всех функций $ f \in W_p^{\l}(D), $ которые для любого
$ \lambda \in \Z_{+ \l}^{\prime d} $ удовлетворяют условию
$$
\int_0^\infty t^{-1 -\theta (\alpha -\l)}
(\Omega^{l -\l}(\D^{\lambda}f, t)_{L_p(D)})^{\theta}\,dt < \infty,
$$
а через $ (\mathcal B_{p,\theta}^\alpha)^{\l}(D) $ обозначим множество
всех функций $ f \in W_p^{\l}(D), $ для которых при любом
$ \lambda \in \Z_{+ \l}^{\prime d} $ соблюдается неравенство
$$
\int_0^\infty t^{-1 -\theta (\alpha -\l)}
(\Omega^{l -\l}(\D^{\lambda}f, t)_{L_p(D)})^{\theta}\,dt \le 1.
$$

В пространстве $ (B_{p,\theta}^\alpha)^{\l}(D) $ фиксируется норма
$$
\| f \|_{(B_{p,\theta}^\alpha)^{\l}(D)} = \max\biggl(\| f \|_{W_p^{\l}(D)},
\max_{ \lambda \in \Z_{+ \l}^{\prime d} } \left(\int_0^\infty
t^{-1 -\theta (\alpha -\l)} (\Omega^{l -\l}(\D^{\lambda}f,
t)_{L_p(D)})^{\theta}\,dt \right)^{1/\theta} \biggr).
$$

При $ \theta = \infty $ положим $ (B_{p,\infty}^\alpha)^{\l}(D) =
(H_p^\alpha)^{\l}(D). $

При тех же условиях на параметры  обозначим через
$ (B_{p,\theta}^\alpha)^\prime(D) $ пространство всех функций
$ f \in L_p(D), $ которые при $ l = l(\alpha) $ подчинены условию
$$
\int_0^\infty t^{-1 -\theta \alpha}
(\Omega^{\prime l}(f, t)_{L_p(D)})^{\theta}\,dt < \infty,
$$
а через $ (\mathcal B_{p,\theta}^\alpha)^\prime(D) $ -- множество функций
$ f \in L_p(D), $ которые удовлетворяют неравенству
$$
\int_0^\infty t^{-1 -\theta \alpha}
(\Omega^{\prime l}(f, t)_{L_p(D)})^{\theta}\,dt \le 1.
$$

В пространстве $ (B_{p,\theta}^\alpha)^\prime(D) $ определяется норма
$$
\| f \|_{(B_{p,\theta}^\alpha)^\prime(D)} = \max\biggl(\| f \|_{L_p(D)},
\left(\int_0^\infty
t^{-1 -\theta \alpha} (\Omega^{\prime l}(f, t)_{L_p(D)})^{\theta}\,dt
\right)^{1/\theta} \biggr).
$$

При $ \theta = \infty $ положим $ (B_{p,\infty}^\alpha)^\prime(D) =
(H_p^\alpha)^\prime(D). $

В случае, когда число $ \l = \l(\alpha) \in \Z_+ $ задаётся равенством
$ \l(\alpha) = \max \{m \in \Z_+: m < \alpha\}, $
пространство $ (B_{p,\theta}^\alpha)^{\l}(D) ((H_p^\alpha)^{\l}(D))$ обычно
обозначается $ B_{p,\theta}^\alpha(D) (H_p^\alpha(D)).$

Неоднократно будет использоваться применявшаяся в [5]

Лемма 1.1.1

Пусть $ l, m \in \Z_+,\ d \in \N,\ 1 \le p \le \infty $ и $ D $ ---
область в $ \R^d. $ Тогда для $ f \in W_p^l(D) $
при $ h \in \R^d $ справедлива оценка
\begin{equation*} \tag{1.1.3}
\| \Delta_h^{l +m} f\|_{L_p(D_{(l +m)h})} \le l! \cdot \|h\|^l
\sum_{\lambda \in \Z_{+ l}^{\prime \,d}} \| \Delta_h^m \D^\lambda f\|_{L_P(D_{m h})}.
\end{equation*}

В тех же обозначениях, что и выше, принимая во внимание то обстоятельство,
что для $ f \in (B_{p,\theta}^\alpha)^{\l}(D) $ при
$ \lambda \in \Z_{+ \l}^{\prime d} $ и $ t >0 $ справедлива оценка
\begin{multline*}
t^{-(\alpha -\l)} \Omega^{l -\l}(\D^\lambda f, t)_{L_p(D)} =
(t^{-1 -\theta(\alpha -\l)} \cdot
(\Omega^{l -\l}(\D^\lambda f, t)_{L_p(D)})^\theta t)^{1/\theta} \\
 \le \left(2^{1 +\theta(\alpha -\l)} \cdot \int_t^{2t} \tau^{-1 -\theta(\alpha -\l)}
(\Omega^{l -\l}(\D^\lambda f, \tau)_{L_p(D)})^\theta
\,d\tau\right)^{1/\theta} \\
 \le 2^{1/\theta +\alpha -\l}
\left(\int_0^\infty \tau^{-1 -\theta(\alpha -\l)}
(\Omega^{l -\l}(\D^\lambda f, \tau)_{L_p(D)})^\theta \,d\tau\right)^{1/\theta},
\end{multline*}
заключаем, что
\begin{multline*}
(B_{p,\theta}^\alpha)^{\l}(D) \subset (H_p^\alpha)^{\l}(D), \\
D \text{ -- область в } \R^d, \alpha \in \R_+, 1 \le p < \infty,
1 \le \theta < \infty, \l \in \Z_+: \l < \alpha.
\end{multline*}

Учитывая, что для $ f \in (B_{p, \theta}^\alpha)^\prime(D), \alpha \in
\R_+, 1 \le p, \theta < \infty, D $ -- область в $ \R^d, $ при
$ t \in \R_+, l = l(\alpha) $ выполняется неравенство
\begin{multline*}
t^{-\alpha} \Omega^{\prime l}(f, t)_{L_p(D)} = (t^{-1 -\theta \alpha} \cdot
(\Omega^{\prime l}(f, t)_{L_p(D)})^\theta t)^{1/\theta} \\
 \le \left(\int_t^{2t} 2^{1 +\theta \alpha} \cdot \tau^{-1 -\theta \alpha}
(2^{d /p} \Omega^{\prime l}(f, \tau)_{L_p(D)})^\theta
\,d \tau \right)^{1/\theta} \\
\le 2^{\alpha +1/\theta +d /p} \biggl(\int_{\R_+} \tau^{-1 -\theta \alpha}
(\Omega^{\prime l}(f,\tau)_{L_p(D)})^\theta d\tau \biggr)^{1/\theta},
\end{multline*}
заключаем, что
\begin{multline*} \tag{1.1.4}
(B_{p, \theta}^\alpha)^\prime(D) \subset
(H_p^\alpha)^\prime(D);
(\mathcal B_{p, \theta}^\alpha)^\prime(D) \subset
c_1(\alpha,d) (\mathcal H_p^\alpha)^\prime(D); \\
\| f\|_{(H_p^\alpha)^\prime(D)} \le c_1(\alpha,d)
\| f\|_{(B_{p, \theta}^\alpha)^\prime(D)},
\end{multline*}
где $ c_1(\alpha,d) = 2^{\alpha +d +1}. $
Из (1.1.2) и (1.1.3) следует, что для $ f \in (B_{p, \theta}^\alpha)^{\l}(D),
\l \in \Z_+: \l < \alpha, l = l(\alpha) $ при $ t \in \R_+ $ выполняется
неравенство
\begin{multline*}
\Omega^{\prime l}(f, t)_{L_p(D)} \le
\Omega^l(f, t)_{L_p(D)} = \supvrai_{ \xi \in t B^d}
\| \Delta_{\xi}^{\l +(l -\l)} f\|_{L_p(D_{(\l +(l -\l)) \xi })} \le \\
\supvrai_{ \xi \in t B^d} \l! \cdot \|\xi\|^{\l}
\sum_{\lambda \in \Z_{+ \l}^{\prime \,d}} \| \Delta_\xi^{l -\l}
\D^\lambda f\|_{L_P(D_{(l -\l) \xi})} \le \\
\l! t^{\l} \cdot \sum_{\lambda \in \Z_{+ \l}^{\prime d}}
\supvrai_{ \xi \in t B^d}
\| \Delta_\xi^{l -\l} \D^\lambda f\|_{L_P(D_{(l -\l) \xi})} =
\l! t^{\l} \sum_{\lambda \in \Z_{+ \l}^{\prime \,d}}
\Omega^{l -\l}(\D^{\lambda} f, t)_{L_p(D)},
\end{multline*}
и, значит,
\begin{equation*} \tag{1.1.5}
(B_{p, \theta}^\alpha)^{\l}(D) \subset (B_{p, \theta}^\alpha)^\prime(D), \\
(\mathcal B_{p, \theta}^\alpha)^{\l}(D) \subset c_2(\l,d)
(\mathcal B_{p, \theta}^\alpha)^\prime(D)
\end{equation*}
и
\begin{multline*} \tag{1.1.6}
\| f\|_{(B_{p, \theta}^\alpha)^\prime(D)} \le c_2
\| f\|_{(B_{p, \theta}^\alpha)^{\l}(D)}, \
\alpha \in \R_+, 1 \le p < \infty, 1 \le \theta \le \infty, \\
D \text{ -- произвольная область в } \R^d.
\end{multline*}

Обозначим через $ C^\infty(D) $ пространство бесконечно дифференцируемых
функций в области $ D \subset \R^d, $ а через $ C_0^\infty(D) $ -- пространство
функций $ f \in C^\infty(\R^d), $ у которых носитель $ \supp f \subset D. $

В заключение этого пункта введём ещё несколько обозначений.
Для банахова пространства $ X $ (над $ \R$) обозначим $ B(X) = \{x \in X:
\|x\|_X \le 1\}. $

Для банаховых пространств $ X,Y $ через $ \mathcal B(X,Y) $ обозначим банахово
пространство, состоящее из непрерывных линейных операторов $ T: X \mapsto Y, $
с нормой
$$
\|T\|_{\mathcal B(X,Y)} = \sup_{x \in B(X)} \|Tx\|_Y.
$$
Отметим, что если $ X=Y,$ то $ \mathcal B(X,Y) $ является банаховой алгеброй.
\bigskip

1.2. В этом пункте рассматриваются некоторые классы областей, в которых определены функции из пространств, изучаемых в работе.

Сначала введём обозначения.
Для $ t \in \R^d $ через $ 2^t $ обозначим вектор
$ 2^t = (2^{t_1}, \ldots, 2^{t_d}). $
При $ d \in \N, \kappa \in \Z_+^d, \nu \in \Z^d $ обозначим
$ Q_{\kappa, \nu}^d = 2^{-\kappa} \nu +2^{-\kappa} I^d, $
$ \overline Q_{\kappa, \nu}^d = 2^{-\kappa} \nu +2^{-\kappa} \overline I^d, $
а при $ k \in \Z_+, \nu \in \Z^d $ положим
$ Q_{k, \nu}^d = Q_{k \e, \nu}^d = 2^{-k} \nu +2^{-k} I^d, $
$ \overline Q_{k, \nu}^d = \overline Q_{k \e, \nu}^d = 2^{-k} \nu +
2^{-k} \overline I^d. $

Определим семейство классов областей, зависящее от параметра
$ \alpha \in \R_+^d, $ одним из элементов которого (при $ \alpha = \e $ )
является класс областей, с которыми мы будем иметь дело ниже.

Определение 1

При $ d \in \n, \alpha \in \R_+^d $ будем говорить, что область $
D \subset \R^d $ является областью $ \alpha $-типа, если
существуют константы $ K^0(d,D,\alpha) \in \Z_+, \mathpzc
k^0(d,D,\alpha) \in \Z_+, \Gamma^0(d,D,\alpha) \in \R_+,
c_0(d,D,\alpha) >0 $ такие, что для $ k \in \Z_+: k \ge K^0, $
при $ \kappa \in \Z_+^d $ с компонентами $ \kappa_j = [k
/\alpha_j], \ j =1,\ldots,d, $ выполняются следующие условия:

1) для каждого $ \boldsymbol{\nu} \in \Z^d: Q_{\kappa,\boldsymbol{\nu}}^d \cap
D \ne \emptyset, $ существует $ \boldsymbol{\nu^\prime} \in \Z^d $ такой, что соблюдается включение
\begin{equation*}
\overline Q_{\kappa,\boldsymbol{\nu^\prime}}^d \subset
D \cap (2^{-\kappa} \boldsymbol{\nu} +\Gamma^0 2^{-\kappa} B^d);
\end{equation*}

2) для любых $ \boldsymbol{\nu}, \boldsymbol{\nu^\prime} \in \Z^d $ таких,
что $ \overline Q_{\kappa,\boldsymbol{\nu}}^d \subset D, \
\overline Q_{\kappa,\boldsymbol{\nu^\prime}}^d \subset D, $ существуют
последовательнности $ \nu^\iota \in \Z^d, \ \iota =0,\ldots,\Iota; \
j^\iota \in \{1,\ldots,d\}, \ \epsilon^\iota \in \{-1,1\}, \
\iota =0,\ldots,\Iota -1, $ со следующими свойствами:

\begin{equation*} \tag{1.2.1}
\Iota \le c_0 \|\boldsymbol{\nu} -\boldsymbol{\nu^\prime}\|;
\end{equation*}

\begin{equation*} \tag{1.2.2}
\overline Q_{\mathpzc k^0 \e +\kappa,\nu^0}^d \subset
\overline Q_{\kappa,\boldsymbol{\nu}}^d, \\
\overline Q_{\mathpzc k^0 \e +\kappa,\nu^{\Iota}}^d \subset
\overline Q_{\kappa,\boldsymbol{\nu^\prime}}^d, \\
\nu^{\iota +1} = \nu^\iota +\epsilon^\iota e_{j^\iota}, \ \iota =0,\ldots,\Iota -1;
\end{equation*}

при $ \iota =0,\ldots,\Iota $ справедливо включение
\begin{equation*} \tag{1.2.3}
\overline Q_{\mathpzc k^0 \e +\kappa,\nu^\iota}^d \subset D.
\end{equation*}

Отметим, что класс областей $ \alpha $-типа в смысле определения 1 шире, чем
класс областей $ \alpha $-типа в смысле определения из [9]. Однако легко
проверить (см. [10]), что все утверждения из [9], относящиеся к функциям,
определённым в областях $\alpha$-типа в указанном там смысле,
без изменения формулировок справедливы  и в случае, когда область $\alpha$-типа
понимается в смысле определения 1.

Напомним определение области с липшицевой границей.

Определение 2

Будем говорить, что область $ D \subset \R^d $ является областью с
липшицевой границей, если для каждой  точки $ x^0, $ принадлежащей границе
$ \mathcal D D $ области $ D $ существуют ортогональное аффинное преобразование
$ h: \R_y^{d-1} \times \R_z^1 \mapsto \R_x^d, \ x = h(y,z), $ окрестность $ V $
точки $ 0 $ в $ \R_y^{d-1}, $ окрестность $ W $ точки $ 0 $ в
$ \R_z^1, $ функция $ \phi: V \ni y \mapsto \phi(y) = z \in W $ и константа
$ M > 0 $ такие, что

1)
\begin{equation*}
h(0) = x^0;
\end{equation*}

2)
\begin{equation*} \tag{1.2.4}
h^{-1}(D) \cap (V \times W) = \{(y,z) \in V \times W: z > \phi(y)\}, \ \phi(0) =0;
\end{equation*}

3) для любых $ y_1, y_2 \in V $ соблюдается неравенство
\begin{equation*} \tag{1.2.5}
| \phi(y_1) -\phi(y_2)| \le M \| y_1 -y_2\|.
\end{equation*}

Например, область
$$
D = \inter (\cup_{k \in \Z_+} (2^{-k} e_1 +2^{-k} \overline I^2)) \subset \R^2
$$
является областью $ \e $-типа, но не является областью с липшицевой границей.

Однако имеет место следующая теорема.

Теорема 1.2.1

Пусть $ D $ -- ограниченная область в $ \R^d $ с липшицевой границей. Тогда
$ D $ является областью $\e$-типа.

Для доказательства теоремы потребуются вспомогательные утверждения.

Лемма 1.2.2

При $ d \in \N $ существует константа $ c_1(d) > 0 $ такая, что для любого
$ k \in \Z_+ $ и любых $ x^0, \xi \in \R^d $ существуют последовательности
$ \nu^\iota \in \Z^d, \ \iota =0,\ldots,\Iota; \ j^\iota \in \Nu_{1,d}^1, \
\epsilon^\iota \in \{-1,1\}, \ \iota =0,\ldots,\Iota -1, $ со следующими
свойствами:

1) \begin{equation*}
x^0 \in \overline Q_{k,\nu^0}^d, \\
(x^0 +\xi) \in \overline Q_{k,\nu^{\Iota}}^d;
\end{equation*}

2) для каждого $ \iota = 0,\ldots,\Iota $ существует $ n \in \Z^d, $ для
которого существует $ t \in [0,1]: x_t = (x^0 +t \xi) \in \overline Q_{k,n}^d $
и $ \| \nu^\iota -n \| \le 1; $

3) \begin{equation*}
\nu^{\iota +1} = \nu^\iota +\epsilon^\iota e_{j^\iota}, \ \iota =0,\ldots,\Iota -1;
\end{equation*}

4) \begin{equation*}
\Iota \le c_1 (2^k \| \xi \| +1).
\end{equation*}

Доказательство.

Для $ k \in \Z_+, x^0, \xi \in \R^d $ обозначим
$$
\mathfrak N_k(x^0,\xi) = \{n \in \Z^d \mid \ \exists t \in [0,1]: x_t =
(x^0 +t \xi) \in \overline Q_{k,n}^d\}
$$
и оценим сверху $ \card \mathfrak N_k(x^0,\xi). $
Для этого обозначим через $ L $ линейную оболочку вектора $ \xi, $ а через
$ L^\perp $ -- ортогональное дополнение к $ L $ в $ l_2^d. $ Фиксируем вектор
$ f_1 = \| \xi \|_{l_2^d}^{-1} \xi \in L $ и ортонормированную систему векторов
$ \{f_2,\ldots,f_d\} \subset L^\perp. $  Пусть
$ x \in \cup_{n \in \mathfrak N_k(x^0,\xi)} Q_{k,n}^d. $
Тогда, выбирая $ n \in \mathfrak N_k(x^0,\xi) $ и $ t \in [0,1] $ такие,
что $ x \in Q_{k,n}^d, x_t \in \overline Q_{k,n}^d, $ имеем
\begin{equation*}
x -x_t = \sum_{j =1}^d z_j f_j
\end{equation*}
или
\begin{equation*}
x = x_t +\sum_{j =1}^d z_j f_j = x^0 +t \| \xi \|_{l_2^d} f_1 +\sum_{j =1}^d z_j f_j =
x^0 +(t \| \xi \|_{l_2^d} +z_1) f_1 +\sum_{j =2}^d z_j f_j,
\end{equation*}
причём вследствие теоремы Пифагора справедливы неравенства
\begin{equation*}
| z_1| = \| z_1 f_1 \|_{l_2^d} \le \| x -x_t \|_{l_2^d}, \\
\| \sum_{j =2}^d z_j f_j \|_{l_2^d} \le \| x -x_t \|_{l_2^d},
\end{equation*}
а для $ x = 2^{-k} n +2^{-k} \eta, \ \eta \in I^d; \ x_t = 2^{-k} n +2^{-k} \eta_t, \
\eta_t \in \overline I^d, $ выполняется неравенство
\begin{equation*}
\| x -x_t \|_{l_2^d} = \| 2^{-k} \eta -2^{-k} \eta_t \|_{l_2^d} =
2^{-k} \| \eta -\eta_t \|_{l_2^d} \le d^{1/2} 2^{-k},
\end{equation*}
что в соединении с предыдущим приводит к соотношению
\begin{multline*} \tag{1.2.6}
x = x^0 +y_1 f_1 +\sum_{j =2}^d y_j f_j, \text{ где } y_1 =
(t \| \xi \|_{l_2^d} +z_1) \in [-d^{1/2} 2^{-k}, \| \xi \|_{l_2^d} +d^{1/2} 2^{-k}], \\
\| (y_2,\ldots,y_d) \|_{l_2^{d -1}} = \| \sum_{j =2}^d y_j f_j \|_{l_2^d} \le
d^{1/2} 2^{-k}, \ y_j = z_j, \ j =2,\ldots,d.
\end{multline*}

Определяя ортогональный линейный оператор $ F: l_2^d \mapsto l_2^d $
равенством $ F y = \sum_{j =1}^d y_j f_j, \ y = (y_1,\ldots,y_d) \in \R^d, $
и обозначая через $ \mathcal F = \{\mathcal F_{i,j}, \ i,j =1,\ldots,d\}, $
соответствующую ему ортогональную матрицу, т.е. матрицу, обладающую тем
свойством, что для $ y \in \R^d $ соблюдаются равенства $ (F y)_i =
\sum_{j =1}^d \mathcal F_{i,j} y_j, \ i =1,\ldots,d, $ из (1.2.6) видим, что
имеет место включение
\begin{equation*} \tag{1.2.7}
\cup_{n \in \mathfrak N_k(x^0,\xi)} Q_{k,n}^d \subset
(x^0 +F(\mathfrak C_k)), \
\end{equation*}
где $ \mathfrak C_k = ([-d^{1/2} 2^{-k}, \| \xi \|_{l_2^d}
+d^{1/2} 2^{-k}] \times (d^{1/2} 2^{-k} B(l_2^{d -1}))). $

Поскольку кубы в левой части (1.2.7) попарно не пересекаются, то
из (1.2.7) следует, что
\begin{multline*}
2^{-kd} \card \mathfrak N_k(x^0,\xi) = \sum_{n \in \mathfrak N_k(x^0,\xi)}
\mes Q_{k,n}^d = \mes (\cup_{n \in \mathfrak N_k(x^0,\xi)} Q_{k,n}^d) \\
\le \mes(x^0 +F(\mathfrak C_k)) = \mes(F(\mathfrak C_k)) =
\int_{F(\mathfrak C_k)} dx = \int_{\mathfrak C_k} | \det \mathcal F| dy = \\
\int_{\mathfrak C_k} dy = \int_{[-d^{1/2} 2^{-k}, \| \xi \|_{l_2^d} +d^{1/2} 2^{-k}]
\times (d^{1/2} 2^{-k} B(l_2^{d -1}))} dy = \\
\int_{[-d^{1/2} 2^{-k}, \| \xi \|_{l_2^d} +d^{1/2} 2^{-k}]}
(\int_{ d^{1/2} 2^{-k} B(l_2^{d -1})} dy_2 \ldots dy_d) dy_1 = \\
\mes([-d^{1/2} 2^{-k}, \| \xi \|_{l_2^d} +d^{1/2} 2^{-k}])
\cdot \mes(d^{1/2} 2^{-k} B(l_2^{d -1})) = \\
(\| \xi \|_{l_2^d} +2 d^{1/2} 2^{-k})
c_2(d) 2^{-k (d -1)},
\end{multline*}
откуда
\begin{multline*} \tag{1.2.8}
\card \mathfrak N_k(x^0,\xi) \le 2^{kd} (\| \xi \|_{l_2^d} +2 d^{1/2} 2^{-k})
c_2(d) 2^{-k (d -1)} = \\ c_2 ( 2^k \| \xi \|_{l_2^d} +2 d^{1/2}) \le
c_3(d) ( 2^k \| \xi \| +1).
\end{multline*}

Отметим ещё некоторые факты, используемые ниже.

Замечание 1

При $ k \in \Z_+, \ x^0, \xi \in \R^d, \ t^0 \in \R $ существует $
n \in \Z^d, $ для которого $ x_{t^0} \in \overline Q_{k,n}^d $ и
существует $ t > t^0: x_t \in \overline Q_{k,n}^d. $ В самом деле,
исходя из того, что $ \R^d = \cup_{n^\prime \in \Z^d} \overline
Q_{k,n^\prime}^d, $ выберем $ n^\prime \in \Z^d, $ для которого $
x_{t^0} \in \overline Q_{k,n^\prime}^d $ и определим при $ j
=1,\ldots,d $ компоненты $ n \in \Z^d $ соотношениями:

если $ (x_{t^0})_j = 2^{-k} n_j^\prime, $ то
\begin{equation*}
n_j = \begin{cases} n_j^\prime \text{ при } \xi_j \ge 0; \\
n_j^\prime -1 \text{ при } \xi_j < 0,
\end{cases}
\end{equation*}
если $ (x_{t^0})_j = 2^{-k} (n_j^\prime +1), $ то
\begin{equation*}
n_j = \begin{cases} n_j^\prime \text{ при } \xi_j < 0; \\
n_j^\prime +1 \text{ при } \xi_j \ge 0,
\end{cases}
\end{equation*}
если $ 2^{-k} n_j^\prime < (x_{t^0})_j < 2^{-k} (n_j^\prime +1), $ то
$ n_j = n_j^\prime. $

Замечание 2

Если при $ k \in \Z_+ $ для $ n, n^\prime \in \Z^d $ существует
$ x \in \overline Q_{k,n}^d \cap \overline Q_{k,n^\prime}^d, $ то
$ \| n -n^\prime \| \le 1. $ Действительно, в условиях замечания
при $ j =1,\ldots,d $ имеем
$$
2^{-k} n_j \le x_j \le 2^{-k} n_j^\prime +2^{-k},
$$
откуда
$$
n_j \le n_j^\prime +1,
$$
точно также
$$
n_j^\prime \le n_j +1,
$$
или
$$
n_j^\prime -1 \le n_j,
$$
т.е.
$$
-1 \le n_j -n_j^\prime \le 1
$$
и
$$
| n_j -n_j^\prime | \le 1.
$$

Пусть теперь $ k \in \Z_+, \ x^0, \xi \in \R^d. $ Построим последовательности
$ \{ t^s \ge 0, \ s \in \Z_+\} $ и $ \{ n^s \in \Z^d, \ s \in \Z_+\} $
следующим образом. Положим $ t^0 = 0 $ и фиксируем $ n^0 \in \Z^d, $ для
которого $ x^0 = x_{t^0} \in \overline Q_{k,n^0}^d $ и существует $ t > t^0:
x_t \in \overline Q_{k,n^0}^d $ (см. замечание 1). При $ s \in \Z_+ $ положим
$$
t^{s +1} = \max \{t \in \R: x_t \in \overline Q_{k,n^s}^d\},
$$
и ввиду замечания 1 выберем $ n^{s +1} \in \Z^d $ так, что $ x_{t^{s +1}} \in
\overline Q_{k,n^{s +1}}^d $ и существует $ t > t^{s +1}: x_t \in
\overline Q_{k,n^{s +1}}^d. $ Ясно, что $ t^s < t^{s +1}, \ s \in \Z_+; \
n^s \ne n^{s^\prime} $ при $ s,s^\prime \in \Z_+: s \ne s^\prime. $
Учитывая эти обстоятельства, а также то, что для $ s \in \Z_+: t^s \in [0,1], $
мультииндекс $ n^s $ принадлежит $ \mathfrak N_k(x^0,\xi), $ с учётом
(1.2.8) заключаем, что
\begin{equation*} \tag{1.2.9}
\card \{ s \in \Z_+: t^s \in [0,1]\} \le \card \mathfrak N_k(x^0,\xi) < \infty.
\end{equation*}
Принимая во внимание сказанное и полагая $ S = \max_{ s \in \Z_+: t^s \in [0,1]} s, $
получаем, что $ t^s \in [0,1], \ s =0,\ldots,S, \ t^s > 1 $ при $ s > S. $
С учётом этого из (1.2.9) и (1.2.8) выводим
\begin{equation*} \tag{1.2.10}
S +1 \le c_3 ( 2^k \| \xi \| +1).
\end{equation*}

Отметим ещё, что вследствие принадлежности $ x_{t^S}, x_{t^{S +1}}
\in \overline Q_{k,n^S}^d, $ а также выпуклости этого куба и соотношения
$ t^S \le 1 < t^{S +1} $ точка
\begin{equation*} \tag{1.2.11}
x^0 +\xi = x_1 \in \overline Q_{k,n^S}^d.
\end{equation*}

Далее, поскольку $ n^s \ne n^{s +1} $ и
\begin{equation*} \tag{1.2.12}
x_{t^{s +1}} \in \overline Q_{k,n^s}^d \cap \overline Q_{k,n^{s +1}}^d,
\end{equation*}
то в силу замечания 2 заключаем, что
$ \| n^s -n^{s +1} \| =1, \ s \in \Z_+. $ Полагая
\begin{equation*} \tag{1.2.13}
\Iota_s = \card \{j =1,\ldots,d: n_j^s -n_j^{s +1} \ne 0\} \le d,
\end{equation*}
видим, что существуют последовательности $ \{j_s^{\iota_s} \in \Nu_{1,d}^1, \
\iota_s =1,\ldots,\Iota_s\}, \
\{\epsilon_s^{\iota_s} \in \{-1,1\}, \ \iota_s =1,\ldots,\Iota_s\} $ такие, что
\begin{equation*} \tag{1.2.14}
n_{j_s^{\iota_s}}^{s +1} -n_{j_s^{\iota_s}}^s = \epsilon_s^{\iota_s}, \
\iota_s =1,\ldots,\Iota_s, \\
n_j^{s +1} = n_j^s, \ j \in \Nu_{1,d}^1: j \ne j_s^{\iota_s}, \
\iota_s =1,\ldots,\Iota_s, \ s \in \Z_+.
\end{equation*}

Теперь построим последовательности
$ \nu^\iota \in \Z^d, \ \iota =0,\ldots,\Iota; \ j^\iota \in \Nu_{1,d}^1,
\epsilon^\iota \in \{-1,1\}, \ \iota =0,\ldots,\Iota -1, $ с требуемыми
свойствами. Сначала положим $ \nu^0 =
n^0. $ А для $ \iota = \Iota_0 + \ldots +\Iota_{s -1} +\iota_s, \
\iota_s =1,\ldots,\Iota_s, \ s =0,\ldots,S -1, $
с учётом (1.2.14) определим
\begin{equation*} \tag{1.2.15}
\nu^\iota = n^s +\sum_{i =1}^{\iota_s} \epsilon_s^i e_{j_s^i}, \\
j^{\iota -1} = j_s^{\iota_s},
\epsilon^{\iota -1} = \epsilon_s^{\iota_s}.
\end{equation*}
При этом для $ \Iota = \sum_{s =0}^{S -1} \Iota_s $ имеет место равенство
$ \nu^{\Iota} = n^S, $ и ввиду (1.2.13), (1.2.10) выполняется неравенство
\begin{equation*} \tag{1.2.16}
\Iota \le \sum_{s =0}^{S -1} d = d S < d (S +1) \le d c_3 ( 2^k \| \xi \| +1) =
c_1 ( 2^k \| \xi \| +1).
\end{equation*}
Сопоставляя сказанное, в частности, (1.2.11), (1.2.12), (1.2.15), (1.2.16), убеждаемся в
справедливости пп. 1) -- 4) леммы 1.2.2. $ \square $

Лемма 1.2.3

Пусть $ D \subset \R^d $ -- область с липшицевой границей. Тогда для каждой
точки $ x^0 \in \overline D $ существуют константы $ \delta^0(d,D,x^0) > 0,
K^1(d,D,x^0) \in \Z_+, \mathpzc k^1(d,D,x^0) \in \Z_+, \gamma^1(d,D,x^0) \in \R_+,
c_4(d,D,x^0) >0 $ такие, что для $ U_{x^0} = x^0 +2 \delta^0 B_0^d $ при
$ k \in \Z_+: k \ge K^1, $ выполняются следующие условия:

1) для каждого $ \nu \in \Z^d: Q_{k,\nu}^d \cap (D \cap U_{x^0}) \ne \emptyset, $
существует $ \nu^\prime \in \Z^d $ такой, что соблюдается включение
\begin{equation*} \tag{1.2.17}
\overline Q_{k,\nu^\prime}^d \subset
D \cap (2^{-k} \nu +2^{-k} \gamma^1 B^d);
\end{equation*}

2) для любых $ \nu, \nu^\prime \in \Z^d $ таких,
что $ \overline Q_{k,\nu}^d \subset D \cap U_{x^0},
\overline Q_{k,\nu^\prime}^d \subset D \cap U_{x^0}, $
при $ \mathpzc k \in \Z_+: \mathpzc k \ge \mathpzc k^1, $ существуют
последовательности $ \nu^\iota \in \Z^d, \ \iota =0,\ldots,\Iota; \
j^\iota \in \Nu_{1,d}^1, \epsilon^\iota \in \{-1,1\}, \ \iota =0,\ldots,\Iota -1, $
со следующими свойствами:

\begin{equation*} \tag{1.2.18}
\Iota \le c_4 2^{\mathpzc k} \|\nu -\nu^\prime\|;
\end{equation*}

\begin{equation*} \tag{1.2.19}
\overline Q_{\mathpzc k +k,\nu^0}^d \subset \overline Q_{k,\nu}^d, \\
\overline Q_{\mathpzc k +k,\nu^{\Iota}}^d \subset \overline Q_{k,\nu^\prime}^d, \\
\nu^{\iota +1} = \nu^\iota +\epsilon^\iota e_{j^\iota}, \ \iota =0,\ldots,\Iota -1;
\end{equation*}

при $ \iota =0,\ldots,\Iota $ справедливо включение
\begin{equation*} \tag{1.2.20}
\overline Q_{\mathpzc k +k,\nu^\iota}^d \subset D.
\end{equation*}

Доказательство.

Фиксируем точку $ x^0 \in \overline D. $ Сначала рассмотрим случай, когда
$ x^0 \in \mathcal D D. $ В этом случае возьмём отображения и окрестности,
соответствующие точке $ x^0, $ со свойствами, приведенными в определении 2.
Тогда существует линейный ортогональный оператор $ \mathcal A: l_2^d \mapsto l_2^d $
такой, что для $ x \in \R^d $ выполняется равенство
\begin{equation*} \tag{1.2.21}
h^{-1} x = \mathcal A(x -x^0).
\end{equation*}
При этом существует ортогональная матрица $ \{a_{i,j} \in \R: \ i,j =1,\ldots,d\} $
такая, что для $ x \in \R^d $ и $ (y,z) = \mathcal A x \in \R^{d -1} \times \R  $
имеют место равенства
\begin{equation*}
y_i = \sum_{j =1}^d a_{i,j} x_j, \ i =1,\ldots,d -1; \\
z = \sum_{j =1}^d a_{d,j} x_j,
\end{equation*}
или
\begin{equation*} \tag{1.2.22}
y = A x, z = (a,x), \
\end{equation*}
где $ A: \R^d \mapsto \R^{d -1} $ -- линейный оператор, определяемый
равенствами
\begin{equation*}
(A x)_i = \sum_{j =1}^d a_{i,j} x_j, \ i =1,\ldots,d -1, \ x \in \R^d,
\end{equation*}
а вектор $ a \in \R^d $ определяется равенством
$ a = (a_{d,1},\ldots,a_{d,d}). $
Из ортогональности матрицы $ \{a_{i,j}: i,j =1,\ldots,d\} $ следует, что
\begin{equation*} \tag{1.2.23}
(a,a) = \|a\|_{l_2^d}^2 =1, \\
A a = 0.
\end{equation*}

Далее, выберем $ \rho \in \R_+, $ для которого справедливо включение
\begin{equation*} \tag{1.2.24}
2 \rho B_0^1 \subset W,
\end{equation*}
и вследствие непрерывности функции $ \phi $ подберём $ r \in \R_+ $ так,
чтобы соблюдались соотношения
\begin{equation*} \tag{1.2.25}
2 r B_0^{d -1} \subset V, \\
| \phi(y) | < \rho \text{ при } y \in 2 r B_0^{d -1}.
\end{equation*}

Фиксировав $ \rho, r \in \R_+, $ для которых соблюдаются соотношения (1.2.24), (1.2.25),
возьмём произвольное $ \delta \in \R_+ $ такое, что для $ u_{x^0} =
x^0 +2 \delta B_0^d $ справедливо включение
\begin{equation*} \tag{1.2.26}
\{ \mathcal A(x -x^0): x \in u_{x^0} \} \subset (r B_0^{d -1}) \times (\rho B_0^1).
\end{equation*}

Тогда для $ x \in (D \cap u_{x^0}) $ в силу (1.2.26), (1.2.22), а также
благодаря (1.2.24), (1.2.25), (1.2.21), (1.2.4) имеем
\begin{equation*} \tag{1.2.27}
\| A(x -x^0) \| < r, \ | (a, x -x^0)| < \rho,
\end{equation*}
и
\begin{equation*} \tag{1.2.28}
(a, x -x^0) > \phi(A(x -x^0)).
\end{equation*}

Пусть $ k \in \Z_+ $ и $ \nu \in \Z^d $ таковы, что $ Q_{k,\nu}^d \cap
(D \cap u_{x^0}) \ne \emptyset. $ Беря $ x \in Q_{k,\nu}^d \cap
(D \cap u_{x^0}), $ для $ t \in \R_+ $ и $ \nu^\prime \in \Z^d $ такого, что
$ x +2^{-k} t a \in \overline Q_{k,\nu^\prime}^d, $ для $ x^\prime \in
\overline Q_{k,\nu^\prime}^d $ с учётом (1.2.27), (1.2.23) выводим
\begin{multline*} \tag{1.2.29}
\| A(x^\prime -x^0)\| = \| A(x^\prime -(x +2^{-k} t a) + (x
+2^{-k} t a) -x^0) \| = \\
\| A(x^\prime -(x +2^{-k} t a)) + A(x
+2^{-k} t a -x^0) \| = \\ \| A(x^\prime -(x +2^{-k} t a)) + A(x
-x^0) +2^{-k} t A a \| = \\
 \| A(x^\prime -(x +2^{-k} t a)) + A(x
-x^0) \| \le  \| A(x^\prime -(x +2^{-k} t a)) \| + \| A(x -x^0)\| \le \\
\| A \|_{\mathcal B(l_\infty^d, l_\infty^{d -1})} \|
x^\prime -(x +2^{-k} t a) \| +\| A(x -x^0) \| \le \\ 
c_5(d,D,x^0)
\max_{\xi, \xi^\prime \in \overline Q_{k,\nu^\prime}^d} \|
\xi^\prime -\xi \| +\| A(x -x^0) \| < c_5 2^{-k} +r,
\end{multline*}
\begin{multline*} \tag{1.2.30}
| (a, x^\prime -x^0)| = | (a, x^\prime -(x +2^{-k} t a) + (x
+2^{-k} t a) -x^0) | = \\
| (a, x^\prime -(x +2^{-k} t a)) + (a, x
+2^{-k} t a -x^0) | = \\ 
| (a, x^\prime -(x +2^{-k} t a)) + (a, x
-x^0) +2^{-k} t (a,a) | \le \\
| (a, x^\prime -(x +2^{-k} t a)) | + |
(a, x -x^0)| +2^{-k} t \le \| a \|_{l_1^d} \cdot \| x^\prime -(x
+2^{-k} t a) \| + | (a, x -x^0)| +2^{-k} t \le \\ 
c_6(d,D,x^0)
\max_{\xi, \xi^\prime \in \overline Q_{k,\nu^\prime}^d} \|
\xi^\prime -\xi \| +| (a, x -x^0)| +2^{-k} t < c_6 2^{-k} +\rho
+2^{-k} t.
\end{multline*}
Принимая во внимание (1.2.29), (1.2.30), (1.2.24), (1.2.25), видим, что в
рассматриваемой ситуации при соблюдении условий
\begin{equation*} \tag{1.2.31}
c_5 2^{-k} < r, \ c_6 2^{-k} +2^{-k} t < \rho,
\end{equation*}
имеют место включения
\begin{equation*} \tag{1.2.32}
A(x^\prime -x^0) = y^\prime \in V, \ (a, x^\prime -x^0) = z^\prime \in W,
\end{equation*}
и, благодаря (1.2.23), (1.2.28), (1.2.5), выполняется неравенство
\begin{multline*} \tag{1.2.33}
(a, x^\prime -x^0) -\phi(A(x^\prime -x^0)) =
(a, x^\prime -(x +2^{-k} t a) +
(x +2^{-k} t a) -x^0) -\phi(A(x^\prime -x^0)) = \\
(a, x^\prime -(x +2^{-k} t a)) +
(a, x +2^{-k} t a -x^0) -\phi(A(x^\prime -x^0)) = \\
(a, x^\prime -(x +2^{-k} t a)) +
(a, x -x^0) +2^{-k} t (a,a) -\phi(A(x^\prime -x^0)) = \\
2^{-k} t +(a, x^\prime -(x +2^{-k} t a)) +(a, x -x^0) -\phi(A(x +2^{-k} t a -x^0)) +\\
\phi(A(x +2^{-k} t a -x^0)) -
\phi(A(x^\prime -x^0)) = \\
2^{-k} t +(a, x^\prime -(x +2^{-k} t a)) +(a, x -x^0) -\phi(A(x -x^0)) +\\
\phi(A(x +2^{-k} t a -x^0)) -
\phi(A(x^\prime -x^0)) > \\
2^{-k} t -| (a, x^\prime -(x +2^{-k} t a))| -| \phi(A(x +2^{-k} t a -x^0))
-\phi(A(x^\prime -x^0))| \ge \\
2^{-k} t -\| a \|_{l_1^d} \| x^\prime -(x +2^{-k} t a) \| -
M \| A(x +2^{-k} t a -x^\prime)\| \ge \\
2^{-k} t -\| a \|_{l_1^d} \| x^\prime -(x +2^{-k} t a) \| -
M \| A\|_{\mathcal B(l_\infty^d,l_\infty^{d -1})}
\| x +2^{-k} t a -x^\prime\| = \\
2^{-k} t -(c_6 +M c_5) \| x +2^{-k} t a -x^\prime \| \ge
2^{-k} t -(c_6 +M c_5) 2^{-k} > 0
\end{multline*}
при условии
\begin{equation*} \tag{1.2.34}
t > c_6 +M c_5.
\end{equation*}

Сопоставляя (1.2.32), (1.2.33), (1.2.22), (1.2.21), (1.2.4),
заключаем, что для фиксированного $ t, $ удовлетворяющего
(1.2.34), при $ \delta > 0, $ для которого соблюдается (1.2.26),
при $ k \in \Z_+, $ подчинённого неравенствам (1.2.31), для $
\nu \in \Z^d $ такого, что $ Q_{k,\nu}^d \cap (D \cap u_{x^0}) \ne
\emptysetб\, x \in Q_{k,\nu}^d \cap (D \cap u_{x^0}), $ и $
\nu^\prime \in \Z^d: (x +2^{-k} t a) \in \overline
Q_{k,\nu^\prime}^d, $ верно включение $ \overline
Q_{k,\nu^\prime}^d \subset D, $ причём для $ x^\prime \in
\overline Q_{k,\nu^\prime}^d $ справедливо неравенство
\begin{multline*}
\| x^\prime -2^{-k} \nu \| \le \| x^\prime -(x +2^{-k} t a) \| +
\| x +2^{-k} t a -2^{-k} \nu \| \le \\
 \| x^\prime -(x +2^{-k} t a) \| +
\| x -2^{-k} \nu \| +2^{-k} t \|a\| \le \\
\max_{\xi, \xi^\prime \in \overline Q_{k,\nu^\prime}^d} \| \xi^\prime -\xi \| +
\max_{\xi, \xi^\prime \in \overline Q_{k,\nu}^d} \| \xi^\prime -\xi \| +
2^{-k} t = 2^{-k} +2^{-k} +2^{-k} t = (t +2) 2^{-k},
\end{multline*}
т.е. имеет место (1.2.17) при $ \gamma^1 = t +2.$

Перейдём к доказательству п. 2). Обозначим через $ \epsilon $ вектор с
координатами
\begin{equation*}
\epsilon_j = \begin{cases} 1, \text{ если } a_{d,j} < 0; \\
0, \text{ если } a_{d,j} \ge 0, \
\end{cases}
j =1,\ldots,d.
\end{equation*}
Тогда для $ j =1,\ldots,d: a_{d,j} < 0, $ имеем $ \epsilon_j +t a_{d,j} =
1 +t a_{d,j}, $ причём
$$
1 +t a_{d,j} \ge 0 \text{ при } t \le 1 /| a_{d,j}|; \\
1 +t a_{d,j} \le 1 \text{ при } t \ge 0,
$$
а для $ j =1,\ldots,d: a_{d,j} \ge 0, $ имеем $ \epsilon_j +t a_{d,j} =
t a_{d,j}, $ причём
$$
t a_{d,j} \ge 0 \text{ при } t \ge 0; \\
t a_{d,j} \le 1 \text{ при } t \le 1 / a_{d,j}.
$$
Таким образом, при $ 0 \le t \le \min_{j =1,\ldots,d} 1 / | a_{d,j}| $ вектор
$ \epsilon +t a \in \overline I^d. $ Положим
$$
t^0 = \max\{ t \in \R_+: \epsilon +t a \in \overline I^d\}.
$$

Теперь снова возьмём произвольное $ \delta \in \R_+ $ такое, что для
$ u_{x^0} = x^0 +2 \delta B_0^d $ соблюдается (1.2.26). Пусть $ k \in \Z_+ $ и
$ \nu, \nu^\prime \in \Z^d $ таковы, что
\begin{equation*} \tag{1.2.35}
\overline Q_{k,\nu}^d \subset D \cap u_{x^0}, \\
\overline Q_{k,\nu^\prime}^d \subset D \cap u_{x^0}.
\end{equation*}
Выберем точки
$$
x_\epsilon = 2^{-k} \nu +2^{-k} \epsilon +2^{-k} t^0 a \in
\overline Q_{k,\nu}^d, \\
x_\epsilon^\prime = 2^{-k} \nu^\prime +2^{-k} \epsilon +2^{-k} t^0 a \in
\overline Q_{k,\nu^\prime}^d.
$$
Тогда, учитывая, что для $ x \in Q_{k,\nu}^d (x \in Q_{k,\nu^\prime}^d) $ и
$ t \in \R_+ $ в силу (1.2.29), (1.2.30) (при $ x^\prime = x +2^{-k} t a $)
имеют место неравенства
\begin{equation*}
\| A(x +2^{-k} t a -x^0)\| < c_5 2^{-k} +r, \\
| (a, x +2^{-k} t a -x^0)| < c_6 2^{-k} +\rho +2^{-k} t,
\end{equation*}
вследствие непрерывности функций $ \| A(\cdot) \|, | (a, \cdot)| $ получаем, что
\begin{multline*} \tag{1.2.36}
\| A(x_\epsilon -x^0)\| = \| A(2^{-k} \nu +2^{-k} \epsilon +2^{-k} t^0 a -x^0)\|
\le c_5 2^{-k} +r, \\
| (a, x_\epsilon -x^0)| = | (a, 2^{-k} \nu +2^{-k} \epsilon +2^{-k} t^0 a -x^0)|
\le c_6 2^{-k} +\rho +2^{-k} t^0 \
\end{multline*}
и при $ t \ge 0 $ --
\begin{multline*} \tag{1.2.37}
\| A(x_\epsilon^\prime +2^{-k} t a -x^0)\| = \\
\| A(2^{-k} \nu^\prime +2^{-k} \epsilon +2^{-k} t^0 a +2^{-k} t a -x^0)\| \\
\le c_5 2^{-k} +r, 
| (a, x_\epsilon^\prime +2^{-k} t a -x^0)| =\\
 | (a, 2^{-k} \nu^\prime
+2^{-k} \epsilon +2^{-k} t^0 a +2^{-k} t a -x^0)|
\le c_6 2^{-k} +\rho +2^{-k} (t^0 +t).
\end{multline*}

При $ k \in \Z_+, \ \nu,\nu^\prime \in \Z^d, $ удовлетворяющих (1.2.35), для
$ \mathpzc k \in \Z_+, \ t^1 \in \R_+ $ и $ n,n^\prime \in \Z^d: \| n -n^\prime \|
\le 1 $ и существует $ t \in [0,1], $ для которого $ \mathfrak x_t =
x_\epsilon^\prime +2^{-k} t t^1 a \in \overline Q_{\mathpzc k +k,n}^d, $
ввиду (1.2.37) для $ x^\prime \in \overline Q_{\mathpzc k +k,n^\prime}^d $ выполняются неравенства
\begin{multline*} \tag{1.2.38}
\| A(x^\prime -x^0)\| = \| A(x^\prime -2^{-\mathpzc k -k} n^\prime
+2^{-\mathpzc k -k} n^\prime -2^{-\mathpzc k -k} n +2^{-\mathpzc k -k} n
-\mathfrak x_t +\mathfrak x_t -x^0)\| \le \\
\| A(x^\prime -2^{-\mathpzc k -k} n^\prime)\|
+\| A(2^{-\mathpzc k -k} n^\prime -2^{-\mathpzc k -k} n)\|
+\| A(2^{-\mathpzc k -k} n
-\mathfrak x_t)\|
+\| A(\mathfrak x_t -x^0)\| \le \\
c_5 \| x^\prime -2^{-\mathpzc k -k} n^\prime\|
+c_5 \| 2^{-\mathpzc k -k} n^\prime -2^{-\mathpzc k -k} n\|+ \\
c_5 \| 2^{-\mathpzc k -k} n -\mathfrak x_t\|
+\| A(x_\epsilon^\prime +2^{-k} t t^1 a -x^0)\| \le \\
c_5 2^{-\mathpzc k -k} +c_5 2^{-\mathpzc k -k} \| n^\prime -n\|
+c_5 2^{-\mathpzc k -k} +c_5 2^{-k} +r \le \\ c_5 2^{-k} +c_5 2^{-k} +c_5 2^{-k}
+c_5 2^{-k} +r = 4 c_5 2^{-k} +r, \\
| (a, x^\prime -x^0)| = | (a, x^\prime -2^{-\mathpzc k -k} n^\prime
+2^{-\mathpzc k -k} n^\prime -2^{-\mathpzc k -k} n +2^{-\mathpzc k -k} n
-\mathfrak x_t +\mathfrak x_t -x^0)| \le \\
| (a, x^\prime -2^{-\mathpzc k -k} n^\prime)|
+| (a, 2^{-\mathpzc k -k} n^\prime -2^{-\mathpzc k -k} n)|
+| (a, 2^{-\mathpzc k -k} n -\mathfrak x_t)|
+| (a, \mathfrak x_t -x^0)| \le \\
c_6 \| x^\prime -2^{-\mathpzc k -k} n^\prime\|
+c_6 \| 2^{-\mathpzc k -k} n^\prime -2^{-\mathpzc k -k} n\|
+c_6 \| 2^{-\mathpzc k -k} n -\mathfrak x_t\|
+| (a, \mathfrak x_t -x^0)| \le \\
c_6 2^{-\mathpzc k -k} +c_6 2^{-\mathpzc k -k} \| n^\prime -n\|
+c_6 2^{-\mathpzc k -k} +| (a, x_\epsilon^\prime +2^{-k} t t^1 a -x^0)| \le \\
c_6 2^{-k} +c_6 2^{-k} +c_6 2^{-k} +c_6 2^{-k} +\rho +2^{-k} (t^0 +t t^1) \le
4 c_6 2^{-k} +\rho +2^{-k} (t^0 +t^1).
\end{multline*}
Принимая во внимание (1.2.38), получаем, что для $ \mathcal M \in \R_+, $
при $ \delta \in \R_+, $ для которого имеет место (1.2.26), и $ k \in \Z_+ $
таких, что соблюдаются условия
\begin{equation*} \tag{1.2.39}
4 c_5 2^{-k} < r, \\
4 c_6 2^{-k} +2^{-k} t^0 +4 \mathcal M \delta < \rho,
\end{equation*}
для $ \nu, \nu^\prime \in \Z^d $ таких, что выполняется (1.2.35), при $ t^1 =
\mathcal M \| \nu -\nu^\prime \| $ и $ \mathpzc k \in \Z_+ $ для $ n,n^\prime
\in \Z^d: \| n -n^\prime \| \le 1 $ и существует $ t \in [0,1], $ для
которого $ \mathfrak x_t \in \overline Q_{\mathpzc k +k,n}^d, $
для $ x^\prime \in \overline Q_{\mathpzc k +k,n^\prime}^d $ справедливы неравенства
\begin{multline*}
\| A(x^\prime -x^0)\| \le 4 c_5 2^{-k} +r < r +r = 2 r, \\
| (a, x^\prime -x^0)| \le 4 c_6 2^{-k} +\rho +2^{-k} (t^0 +t^1) =
4 c_6 2^{-k} +2^{-k} t^0 +2^{-k} t^1 +\rho = \\
4 c_6 2^{-k} +2^{-k} t^0 +2^{-k} \mathcal M \| \nu -\nu^\prime \| +\rho =
4 c_6 2^{-k} +2^{-k} t^0 +\mathcal M \| 2^{-k} \nu -2^{-k} \nu^\prime \| +\rho \le \\
4 c_6 2^{-k} +2^{-k} t^0 +\mathcal M (\| 2^{-k} \nu -x^0 \| +\| x^0 -2^{-k} \nu^\prime \|) +\rho < \\
4 c_6 2^{-k} +2^{-k} t^0 +\mathcal M ( 2 \delta +2 \delta) +\rho =
4 c_6 2^{-k} +2^{-k} t^0 +4 \mathcal M \delta +\rho < \rho +\rho = 2 \rho,
\end{multline*}
что ввиду (1.2.24), (1.2.25), (1.2.22), (1.2.21) влечёт включение
\begin{equation*} \tag{1.2.40}
h^{-1} (x^\prime) \in V \times W,
\end{equation*}
причём в силу (1.2.23), (1.2.28) (при $ x = 2^{-k} \nu^\prime +2^{-k} \epsilon $),
(1.2.5) справедливо неравенство
\begin{multline*} \tag{1.2.41}
(a, x^\prime -x^0) -\phi(A(x^\prime -x^0)) =
(a, x^\prime -\mathfrak x_t +\mathfrak x_t -x^0) -
\phi( A( x^\prime -\mathfrak x_t +\mathfrak x_t -x^0)) = \\
(a, x^\prime -\mathfrak x_t)
+(a, \mathfrak x_t -x^0) -
\phi( A(x^\prime -\mathfrak x_t)
+A( \mathfrak x_t -x^0)) = \\
(a, x^\prime -\mathfrak x_t)
+(a, 2^{-k} \nu^\prime +2^{-k} \epsilon -x^0) +2^{-k} t^0 (a,a) +2^{-k} t t^1 (a,a) -\\
\phi( A(x^\prime -\mathfrak x_t)
+A( 2^{-k} \nu^\prime +2^{-k} \epsilon -x^0) +2^{-k} t^0 A a
+2^{-k} t t^1 A a ) = \\
(a, x^\prime -\mathfrak x_t)
+(a, 2^{-k} \nu^\prime +2^{-k} \epsilon -x^0) +2^{-k} t^0 +2^{-k} t t^1 -\\
\phi(A( 2^{-k} \nu^\prime +2^{-k} \epsilon -x^0)) +
\phi(A( 2^{-k} \nu^\prime +2^{-k} \epsilon -x^0)) - \\
\phi( A(x^\prime -\mathfrak x_t)
+A( 2^{-k} \nu^\prime +2^{-k} \epsilon -x^0)) = \\
2^{-k} t^0 +2^{-k} t t^1 +(a, x^\prime -\mathfrak x_t)
+(a, 2^{-k} \nu^\prime +2^{-k} \epsilon -x^0) -
\phi(A( 2^{-k} \nu^\prime +2^{-k} \epsilon -x^0)) + \\
\phi(A( 2^{-k} \nu^\prime +2^{-k} \epsilon -x^0)) -
\phi( A(x^\prime -\mathfrak x_t)
+A( 2^{-k} \nu^\prime +2^{-k} \epsilon -x^0)) > \\
2^{-k} t^0 +(a, x^\prime -\mathfrak x_t) +
\phi(A( 2^{-k} \nu^\prime +2^{-k} \epsilon -x^0)) -
\phi( A(x^\prime -\mathfrak x_t)
+A( 2^{-k} \nu^\prime +2^{-k} \epsilon -x^0)) \ge \\
2^{-k} t^0 -| (a, x^\prime -\mathfrak x_t)| -
| \phi(A( 2^{-k} \nu^\prime +2^{-k} \epsilon -x^0)) -
\phi( A(x^\prime -\mathfrak x_t)
+A( 2^{-k} \nu^\prime +2^{-k} \epsilon -x^0))| \ge \\
2^{-k} t^0 -| (a, x^\prime -\mathfrak x_t)| -
M \| A(x^\prime -\mathfrak x_t) \| \ge \\
2^{-k} t^0 -c_6 \| x^\prime -2^{-\mathpzc k -k} n^\prime
+2^{-\mathpzc k -k} n^\prime -2^{-\mathpzc k -k} n
+2^{-\mathpzc k -k} n -\mathfrak x_t\| -\\
M c_5 \| x^\prime -2^{-\mathpzc k -k} n^\prime
+2^{-\mathpzc k -k} n^\prime -2^{-\mathpzc k -k} n
+2^{-\mathpzc k -k} n -\mathfrak x_t \| \ge \\
2^{-k} t^0 -c_6 (\| x^\prime -2^{-\mathpzc k -k} n^\prime\|
+\| 2^{-\mathpzc k -k} n^\prime -2^{-\mathpzc k -k} n\|
+\| 2^{-\mathpzc k -k} n -\mathfrak x_t\|) - \\
M c_5 (\| x^\prime -2^{-\mathpzc k -k} n^\prime\|
+\| 2^{-\mathpzc k -k} n^\prime -2^{-\mathpzc k -k} n\|
+\| 2^{-\mathpzc k -k} n -\mathfrak x_t \|) \ge \\
2^{-k} t^0 -c_6 (2^{-\mathpzc k -k}
+2^{-\mathpzc k -k} \| n^\prime -n\|
+2^{-\mathpzc k -k})-\\
M c_5 (2^{-\mathpzc k -k}
+2^{-\mathpzc k -k} \| n^\prime -n\|
+2^{-\mathpzc k -k}) \ge 
2^{-k} t^0 -3 c_6 2^{-\mathpzc k -k}
-M 3 c_5 2^{-\mathpzc k -k} = \\
2^{-k} (t^0 -3 c_6 2^{-\mathpzc k} -3 M c_5 2^{-\mathpzc k}) > 0
\end{multline*}
при условии, что
\begin{equation*} \tag{1.2.42}
3 c_6 2^{-\mathpzc k} +3 M c_5 2^{-\mathpzc k} < t^0.
\end{equation*}

На основании (1.2.40), (1.2.41), (1.2.22), (1.2.21), (1.2.4), заключаем, что для
$ \mathcal M \in \R_+, $ при $ \delta \in \R_+, k \in \Z_+, \mathpzc k \in \Z_+, $
для которых соблюдаются (1.2.26), (1.2.39), (1.2.42), для
$ \nu, \nu^\prime \in \Z^d $ таких, что имеет место (1.2.35), при
$ t^1 = \mathcal M \| \nu -\nu^\prime \| $

для $ n,n^\prime \in \Z^d: \| n -n^\prime \| \le 1 $ и существует
$ t \in [0,1], $ для которого $ \mathfrak x_t \in
\overline Q_{\mathpzc k +k,n}^d, $
справедливо включение
\begin{equation*} \tag{1.2.43}
\overline Q_{\mathpzc k +k,n^\prime}^d \subset D.
\end{equation*}

Далее, при $ k \in \Z_+ $ для $ \nu, \nu^\prime \in \Z^d, $ удовлетворяющих
(1.2.35), $ t^1 \in \R_+, t \in [0,1] $ для $ \mathfrak X_t = (1 -t) x_\epsilon +
t (x_\epsilon^\prime +2^{-k} t^1 a) $ вследствие (1.2.36), (1.2.37)
выполняется неравенство
\begin{multline*} \tag{1.2.44}
\| A(\mathfrak X_t -x^0)\| = \| A((1 -t) x_\epsilon +
t (x_\epsilon^\prime +2^{-k} t^1 a) -x^0) \| = \\
\| A((1 -t)(x_\epsilon -x^0) +t (x_\epsilon^\prime +2^{-k} t^1 a -x^0)) \| = \\
\| (1 -t) A(x_\epsilon -x^0) +t A(x_\epsilon^\prime +2^{-k} t^1 a -x^0)\| \le \\
(1 -t) \| A(x_\epsilon -x^0)\| +
t \| A(x_\epsilon^\prime +2^{-k} t^1 a -x^0)\| \le \\
(1 -t) (c_5 2^{-k} +r) +t (c_5 2^{-k} +r) = c_5 2^{-k} +r, \\
| (a, \mathfrak X_t -x^0)| = | (a, (1 -t) x_\epsilon +
t (x_\epsilon^\prime +2^{-k} t^1 a) -x^0) | = \\
| (a, (1 -t)(x_\epsilon -x^0) +t (x_\epsilon^\prime +2^{-k} t^1 a -x^0)) | = \\
| (1 -t) (a, x_\epsilon -x^0) +t (a, x_\epsilon^\prime +2^{-k} t^1 a -x^0)| \le \\
(1 -t) | (a, x_\epsilon -x^0)| +
t | (a, x_\epsilon^\prime +2^{-k} t^1 a -x^0)| \le \\
(1 -t) (c_6 2^{-k} +\rho +2^{-k} t^0) +t (c_6 2^{-k} +\rho +2^{-k} (t^0 +t^1)) = \\
c_6 2^{-k} +\rho +2^{-k} t^0 +t 2^{-k} t^1 \le
c_6 2^{-k} +2^{-k} (t^0 +t^1) +\rho.
\end{multline*}

При $ k \in \Z_+, \ \nu,\nu^\prime \in \Z^d, $ удовлетворяющих (1.2.35), для
$ \mathpzc k \in \Z_+, \ t^1 \in \R_+ $ для $ n,n^\prime \in \Z^d:
\| n -n^\prime \| \le 1 $ и существует $ t \in [0,1], $ для которого
$ \mathfrak X_t \in \overline Q_{\mathpzc k +k,n}^d, $
для $ x^\prime \in \overline Q_{\mathpzc k +k,n^\prime}^d $ с учётом (1.2.44)
соблюдаются неравенства
\begin{multline*} \tag{1.2.45}
\| A(x^\prime -x^0)\| = \| A(x^\prime -2^{-\mathpzc k -k} n^\prime
+2^{-\mathpzc k -k} n^\prime -2^{-\mathpzc k -k} n +2^{-\mathpzc k -k} n
-\mathfrak X_t +\mathfrak X_t -x^0)\| \le \\
\| A(x^\prime -2^{-\mathpzc k -k} n^\prime)\|
+\| A(2^{-\mathpzc k -k} n^\prime -2^{-\mathpzc k -k} n)\|
+\| A(2^{-\mathpzc k -k} n -\mathfrak X_t)\|
+\| A(\mathfrak X_t -x^0)\| \le \\
c_5 \| x^\prime -2^{-\mathpzc k -k} n^\prime\|
+c_5 \| 2^{-\mathpzc k -k} n^\prime -2^{-\mathpzc k -k} n\|
+c_5 \| 2^{-\mathpzc k -k} n -\mathfrak X_t\|
+\| A(\mathfrak X_t -x^0)\| \le \\
c_5 2^{-\mathpzc k -k} +c_5 2^{-\mathpzc k -k} \| n^\prime -n\|
+c_5 2^{-\mathpzc k -k} +c_5 2^{-k} +r \le \\
c_5 2^{-k} +c_5 2^{-k} +c_5 2^{-k}
+c_5 2^{-k} +r = 4 c_5 2^{-k} +r, \\
| (a, x^\prime -x^0)| = | (a, x^\prime -2^{-\mathpzc k -k} n^\prime
+2^{-\mathpzc k -k} n^\prime -2^{-\mathpzc k -k} n +2^{-\mathpzc k -k} n
-\mathfrak X_t +\mathfrak X_t -x^0)| \le \\
| (a, x^\prime -2^{-\mathpzc k -k} n^\prime)|
+| (a, 2^{-\mathpzc k -k} n^\prime -2^{-\mathpzc k -k} n)|
+| (a, 2^{-\mathpzc k -k} n -\mathfrak X_t)|
+| (a, \mathfrak X_t -x^0)| \le \\
c_6 \| x^\prime -2^{-\mathpzc k -k} n^\prime\|
+c_6 \| 2^{-\mathpzc k -k} n^\prime -2^{-\mathpzc k -k} n\|
+c_6 \| 2^{-\mathpzc k -k} n -\mathfrak X_t\|
+| (a, \mathfrak X_t -x^0)| \le \\
c_6 2^{-\mathpzc k -k} +c_6 2^{-\mathpzc k -k} \| n^\prime -n\|
+c_6 2^{-\mathpzc k -k} +| (a, \mathfrak X_t -x^0)| \le \\
c_6 2^{-k} +c_6 2^{-k} +c_6 2^{-k} +c_6 2^{-k} +\rho +2^{-k} (t^0 +t^1) =
4 c_6 2^{-k} +2^{-k} (t^0 +t^1) +\rho.
\end{multline*}

Учитывая (1.2.45), приходим к выводу, что для $ \mathcal M \in \R_+, $ при
$ \delta \in \R_+, $ для которого имеет место (1.2.26), и $ k \in \Z_+ $ таких, что
выполняются условия (1.2.39), для $ \nu, \nu^\prime \in \Z^d $ таких, что
соблюдается (1.2.35), при $ t^1 = \mathcal M \| \nu -\nu^\prime \| $ и
$ \mathpzc k \in \Z_+ $ для $ n,n^\prime \in \Z^d: \| n -n^\prime \| \le 1 $ и
существует $ t \in [0,1], $ для которого $ \mathfrak X_t \in
\overline Q_{\mathpzc k +k,n}^d, $ для $ x^\prime \in
\overline Q_{\mathpzc k +k,n^\prime}^d $ справедливы неравенства
\begin{multline*}
\| A(x^\prime -x^0)\| \le 4 c_5 2^{-k} +r < r +r = 2 r, \\
| (a, x^\prime -x^0)| \le 4 c_6 2^{-k} +2^{-k} (t^0 +t^1) +\rho = \\
4 c_6 2^{-k} +2^{-k} t^0 +2^{-k} t^1 +\rho =
4 c_6 2^{-k} +2^{-k} t^0 +2^{-k} \mathcal M \| \nu -\nu^\prime \| +\rho = \\
4 c_6 2^{-k} +2^{-k} t^0 +\mathcal M \| 2^{-k} \nu -2^{-k} \nu^\prime \| +\rho \le \\
4 c_6 2^{-k} +2^{-k} t^0 +\mathcal M (\| 2^{-k} \nu -x^0 \| +\| x^0 -2^{-k} \nu^\prime \|) +\rho < \\
4 c_6 2^{-k} +2^{-k} t^0 +\mathcal M ( 2 \delta +2 \delta) +\rho = \\
4 c_6 2^{-k} +2^{-k} t^0 +4 \mathcal M \delta +\rho < \rho +\rho = 2 \rho,
\end{multline*}
что благодаря (1.2.24), (1.2.25), (1.2.22), (1.2.21) приводит к
соотношению (1.2.40), причём в силу (1.2.23), (1.2.28)
(при $ x = 2^{-k} \nu +2^{-k} \epsilon $), (1.2.5) имеет место неравенство
\begin{multline*} \tag{1.2.46}
(a, x^\prime -x^0) -\phi( A(x^\prime -x^0)) =
(a, x^\prime -\mathfrak X_t +\mathfrak X_t -x^0) -
\phi( A( x^\prime -\mathfrak X_t +\mathfrak X_t -x^0)) = \\
(a, x^\prime -\mathfrak X_t)
+(a, \mathfrak X_t -x^0) -
\phi( A(x^\prime -\mathfrak X_t)
+A( \mathfrak X_t -x^0)) = \\
(a, x^\prime -\mathfrak X_t)
+(a, (1 -t) x_\epsilon +t (x_\epsilon^\prime +2^{-k} t^1 a) -x^0) -\\
\phi( A(x^\prime -\mathfrak X_t)
+A( (1 -t) x_\epsilon +t (x_\epsilon^\prime +2^{-k} t^1 a) -x^0)) = \\
(a, x^\prime -\mathfrak X_t)
+(a, x_\epsilon +t (x_\epsilon^\prime -x_\epsilon +2^{-k} t^1 a) -x^0) -\\
\phi( A(x^\prime -\mathfrak X_t)
+A( x_\epsilon +t (x_\epsilon^\prime -x_\epsilon +2^{-k} t^1 a) -x^0)) = \\
(a, x^\prime -\mathfrak X_t)
+(a, x_\epsilon -x^0) +t (a, x_\epsilon^\prime -x_\epsilon) +2^{-k} t t^1(a, a) -\\
\phi( A(x^\prime -\mathfrak X_t)
+A( x_\epsilon -x^0) +t A(x_\epsilon^\prime -x_\epsilon)
+2^{-k} t t^1 A a) = \\
(a, x^\prime -\mathfrak X_t)
+(a, 2^{-k} \nu +2^{-k} \epsilon -x^0) +2^{-k} t^0 (a,a) +
t (a, 2^{-k} \nu^\prime -2^{-k} \nu) +2^{-k} t t^1(a, a) - \\
\phi( A(x^\prime -\mathfrak X_t)
+A( 2^{-k} \nu +2^{-k} \epsilon -x^0) +2^{-k} t^0 A a +
t A(2^{-k} \nu^\prime -2^{-k} \nu)
+2^{-k} t t^1 A a) = \\
(a, x^\prime -\mathfrak X_t)
+(a, 2^{-k} \nu +2^{-k} \epsilon -x^0) +2^{-k} t^0 +
t (a, 2^{-k} \nu^\prime -2^{-k} \nu) +2^{-k} t t^1 - \\
\phi(A( 2^{-k} \nu +2^{-k} \epsilon -x^0)) +
\phi(A( 2^{-k} \nu +2^{-k} \epsilon -x^0)) -\\
\phi( A(x^\prime -\mathfrak X_t)
+A( 2^{-k} \nu +2^{-k} \epsilon -x^0) +
t A(2^{-k} \nu^\prime -2^{-k} \nu)) > \\
2^{-k} t^0 -| (a, x^\prime -\mathfrak X_t)|
-t | (a, 2^{-k} \nu^\prime -2^{-k} \nu)| +2^{-k} t t^1
-| \phi(A( 2^{-k} \nu +2^{-k} \epsilon -x^0)) -\\
\phi( A(x^\prime -\mathfrak X_t)
+A( 2^{-k} \nu +2^{-k} \epsilon -x^0) +
t A(2^{-k} \nu^\prime -2^{-k} \nu))| \ge \\
2^{-k} t^0 -| (a, x^\prime -\mathfrak X_t)|
-t | (a, 2^{-k} \nu^\prime -2^{-k} \nu)| +\\
2^{-k} t t^1-M \| A(x^\prime -\mathfrak X_t)
+t A(2^{-k} \nu^\prime -2^{-k} \nu)\| \ge \\
2^{-k} t^0 -| (a, x^\prime -2^{-\mathpzc k -k} n^\prime)|
-| (a, 2^{-\mathpzc k -k} n^\prime -2^{-\mathpzc k -k} n)|
-| (a, 2^{-\mathpzc k -k} n -\mathfrak X_t)| \\
-t | (a, 2^{-k} \nu^\prime -2^{-k} \nu)| +2^{-k} t t^1
-M (\| A(x^\prime -2^{-\mathpzc k -k} n^\prime)\|
+\| A( 2^{-\mathpzc k -k} n^\prime -2^{-\mathpzc k -k} n)\| \\
+\| A( 2^{-\mathpzc k -k} n -\mathfrak X_t)\|
+t \| A(2^{-k} \nu^\prime -2^{-k} \nu)\|) \ge \\
2^{-k} t^0 -c_6 \| x^\prime -2^{-\mathpzc k -k} n^\prime\|
-c_6 \| 2^{-\mathpzc k -k} n^\prime -2^{-\mathpzc k -k} n\|
-c_6 \| 2^{-\mathpzc k -k} n -\mathfrak X_t\| \\
-t c_6 \| 2^{-k} \nu^\prime -2^{-k} \nu\| +2^{-k} t t^1
-M (c_5 \| x^\prime -2^{-\mathpzc k -k} n^\prime\|
+c_5 \| 2^{-\mathpzc k -k} n^\prime -2^{-\mathpzc k -k} n\| \\
+c_5 \| 2^{-\mathpzc k -k} n -\mathfrak X_t\|
+t c_5 \| 2^{-k} \nu^\prime -2^{-k} \nu\|) \ge \\
2^{-k} t^0 -c_6 2^{-\mathpzc k -k}
-c_6 2^{-\mathpzc k -k} \| n^\prime -n\|
-c_6 2^{-\mathpzc k -k}
-t c_6 2^{-k} \| \nu^\prime -\nu\| +2^{-k} t t^1 \\
-M (c_5 2^{-\mathpzc k -k}
+c_5 2^{-\mathpzc k -k} \| n^\prime -n\|
+c_5 2^{-\mathpzc k -k}
+t c_5 2^{-k} \| \nu^\prime -\nu\|) \ge \\
2^{-k} t^0 -c_6 2^{-\mathpzc k -k}
-c_6 2^{-\mathpzc k -k}
-c_6 2^{-\mathpzc k -k}
-t c_6 2^{-k} \| \nu^\prime -\nu\| +2^{-k} t t^1 \\
-M (c_5 2^{-\mathpzc k -k}
+c_5 2^{-\mathpzc k -k}
+c_5 2^{-\mathpzc k -k}
+t c_5 2^{-k} \| \nu^\prime -\nu\|) = \\
2^{-k} t^0 -3 c_6 2^{-\mathpzc k -k} -t c_6 2^{-k} \| \nu^\prime
-\nu\| +t 2^{-k} \mathcal M \| \nu^\prime -\nu\|
-M (3 c_5 2^{-\mathpzc k -k} +t c_5 2^{-k} \| \nu^\prime -\nu\|) = \\
2^{-k} t^0 -3 c_6 2^{-\mathpzc k -k} -t c_6 2^{-k} \| \nu^\prime
-\nu\| +t 2^{-k} \mathcal M \| \nu^\prime -\nu\|
-3 M c_5 2^{-\mathpzc k -k} -t M c_5 2^{-k} \| \nu^\prime -\nu\| = \\
2^{-k} (t^0 -3 c_6 2^{-\mathpzc k} -3 M c_5 2^{-\mathpzc k}) +t
2^{-k} \| \nu^\prime -\nu\| (\mathcal M -c_6 -M c_5) > 0
\end{multline*}
при соблюдении условий (1.2.42) и
\begin{equation*} \tag{1.2.47}
\mathcal M -c_6 -M c_5 \ge 0.
\end{equation*}

Опираясь на (1.2.40), (1.2.46), (1.2.22), (1.2.21), (1.2.4), заключаем,
что для $ \mathcal M \in \R_+, $ удовлетворяющего (1.2.47), при
$ \delta \in \R_+, k \in \Z_+, \mathpzc k \in \Z_+, $ для которых соблюдаются
(1.2.26), (1.2.39), (1.2.42), для $ \nu, \nu^\prime \in \Z^d $ таких, что имеет
место (1.2.35), при $ t^1 = \mathcal M \| \nu -\nu^\prime \| $ для $ n,n^\prime
\in \Z^d: \| n -n^\prime \| \le 1 $ и существует $ t \in [0,1], $ для которого
$ \mathfrak X_t \in \overline Q_{\mathpzc k +k,n}^d, $ верно включение
(1.2.43).

Фиксируя $ \mathcal M \in \R_+, $ для которого имеет место (1.2.47),
при $ \delta \in \R_+, k, \mathpzc k \in \Z_+, $ подчинённых условиям (1.2.26),
(1.2.39), (1.2.42), для $ \nu, \nu^\prime \in \Z^d, $ таких, что выполняется (1.2.35),
построим последовательности $ \{\nu^\iota, \ \iota = 0,\ldots,\Iota\}, \
\{j^\iota, \ \iota =0,\ldots,\Iota -1\}, \ \{ \epsilon^\iota, \ \iota =0,\ldots,\Iota -1\}, $
удовлетворяющие (1.2.18), (1.2.19), (1.2.20). Для этого, пользуясь леммой 1.2.2, построим
последовательности
\begin{multline*} \tag{1.2.48}
\{\nu_s^{\iota_s} \in \Z^d, \ \iota_s =0,\ldots,\Iota_s\}, \\
\{j_s^{\iota_s} \in \Nu_{1,d}^1, \ \iota_s =0,\ldots,\Iota_s -1\}, \\
\{\epsilon_s^{\iota_s} \in \{-1,1\}, \ \iota_s =0,\ldots,\Iota_s -1\}, \ s =1,3,
\end{multline*}
для которых соблюдаются условия пп. 1) -- 4) леммы 1.2.2 при $ \mathpzc k +k $
вместо $ k $ и
\begin{multline*} \tag{1.2.49}
\text{ при } x^0 = \mathfrak X_0 = x_\epsilon, \\
\xi = \mathfrak X_1 -\mathfrak X_0 =
(x_\epsilon^\prime +2^{-k} t^1 a) -x_\epsilon, \ s =1;
\end{multline*}
а также
\begin{multline*} \tag{1.2.50}
\text{ при } x^0 = \mathfrak x_1 = x_\epsilon^\prime +2^{-k} t^1 a =
\mathfrak X_1, \
\xi = \mathfrak x_0 -\mathfrak x_1 = \\
x_\epsilon^\prime -
(x_\epsilon^\prime +2^{-k} t^1 a) = -2^{-k} t^1 a, \ s =3,
\end{multline*}
где $ t^1 = \mathcal M \| \nu -\nu^\prime \| $ и соблюдается (1.2.47).
Учитывая, что
\begin{equation*}
\mathfrak X_0 \in \overline Q_{k,\nu}^d = \cup_{\mu \in \Nu_{0, 2^{\mathpzc k} -1}^d}
\overline Q_{\mathpzc k +k, 2^{\mathpzc k} \nu +\mu}^d,
\end{equation*}
выберем $ n \in \Z^d $ такое, что
\begin{equation*} \tag{1.2.51}
\mathfrak X_0 \in \overline Q_{\mathpzc k +k, n}^d \subset \overline Q_{k,\nu}^d.
\end{equation*}
Точно так же возьмём $ n^\prime \in \Z^d, $ для которого
\begin{equation*} \tag{1.2.52}
\mathfrak x_1 +(\mathfrak x_0 -\mathfrak x_1) = \mathfrak x_0 \in
\overline Q_{\mathpzc k +k, n^\prime}^d \subset \overline Q_{k,\nu^\prime}^d.
\end{equation*}

Опираясь на замечание 2 (при $ \mathpzc k +k $ вместо $ k $ ) и соотношения
(1.2.14), (1.2.15), с учётом (1.2.51), (1.2.49), (1.2.50), (1.2.52), а также п. 1) леммы 1.2.2 построим
последовательности (1.2.48) при $ s =0,2,4, $ для которых $ \nu_0^0 = n, \
\nu_4^{\Iota_4} = n^\prime, \ \nu_s^{\Iota_s} = \nu_{s +1}^0, \ s =0,1,2,3; \
\| \nu_s^{\iota_s} -\nu_s^0 \| \le 1 \text{ при } \iota_s =1,\ldots,\Iota_s, \
s =0,2,4; \ \Iota_s $ подчинено (1.2.13) при $ s =0,2,4, $ и выполняются
равенства п. 3) леммы 1.2.2.
Исходя из всего сказанного, положим $ \nu^0 = \nu_0^0 $ и определим при
$ \iota = \Iota_0 + \ldots +\Iota_{s -1} +\iota_s, \ \iota_s =1,\ldots,\Iota_s, \
s =0,1,2,3,4, $ значения
\begin{equation*}
\nu^\iota = \nu_s^{\iota_s}, \\
j^{\iota -1} = j_s^{\iota_s -1}, \\
\epsilon^{\iota -1} = \epsilon_s^{\iota_s -1}.
\end{equation*}

Тогда в силу построений, а также п. 4) леммы 1.2.2, (1.2.49), (1.2.50) и
(1.2.13) для $ \Iota = \sum_{s =0}^4 \Iota_s $ справедливо неравенство
\begin{multline*}
\Iota \le d +c_1 (2^{\mathpzc k +k} \| x_\epsilon^\prime +2^{-k} t^1 a -
x_\epsilon \| +1) +d +c_1 (2^{\mathpzc k +k} \| -2^{-k} t^1 a\| +1) +d = \\
3d +c_1 (2^{\mathpzc k +k} \| 2^{-k} \nu^\prime +2^{-k} \epsilon
+2^{-k} t^0 a +2^{-k} t^1 a -
(2^{-k} \nu +2^{-k} \epsilon +2^{-k} t^0 a) \| +1) +\\
c_1 (2^{\mathpzc k +k} \| 2^{-k} t^1 a\| +1) = \\
3d +c_1 (2^{\mathpzc k +k} \| 2^{-k} \nu^\prime -2^{-k} \nu
+2^{-k} \mathcal M \| \nu -\nu^\prime \| a \| +1) +\\
c_1 (2^{\mathpzc k +k}
\| 2^{-k} \mathcal M \| \nu -\nu^\prime \| a\| +1) \le \\
3d +c_1 (2^{\mathpzc k +k} (2^{-k} \| \nu^\prime -\nu\|
+2^{-k} \mathcal M \| \nu -\nu^\prime \| \cdot \| a \|) +1) +\\
c_1 (2^{\mathpzc k +k} 2^{-k} \mathcal M \| \nu -\nu^\prime \| \cdot \| a \| +1) \le \\
3d +c_1 (2^{\mathpzc k} ( \| \nu^\prime -\nu\|
+\mathcal M \| \nu -\nu^\prime \| ) +1) +\\
c_1 (2^{\mathpzc k} \mathcal M \| \nu -\nu^\prime \| +1) = \\
3d +c_1 ((2 \mathcal M +1) 2^{\mathpzc k} \| \nu -\nu^\prime \| +2) \le
c_4 2^{\mathpzc k} \| \nu -\nu^\prime\|
\text{ при } \nu \ne \nu^\prime,
\end{multline*}
что совпадает с (1.2.18). Ввиду построений и пп. 1), 3) леммы 1.2.2
выполняется (1.2.19), а согласно построениям, п. 2) леммы 1.2.2 и
включению (1.2.43) имеет место (1.2.20).

Итак, фиксируя $ t, $ удовлетворяющее (1.2.34), а $ \mathcal M, $ подчинённое
(1.2.47), и подбирая $ \delta^0 \in \R_+, K^1 \in \Z_+, \mathpzc k^1 \in \Z_+, $
для которых соблюдаются (1.2.26), (1.2.31), (1.2.39), (1.2.42) при
$ \delta^0 $ вместо $ \delta, \ K^1 $ вместо $ k, \ \mathpzc k^1 $ вместо
$ \mathpzc k, $ получаем, что при $ k \ge K^1, \ \mathpzc k \ge \mathpzc k^1 $
справедливы пп. 1) и 2), что завершает доказательство леммы для $ x^0 \in \D D. $

Если $ x^0 \in D, $ выберем $ \delta^0 > 0 $ так, чтобы $ U_{x^0} =
(x^0 +2 \delta^0 B_0^d) \subset D. $ Тогда при $ k \in \Z_+ $ для $ \nu \in
\Z^d: Q_{k,\nu}^d \cap U_{x^0} \ne \emptyset, $ при $ j =1,\ldots,d $
выполняются неравенства
\begin{equation*} \tag{1.2.53}
x_j^0 -2 \delta^0 < 2^{-k} \nu_j +2^{-k}, \\
2^{-k} \nu_j < x_j^0 +2 \delta^0.
\end{equation*}
При $ k \in \Z_+: 3 2^{-k} < 4 \delta^0, $ и $ \nu \in \Z^d: Q_{k,\nu}^d \cap
U_{x^0} \ne \emptyset, $ определим координаты $ \nu^\prime \in \Z^d, $
полагая при $ j =1,\ldots,d $ значение
\begin{equation*} \tag{1.2.54}
\nu_j^\prime = \begin{cases} \nu_j -1, \text{ если } 2^{-k} (\nu_j -1) > x_j^0 -2 \delta^0; \\
\nu_j +1, \text{ если } 2^{-k} (\nu_j +2) < x_j^0 +2 \delta^0,
\end{cases}
\end{equation*}
при этом для каждого $ j =1,\ldots,d $ ситуация, когда $ 2^{-k} (\nu_j -1) \le
x_j^0 -2 \delta^0, \ 2^{-k} (\nu_j +2) \ge x_j^0 +2 \delta^0, $ невозможна,
ибо в таком случае
$$
2^{-k} (\nu_j +2) -2^{-k} (\nu_j -1) \ge x_j^0 +2 \delta^0 -(x_j^0 -2 \delta^0)
$$
или $ 3 2^{-k} \ge 4 \delta^0, $ что не верно.
Из (1.2.53), (1.2.54) следует (1.2.17) при $ \gamma^1 =2. $

Наконец, при $ k \in \Z_+ $ для $ \nu, \nu^\prime \in \Z^d $ таких, что
\begin{equation*}
\overline Q_{k,\nu}^d \subset U_{x^0}, \\
\overline Q_{k,\nu^\prime}^d \subset U_{x^0},
\end{equation*}
при $ \mathpzc k \in \Z_+ $ положим
\begin{multline*}
\Iota = \sum_{j =1}^d | 2^{\mathpzc k} \nu_j^\prime -2^{\mathpzc k} \nu_j |, \\
\nu^{\Iota} = 2^{\mathpzc k} \nu^\prime
\end{multline*}
и для $ \iota =0,\ldots,\Iota -1 $ такого, что
\begin{multline*}
\iota = \sum_{j =1,\ldots,d: j < j^\iota} 2^{\mathpzc k} | \nu_j^\prime -
\nu_j| +s^\iota, \\
s^\iota =0,\ldots, 2^{\mathpzc k} | \nu_{j^\iota}^\prime -\nu_{j^\iota}| -1, \\
j^\iota \in \Nu_{1,d}^1,
\end{multline*}
определим $ \nu^\iota, \ \epsilon^\iota $ равенствами
\begin{multline*} \tag{1.2.55}
\nu_j^\iota = \begin{cases} 2^{\mathpzc k} \nu_j^\prime \text{ при } 1 \le j < j^\iota; \\
2^{\mathpzc k} \nu_{j^\iota} +\sgn (\nu_{j^\iota}^\prime -\nu_{j^\iota}) s^\iota
\text{ при } j = j^\iota; \\
2^{\mathpzc k} \nu_j \text{ при } j^\iota < j \le d, \\
\end{cases}
\epsilon^\iota = \sgn (\nu_{j^\iota}^\prime -\nu_{j^\iota}).
\end{multline*}
Опираясь на (1.2.55), несложно проверить, что выполняются соотношения
(1.2.19) и (1.2.20). А
$$
\Iota = \sum_{j =1}^d 2^{\mathpzc k} | \nu_j^\prime -\nu_j | = 2^{\mathpzc k}
\sum_{j =1}^d | \nu_j^\prime -\nu_j| \le d 2^{\mathpzc k} \| \nu^\prime -\nu\|,
$$
т.е. имеет место (1.2.18). $ \square $

Лемма 1.2.4

Пусть $ D $ -- область в $ \R^d $ и $ \{U_i \subset D, \ i \in \mathcal I\} $ --
конечный набор непустых открытых множеств таких, что
$ D = \cup_{i \in \mathcal I} U_i. $ Тогда для любых точек $ x^0, x^1 \in D $
существует последовательность различных между собой элементов
$ \{i_s \in \mathcal I, \ s =0,\ldots,S\} $ такая, что $ x^0 \in U_{i_0}, \
x^1 \in U_{i_S} $ и при $ s =0,\ldots, S -1 $ пересечение
$ U_{i_s} \cap U_{i_{s +1}} \ne \emptyset. $

Доказательство.

Исходя из условий леммы, для $ x^0, x^1 \in D $ выберем $ i_0 \in \mathcal I, $
для которого $ x^0 \in U_{i_0}, $ и положим $ \mathcal I_0 = \{i_0\}. $
Далее, при $ s \in \Z_+, $ предполагая, что уже построены множества
$ \mathcal I_0,\ldots, \mathcal I_s $ и множество
$ \mathcal I \setminus (\cup_{s^\prime =0}^s \mathcal I_{s^\prime}) \ne \emptyset, $
определим множество $ \mathcal I_{s +1} $ равенством
\begin{equation*} \tag{1.2.56}
\mathcal I_{s +1} = \{ i \in \mathcal I \setminus (\cup_{s^\prime =0}^s
\mathcal I_{s^\prime}): U_i \cap (\cup_{j \in (\cup_{s^\prime =0}^s
\mathcal I_{s^\prime})} U_j) \ne \emptyset.
\end{equation*}
При этом множество $ \mathcal I_{s +1} \ne \emptyset, $ ибо если это не так, то
$ D $ представляется в виде объединения непустых непересекающихся открытых
множеств
$$
(\cup_{ i \in \mathcal I \setminus (\cup_{s^\prime =0}^s
\mathcal I_{s^\prime})} U_i)
$$
и
$$
(\cup_{j \in (\cup_{s^\prime =0}^s \mathcal I_{s^\prime})} U_j),
$$
что противоречит связности $ D. $
Кроме того, по построению $ \mathcal I_{s +1} \cap \mathcal I_{s^\prime} =
\emptyset, \ s^\prime =0,\ldots. $ Учитывая последнее обстоятельство,
заключаем, что для $ i \in \mathcal I_{s +1} $ верно соотношение
\begin{equation*} \tag{1.2.57}
U_i \cap (\cup_{j \in (\cup_{s^\prime =0}^{s -1} \mathcal I_{s^\prime})} U_j)
= \emptyset,
\end{equation*}
поскольку в противном случае ввиду (1.2.56) также $ i \in \mathcal I_s, $
что невозможно.
И ещё
\begin{equation*}
s +2 \le \sum_{s^\prime =0}^{s +1} \card \mathcal I_{s^\prime} =
\card(\cup_{s^\prime =0}^{s +1} \mathcal I_{s^\prime}) \le \card \mathcal I.
\end{equation*}

Продолжая построение таких множеств, в силу отмеченных обстоятельств через
конечное число шагов прийдём к тому, что при некотором $ S^0 \in \Z_+ $ будет
$ \mathcal I_s \ne \emptyset, \ s =0,\ldots,S^0, $ а
$ \mathcal I \setminus \cup_{s =0}^{S^0} \mathcal I_s = \emptyset $ или
$ \mathcal I = \cup_{s =0}^{S^0} \mathcal I_s. $
Отсюда получаем, что $ D = \cup_{s =0}^{S^0} (\cup_{i \in \mathcal I_s} U_i). $
Поэтому существует $ S \in \Nu_{0,S^0}^1, $ для которого
$ x^1 \in (\cup_{i \in \mathcal I_S} U_i), $ а, следовательно, существует
$ i_S \in \mathcal I_S $ такое, что $ x^1 \in U_{i_S}. $ Из (1.2.56), (1.2.57)
вытекает, что $ U_{i_S} \cap (\cup_{j \in \mathcal I_{S -1}} U_j) \ne
\emptyset, $ и, значит, существует $ i_{S -1} \in \mathcal I_{S -1}, $ для
которого $ U_{i_S} \cap U_{i_{S -1}} \ne \emptyset. $ Продолжая действовать
таким образом, через $ S $ шагов получим последовательность элементов
$ \{i_s \in \mathcal I_s, \ s =0,\ldots,S\}, $ обладающую требуемыми
свойствами. $ \square $

Доказательство теоремы 1.2.1.

В условиях теоремы для каждой точки $ x^0 \in \overline D $ построим её
окрестность $ x^0 +\frac{1}{2} \delta^0(d,D,x^0) B_0^d $ с константой
$ \delta^0(d,D,x^0) $ из леммы 1.2.3. В силу ограниченности и замкнутости
(а, следовательно, компактности) $ \overline D $ выберем конечное число
точек $ x_i \in \overline D, \ i =1,\ldots,N, $ для которых $ \overline D
\subset (\cup_{i =1}^N (x_i +\frac{1}{2} \delta^0(d,D,x_i) B_0^d)). $ Отсюда,
обозначая через $ \mathcal U_i = D \cap (x_i +\frac{1}{2} \delta^0(d,D,x_i)
B_0^d) \ne \emptyset, \ i =1,\ldots,N, $ имеем
\begin{equation*} \tag{1.2.58}
D = \cup_{i =1}^N \mathcal U_i.
\end{equation*}

Тогда при
\begin{equation*} \tag{1.2.59}
k \in \Z_+: k \ge \max_{i =1,\ldots,N} K^1(d,D,x_i),
\end{equation*}
для $ \nu \in \Z^d: Q_{k,\nu}^d \cap D \ne \emptyset, $ ввиду (1.2.58)
существует $ i =1,\ldots,N, $ для которого $ Q_{k,\nu}^d \cap \mathcal U_i
\ne \emptyset, $ а, следовательно,
$ Q_{k,\nu}^d \cap (D \cap U_{x_i}) \ne \emptyset $ (см. лемму 1.2.3). Отсюда на
основании п. 1) леммы 1.2.3 (при $ x^0 = x_i $ ) заключаем,
что существует $ \nu^\prime \in \Z^d $ такое, что
$ \overline Q_{k,\nu^\prime}^d \subset D \cap (2^{-k} \nu +2^{-k} \Gamma^0 B^d) $
при $ \Gamma^0 = \max_{i =1,\ldots,N} \gamma^1(d,D,x_i), $ т.е. имеет место
п. 1) определения 1.

Теперь проверим соблюдение п. 2) определения 1. Для этого сначала
заметим, что вследствие открытости множеств $ \mathcal U_i, \ i =1,\ldots,N $
существует константа $ K^2(d,D) \in \Z_+ $ такая, что для любых
$ i,j =1,\ldots,N: \mathcal U_i \cap \mathcal U_j \ne \emptyset, $
при $ k \in \Z_+: k \ge K^2, $ существует мультииндекс $ n \in \Z^d $
такой, что куб $ \overline Q_{k,n}^d \subset \mathcal U_i \cap \mathcal U_j. $

Далее, положим $ \delta = \min{i =1,\ldots,N} \delta^0(d,D,x_i) $ (см. лемму 1.2.3).
Пусть $ k \in \Z_+: k \ge K^2, $ и помимо (1.2.59) подчинено неравенству

\begin{equation*} \tag{1.2.60}
2^{-k} < \delta /2,
\end{equation*}
и $ \nu, \nu^\prime \in \Z^d $ таковы, что
$$
\overline Q_{k,\nu}^d \subset D, \\
\overline Q_{k,\nu^\prime}^d \subset D.
$$
Опираясь на лемму 1.2.4 с учётом (1.2.58), выберем последовательность
различных между собой чисел $ \{i_s \in \Nu_{1,N}^1: \ s =0,\ldots,S\} $ такую, что
$$
2^{-k} \nu \in \mathcal U_{i_0}, \ 2^{-k} \nu^\prime \in \mathcal U_{i_S}, \\
\mathcal U_{i_s} \cap \mathcal U_{i_{s +1}} \ne \emptyset, \ s =0,\ldots,S -1.
$$
При этом для $ x \in \overline Q_{k,\nu}^d $ с учётом (1.2.60) выполняется
неравенство
\begin{multline*} \tag{1.2.61}
\| x -x_{i_0} \| \le \| x -2^{-k} \nu \| +\| 2^{-k} \nu -x_{i_0} \| \le 2^{-k} +
\frac{1}{2} \delta^0(d,D,x_{i_0}) < \delta /2 +\frac{1}{2} \delta^0(d,D,x_{i_0}) \\
\le \delta^0(d,D,x_{i_0}) < 2 \delta^0(d,D,x_{i_0}), \\
\text{ т.е. } \overline Q_{k,\nu}^d \subset D \cap U_{x_{i_0}}.
\end{multline*}
Точно так же для $ x \in \overline Q_{k,\nu^\prime}^d $ с учётом (1.2.60)
справедливо неравенство
\begin{multline*} \tag{1.2.62}
\| x -x_{i_S} \| \le \| x -2^{-k} \nu^\prime \| +\| 2^{-k} \nu^\prime -x_{i_S} \| \le \\
2^{-k} +\frac{1}{2} \delta^0(d,D,x_{i_S}) < \delta /2 +
\frac{1}{2} \delta^0(d,D,x_{i_S}) \\
\le \delta^0(d,D,x_{i_S}) < 2 \delta^0(d,D,x_{i_S}), \\
\text{ т.е. } \overline Q_{k,\nu^\prime}^d \subset D \cap U_{x_{i_S}}.
\end{multline*}

Если $ \| 2^{-k} \nu -2^{-k} \nu^\prime \| < \delta, $ то для $ x \in
\overline Q_{k,\nu^\prime}^d $ ввиду (1.2.60) соблюдается неравенство
\begin{multline*} \tag{1.2.63}
\| x -x_{i_0} \| \le \| x -2^{-k} \nu^\prime \| +\| 2^{-k} \nu^\prime -2^{-k} \nu \|
+\| 2^{-k} \nu -x_{i_0} \| \le \\
2^{-k} +\delta +
\frac{1}{2} \delta^0(d,D,x_{i_0}) < \delta /2 +\delta +\frac{1}{2}
\delta^0(d,D,x_{i_0}) \le 
2 \delta^0(d,D,x_{i_0}), \\
\text{ т.е. } \overline Q_{k,\nu^\prime}^d \subset D \cap U_{x_{i_0}}.
\end{multline*}
Полагая $ \mathpzc k^0 = \max_{i =1,\ldots,N} \mathpzc k^1(d,D,x_i), $ в силу
(1.2.59), (1.2.61), (1.2.63) в соответствии с п. 2) леммы 1.2.3
(при $ x^0 = x_{i_0} $) можно построить последовательности
$ \nu^\iota \in \Z^d, \ \iota =0,\ldots,\Iota; \
j^\iota \in \Nu_{1,d}^1, \ \epsilon^\iota \in \{-1,1\}, \ \iota =0,\ldots,\Iota -1, $
для которых имеют место соотношения (1.2.18), (1.2.19), (1.2.20) при
$ x^0 = x_{i_0}, \mathpzc k = \mathpzc k^0, $ а, значит, выполняются условия
(1.2.1), (1.2.2), (1.2.3) при $ \kappa = k \e, \ c_0 \ge
(\max_{i =1,\ldots,N} c_4(d,D,x_i)) 2^{\mathpzc k^0}. $

Если же $ \| 2^{-k} \nu -2^{-k} \nu_\prime \| \ge \delta, $ то при
$ s =1,\ldots,S $ выберем $ n_s \in \Z^d $ таким образом, что куб
$ \overline Q_{k,n_s}^d \subset \mathcal U_{i_{s -1}} \cap \mathcal U_{i_s}
\subset (D \cap U_{x_{i_{s -1}}}) \cap (D \cap U_{x_{i_s}}). $
Отсюда, полагая $ n_0 = \nu, n_{S +1} = \nu^\prime $ и учитывая (1.2.61),
(1.2.62), на основании п. 2) леммы 1.2.3 при $ s =0,\ldots,S $ построим
последовательности $ \nu_s^{\iota_s} \in \Z^d, \ \iota_s =0,\ldots,\Iota_s; \
j_s^{\iota_s} \in \Nu_{1,d}^1, \ \epsilon_s^{\iota_s} \in \{-1,1\}, \
\iota_s =0,\ldots,\Iota_s -1, $ для которых справедливы соотношения
(1.2.18), (1.2.19), (1.2.20) при $ x^0 = x_{i_s}, \
\nu = n_s, \ \nu^\prime = n_{s +1}, \ \mathpzc k = \mathpzc k^0, \ \Iota =
\Iota_s, \ \iota = \iota_s. $
Далее, имея в виду (1.2.55), при $ s =1,\ldots,S $ построим последовательности
$ \nu_s^{\prime \iota_s^\prime} \in \Z^d, \
\iota_s^\prime =0,\ldots,\Iota_s^\prime; \ j_s^{\prime \iota_s^\prime} \in
\Nu_{1,d}^1, \ \epsilon_s^{\prime \iota_s^\prime} \in \{-1,1\}, \
\iota_s^\prime =0,\ldots,\Iota_s^\prime -1, $ обладающие следующими свойствами:
\begin{gather*} \tag{1.2.64}
\Iota_s^\prime \le d 2^{\mathpzc k^0}, 
\nu_s^{\prime 0} = \nu_{s -1}^{\Iota_{s -1}}, 
\nu_s^{\prime \Iota_s^\prime} = \nu_s^0, \\
\nu_s^{\prime \iota_s^\prime +1} = \nu_s^{\prime \iota_s^\prime} +
\epsilon_s^{\prime \iota_s^\prime} e_{j_s^{\prime \iota_s^\prime}}, 
\iota_s^\prime =0,\ldots,\Iota_s^\prime -1; \\
\overline Q_{\mathpzc k^0 +k,\nu_s^{\prime \iota_s^\prime}}^d \subset
\overline Q_{k,n_s}^d, \ \iota_s^\prime =0,\ldots,\Iota_s^\prime.
\end{gather*}

Теперь определим интересующие нас последовательности так. Задавая
$ \Iota_0^\prime = 0, $ положим
\begin{equation*} tag{1.2.65}
\Iota = \sum_{s =0}^S (\Iota_s^\prime +\Iota_s).
\end{equation*}
При $ \iota =0,\ldots,\Iota $ для
$$
\iota = \sum_{s^\prime =0}^{s -1} (\Iota_{s^\prime}^\prime +
\Iota_{s^\prime}) +\iota_s^\prime,
$$
зададим значения
\begin{equation*}
\nu^\iota = \nu_s^{\prime \iota_s^\prime}, \
\iota_s^\prime =0,\ldots,\Iota_s^\prime, 
j^\iota = j_s^{\prime \iota_s^\prime}, \
\epsilon^\iota = \epsilon_s^{\prime \iota_s^\prime}, \
\iota_s^\prime =0,\ldots,\Iota_s^\prime -1,
\end{equation*}
а для
$$
\iota = \sum_{s^\prime =0}^{s -1} (\Iota_{s^\prime}^\prime +
\Iota_{s^\prime}) +\Iota_s^\prime +\iota_s,
$$
определим значения
\begin{equation*}
\nu^\iota = \nu_s^{\iota_s}, \ \iota_s =0,\ldots,\Iota_s, 
j^\iota = j_s^{\iota_s}, \ \epsilon^\iota = \epsilon_s^{\iota_s}, \
\iota_s =0,\ldots,\Iota_s -1.
\end{equation*}

По построению и из (1.2.19), (1.2.20) (при $ \mathpzc k^0 $ вместо
$ \mathpzc k $) следует справедливость (1.2.2), (1.2.3) (при $ \kappa = k \e $).
Для проверки соблюдения (1.2.1), используя (1.2.65), (1.2.64), (1.2.18) и
учитывая, что $ S +1 \le N, $ выводим
\begin{multline*}
\Iota = \sum_{s =1}^S \Iota_s^\prime +\sum_{s =0}^S \Iota_s \le \sum_{s =1}^S
d 2^{\mathpzc k^0} +\sum_{s =0}^S c_4(d,D,x_{i_s}) 2^{\mathpzc k^0}
\| n_s -n_{s +1}\| \le \\
S d 2^{\mathpzc k^0} +(\max_{i =1,\ldots,N} c_4(d,D,x_i)) 2^{\mathpzc k^0}
\sum_{s =0}^S \| n_s -n_{s +1}\| \le \\
c_7(d,D) +c_8(d,D)
\sum_{s =0}^S \| n_s -n_{s +1}\| = \\
c_7 +c_8 \sum_{s =0}^S (2^{-k} \| n_s -n_{s +1}\|
2^k \| \nu -\nu^\prime\|^{-1} \| \nu -\nu^\prime\|) = \\
c_7 +c_8 \sum_{s =0}^S (\| 2^{-k} n_s -2^{-k} n_{s +1}\| \cdot
\| 2^{-k} \nu -2^{-k} \nu^\prime\|^{-1} \cdot \| \nu -\nu^\prime\|) \le \\
c_7 +c_8 (S +1) (\max_{i =1,\ldots,N} \diam(\mathcal U_i)) \delta^{-1}
\| \nu -\nu^\prime\| \le \\
c_7 \| \nu -\nu^\prime \| +c_8 N
(\max_{i =1,\ldots,N} \delta^0(d,D,x_i)) \delta^{-1} \| \nu -\nu^\prime\|,
\end{multline*}
т.е. имеет место (1.2.1) при соответствующем выборе $ c_0. \square $
\bigskip

1.3. В этом пункте вводятся в рассмотрение пространства кусочно-
полиномиальных функций и операторы в них, которые используются для построения
средств приближения для функций из изучаемых нами пространств.
Но сначала приведём некоторые вспомогательные сведения, касающиеся

системы разбиений единицы, используемой для построения
кусочно-полиномиальных функций из интересующих нас пространств. Начнём с
обозначений.

Обозначим через $ \psi^{0,1} $ характеристическую функцию интервала
$ I, $ т.е. функцию, определяемую равенством
$$
\psi^{0, 1}(x) = \begin{cases} 1,  \text{ для } x \in I; \\
0,  \text{ для } x \in \R \setminus I.
\end{cases}
$$
При $ m \in \N $ положим
$$
\psi^{m,1}(x) = \int_I \psi^{m -1,1}(x -y) dy (\text{см. } [11]),
$$
а для $ d \in \N, m \in \Z_+ $ определим
$$
\psi^{m, d}(x) = \prod_{j=1}^d \psi^{m, 1}(x_j), \ x = (x_1,\ldots,x_d) \in \R^d.
$$

Напомним свойства функций $ \psi^{m,d}, $ на которые мы будем
опираться в наших рассмотрениях.

При $ m \in \Z_+, d \in \N $ функция $ \psi^{m,d} $
обладает следующими свойствами (см. [11]):

1) имеет место равенство
\begin{equation*} \tag{1.3.1}
\sgn \psi^{m,d}(x) = \begin{cases} 1,  \text{ для $ x  \in
(m +1) I^d; $ }  \\ 0,   \text{ для $ x \in \R^d \setminus (m +1) I^d; $}
\end{cases}
\end{equation*}

2) для каждого $ \lambda \in \Z_{+ m}^d $ (обобщённая) производная
\begin{equation*} \tag{1.3.2}
\D^\lambda \psi^{m, d} \in L_\infty(\R^d);
\end{equation*}

3) для $ x \in \R^d $ (при $ m =0 $ почти для всех $ x \in \R^d $)
справедливо равенство
\begin{equation*} \tag{1.3.3}
\sum_{\nu \in \Z^d} \psi^{m,d}(x -\nu) =1;
\end{equation*}

4) при $ m \in \N $ для всех $ x \in \R $ (при $ m =0 $ почти для всех $ x \in \R $)
соблюдается равенство
\begin{equation*} \tag{1.3.4}
\psi^{m, 1}(x) = \sum_{\mu \in \Nu_{0, m +1}^1} a_{\mu}^m
\psi^{m, 1}(2x -\mu),
\end{equation*}
где $ a_\mu^m = 2^{-m} C_{m+1}^\mu. $
Используя разложение Ньютона для $ (1 +1)^{m +1} $ и $ (-1 +1)^{m +1}, $ легко
проверить, что при $ m \in \Z_+ $ выполняются равенства
\begin{equation*} \tag{1.3.5}
\sum_{\mu \in \Nu_{0,m +1}^1 \cap (2 \Z)} a_\mu^m =1,
\sum_{\mu \in \Nu_{0,m +1}^1 \cap (2 \Z +1)} a_\mu^m =1.
\end{equation*}

Из (1.3.4) вытекает, что для $ x \in \R^d $ (при $ m =0 $ почти для всех
$ x \in \R^d $) имеет место равенство
\begin{multline*} \tag{1.3.6}
\psi^{m,d}(x) = \prod_{j =1}^d \psi^{m, 1}(x_j) = \prod_{j =1}^d
(\sum_{\mu_j \in \Nu_{0, m +1}^1} a_{\mu_j}^m \psi^{m, 1}(2x_j -\mu_j)) \\
= \sum_{\mu \in \prod_{j =1}^d \{\mu_j \in \Nu_{0, m +1}^1\}}
(\prod_{j =1}^d a_{\mu_j}^m) (\prod_{j =1}^d \psi^{m, 1}(2x_j
-\mu_j)) = \sum_{\mu \in \Nu_{0, m +1}^d} A_{\mu}^{m,d}
\psi^{m,d}(2x -\mu),
\end{multline*}
где $ A_{\mu}^{m,d} = \prod_{j =1}^d a_{\mu_j}^m =
2^{-md} \prod_{j =1}^d C_{m +1}^{\mu_j}. $

Для $ k,m \in \Z_+, d \in \N, \nu \in \Z^d $ обозначим
$$
g_{k,\nu}^{m,d}(x) = \psi^{m,d}(2^k x -\nu) =
\prod_{j =1}^d \psi^{m, 1}( 2^k x_j -\nu_j), x \in \R^d.
$$

Из (1.3.1) следует, что при $ d \in \N, m, k \in \Z_+, \nu \in \Z^d $ носитель
\begin{equation*} \tag{1.3.7}
\supp g_{k, \nu}^{m,d} = 2^{-k} \nu +2^{-k} (m +1) \overline I^d.
\end{equation*}

Отметим некоторые полезные для нас свойства функций $ g_{k, \nu}^{m,d}. $

При $ d \in \N, m, k \in \Z_+ $ для каждого $ \nu^\prime \in \Z^d $ ввиду
(1.3.7) имеет место равенство
\begin{equation*} \tag{1.3.8}
\{ \nu \in \Z^d: Q_{k, \nu^\prime}^d \cap
\supp g_{k, \nu}^{m,d} \ne \emptyset\} = \nu^\prime +\Nu_{-m,0}^d.
\end{equation*}

Из (1.3.3) вытекает, что при $ d \in \N, m, k \in \Z_+ $ для любой области
$ D \subset \R^d $ почти для всех
$ x \in D $ соблюдается равенство
\begin{equation*} \tag{1.3.9}
\sum_{ \nu \in \Z^d: \supp g_{k, \nu}^{m,d} \cap D \ne \emptyset}
g_{k, \nu}^{m,d}(x) =1.
\end{equation*}

Имея в виду (1.3.2), отметим, что при $ d \in \N, m, k \in \Z_+, \nu \in \Z^d,
\lambda \in \Z_{+ m}^d $ выполняется равенство
\begin{multline*} \tag{1.3.10}
\| \D^\lambda g_{k, \nu}^{m,d} \|_{L_\infty (\R^d)} =
2^{k | \lambda|} \| \D^\lambda \psi^{m,d} \|_{L_\infty(\R^d)} =
c_1(d,m,\lambda) 2^{k | \lambda|}.
\end{multline*}

Введём в рассмотрение следующие пространства кусочно-полиномиальных
функций.

При $ d \in \N, l \in \Z_+, m \in \N, k \in \Z_+ $ и области
$ D \subset \R^d, $ полагая
\begin{equation*}
N_k^{d,m,D} = \{\nu \in \Z^d: \supp g_{k, \nu}^{m,d}
\cap D \ne \emptyset\},
\end{equation*}
через $ \mathcal P_k^{l,d,m,D} $ обозначим линейное пространство,
состоящее из функций $ f: \R^d \mapsto \R, $ для каждой из которых существует
набор полиномов
$ \{f_\nu \in \mathcal P^{l,d}, \ \nu \in N_k^{d,m,D}\} $ такой, что
для $ x \in \R^d $ выполняется равенство
\begin{equation*} \tag{1.3.11}
f(x) = \sum_{\nu \in N_k^{d,m,D}} f_\nu(x) g_{k, \nu}^{m,d}(x).
\end{equation*}

Нетрудно проверить, что при $ d \in \N, l \in \Z_+, m \in \N,
k \in \Z_+ $ и ограниченной области $ D \subset \R^d $ отображение,
которое каждому набору полиномов $ \{f_\nu \in \mathcal P^{l,d},
\nu \in N_k^{d,m,D} \} $ ставит в соответствие функцию $ f, $ задаваемую
равенством (1.3.11), является изоморфизмом прямого произведения
$ \card N_k^{d,m,D} $ экземпляров пространства $ \mathcal P^{l,d} $
на пространство $ \mathcal P_k^{l,d,m,D}. $

Предложение 1.3.1

Пусть $ d \in \N, l \in \Z_+, m \in \N, k \in \Z_+, D $ -- ограниченная область
в $ \R^d. $ Тогда линейный оператор $ H_k^{l,d,m,D}:
\mathcal P_k^{l,d,m,D} \mapsto \mathcal P_{k +1}^{l,d,m,D}, $
значение которого на функции $ f \in \mathcal P_k^{l,d,m,D}, $ задаваемой
равенством (1.3.11), определяется соотношением
\begin{multline*} \tag{1.3.12}
(H_k^{l,d,m,D} f)(x) = \\
\sum_{\nu \in N_{k +1}^{d,m,D}}
\biggl(\sum_{\substack{\nu^\prime \in N_k^{d,m,D}, \m \in \Nu_{0, m +1}^d: \\
2 \nu^\prime +\m = \nu }} A_{\m}^{m,d}
f_{\nu^\prime}(x) \biggr) g_{k +1, \nu}^{m,d}(x), \ x \in \R^d,
\end{multline*}
обладает тем свойством, что для $ f \in \mathcal P_k^{l,d,m,D} $ выполняется равенство
\begin{equation*} \tag{1.3.13}
(H_k^{l,d,m,D} f) \mid_{D} = f \mid_{D}.
\end{equation*}

Доказательство.

Прежде всего, заметим, что в условиях предложения для функции
$ f \in \mathcal P_k^{l,d,m,D}, $ задаваемой равенством (1.3.11), в силу
(1.3.6) при $ x \in D $ имеет место равенство
\begin{multline*}
f(x) = \sum_{\nu^\prime \in N_k^{d,m,D}} f_{\nu^\prime}(x)
\psi^{m,d}(2^k x -\nu^\prime) = \\
\sum_{\nu^\prime \in N_k^{d,m,D}} f_{\nu^\prime}(x)
\biggl(\sum_{\m \in \Nu_{0, m +1}^d}
A_{\m}^{m,d} \psi^{m,d}(2^{k +1} x -2 \nu^\prime -\m)\biggr) = \\
\sum_{\nu^\prime \in N_k^{d,m,D}} \sum_{\m \in \Nu_{0, m +1}^d}
A_{\m}^{m,d} f_{\nu^\prime}(x) \psi^{m,d}(2^{k +1} x -2 \nu^\prime -\m) = \\
\sum_{\nu^\prime \in N_k^{d,m,D}, \ \m \in \Nu_{0, m +1}^d}
A_{\m}^{m,d} f_{\nu^\prime}(x) \psi^{m,d}(2^{k +1} x -2 \nu^\prime -\m)
= \\
\sum_{\nu \in \Z^d} \biggl(\sum_{\substack{\nu^\prime \in N_k^{d,m,D}, \m \in
\Nu_{0, m +1}^d: \\ 2 \nu^\prime +\m =  \nu}} A_{\m}^{m,d} f_{\nu^\prime}(x)
\psi^{m,d}(2^{k +1} x -2 \nu^\prime -\m) \biggr)  = \\
\sum_{\nu \in \Z^d} \biggl(\sum_{\substack{\nu^\prime \in N_k^{d,m,D}, \m \in
\Nu_{0, m +1}^d: \\ 2 \nu^\prime +\m = \nu}} A_{\m}^{m,d} f_{\nu^\prime}(x)
\psi^{m,d}(2^{k +1} x -\nu) \biggr) = \\
\sum_{\nu \in \Z^d} \biggl(\sum_{\substack{\nu^\prime \in N_k^{d,m,D}, \m \in
\Nu_{0, m +1}^d: \\ 2 \nu^\prime +\m = \nu}} A_{\m}^{m,d} f_{\nu^\prime}(x) \biggr)
g_{k +1, \nu}^{m,d}(x) = \\
\sum_{\nu \in \Z^d: \supp g_{k +1, \nu}^{m,d} \cap D \ne \emptyset}
\biggl(\sum_{\substack{\nu^\prime \in N_k^{d,m,D}, \m \in \Nu_{0, m +1}^d: \\
2 \nu^\prime +\m = \nu}} A_{\m}^{m,d} f_{\nu^\prime}(x) \biggr)
g_{k +1, \nu}^{m,d}(x) = \\
\sum_{\nu \in N_{k +1}^{d,m,D}}
\biggl(\sum_{\substack{\nu^\prime \in N_k^{d,m,D}, \m \in \Nu_{0, m +1}^d: \\
2 \nu^\prime +\m = \nu}} A_{\m}^{m,d} f_{\nu^\prime}(x)\biggr)
g_{k +1, \nu}^{m,d}(x).
\end{multline*}

Сопоставляя последнее равенство с (1.3.12), приходим к (1.3.13). $ \square $

Замечание.

В условиях предложения 1.3.1, если для $ \nu \in N_{k +1}^{d,m,D}, \
\m \in \Nu_{0, m +1}^d $ существует
$ \nu^\prime \in \Z^d, $ для которого $ 2 \nu^\prime +\m = \nu, $
то для такого $ \nu^\prime $ справедливо включение $ \nu^\prime \in N_k^{d,m,D}. $

В самом деле, при соблюдении условий замечания, выбирая $ x \in D \cap
\supp g_{k +1, \nu}^{m,d}, $ с учётом (1.3.7) получаем, что
\begin{equation*}
2^{-(k +1)} \nu_j \le x_j \le 2^{-(k +1)} \nu_j +2^{-(k +1)} (m +1), \ j =1, \ldots,d,
\end{equation*}
откуда
\begin{multline*}
2^{-k} \nu^\prime_j = 2^{-k} (\nu_j -\m_j) /2 =
2^{-(k +1)} \nu_j -2^{-(k +1)} \m_j \le \\ 2^{-(k +1)} \nu_j
\le x_j \le 2^{-(k +1)} \nu_j +2^{-(k +1)} (m +1) = \\
 2^{-k} (\nu_j -\m_j) /2 +
2^{-(k +1)} \m_j +2^{-(k +1)} (m +1) = \\
2^{-k} \nu^\prime_j +2^{-k} (\m_j +m +1) /2 \le
2^{-k} \nu^\prime_j +2^{-k} (m +1), \ j =1,\ldots,d,
\end{multline*}
т.е. $ x \in D \cap \supp g_{k, \nu^\prime}^{m,d}, $ а, значит,
$ \nu^\prime \in N_k^{d,m,D}. $

При $ d, m \in \N, \nu \in \Z^d $ обозначим через
$ \M^{m,d}(\nu) $ множество
\begin{equation*} \tag{1.3.14}
\M^{m,d}(\nu) = \{ \m \in \Nu_{0, m +1}^d: (\nu_j -\m_j) /2 \in \Z \
\forall j = 1,\ldots,d\}
\end{equation*}
и каждой паре $ \nu \in \Z^d, \m \in \M^{m,d}(\nu) $ сопоставим
$ \n(\nu,\m) \in \Z^d, $ полагая
\begin{equation*} \tag{1.3.15}
\n(\nu,\m) = 2^{-1} (\nu -\m).
\end{equation*}

Предложение 1.3.2

В условиях предложения 1.3.1 для $ f \in \mathcal P_k^{l,d,m,D}, $ задаваемой
равенством (1.3.11), имеет место представление
\begin{multline*} \tag{1.3.16}
(H_k^{l,d,m,D} f)(x) = \\
\sum_{\nu \in N_{k +1}^{d,m,D}}
\biggl(\sum_{\m \in \M^{m,d}(\nu)} A_{\m}^{m,d} f_{\n(\nu,\m)}(x) \biggr)
g_{k +1, \nu}^{m,d}(x), \ x \in \R^d.
\end{multline*}

Доказательство.

Исходя из (1.3.12), принимая во внимание замечание после предложения 1.3.1, с
учётом (1.3.14), (1.3.15) получаем
\begin{multline*}
(H_k^{l,d,m,D} f)(x) = \\
\sum_{\nu \in N_{k +1}^{d,m,D}}
\biggl(\sum_{\substack{\nu^\prime \in N_k^{d,m,D}, \m \in \Nu_{0, m +1}^d: \\
2 \nu^\prime +\m = \nu }} A_{\m}^{m,d}
f_{\nu^\prime}(x) \biggr) g_{k +1, \nu}^{m,d}(x) = \\
\sum_{\nu \in N_{k +1}^{d,m,D}}
\biggl(\sum_{\substack{\m \in \Nu_{0, m +1}^d, \nu^\prime \in \Z^d: 2^{-1} (\nu -\m) \in N_k^{d,m,D}, \\
\nu^\prime = 2^{-1} (\nu -\m)}} A_{\m}^{m,d}
f_{\nu^\prime}(x) \biggr) g_{k +1, \nu}^{m,d}(x) = \\
\sum_{\nu \in N_{k +1}^{d,m,D}}
\biggl(\sum_{\substack{\m \in \Nu_{0, m +1}^d, \nu^\prime \in \Z^d: 2^{-1} (\nu -\m) \in \Z^d, \\
\nu^\prime = 2^{-1} (\nu -\m)}} A_{\m}^{m,d}
f_{\nu^\prime}(x) \biggr) g_{k +1, \nu}^{m,d}(x) = \\
\sum_{\nu \in N_{k +1}^{d,m,D}}
\biggl(\sum_{\substack{\m \in \M^{m,d}(\nu), \nu^\prime \in \Z^d:
\nu^\prime = \n(\nu,\m)}} A_{\m}^{m,d}
f_{\nu^\prime}(x) \biggr) g_{k +1, \nu}^{m,d}(x) = \\
\sum_{\nu \in N_{k +1}^{d,m,D}}
\biggl(\sum_{\m \in \M^{m,d}(\nu)} A_{\m}^{m,d} f_{\n(\nu,\m)}(x) \biggr)
g_{k +1, \nu}^{m,d}(x), \ x \in \R^d. \square
\end{multline*}

Лемма 1.3.3

При $ d \in \N, \nu \in \Z^d, m \in \N $ имеет место равенство
\begin{equation*} \tag{1.3.17}
\sum_{\m \in \M^{m,d}(\nu)} A_{\m}^{m,d} =1.
\end{equation*}

Доказательство.

Используя (1.3.14) и (1.3.5), имеем
\begin{multline*}
\sum_{\m \in \M^{m,d}(\nu)} A_{\m}^{m,d} =
\sum_{\mu \in \Nu_{0,m +1}^d: (\nu_j -\mu_j) /2 \in \Z, \ j =1, \ldots,d}
\prod_{j =1}^d a_{\mu_j}^m = \\
\sum_{\mu \in \prod_{j =1}^d \{\mu_j \in \Nu_{0,m +1}^1: (\nu_j -\mu_j) /2 \in \Z\}}
\prod_{j =1}^d a_{\mu_j}^m =
\prod_{j =1}^d \biggl(\sum_{\mu_j \in \Nu_{0,m +1}^1: (\nu_j -\mu_j) /2 \in \Z}
a_{\mu_j}^m\biggr) = \\
\prod_{j =1}^d \biggl(\sum_{\mu_j \in \Nu_{0,m +1}^1 \cap (2 \Z +\nu_j)}
a_{\mu_j}^m\biggr) = \prod_{j =1}^d 1 = 1. \square
\end{multline*}
\bigskip

1.4. В этом пункте определяются средства приближения функций из
рассматриваемых нами пространств.

Сначала приведём некоторые факты, относящиеся к полиномам, которыми
мы будем пользоваться ниже.

Сформулируем лемму 1.4.1, доказательство которой можно найти в [12].

Лемма 1.4.1

Пусть $ l \in \Z_+,\ d \in \N,\ 1 \le p, q \le \infty,\ \lambda \in \Z_{+ l}^d $
и $ \rho, \sigma >0. $ Тогда существует константа $ c_1(l,d,\lambda,\rho,\sigma) >0 $
такая, что для любых измеримых по
Лебегу множеств $ D, Q \subset \R^d, $ для которых можно найти
$ \delta >0 $ и $ x^0 \in \R^d $ такие, что
$ \sup_{x \in D} \| x -x^0 \| \le \rho \delta $ и $ x^0 +\sigma \delta I^d \subset Q, $
для любого полинома $ f \in \mathcal P^{l,d} $
выполняется неравенство
\begin{equation*} \tag{1.4.1}
\| \D^\lambda f \|_{L_q(D)} \le c_1
\delta^{-|\lambda| -d /p +d /q} \|f\|_{L_p(Q)}.
\end{equation*}

Справедлива также следующая лемма.

Лемма 1.4.2

Пусть $ l \ \in \Z_+, \ d \in \N, \ \rho \in \R_+^d. $
Тогда существует константа $ c_2(l,d,\rho) >0 $ такая, что
для любого $ \delta >0 $ и любого $ x^0 \in \R^d $ для $ Q = x^0 +\delta \rho I^d $
можно построить систему линейных операторов
$$
P_{\delta,x^0,\rho}^{l,d,\lambda}: L_1(Q) \mapsto \mathcal P^{l,d}, \
\lambda \in \Z_{+ l}^d,
$$
обладающих следующими свойствами:

1) для $ f \in \mathcal P^{l,d} $ при $ \lambda \in \Z_{+ l}^d $
имеет место равенство
\begin{equation*} \tag{1.4.2}
P_{\delta,x^0,\rho}^{l,d,\lambda}(f \mid_Q) = f;
\end{equation*}

2) при $ 1 \le p \le \infty $ для $ f \in L_p(Q), \lambda \in \Z_{+ l}^d $
справедливо неравенство
\begin{equation*} \tag{1.4.3}
\| P_{\delta,x^0,\rho}^{l,d,\lambda} f \|_{L_p(Q)} \le
c_2 \|f\|_{L_p(Q)};
\end{equation*}

3) если для функции $ f \in L_1(Q) $ при $ \lambda \in \Z_{+ l}^d $ её
обобщённая производная $ \D^\lambda f \in L_1(Q), $ то
\begin{equation*}
\D^\lambda P_{\delta,x^0,\rho}^{l,d,\lambda} f =
P_{\delta,x^0,\rho}^{l -|\lambda|,d,0}  \D^\lambda f.
\end{equation*}

Для доказательства леммы 1.4.2 достаточно повторить доказательство
теоремы 1.3.1 из [5], заменив в нём $ I^d $ на $ \rho I^d. $

Прежде чем формулировать теорему 1.4.3, введём следующие обозначения.
Для области $ D \subset \R^d $ и $ x \in D $ обозначим через $ G(D,x) $
множество
$$
G(D,x) = \{y \in \R^d: x +ty \in D\ \forall t \in \overline I\}.
$$
Обозначим ещё для области $ D \subset \R^d $ и $ x \in D $ при $ l \in \N $ и
$ k =1,\ldots,l $ через
$$
F_k^l(D,x) = \{y \in \R^d: (x -y)/k \in B^d,\ y +tl(x -y)/k \in D \
\forall \ t \in \overline I\}.
$$

Теперь мы можем сформулировать теорему  1.4.3.

Теорема 1.4.3

Пусть $ l,d \in \N,\ 1<p<\infty $ и $ D $ --- ограниченная область в
$ \R^d, $ обладающая следующими свойствами:

1) существует константа $ \sigma >0 $ такая, что для любого
$ x \in D $ имеет место неравенство
$$
\mes G(D,x) \ge \sigma;
$$

2) при $ k =1,\ldots,l $ отображение $ F_k^l(D,\cdot), $ которое каждой точке
$ x \in D $ ставит в соответствие измеримое множество $ F_k^l(D,x), $
равномерно непрерывно на $ D $ в том смысле, что для любого
$ \epsilon >0 $ существует $ \delta >0 $ такое, что
для $ x, x^\prime \in D $ таких, что $ \|x -x^\prime\| < \delta, $
выполняется неравенство
$$
\mes(F_k^l(D,x) \triangle F_k^l(D,x^\prime)) < \epsilon.
$$
Тогда существует константа $ c_3(l,d,p,D) >0 $ такая, что для любой
функции $ f \in L_p(D) $ существует полином $ g \in \mathcal P^{l -1,d} $
такой, что
\begin{equation*} \tag{1.4.4}
\| f -g \|_{L_p(D)} \le c_3 \cdot \biggl(\int_{B^d} \int_{D_{l \xi}}
| \Delta_\xi^l f(x)|^p \,dx \,d\xi \biggr)^{1/p}.
\end{equation*}

Доказательство теоремы 1.4.3 содержится в [5].

Замечание.

Несложно проверить, что при $ d \in \N, \rho \in \R_+^d $ область $ \rho I^d $
удовлетворяет условиям пп. 1), 2) теоремы 1.4.3. Тем самым,
в (1.4.4) в качестве $ D $ можно положить $ \rho I^d. $

Теорема 1.4.4

Пусть $ l,d \in \N, \ \rho \in \R_+^d, \ 1 < p < \infty. $ Тогда существует
константа $ c_4(l,d,\rho,p) >0 $ такая, что для любых
$ x^0 \in \R^d $ и $ \delta >0 $ для $ Q = x^0 +\delta \rho I^d $
для любой функции $ f \in L_p(Q) $ существует полином
$ g \in \mathcal P^{l -1,d}, $ для которого выполняется неравенство
\begin{equation*} \tag{1.4.5}
\| f -g \|_{L_p(Q)} \le c_4 \delta^{-d/p} \biggl(\int_{\delta B^d}
\int_{Q_{l \xi}} | \Delta_\xi^l f(x)|^p \,dx \,d\xi \biggr)^{1/p}.
\end{equation*}

Доказательство.

Для $ f \in L_p(Q) $ положим $ \phi(x) = f(x^0 +\delta x) \in L_p(\rho I^d) $
и выберем $ \psi \in \mathcal P^{l -1,d}, $ для которого согласно
(1.4.4) (см. замечание после теоремы 1.4.3) соблюдается неравенство
\begin{equation*}
\| \phi -\psi \|_{L_p(\rho I^d)} \le c_4\cdot
\biggl(\int_{B^d} \int_{(\rho I^d)_{l \eta}} | \Delta_\eta^l \phi(y)|^p \,dy \,d\eta \biggr)^{1/p}.
\end{equation*}
Тогда отсюда, учитывая, что
\begin{multline*}
\| \phi -\psi \|_{L_p(\rho I^d)} = \delta^{-d/p}
\| f(\cdot) -\psi((\cdot -x^0)/\delta) \|_{L_p(Q)}, \\
\left(\int_{B^d} \int_{(\rho I^d)_{l \eta}} | \Delta_\eta^l \phi(y)|^p \,dy \,d\eta \right)^{1/p}
= \delta^{-2d/p} \left(\int_{\delta B^d}\int_{Q_{l \xi}} | \Delta_\xi^l
f(x)|^p \,dx \,d\xi \right)^{1/p},
\end{multline*}
а также то обстоятельство, что
$$
\psi((\cdot -x^0)/\delta) \in \mathcal P^{l -1,d},
$$
приходим к (1.4.5). $ \square $

Теорема 1.4.5

Пусть $ l, d \in \N, \ \rho \in \R_+^d, \ 1 < p < \infty. $
Тогда существует константа $ c_5(l,d,\rho,p) >0 $ такая, что
для любого $ \delta >0 $ и любого $ x^0 \in \R^d $ для $ Q = x^0 +\delta \rho I^d, $
для $ f \in L_p(Q), \ \lambda \in \Z_{+ l -1}^d $ выполняется неравенство
\begin{equation*} \tag{1.4.6}
\| f -P_{\delta,x^0,\rho}^{l -1,d,\lambda} f \|_{L_p(Q)} \le
c_5 \delta^{-d/p} \left(\int_{\delta B^d} \int_{Q_{l \xi}}
| \Delta_\xi^l f(x)|^p \,dx \,d\xi \right)^{1/p}.
\end{equation*}

Доказательство.

Выбирая для функции $ f \in L_p(Q) $ полином
$ g \in \mathcal P^{l -1,d}, $ для которого соблюдается
(1.4.5), и применяя (1.4.2), (1.4.3), имеем
\begin{multline*}
\| f -P_{\delta,x^0,\rho}^{l -1,d,\lambda} f \|_{L_p(Q)}  \le
\| f -g \|_{L_p(Q)}
+\| P_{\delta,x^0,\rho}^{l -1,d,\lambda}(g -f) \|_{L_p(Q)} \\
 \le (1 +c_2) \| f -g \|_{L_p(Q)}.
\end{multline*}
Подставляя сюда (1.4.5), получаем (1.4.6). $ \square $

При $ d \in \N, l, k \in \Z_+, \nu \in \Z^d $ определим линейный
оператор $ S_{k,\nu}^{l,d}: L_1(Q_{k,\nu}^d) \mapsto \mathcal P^{l,d}, $
полагая $ S_{k,\nu}^{l,d} = P_{\delta, x^0,\e}^{l,d,0} $ при
$ \delta = 2^{-k}, x^0 = 2^{-k} \nu $ (см. лемму 1.4.2).
Отметим, что в ситуации, когда $ Q_{k,\nu}^d \subset D, $ где $ D $ --
область в $ \R^d, $ для $ f \in L_1(D) $ вместо $ S_{k,\nu}^{l,d}( f \mid_{Q_{k,\nu}^d}) $
будем писать $ S_{k,\nu}^{l,d} f.$

Для области $ D \subset \R^d $ и $ k \in \Z_+ $ таких, что множество
\begin{equation*} \tag{1.4.7}
\{\nu^\prime \in \Z^d: \overline Q_{k,\nu^\prime}^d \subset D\} \ne \emptyset,
\end{equation*}
для каждого $ \nu \in \Z^d $ фиксируем некоторый $ \nu_k^D(\nu) \in \Z^d, $
для которого
\begin{equation*} \tag{1.4.8}
\overline Q_{k,\nu_k^D(\nu)}^d \subset D, \\
\text{ а } \\
\|\nu -\nu_k^D(\nu)\| = \min_{\nu^\prime \in \Z^d:
\overline Q_{k,\nu^\prime}^d \subset D} \|\nu -\nu^\prime\|,
\end{equation*}
и при $ d \in \N, l \in \Z_+, m \in \N $ для ограниченной области
$ D \subset \R^d $ и $ k \in \Z_+, $ удовлетворяющих (1.4.7),
определим линейный непрерывный оператор
$ E_k^{l,d,m,D}: L_1(D) \mapsto \mathcal P_k^{l,d,m,D} \cap
W_\infty^m(\R^d) $ равенством
\begin{equation*} \tag{1.4.9}
E_k^{l,d,m,D} f = \sum_{\nu \in N_k^{d,m,D}}
(S_{k,\nu_k^D(\nu)}^{l,d} f ) g_{k, \nu}^{m,d}, \ f \in L_1(D).
\end{equation*}

Обозначим ещё через $ \mathcal I^D $ линейное отображение, которое каждой
функции $ f, $ заданной на области $ D \subset \R^d, $ сопоставляет функцию
$ \mathcal I^D f, $ определяемую
на $ \R^d $ равенством
\begin{equation*}
(\mathcal I^D f)(x) = \begin{cases} f(x), \text{ при } x \in D; \\
0, \text{ при } x \in \R^d \setminus D.
\end{cases}
\end{equation*}

Предложение 1.4.6

Пусть $ d \in \N, \ l \in \N, \ m \in \N, \ 1 < p < \infty $
и ограниченная область $ D \subset \R^d $ такова, что существуют $ k^0(d,D)
\in \Z_+, \ \gamma^0(d,D) \in \R_+, $ для которых при любом $ k \in \Z_+:
k \ge k^0, $ для каждого $ \nu \in \Z^d:
Q_{k,\nu}^d \cap D \ne \emptyset, $ существует $ \nu^\prime \in \Z^d $
такой, что
\begin{equation*} \tag{1.4.10}
\overline Q_{k,\nu^\prime}^d \subset D \cap (2^{-k} \nu +\gamma^0 2^{-k} B^d).
\end{equation*}
Тогда для любой функции $ f \in L_p(D) $ в $ L_p(D) $ имеет место равенство
\begin{equation*} \tag{1.4.11}
f = (E_{k^0}^{l -1,d,m,D} f) \mid_D +\sum_{k = k^0 +1}^\infty
(E_k^{l -1,d,m,D} f -E_{k -1}^{l -1,d,m,D} f) \mid_D.
\end{equation*}

Доказательство.

Сначала отметим свойства некоторых вспомогательных множеств и других
объектов, которые понадобятся для доказательства предложения.

Для $ d \in \N, \ k \in \Z_+, \ m \in \N $ и области $ D \subset \R^d $
обозначим
$$
G_k^{d,m,D} = \cup_{\nu \in N_k^{d,m,D}} \supp g_{k, \nu}^{m,d},
$$
и рассмотрим представление
\begin{multline*} \tag{1.4.12}
G_k^{d,m,D} = G_k^{d,m,D} \cap \R^d = G_k^{d,m,D} \cap
((\cup_{n \in \Z^d} Q_{k,n}^d) \cup \mathcal A_k^d) = \\
(\cup_{n \in \Z^d} (G_k^{d,m,D} \cap Q_{k,n}^d)) \cup
(G_k^{d,m,D} \cap \mathcal A_k^d) = \\
(\cup_{n \in \Z^d: Q_{k,n}^d \cap G_k^{d,m,D} \ne \emptyset}
(G_k^{d,m,D} \cap Q_{k,n}^d)) \cup (G_k^{d,m,D} \cap \mathcal A_k^d), \\
\text{ где } \mes \mathcal A_k^d =0, \ \mathcal A_k^d \cap Q_{k,n}^d = \emptyset, \
Q_{k,n}^d \cap Q_{k,n^\prime}^d = \emptyset, \ n,n^\prime \in
\Z^d: n \ne n^\prime.
\end{multline*}
Заметим, что если при $ n \in \Z^d $ пересечение
\begin{equation*}
Q_{k,n}^d \cap G_k^{d,m,D} \ne \emptyset,
\end{equation*}
т.е.
\begin{equation*}
Q_{k,n}^d \cap (\cup_{\nu \in N_k^{d,m,D}} \supp g_{k, \nu}^{m,d}) =
\cup_{\nu \in N_k^{d,m,D}} (Q_{k,n}^d \cap \supp g_{k, \nu}^{m,d})
\ne \emptyset,
\end{equation*}
то существует $ \nu \in N_k^{d,m,D}, $ для которого
\begin{equation*}
Q_{k,n}^d \cap \supp g_{k, \nu}^{m,d} \ne \emptyset.
\end{equation*}
Для $ n \in \Z^d, \ \nu \in N_k^{d,m,D}: Q_{k,n}^d \cap
\supp g_{k, \nu}^{m,d} \ne \emptyset, $ выбирая $ x \in
Q_{k,n}^d \cap \supp g_{k, \nu}^{m,d}, $ с учётом (1.3.7) при $ j =1,\ldots,d $
имеем
$$
2^{-k} n_j < x_j \le 2^{-k} \nu_j +2^{-k} (m +1); \
2^{-k} \nu_j \le x_j < 2^{-k} n_j +2^{-k},
$$
или
$$
n_j < \nu_j +m +1; \
\nu_j < n_j +1,
$$
откуда
$$
n_j +1 \le \nu_j +m +1; \
\nu_j \le n_j,
$$
следовательно,
$$
2^{-k} n_j +2^{-k} \le 2^{-k} \nu_j +2^{-k} (m +1); \
2^{-k} \nu_j \le 2^{-k} n_j,
$$
а это значит, что
\begin{equation*} \tag{1.4.13}
Q_{k,n}^d \subset \overline Q_{k,n}^d \subset
\supp g_{k, \nu}^{m,d} \subset G_k^{d,m,D}, \
n \in \Z^d, \ \nu \in N_k^{d,m,D}: Q_{k,n}^d \cap
\supp g_{k, \nu}^{m,d} \ne \emptyset.
\end{equation*}
Из (1.4.12) и (1.4.13) получаем
\begin{multline*} \tag{1.4.14}
G_k^{d,m,D} = (\cup_{n \in \Z^d: Q_{k,n}^d \cap G_k^{d,m,D} \ne
\emptyset} Q_{k,n}^d) \cup (G_k^{d,m,D} \cap \mathcal A_k^d), \\
d \in \N, \ m \in \N, \ k \in \Z_+, \ D  \text{ -- область в } \R^d.
\end{multline*}

В условиях предложения при $ k \ge k^0 $ для $ \nu \in N_k^{d,m,D}, $
выбирая точку $ z \in (2^{-k} \nu +2^{-k} (m +1) \overline I^d) \cap D, $ и
учитывая открытость множества $ D, $ найдём точку
$ y \in (2^{-k} \nu +2^{-k} (m +1) I^d) \cap D. $
Затем возьмём куб $ \overline Q_{k,\rho}^d, \ \rho \in \Z^d, $
содержащий точку $ y. $ Тогда поскольку
$ (2^{-k} \nu +2^{-k} (m +1) I^d) \cap D $ является
окрестностью точки $ y, $ принадлежащей замыканию $ Q_{k,\rho}^d, $ то
множество $ (2^{-k} \nu +2^{-k} (m +1) I^d) \cap D \cap
Q_{k,\rho}^d \ne \emptyset. $
Принимая во внимание это обстоятельство и (1.4.13), получаем, что
\begin{equation*} \tag{1.4.15}
Q_{k,\rho}^d \subset \overline Q_{k,\rho}^d \subset
\supp g_{k, \nu}^{m,d},
\end{equation*}
и $ Q_{k,\rho}^d \cap D \ne \emptyset. $
Исходя из условий предложения, согласно  (1.4.10) выберем $ \nu^\prime \in \Z^d, $
для которого
\begin{equation*} \tag{1.4.16}
\overline Q_{k,\nu^\prime}^d \subset D \cap
(2^{-k} \rho +\gamma^0 2^{-k} B^d).
\end{equation*}
Из (1.4.16) и определения $ \nu_k^D(\nu) $ (см. (1.4.8)) имеем
\begin{equation*} \tag{1.4.17}
\|\nu -\nu_k^D(\nu)\| \le \|\nu -\nu^\prime\| \le
\|\nu -\rho\| +\|\\rho -\nu^\prime\|.
\end{equation*}
Ввиду (1.4.15), (1.3.7) получаем
$$
| 2^{-k} \nu_j -2^{-k} \rho_j | \le 2^{-k} (m +1),
$$
или
$$
| \nu_j -\rho_j | \le (m +1), \ j =1,\ldots,d,
$$
т.е.
\begin{equation*} \tag{1.4.18}
\|\nu -\rho\| \le m +1.
\end{equation*}
Кроме того, из (1.4.16) следует, что
$$
| 2^{-k} \nu_j^\prime -2^{-k} \rho_j | \le 2^{-k} \gamma^0,
$$
или
$$
| \nu_j^\prime -\rho_j | \le \gamma^0, \ j =1,\ldots,d,
$$
т.е.
\begin{equation*} \tag{1.4.19}
\| \nu^\prime -\rho\| \le \gamma^0.
\end{equation*}
Объединяя (1.4.17), (1.4.18), (1.4.19), приходим к выводу, что
\begin{equation*} \tag{1.4.20}
\|\nu -\nu_k^D(\nu)\| \le m +1 +\gamma^0 = c_6(d,m,D), \
\nu \in N_k^{d,m,D}, \ k \ge k^0.
\end{equation*}

Отметим ещё, что в условиях предложения при $ k \ge k^0 $
для $ n \in \Z^d: Q_{k,n}^d \cap G_k^{d,m,D} \ne \emptyset, \
\nu \in N_k^{d,m,D}: \supp g_{k, \nu}^{m,d} \cap Q_{k,n}^d
\ne \emptyset, $ для $ x \in \overline Q_{k,\nu_k^D(\nu)}^d $ в силу
(1.4.13), (1.3.7) и (1.4.20) при $ j =1,\ldots, d $ выполняется неравенство
\begin{multline*}
| x_j -2^{-k} n_j | \le \\
| x_j -2^{-k} (\nu_k^D(\nu))_j |
+| 2^{-k} (\nu_k^D(\nu))_j -2^{-k} \nu_j |
+| 2^{-k} \nu_j -2^{-k} n_j | \le \\
2^{-k} +2^{k} \|\nu -\nu_k^D(\nu)\| +2^{-k} (m +1) \le
2^{-k} (1 +c_6 +m +1) = \\
\gamma^1(d,m,D) 2^{-k},
\end{multline*}
т.е.
\begin{equation*} \tag{1.4.21}
Q_{k,\nu_k^D(\nu)}^d \subset \overline Q_{k,\nu_k^D(\nu)}^d
\subset (2^{-k} n +\gamma^1 2^{-k} B^d).
\end{equation*}

В условиях предложения при $ k \ge k^0, \ n \in \Z^d: Q_{k,n}^d \cap G_k^{d,m,D}
\ne \emptyset, $ зададим
\begin{equation*} \tag{1.4.22}
x_{k,n}^{d,m,D} = 2^{-k} n -\gamma^1 2^{-k} \e; \
\delta_{k,n}^{d,m,D} = 2 \gamma^1 2^{-k}
\end{equation*}
и определим куб $ D_{k,n}^{d,m,D} $ равенством
\begin{equation*} \tag{1.4.23}
D_{k,n}^{d,m,D} = x_{k,n}^{d,m,D} +
\delta_{k,n}^{d,m,D} I^d = \inter (2^{-k} n +\gamma^1 2^{-k} B^d).
\end{equation*}

Из (1.4.22), (1.4.23) с учётом того, что $ \gamma^1 > 1, $ видно, что
справедливо включение
\begin{equation*} \tag{1.4.24}
Q_{k, n}^d \subset D_{k,n}^{d,m,D}, \ n \in \Z^d:
Q_{k,n}^d \cap G_k^{d,m,D} \ne \emptyset, \ k \ge k^0.
\end{equation*}

При $ k \ge k^0, \ n \in \Z^d: Q_{k,n}^d \cap G_k^{d,m,D} \ne \emptyset, $
введём в рассмотрение линейный оператор
$ \mathcal S_{k,n}^{l -1,d,m,D}, $ полагая
\begin{equation*}
\mathcal S_{k,n}^{l -1,d,m,D} = P_{\delta,x^0,\e}^{l -1,d,0}
\end{equation*}
при $ x^0 = x_{k,n}^{d,m,D}, \ \delta = \delta_{k,n}^{d,m,D}. $

Учитывая (1.4.24), (1.4.23), (1.4.22),  нетрудно видеть, что в условиях
предложения существует константа $ c_7(d,m,D) >0 $ такая, что для любого
$ k \in \Z_+: k \ge k^0, $ для каждого $ x \in \R^d $ число
\begin{equation*} \tag{1.4.25}
\card \{ n \in \Z^d: Q_{k,n}^d \cap G_k^{d,m,D} \ne \emptyset, \
x \in D_{k,n}^{d,m,D} \} \le c_7.
\end{equation*}

Из (1.4.21) и (1.4.23) следует, что при $ k \ge k^0 $ для $ n \in \Z^d:
Q_{k,n}^d \cap G_k^{d,m,D} \ne \emptyset, \ \nu \in N_k^{d,m,D}:
\supp g_{k, \nu}^{m,d} \cap Q_{k,n}^d \ne \emptyset, $ имеет место включение
\begin{equation*} \tag{1.4.26}
Q_{k,\nu_k^D(\nu)}^d \subset D_{k,n}^{d,m,D}.
\end{equation*}

Из (1.3.8) вытекает, что при $ k \in \Z_+ $ для $ n \in \Z^d:
Q_{k,n}^d \cap G_k^{d,m,D} \ne \emptyset, $ верно неравенство
\begin{equation*} \tag{1.4.27}
\card \{ \nu \in N_k^{d,m,D}: \supp g_{k, \nu}^{m,d} \cap
Q_{k,n}^d \ne \emptyset \} \le c_8(d,m).
\end{equation*}

Пусть теперь в условиях предложения $ f \in L_p(D) $ и $ k \in \Z_+: k \ge k^0. $
Тогда, полагая $ F = \mathcal I^D f, $ ввиду (1.3.9) имеем
\begin{multline*} \tag{1.4.28}
\biggl\|f -\biggl((E_{k^0}^{l -1,d,m,D} f) \mid_D
+\sum_{\k = k^0 +1}^k (E_{\k}^{l -1,d,m,D} f
-E_{\k -1}^{l -1,d,m,D} f) \mid_D\biggr) \biggr\|_{L_p(D)}^p = \\
\| F \mid_D -(E_{k}^{l -1,d,m,D} f) \mid_D \|_{L_p(D)}^p = \\
\biggl\| \biggl(F \biggl(\sum_{ \nu \in \Z^d: \supp g_{k, \nu}^{m,d} \cap D \ne \emptyset}
g_{k, \nu}^{m,d}\biggr)\biggr) \big|_D -(E_k^{l -1,d,m,D} f) \mid_D \biggr\|_{L_p(D)}^p \le \\
\biggl\| F \biggl(\sum_{ \nu \in N_k^{d,m,D}}
g_{k, \nu}^{m,d}\biggr) -(E_k^{l -1,d,m,D} f) \biggr\|_{L_p(\R^d)}^p = \\
\biggl\| F \biggl(\sum_{ \nu \in N_k^{d,m,D}}
g_{k, \nu}^{m,d}\biggr) -\sum_{\nu \in N_k^{d,m,D}}
(S_{k, \nu_k^D(\nu)}^{l -1,d} f) g_{k, \nu}^{m,d} \biggr\|_{L_p(\R^d)}^p = \\
\biggl\| F \biggl(\sum_{ \nu \in N_k^{d,m,D}}
g_{k, \nu}^{m,d}\biggr) -\sum_{\nu \in N_k^{d,m,D}}
(S_{k, \nu_k^D(\nu)}^{l -1,d} F) g_{k, \nu}^{m,d} \biggr\|_{L_p(\R^d)}^p = \\
\biggl\| \sum_{ \nu \in N_k^{d,m,D}} (F -
(S_{k, \nu_k^D(\nu)}^{l -1,d} F)) g_{k, \nu}^{m,d} \biggr\|_{L_p(\R^d)}^p = \\
\int_{\R^d} \biggl| \sum_{\nu \in N_k^{d,m,D}} (F(x) -
(S_{k, \nu_k^D(\nu)}^{l -1,d} F)(x))
g_{k, \nu}^{m,d}(x)\biggr|^p dx = \\
\int_{G_k^{d,m,D}} \biggl| \sum_{\nu \in N_k^{d,m,D}} (F(x) -
(S_{k, \nu_k^D(\nu)}^{l -1,d} F)(x))
g_{k, \nu}^{m,d}(x)\biggr|^p dx.
\end{multline*}
В силу (1.4.14) (см. также (1.4.12)) выводим
\begin{multline*} \tag{1.4.29}
\int_{G_k^{d,m,D}} \biggl| \sum_{\nu \in N_k^{d,m,D}} (F(x) -
(S_{k, \nu_k^D(\nu)}^{l -1,d} F)(x))
g_{k, \nu}^{m,d}(x)\biggr|^p dx = \\
\int_{\cup_{n \in \Z^d: Q_{k,n}^d \cap G_k^{d,m,D} \ne \emptyset}
Q_{k,n}^d} \biggl| \sum_{\nu \in N_k^{d,m,D}} (F(x) -
(S_{k, \nu_k^D(\nu)}^{l -1,d} F)(x))
g_{k, \nu}^{m,d}(x)\biggr|^p dx = \\
\sum_{n \in \Z^d: Q_{k,n}^d \cap G_k^{d,m,D} \ne \emptyset}
\int_{Q_{k,n}^d} \biggl| \sum_{\nu \in N_k^{d,m,D}} (F(x) -
(S_{k, \nu_k^D(\nu)}^{l -1,d} F)(x))
g_{k, \nu}^{m,d}(x)\biggr|^p dx = \\
\sum_{n \in \Z^d: Q_{k,n}^d \cap G_k^{d,m,D} \ne \emptyset}
\int_{Q_{k,n}^d} \biggl| \sum_{\nu \in N_k^{d,m,D}:
\supp g_{k, \nu}^{m,d} \cap Q_{k,n}^d \ne \emptyset} (F(x) -
(S_{k, \nu_k^D(\nu)}^{l -1,d} F)(x))
g_{k, \nu}^{m,d}(x)\biggr|^p dx = \\
\sum_{n \in \Z^d: Q_{k,n}^d \cap G_k^{d,m,D} \ne \emptyset}
\biggl\| \sum_{\nu \in N_k^{d,m,D}:
\supp g_{k, \nu}^{m,d} \cap Q_{k,n}^d \ne \emptyset} (F -
(S_{k, \nu_k^D(\nu)}^{l -1,d} F))
g_{k, \nu}^{m,d} \biggr\|_{L_p(Q_{k,n}^d)}^p \le \\
\sum_{n \in \Z^d: Q_{k,n}^d \cap G_k^{d,m,D} \ne \emptyset}
\biggl( \sum_{\nu \in N_k^{d,m,D}:
\supp g_{k, \nu}^{m,d} \cap Q_{k,n}^d \ne \emptyset}
\| (F -(S_{k, \nu_k^D(\nu)}^{l -1,d} F))
g_{k, \nu}^{m,d} \|_{L_p(Q_{k,n}^d)}\biggr)^p.
\end{multline*}

Далее, для $ n \in \Z^d: Q_{k,n}^d \cap G_k^{d,m,D} \ne \emptyset, \
\nu \in N_k^{d,m,D}: \supp g_{k, \nu}^{m,d} \cap Q_{k,n}^d
\ne \emptyset, $ с учётом (1.3.10) получаем
\begin{multline*} \tag{1.4.30}
\| (F -(S_{k, \nu_k^D(\nu)}^{l -1,d} F))
g_{k, \nu}^{m,d} \|_{L_p(Q_{k,n}^d)} \le
\| F -S_{k, \nu_k^D(\nu)}^{l -1,d} F \|_{L_p(Q_{k,n}^d)} \le \\
\| F -\mathcal S_{k,n}^{l -1,d,m,D} F \|_{L_p(Q_{k,n}^d)} +
\| \mathcal S_{k,n}^{l -1,d,m,D} F
-S_{k, \nu_k^D(\nu)}^{l -1,d} F \|_{L_p(Q_{k,n}^d)}.
\end{multline*}

Учитывая (1.4.24), на основании (1.4.6) для $ n \in \Z^d: Q_{k,n}^d
\cap G_k^{d,m,D} \ne \emptyset, $ заключаем, что
\begin{multline*} \tag{1.4.31}
\| F -\mathcal S_{k,n}^{l -1,d,m,D} F \|_{L_p(Q_{k,n}^d)} \le
\| F -\mathcal S_{k,n}^{l -1,d,m,D} F \|_{L_p(D_{k,n}^{d,m,D})} \\
\le c_9 2^{kd /p} \biggl(\int_{ c_{10} 2^{-k} B^d}
\int_{ (D_{k,n}^{d,m,D})_{l \xi}}
|\Delta_{\xi}^{l} F(x)|^p dx d\xi\biggr)^{1/p}.
\end{multline*}

Принимая во внимание (1.4.24), (1.4.26), (1.4.1), (1.4.2), (1.4.3) и снова
(1.4.26), а также, опираясь на (1.4.6), для $ n \in \Z^d: Q_{k,n}^d \cap
G_k^{d,m,D} \ne \emptyset, \ \nu \in N_k^{d,m,D}:
\supp g_{k, \nu}^{m,d} \cap Q_{k,n}^d \ne \emptyset, $
находим, что
\begin{multline*} \tag{1.4.32}
\|\mathcal S_{k,n}^{l -1,d,m,D} F
-S_{k,\nu_k^D(\nu)}^{l -1,d} F \|_{L_p(Q_{k,n}^d)} \le \\
c_{11} \|S_{k,\nu_k^D(\nu)}^{l -1,d}
(\mathcal S_{k,n}^{l -1,d,m,D} F -F) \|_{L_p(Q_{k, \nu_k^D(\nu)}^d)} \le
c_{12} \| F -\mathcal S_{k,n}^{l -1,d,m,D} F \|_{L_p(D_{k,n}^{d,m,D})} \le \\
c_{13} 2^{kd /p} \biggl(\int_{ c_{10} 2^{-k} B^d}
\int_{ (D_{k,n}^{d,m,D})_{l \xi}}
|\Delta_\xi^l F(x)|^p dx d\xi\biggr)^{1/p}.
\end{multline*}

Соединяя (1.4.30), (1.4.31) и (1.4.32), для $ n \in \Z^d: Q_{k,n}^d
\cap G_k^{d,m,D} \ne \emptyset, \ \nu \in N_k^{d,m,D}:
\supp g_{k, \nu}^{m,d} \cap Q_{k,n}^d \ne \emptyset, $ получаем
\begin{equation*}
\|(F -S_{k, \nu_k^D(\nu)}^{l -1,d} F)
g_{k, \nu}^{m,d} \|_{L_p(Q_{k,n}^d)} \le \\
c_{14} 2^{kd /p} \biggl(\int_{ c_{10} 2^{-k} B^d}
\int_{ (D_{k,n}^{d,m,D})_{l \xi}}
|\Delta_\xi^l F(x)|^p dx d\xi\biggr)^{1/p}.
\end{equation*}

Подставляя эту оценку в (1.4.29) и применяя (1.4.27), а затем используя
(1.4.25), выводим
\begin{multline*}
\int_{G_k^{d,m,D}} \biggl| \sum_{\nu \in N_k^{d,m,D}} (F(x) -
(S_{k, \nu_k^D(\nu)}^{l -1,d} F)(x))
g_{k, \nu}^{m,d}(x)\biggr|^p dx \le \\
\sum_{\substack{n \in \Z^d: \\ Q_{k,n}^d \cap G_k^{d,m,D} \ne \emptyset}}
\biggl( \sum_{\substack{\nu \in N_k^{d,m,D}: \\
\supp g_{k, \nu}^{m,d} \cap Q_{k,n}^d \ne \emptyset}}
c_{14} 2^{kd /p} \biggl(\int_{ c_{10} 2^{-k} B^d}
\int_{ (D_{k,n}^{d,m,D})_{l \xi}}
|\Delta_\xi^l F(x)|^p dx d\xi\biggr)^{1/p}\biggr)^p \le \\
\sum_{\substack{n \in \Z^d: \\ Q_{k,n}^d \cap G_k^{d,m,D} \ne \emptyset}}
c_{15}^p 2^{kd} \int_{ c_{10} 2^{-k} B^d}
\int_{D_{k,n}^{d,m,D}} |\Delta_\xi^l F(x)|^p dx d\xi \le \\
c_{15}^p 2^{kd} \int_{ c_{10} 2^{-k} B^d}
\int_{ \R^d} \biggl(\sum_{n \in \Z^d: Q_{k,n}^d \cap G_k^{d,m,D} \ne
\emptyset} \chi_{D_{k,n}^{d,m,D}}(x)\biggr)
|\Delta_\xi^l F(x)|^p dx d\xi \\
\le \biggl(c_{16} 2^{kd /p} \biggl(\int_{c_{10} 2^{-k} B^d}
\int_{ \R^d} |\Delta_\xi^l F(x)|^p dx d\xi\biggr)^{1/p}\biggr)^p.
\end{multline*}
Объединяя последнее неравенство с (1.4.28), для $ f \in L_p(D) $ при
$ k \in \Z_+: k \ge k^0, $ приходим к соотношению
\begin{multline*}
\| f -((E_{k^0}^{l -1,d,m,D} f) \mid_D
+\sum_{\k = k^0 +1}^k (E_{\k}^{l -1,d,m,D} f
-E_{\k -1}^{l -1,d,m,D} f) \mid_D) \|_{L_p(D)} \le \\
c_{16} 2^{kd /p} \biggl(\int_{c_{10} 2^{-k} B^d}
\int_{ \R^d} |\Delta_\xi^l F(x)|^p dx d\xi \biggr)^{1/p} \le \\
c_{17} \Omega^l(\mathcal I^D f,
c_{10} 2^{-k})_{L_p(\R^d)} \to 0 \text{ при } k \to \infty,
\end{multline*}
что влечёт (1.4.11). $ \square $

Как известно, имеет место следующая лемма.

Лемма 1.4.7

Пусть $ d \in \N, \lambda \in \Z_+^d, D $ --- область в $ \R^d $ и
функция $ f \in C^\infty(D), $ а $ g \in L_1(D), $ причём
для каждого $ \mu \in \Z_+^d(\lambda) $ обобщённая производная
$ \D^\mu g \in L_1(D). $ Тогда в пространстве обобщённых функций в
области $ D $ имеет место равенство
\begin{equation*} \tag{1.4.33}
\D^\lambda (fg) = \sum_{ \mu \in \Z_+^d(\lambda)} C_\lambda^\mu
\D^{\lambda -\mu} f \D^\mu g.
\end{equation*}

Лемма 1.4.8

Пусть $ d \in \N, l \in \Z_+, m \in \N, \lambda \in \Z_{+ m}^d,
1 \le s, q \le \infty $ и область $ D \subset \R^d $ удовлетворяет
условиям предложения 1.4.6. Тогда существует константа
$ c_{18}(l,d,m,D,\lambda,s,q) > 0 $ такая, что для любой функции
$ f \in L_s(D) $ при $ k \in \Z_+: k \ge k^0(d,D), $ справедливо
неравенство
\begin{equation*} \tag{1.4.34}
\| \D^\lambda E_k^{l,d,m,D} f \|_{L_q(\R^d)} \le
c_{18} 2^{k (|\lambda| +(d /s -d /q)_+)} \| f\|_{L_s(D)}.
\end{equation*}

Доказательство.

При доказательстве леммы будем использовать объекты и связанные с ними
факты из доказательства предложения 1.4.6. В условиях леммы для $ f \in
L_s(D) $ при $ s \le q, k \in \Z_+: k \ge k^0, $ принимая во внимание (1.4.33),
имеем
\begin{multline*} \tag{1.4.35}
\| \D^\lambda (E_k^{l,d,m,D} f) \|_{L_q(\R^d)} = \\
\biggl\| \D^\lambda (\sum_{\nu \in N_k^{d,m,D}}
(S_{k,\nu_k^D(\nu)}^{l,d} f) g_{k, \nu}^{m,d}) \biggr\|_{L_q(\R^d)} =
\biggl\| \sum_{\nu \in N_k^{d,m,D}} \D^\lambda
((S_{k,\nu_k^D(\nu)}^{l,d} f) g_{k, \nu}^{m,d}) \biggr\|_{L_q(\R^d)} = \\
\biggl\| \sum_{\nu \in N_k^{d,m,D}} \sum_{ \mu \in \Z_+^d(\lambda)}
C_\lambda^\mu (\D^\mu S_{k,\nu_k^D(\nu)}^{l,d} f)
\D^{\lambda -\mu} g_{k, \nu}^{m,d} \biggr\|_{L_q(\R^d)} = \\
\biggl\| \sum_{ \mu \in \Z_+^d(\lambda)} C_\lambda^\mu
\sum_{\nu \in N_k^{d,m,D}}
(\D^\mu S_{k,\nu_k^D(\nu)}^{l,d} f)
\D^{\lambda -\mu} g_{k, \nu}^{m,d} \biggr\|_{L_q(\R^d)} \le \\
\sum_{ \mu \in \Z_+^d(\lambda)} C_\lambda^\mu
\biggl\| \sum_{\nu \in N_k^{d,m,D}}
(\D^\mu S_{k,\nu_k^D(\nu)}^{l,d} f)
\D^{\lambda -\mu} g_{k, \nu}^{m,d} \biggr\|_{L_q(\R^d)}.
\end{multline*}

Оценивая правую часть (1.4.35), при $ \mu \in \Z_+^d(\lambda) $ с учётом
(1.4.14), (1.4.12) получаем (ср. с (1.4.29))
\begin{multline*} \tag{1.4.36}
\biggl\| \sum_{\nu \in N_k^{d,m,D}}
(\D^\mu S_{k,\nu_k^D(\nu)}^{l,d} f)
\D^{\lambda -\mu} g_{k, \nu}^{m,d} \biggr\|_{L_q(\R^d)}^q = \\
\int_{\R^d} \biggl| \sum_{\nu \in N_k^{d,m,D}}
(\D^\mu S_{k,\nu_k^D(\nu)}^{l,d} f)
\D^{\lambda -\mu} g_{k, \nu}^{m,d} \biggr|^q dx = \\
\int_{G_k^{d,m,D}} \biggl| \sum_{\nu \in N_k^{d,m,D}}
(\D^\mu S_{k,\nu_k^D(\nu)}^{l,d} f)
\D^{\lambda -\mu} g_{k, \nu}^{m,d} \biggr|^q dx = \\
\sum_{n \in \Z^d: Q_{k,n}^d \cap G_k^{d,m,D} \ne \emptyset}
\int_{Q_{k,n}^d} \biggl| \sum_{\nu \in N_k^{d,m,D}}
(\D^\mu S_{k,\nu_k^D(\nu)}^{l,d} f)
\D^{\lambda -\mu} g_{k, \nu}^{m,d} \biggr|^q dx = \\
\sum_{n \in \Z^d: Q_{k,n}^d \cap G_k^{d,m,D} \ne \emptyset}
\int_{Q_{k,n}^d} \biggl| \sum_{\nu \in N_k^{d,m,D}: Q_{k,n}^d \cap
\supp g_{k, \nu}^{m,d} \ne \emptyset}
(\D^\mu S_{k,\nu_k^D(\nu)}^{l,d} f)
\D^{\lambda -\mu} g_{k, \nu}^{m,d} \biggr|^q dx = \\
\sum_{n \in \Z^d: Q_{k,n}^d \cap G_k^{d,m,D} \ne \emptyset}
\biggl\| \sum_{\nu \in N_k^{d,m,D}: Q_{k,n}^d \cap
\supp g_{k, \nu}^{m,d} \ne \emptyset}
(\D^\mu S_{k,\nu_k^D(\nu)}^{l,d} f)
\D^{\lambda -\mu} g_{k, \nu}^{m,d} \biggr\|_{L_q(Q_{k,n}^d)}^q \le \\
\sum_{n \in \Z^d: Q_{k,n}^d \cap G_k^{d,m,D} \ne \emptyset}
\biggl(\sum_{\nu \in N_k^{d,m,D}: Q_{k,n}^d \cap
\supp g_{k, \nu}^{m,d} \ne \emptyset}
\| (\D^\mu S_{k,\nu_k^D(\nu)}^{l,d} f)
\D^{\lambda -\mu} g_{k, \nu}^{m,d} \|_{L_q(Q_{k,n}^d)}\biggr)^q.
\end{multline*}

Для оценки правой части (1.4.36), используя сначала (1.3.10), а затем
применяя (1.4.1), при $ n \in \Z^d: Q_{k,n}^d \cap G_k^{d,m,D} \ne
\emptyset, \ \nu \in N_k^{d,m,D}: Q_{k,n}^d \cap
\supp g_{k, \nu}^{m,d} \ne \emptyset, $ выводим
\begin{multline*} \tag{1.4.37}
\| (\D^\mu S_{k,\nu_k^D(\nu)}^{l,d} f)
\D^{\lambda -\mu} g_{k, \nu}^{m,d} \|_{L_q(Q_{k,n}^d)} \le
\| \D^{\lambda -\mu} g_{k, \nu}^{m,d} \|_{L_\infty(\R^d)}
\| \D^\mu S_{k,\nu_k^D(\nu)}^{l,d} f \|_{L_q(Q_{k,n}^d)} = \\
c_{19} 2^{k |\lambda -\mu|}
\| \D^\mu S_{k,\nu_k^D(\nu)}^{l,d} f \|_{L_q(Q_{k,n}^d)} \le
c_{20} 2^{k |\lambda -\mu|} 2^{k(|\mu| +d /s -d /q)}
\| S_{k,\nu_k^D(\nu)}^{l,d} f\|_{L_s(Q_{k,n}^d)} = \\
c_{20} 2^{k(|\lambda| +d /s -d /q)}
\| S_{k,\nu_k^D(\nu)}^{l,d} f\|_{L_s(Q_{k,n}^d)}.
\end{multline*}

Далее, принимая во внимание, что в силу (1.4.21) имеет место включение
\begin{multline*}
Q_{k,n}^d \subset (2^{-k} \nu_k^D(\nu) +(\gamma^1 +1) 2^{-k} B^d), \\
n \in \Z^d: Q_{k,n}^d \cap G_k^{d,m,D} \ne \emptyset, \
\nu \in N_k^{d,m,D}: \supp g_{k, \nu}^{m,d} \cap Q_{k,n}^d
\ne \emptyset,
\end{multline*}
на основании (1.4.1), (1.4.3), (1.4.26) заключаем, что для
$ n \in \Z^d: Q_{k,n}^d \cap G_k^{d,m,D} \ne \emptyset, \
\nu \in N_k^{d,m,D}: \supp g_{k, \nu}^{m,d} \cap Q_{k,n}^d
\ne \emptyset, $ выполняется неравенство
\begin{multline*} \tag{1.4.38}
\| S_{k,\nu_k^D(\nu)}^{l,d} f\|_{L_s(Q_{k,n}^d)} \le
c_{21} \| S_{k,\nu_k^D(\nu)}^{l,d} f \|_{L_s(Q_{k, \nu_k^D(\nu)}^d)} \le \\
c_{22} \| f \|_{L_s(Q_{k, \nu_k^D(\nu)}^d)} \le c_{22}
\| f\|_{L_s(D \cap D_{k, n}^{d,m,D})}.
\end{multline*}

Объединяя (1.4.37) и (1.4.38), находим, что при $ n \in \Z^d: Q_{k,n}^d
\cap G_k^{d,m,D} \ne \emptyset, \
\nu \in N_k^{d,m,D}: Q_{k,n}^d \cap \supp g_{k, \nu}^{m,d} \ne
\emptyset, $ имеет место неравенство
$$
\| (\D^\mu S_{k,\nu_k^D(\nu)}^{l,d} f)
\D^{\lambda -\mu} g_{k, \nu}^{m,d} \|_{L_q(Q_{k,n}^d)} \le
c_{23} 2^{k(|\lambda| +d /s -d /q)}
\| f\|_{L_s(D \cap D_{k,n}^{d,m,D})}.
$$

Подставляя эту оценку в (1.4.36) и применяя (1.4.27), а затем используя
неравенство (1.1.1) при $ a = s /q \le 1 $ и оценку (1.4.25), приходим к
неравенству
\begin{multline*} \tag{1.4.39}
\biggl\| \sum_{\nu \in N_k^{d,m,D}}
(\D^\mu S_{k,\nu_k^D(\nu)}^{l,d} f)
\D^{\lambda -\mu} g_{k, \nu}^{m,d} \biggr\|_{L_q(\R^d)}^q \le \\
\sum_{n \in \Z^d: Q_{k,n}^d \cap G_k^{d,m,D} \ne \emptyset}
\biggl(\sum_{\substack{\nu \in N_k^{d,m,D}: \\ Q_{k,n}^d \cap
\supp g_{k, \nu}^{m,d} \ne \emptyset}}
c_{23} 2^{k(|\lambda| +d /s -d /q)}
\| f\|_{L_s(D \cap D_{k,n}^{d,m,D})}\biggr)^q \le \\
(c_{23} 2^{k(|\lambda| +d /s -d /q)})^q
\sum_{n \in \Z^d: Q_{k,n}^d \cap G_k^{d,m,D} \ne \emptyset}
\biggl( c_8 \| f\|_{L_s(D \cap D_{k,n}^{d,m,D})} \biggr)^q = \\
(c_{24} 2^{k(|\lambda| +d /s -d /q)})^q
\sum_{n \in \Z^d: Q_{k,n}^d \cap G_k^{d,m,D} \ne \emptyset}
\biggl( \int_{D \cap D_{k,n}^{d,m,D}} | f(x)|^s dx\biggr)^{q /s} \le \\
(c_{24} 2^{k(|\lambda| +d /s -d /q)})^q
\biggl(\sum_{n \in \Z^d: Q_{k,n}^d \cap G_k^{d,m,D} \ne \emptyset}
\int_{D \cap D_{k,n}^{d,m,D}} | f(x)|^s dx \biggr)^{q /s} = \\
(c_{24} 2^{k(|\lambda| +d /s -d /q)})^q
\biggl(\sum_{n \in \Z^d: Q_{k,n}^d \cap G_k^{d,m,D} \ne \emptyset}
\int_D \chi_{ D_{k,n}^{d,m,D}}(x) | f(x)|^s dx \biggr)^{q /s} = \\
(c_{24} 2^{k(|\lambda| +d /s -d /q)})^q
\biggl(\int_D \biggl(\sum_{n \in \Z^d: Q_{k,n}^d \cap G_k^{d,m,D} \ne \emptyset}
\chi_{ D_{k,n}^{d,m,D}}(x)\biggr)  | f(x)|^s dx\biggr)^{q /s} \le \\
(c_{24} 2^{k(|\lambda| +d /s -d /q)})^q
\biggl(\int_D c_7 | f(x)|^s dx\biggr)^{q /s} =
(c_{25} 2^{k(|\lambda| +d /s -d /q)} \| f\|_{L_s(D)})^q.
\end{multline*}

Соединяя (1.4.35) с (1.4.39). получаем
\begin{multline*}
\| \D^\lambda E_k^{l,d,m,D} f \|_{L_q(\R^d)} \le
\sum_{ \mu \in \Z_+^d(\lambda)} C_\lambda^\mu
c_{25} 2^{k(|\lambda| +d /s -d /q)} \| f\|_{L_s(D)} = \\
c_{18} 2^{k(|\lambda| +d /s -d /q)} \| f\|_{L_s(D)},
\end{multline*}
что совпадает с (1.4.34) при $ s \le q. $ Справедливость (1.4.34) при $ q < s $
вытекает из соблюдения (1.4.34) при $ q = s $ и неравенства
\begin{multline*}
\| \D^\lambda E_k^{l,d,m,D} f \|_{L_q(\R^d)} =
\| \D^\lambda E_k^{l,d,m,D} f \|_{L_q(G)} \le \\
(\mes G)^{1 /q -1 /s} \| \D^\lambda E_k^{l,d,m,D} f \|_{L_s(G)} \le
(\mes G)^{1 /q -1 /s} \| \D^\lambda E_k^{l,d,m,D} f \|_{L_s(\R^d)},
\end{multline*}
где $ G $ -- ограниченная область в $ \R^d $ такая, что $ G_k^{d,m,D}
\subset G $ при $ k \in \Z_+: k \ge k^0. \square $

Отметим, что в лемме 1.4.8 при $ s \le q $ и в предложении 1.4.6 условие
ограниченности области $ D $ можно опустить.

Введём следующее обозначение. Пусть $ d \in \N, \ D $ -- ограниченная область в
$ \R^d $ и число $ k \in \N $ таковы, что существует $ \nu^\prime \in \Z^d, $
для которого $ \overline Q_{k -1, \nu^\prime}^d \subset D, $ а также пусть
$ l \in \Z_+, \ m \in \N. $ Тогда положим
\begin{equation*} \tag{1.4.40}
\mathcal E_k^{l,d,m,D} = E_k^{l,d,m,D} -
H_{k -1}^{l,d,m,D} E_{k -1}^{l,d,m,D} \ (\text{см. } (1.4.9), (1.3.12)).
\end{equation*}

Предложение 1.4.9

Пусть $ d \in \N, \ l \in \N, \ m \in \N, \ \lambda \in \Z_{+ m}^d, \
D \subset \R^d $ -- ограниченная область $ \e $-типа, $ 1 < p < \infty, \
1 \le q \le \infty. $
Тогда существует константа $ c_{26}(d,l,m,D,\lambda,p,q) > 0 $
такая, что при $ k \in \Z_+: k > K^0(d,D,\e) $
(см. определение 1 при $ \alpha = \e $), для $ f \in L_p(D) $ соблюдается
неравенство
\begin{multline*} \tag{1.4.41}
\| \D^\lambda \mathcal E_k^{l -1,d,m,D} f \|_{L_q(\R^d)}
\le c_{26} 2^{k (| \lambda| +(d /p -d /q)_+)}
\Omega^{\prime l}(f, 2^{-k +1})_{L_p(D)}.
\end{multline*}

Доказательство.

В условиях предложения 1.4.9 пусть $ f \in L_p(D), \ k \in \Z_+: k > K^0 $
(см. определение 1 при $ \alpha = \e $), $ p \le q. $
Тогда, принимая во внимание (1.4.40), (1.3.11), (1.3.16), (1.3.17), (1.4.33),
имеем
\begin{multline*} \tag{1.4.42}
\| \D^\lambda \mathcal E_k^{l -1,d,m,D} f \|_{L_q(\R^d)} =
\| \D^\lambda (E_k^{l -1,d,m,D} f -
H_{k -1}^{l -1,d,m,D} (E_{k -1}^{l -1,d,m,D} f)) \|_{L_q(\R^d)} = \\
\biggl\| \D^\lambda \biggl(\sum_{\nu \in N_k^{d,m,D}}
(S_{k, \nu_k^D(\nu)}^{l -1,d} f) g_{k, \nu}^{m,d} -
H_{k -1}^{l -1,d,m,D} \biggl(\sum_{\nu^\prime \in N_{k -1}^{d,m,D}}
(S_{k -1, \nu_{k -1}^D(\nu^\prime)}^{l -1,d} f)
g_{k -1, \nu^\prime}^{m,d}\biggr)\biggr)  \biggr\|_{L_q(\R^d)} = \\
\biggl\| \D^\lambda \biggl(\sum_{\nu \in N_k^{d,m,D}}
(S_{k, \nu_k^D(\nu)}^{l -1,d} f) g_{k, \nu}^{m,d} -
\sum_{\nu \in N_k^{d,m,D}}
\biggl(\sum_{\m \in \M^{m,d}(\nu)} A_{\m}^{m,d}
S_{k -1, \nu_{k -1}^D(\n(\nu,\m))}^{l -1,d} f\biggr)
g_{k, \nu}^{m,d}\biggr) \biggr\|_{L_q(\R^d)} = \\
\biggl\| \D^\lambda \biggl(\sum_{\nu \in N_k^{d,m,D}}
(\sum_{\m \in \M^{m,d}(\nu)} A_{\m}^{m,d})
(S_{k, \nu_k^D(\nu)}^{l -1,d} f) g_{k, \nu}^{m,d} - \\
\sum_{\nu \in N_k^{d,m,D}}
(\sum_{\m \in \M^{m,d}(\nu)} A_{\m}^{m,d}
S_{k -1, \nu_{k -1}^D(\n(\nu,\m))}^{l -1,d} f)
g_{k, \nu}^{m,d}\biggr) \biggr\|_{L_q(\R^d)} = \\
\biggl\| \D^\lambda \biggl(\sum_{\nu \in N_k^{d,m,D}}
\biggl(\biggl(\sum_{\m \in \M^{m,d}(\nu)} A_{\m}^{m,d}
S_{k, \nu_k^D(\nu)}^{l -1,d} f\biggr) -
\biggl(\sum_{\m \in \M^{m,d}(\nu)} A_{\m}^{m,d}
S_{k -1, \nu_{k -1}^D(\n(\nu,\m))}^{l -1,d} f\biggr)\biggr)
g_{k, \nu}^{m,d}\biggr) \biggr\|_{L_q(\R^d)} = \\
\biggl\| \D^\lambda \biggl(\sum_{\nu \in N_k^{d,m,D}}
\biggl(\sum_{\m \in \M^{m,d}(\nu)} A_{\m}^{m,d}
((S_{k, \nu_k^D(\nu)}^{l -1,d} f) -
(S_{k -1, \nu_{k -1}^D(\n(\nu,\m))}^{l -1,d} f))\biggr)
g_{k, \nu}^{m,d}\biggr) \biggr\|_{L_q(\R^d)} = \\
\biggl\| \sum_{\nu \in N_k^{d,m,D}} \D^\lambda
\biggl(\biggl(\sum_{\m \in \M^{m,d}(\nu)} A_{\m}^{m,d}
((S_{k, \nu_k^D(\nu)}^{l -1,d} f) -
(S_{k -1, \nu_{k -1}^D(\n(\nu,\m))}^{l -1,d} f))\biggr)
g_{k, \nu}^{m,d}\biggr) \biggr\|_{L_q(\R^d)} = \\
\biggl\| \sum_{\nu \in N_k^{d,m,D}} \sum_{ \mu \in \Z_+^d(\lambda)}
C_\lambda^\mu \D^\mu \biggl(\sum_{\m \in \M^{m,d}(\nu)} A_{\m}^{m,d}
((S_{k, \nu_k^D(\nu)}^{l -1,d} f) -
(S_{k -1, \nu_{k -1}^D(\n(\nu,\m))}^{l -1,d} f))\biggr)
\D^{\lambda -\mu} g_{k, \nu}^{m,d} \biggr\|_{L_q(\R^d)} = \\
\biggl\| \sum_{ \mu \in \Z_+^d(\lambda)} C_\lambda^\mu
\sum_{\nu \in N_k^{d,m,D}}
\D^\mu \biggl(\sum_{\m \in \M^{m,d}(\nu)} A_{\m}^{m,d}
((S_{k, \nu_k^D(\nu)}^{l -1,d} f) -
(S_{k -1, \nu_{k -1}^D(\n(\nu,\m))}^{l -1,d} f))\biggr)
\D^{\lambda -\mu} g_{k, \nu}^{m,d} \biggr\|_{L_q(\R^d)} = \\
\biggl\| \sum_{ \mu \in \Z_+^d(\lambda)} C_\lambda^\mu
\sum_{\nu \in N_k^{d,m,D}}
\biggl(\sum_{\m \in \M^{m,d}(\nu)} A_{\m}^{m,d}
\D^\mu (S_{k, \nu_k^D(\nu)}^{l -1,d} f -
S_{k -1, \nu_{k -1}^D(\n(\nu,\m))}^{l -1,d} f)\biggr)
\D^{\lambda -\mu} g_{k, \nu}^{m,d} \biggr\|_{L_q(\R^d)} \le \\
\sum_{ \mu \in \Z_+^d(\lambda)} C_\lambda^\mu
\biggl\| \sum_{\nu \in N_k^{d,m,D}}
\biggl(\sum_{\m \in \M^{m,d}(\nu)} A_{\m}^{m,d}
\D^\mu (S_{k, \nu_k^D(\nu)}^{l -1,d} f -
S_{k -1, \nu_{k -1}^D(\n(\nu,\m))}^{l -1,d} f)\biggr)
\D^{\lambda -\mu} g_{k, \nu}^{m,d} \biggr\|_{L_q(\R^d)}.
\end{multline*}

Оценивая правую часть (1.4.42), при $ \mu \in \Z_+^d(\lambda) $ с учётом
(1.4.14), (1.4.12) получаем
\begin{multline*} \tag{1.4.43}
\biggl\| \sum_{\nu \in N_k^{d,m,D}}
\biggl(\sum_{\m \in \M^{m,d}(\nu)} A_{\m}^{m,d}
\D^\mu (S_{k, \nu_k^D(\nu)}^{l -1,d} f -
S_{k -1, \nu_{k -1}^D(\n(\nu,\m))}^{l -1,d} f)\biggr)
\D^{\lambda -\mu} g_{k, \nu}^{m,d} \biggr\|_{L_q(\R^d)}^q = \\
\int_{\R^d} \biggl| \sum_{\nu \in N_k^{d,m,D}}
\biggl(\sum_{\m \in \M^{m,d}(\nu)} A_{\m}^{m,d}
\D^\mu (S_{k, \nu_k^D(\nu)}^{l -1,d} f -
S_{k -1, \nu_{k -1}^D(\n(\nu,\m))}^{l -1,d} f)\biggr)
\D^{\lambda -\mu} g_{k, \nu}^{m,d} \biggr|^q dx = \\
\int_{G_k^{d,m,D}} \biggl| \sum_{\nu \in N_k^{d,m,D}}
\biggl(\sum_{\m \in \M^{m,d}(\nu)} A_{\m}^{m,d}
\D^\mu (S_{k, \nu_k^D(\nu)}^{l -1,d} f -
S_{k -1, \nu_{k -1}^D(\n(\nu,\m))}^{l -1,d} f)\biggr)
\D^{\lambda -\mu} g_{k, \nu}^{m,d} \biggr|^q dx = \\
\sum_{n \in \Z^d: Q_{k,n}^d \cap G_k^{d,m,D} \ne \emptyset}
\int_{Q_{k,n}^d} \biggl| \sum_{\nu \in N_k^{d,m,D}}
\biggl(\sum_{\m \in \M^{m,d}(\nu)} A_{\m}^{m,d}
\D^\mu (S_{k, \nu_k^D(\nu)}^{l -1,d} f -\\
S_{k -1, \nu_{k -1}^D(\n(\nu,\m))}^{l -1,d} f)\biggr)
\D^{\lambda -\mu} g_{k, \nu}^{m,d} \biggr|^q dx = \\
\sum_{n \in \Z^d: Q_{k,n}^d \cap G_k^{d,m,D} \ne \emptyset}
\int_{Q_{k,n}^d} \biggl| \sum_{\substack{\nu \in N_k^{d,m,D}:\\ Q_{k,n}^d \cap
\supp g_{k, \nu}^{m,d} \ne \emptyset}}
\biggl(\sum_{\m \in \M^{m,d}(\nu)} A_{\m}^{m,d}
\D^\mu (S_{k, \nu_k^D(\nu)}^{l -1,d} f -\\
S_{k -1, \nu_{k -1}^D(\n(\nu,\m))}^{l -1,d} f)\biggr)
\D^{\lambda -\mu} g_{k, \nu}^{m,d} \biggr|^q dx = \\
\sum_{n \in \Z^d: Q_{k,n}^d \cap G_k^{d,m,D} \ne \emptyset}
\biggl\| \sum_{\substack{\nu \in N_k^{d,m,D}: \\ Q_{k,n}^d \cap
\supp g_{k, \nu}^{m,d} \ne \emptyset}}
\biggl(\sum_{\m \in \M^{m,d}(\nu)} A_{\m}^{m,d}
\D^\mu (S_{k, \nu_k^D(\nu)}^{l -1,d} f -\\
S_{k -1, \nu_{k -1}^D(\n(\nu,\m))}^{l -1,d} f)\biggr)
\D^{\lambda -\mu} g_{k, \nu}^{m,d} \biggr\|_{L_q(Q_{k,n}^d)}^q \le \\
\sum_{n \in \Z^d: Q_{k,n}^d \cap G_k^{d,m,D} \ne \emptyset}
\biggl(\sum_{\substack{\nu \in N_k^{d,m,D}: \\ Q_{k,n}^d \cap
\supp g_{k, \nu}^{m,d} \ne \emptyset}}
\biggl\|\biggl(\sum_{\m \in \M^{m,d}(\nu)} A_{\m}^{m,d}
\D^\mu (S_{k, \nu_k^D(\nu)}^{l -1,d} f -\\
S_{k -1, \nu_{k -1}^D(\n(\nu,\m))}^{l -1,d} f)\biggr)
\D^{\lambda -\mu} g_{k, \nu}^{m,d} \biggr\|_{L_q(Q_{k,n}^d)}\biggr)^q \le \\
\sum_{n \in \Z^d: Q_{k,n}^d \cap G_k^{d,m,D} \ne \emptyset}
\biggl(\sum_{\substack{\nu \in N_k^{d,m,D}: \\ Q_{k,n}^d \cap
\supp g_{k, \nu}^{m,d} \ne \emptyset}}
\sum_{\m \in \M^{m,d}(\nu)} A_{\m}^{m,d}
\biggl\| \D^\mu (S_{k, \nu_k^D(\nu)}^{l -1,d} f - \\
S_{k -1, \nu_{k -1}^D(\n(\nu,\m))}^{l -1,d} f)
\D^{\lambda -\mu} g_{k, \nu}^{m,d} \biggr\|_{L_q(Q_{k,n}^d)}\biggr)^q.
\end{multline*}

Для оценки правой части (1.4.43), используя сначала (1.3.10), а
затем применяя (1.4.1), при $ \mu \in \Z_+^d(\lambda), \ n \in \Z^d: Q_{k,n}^d \cap
G_k^{d,m,D} \ne \emptyset, \ \nu \in N_k^{d,m,D}: Q_{k,n}^d \cap
\supp g_{k, \nu}^{m,d} \ne \emptyset, \ \m \in \M^{m,d}(\nu) $
выводим
\begin{multline*} \tag{1.4.44}
\| \D^\mu (S_{k, \nu_k^D(\nu)}^{l -1,d} f -
S_{k -1, \nu_{k -1}^D(\n(\nu,\m))}^{l -1,d} f)
\D^{\lambda -\mu} g_{k, \nu}^{m,d} \|_{L_q(Q_{k,n}^d)} \le \\
\| \D^{\lambda -\mu} g_{k, \nu}^{m,d} \|_{L_\infty(\R^d)}
\| \D^\mu (S_{k, \nu_k^D(\nu)}^{l -1,d} f -
S_{k -1, \nu_{k -1}^D(\n(\nu,\m))}^{l -1,d} f)\|_{L_q(Q_{k,n}^d)} = \\
c_{27} 2^{k |\lambda -\mu|}
\| \D^\mu (S_{k, \nu_k^D(\nu)}^{l -1,d} f -
S_{k -1, \nu_{k -1}^D(\n(\nu,\m))}^{l -1,d} f)\|_{L_q(Q_{k,n}^d)} \le \\
c_{28} 2^{k |\lambda -\mu|} 2^{k (|\mu| +d /p -d /q)}
\| S_{k, \nu_k^D(\nu)}^{l -1,d} f -
S_{k -1, \nu_{k -1}^D(\n(\nu,\m))}^{l -1,d} f\|_{L_p(Q_{k,n}^d)} = \\
c_{28} 2^{k (|\lambda| +d /p -d /q)}
\| S_{k, \nu_k^D(\nu)}^{l -1,d} f -
S_{k -1, \nu_{k -1}^D(\n(\nu,\m))}^{l -1,d} f\|_{L_p(Q_{k,n}^d)}.
\end{multline*}

Оценим норму в правой части (1.4.44). Для этого, фиксировав $ n \in \Z^d:
Q_{k,n}^d \cap G_k^{d,m,D} \ne \emptyset, \ \nu \in N_k^{d,m,D}:
Q_{k,n}^d \cap \supp g_{k, \nu}^{m,d} \ne \emptyset, \
\m \in \M^{m,d}(\nu), $ принимая во внимание (1.4.8), заметим, что
\begin{multline*} \tag{1.4.45}
D \supset \overline Q_{k -1, \nu_{k -1}^D(\n(\nu,\m))}^d =
2^{-(k -1)} \nu_{k -1}^D(\n(\nu,\m)) +
2^{-(k -1)} \overline I^d =
2^{-k} 2 \nu_{k -1}^D(\n(\nu,\m)) +\\
2^{-(k -1)} \overline I^d \supset 2^{-k} 2 \nu_{k -1}^D(\n(\nu,\m)) +
2^{-k} \overline I^d = \overline Q_{k, 2 \nu_{k -1}^D(\n(\nu,\m))}^d.
\end{multline*}
Учитывая, что и $ \overline Q_{k, \nu_k^D(\nu)}^d \subset D, $ выберем
последовательности $ \nu^\iota \in \Z^d, \ \iota =0,\ldots,\Iota; \ j^\iota \in
\Nu_{1,d}^1, \ \epsilon^\iota \in \{-1, 1\}, \ \iota =0,\ldots,\Iota -1, $ для
которых соблюдаются соотношения (1.2.1), (1.2.2), (1.2.3) при
\begin{equation*} \tag{1.4.46}
\boldsymbol{\nu} = \nu_k^D(\nu), \\
\boldsymbol{\nu^\prime} = 2 \nu_{k -1}^D(\n(\nu,\m)), \\
\kappa = k \e.
\end{equation*}
Тогда
\begin{multline*} \tag{1.4.47}
\| S_{k, \nu_k^D(\nu)}^{l -1,d} f -
S_{k -1, \nu_{k -1}^D(\n(\nu,\m))}^{l -1,d} f\|_{L_p(Q_{k,n}^d)} = \\
\biggl\| S_{k, \nu_k^D(\nu)}^{l -1,d} f -
\sum_{\iota =0}^{\Iota}
(S_{\mathpzc k^0 +k, \nu^\iota}^{l -1,d} f -S_{\mathpzc k^0 +k, \nu^\iota}^{l -1,d} f) -
S_{k -1, \nu_{k -1}^D(\n(\nu,\m))}^{l -1,d} f\biggr\|_{L_p(Q_{k,n}^d)} = \\
\biggl\| S_{k, \nu_k^D(\nu)}^{l -1,d} f -
S_{\mathpzc k^0 +k, \nu^0}^{l -1,d} f +\\
\sum_{\iota =0}^{\Iota -1}
(S_{\mathpzc k^0 +k, \nu^\iota}^{l -1,d} f -
S_{\mathpzc k^0 +k, \nu^{\iota +1}}^{l -1,d} f) +
S_{\mathpzc k^0 +k, \nu^{\Iota}}^{l -1,d} f -
S_{k -1, \nu_{k -1}^D(\n(\nu,\m))}^{l -1,d} f\biggr\|_{L_p(Q_{k,n}^d)} \le \\
\| S_{k, \nu_k^D(\nu)}^{l -1,d} f -
S_{\mathpzc k^0 +k, \nu^0}^{l -1,d} f \|_{L_p(Q_{k,n}^d)}+\\
\sum_{\iota =0}^{\Iota -1} \| S_{\mathpzc k^0 +k, \nu^\iota}^{l -1,d} f -
S_{\mathpzc k^0 +k, \nu^{\iota +1}}^{l -1,d} f\|_{L_p(Q_{k,n}^d)} +
\| S_{\mathpzc k^0 +k, \nu^{\Iota}}^{l -1,d} f -
S_{k -1, \nu_{k -1}^D(\n(\nu,\m))}^{l -1,d} f\|_{L_p(Q_{k,n}^d)}.
\end{multline*}

Для проведения оценки слагаемых в правой части (1.4.47) отметим некоторые
полезные для нас факты.
При $ \iota =0,\ldots,\Iota, \ j =1,\ldots,d $ для $ x \in
\overline Q_{\mathpzc k^0 +k,\nu^\iota}^d $ ввиду (1.2.2), (1.4.46), (1.4.21), (1.2.1)
справедливо неравенство
\begin{multline*} \tag{1.4.48}
| x_j -2^{-k} n_j | \le | x_j -2^{-\mathpzc k^0 -k} (\nu^\iota)_j | +
| 2^{-\mathpzc k^0 -k} (\nu^\iota)_j -2^{-k} n_j | \le \\
 2^{-\mathpzc k^0 -k} +
| 2^{-\mathpzc k^0 -k} (\nu^\iota)_j -2^{-k} n_j | = \\
2^{-\mathpzc k^0 -k} +| 2^{-\mathpzc k^0 -k} (\nu^\iota)_j -
\sum_{i =0}^{\iota -1} (2^{-\mathpzc k^0 -k} (\nu^i)_j
-2^{-\mathpzc k^0 -k} (\nu^i)_j) -2^{-k} n_j | = \\
2^{-\mathpzc k^0 -k} +| \sum_{i =0}^{\iota -1} (2^{-\mathpzc k^0 -k} (\nu^{i +1})_j
-2^{-\mathpzc k^0 -k} (\nu^i)_j) +\\
2^{-\mathpzc k^0 -k} (\nu^0)_j
-2^{-k} (\nu_k^D(\nu))_j +2^{-k} (\nu_k^D(\nu))_j -2^{-k} n_j | \le \\
2^{-k} +\sum_{i =0}^{\iota -1} | 2^{-\mathpzc k^0 -k} (\nu^{i +1})_j
-2^{-\mathpzc k^0 -k} (\nu^i)_j| +\\
| 2^{-\mathpzc k^0 -k} (\nu^0)_j
-2^{-k} (\nu_k^D(\nu))_j| +| 2^{-k} (\nu_k^D(\nu))_j -2^{-k} n_j | \le \\
2^{-k} +\sum_{i =0}^{\iota -1} 2^{-k} \| \nu^{i +1} -\nu^i \|
+| 2^{-\mathpzc k^0 -k} (\nu^0)_j -2^{-k} (\nu_k^D(\nu))_j| +
| 2^{-k} (\nu_k^D(\nu))_j -2^{-k} n_j | = \\
2^{-k} +\sum_{i =0}^{\iota -1} 2^{-k} \| \epsilon^i e_{j^i} \|
+| 2^{-\mathpzc k^0 -k} (\nu^0)_j -2^{-k} (\nu_k^D(\nu))_j| +
| 2^{-k} (\nu_k^D(\nu))_j -2^{-k} n_j | = \\
2^{-k} +\sum_{i =0}^{\iota -1} 2^{-k}
+| 2^{-\mathpzc k^0 -k} (\nu^0)_j -2^{-k} (\nu_k^D(\nu))_j| +
| 2^{-k} (\nu_k^D(\nu))_j -2^{-k} n_j | = \\
2^{-k} (\iota +1) +| 2^{-\mathpzc k^0 -k} (\nu^0)_j -2^{-k} (\nu_k^D(\nu))_j| +
| 2^{-k} (\nu_k^D(\nu))_j -2^{-k} n_j | \le \\
2^{-k} (\Iota +1) +2^{-k} +\Gamma^1 2^{-k} \le \\
2^{-k} (c_{29} \| \nu_k^D(\nu) -2 \nu_{k -1}^D(\n(\nu,\m)) \| +2) +
\Gamma^1 2^{-k}.
\end{multline*}

Оценивая правую часть (1.4.48), имеем
\begin{multline*} \tag{1.4.49}
\| \nu_k^D(\nu) -2 \nu_{k -1}^D(\n(\nu,\m)) \| = \\
\| \nu_k^D(\nu) -\nu +\nu -2 \n(\nu,\m) +2 \n(\nu,\m)
-2 \nu_{k -1}^D(\n(\nu,\m)) \| \le \\
\| \nu_k^D(\nu) -\nu \| +\| \nu -
2 \n(\nu,\m) \| +\| 2 \n(\nu,\m)
-2 \nu_{k -1}^D(\n(\nu,\m)) \|.
\end{multline*}
Замечая, что для $ \nu \in N_k^{d,m,D}, \ \m \in \M^{m,d}(\nu) $
при $ j \in \Nu_{1,d}^1 $ ввиду (1.3.15), (1.3.14) справедливо
соотношение
\begin{equation*}
| \nu_j -(2 \n(\nu,\m))_j | = | \nu_j -2 (\nu_j -\m_j) /2| = | \m_j| \le (m +1) ,
\end{equation*}
получаем, что
\begin{equation*} \tag{1.4.50}
\| \nu -2 \n(\nu,\m) \| \le (m +1).
\end{equation*}
Учитывая, что для $ \nu \in N_k^{d,m,D}, \ \m \in \M^{m,d}(\nu) $
в силу (1.3.14), (1.3.15) и замечания после предложения 1.3.1
мультииндекс $ \n(\nu,\m) \in N_{k -1}^{d,m,D}, $ согласно
(1.4.20), находим, что
\begin{multline*} \tag{1.4.51}
\| 2 \n(\nu,\m) -2 \nu_{k -1}^D(\n(\nu,\m)) \| =
2 \| \n(\nu,\m) -\nu_{k -1}^D(\n(\nu,\m)) \| \le 2 c_{6}.
\end{multline*}
Используя для оценки правой части (1.4.49) неравенства (1.4.20), (1.4.50),
(1.4.51), выводим
\begin{multline*} \tag{1.4.52}
\| \nu_k^D(\nu) -2 \nu_{k -1}^D(\n(\nu,\m)) \| \le
c_{6} +(m +1) +2 c_{6} = \\
c_{30}(d,m,D), \ \nu \in N_k^{d,m,D}, \ \m \in \M^{m,d}(\nu).
\end{multline*}
Объединяя (1.4.48), (1.4.52). видим, что для $ \iota =0,\ldots,\Iota $ имеет место включение
\begin{equation*}
\overline Q_{\mathpzc k^0 +k,\nu^\iota}^d \subset (2^{-k} n +\Gamma^2 2^{-k} B^d)
\end{equation*}
с $ \Gamma^2 = (c_{29} c_{30} +2) +\Gamma^1, $
из которого с учётом (1.2.3) (при $ \kappa = k \e $) следует, что
\begin{equation*} \tag{1.4.53}
\overline Q_{\mathpzc k^0 +k,\nu^\iota}^d \subset
D \cap (2^{-k} n +\Gamma^2 2^{-k} B^d), \ \iota =0,\ldots,\Iota.
\end{equation*}

Заметим ещё, что для $ x \in
\overline Q_{k -1, \nu_{k -1}^D(\n(\nu,\m))}^d $
при $ j =1,\ldots,d, $ благодаря (1.4.52), (1.4.21), соблюдается
неравенство
\begin{multline*}
| x_j -2^{-k} n_j | = | x_j -2^{-(k -1)}
(\nu_{k -1}^D(\n(\nu,\m)))_j +2^{-(k -1)}
(\nu_{k -1}^D(\n(\nu,\m)))_j -2^{-k} n_j | \le \\
| x_j -2^{-(k -1)} (\nu_{k -1}^D(\n(\nu,\m)))_j |
+| 2^{-(k -1)} (\nu_{k -1}^D(\n(\nu,\m)))_j -
2^{-k} n_j | \le \\
2^{-(k -1)} +| 2^{-(k -1)}
(\nu_{k -1}^D(\n(\nu,\m)))_j -2^{-k} n_j | = \\
2 2^{-k} +| 2^{-k} 2(\nu_{k -1}^D(\n(\nu,\m)))_j -2^{-k} n_j | = \\
2 2^{-k} +2^{-k} | 2(\nu_{k -1}^D(\n(\nu,\m)))_j -n_j | \le \\
2 2^{-k} +2^{-k} | (2 \nu_{k -1}^D(\n(\nu,\m)))_j -(\nu_k^D(\nu))_j |
+2^{-k} | (\nu_k^D(\nu))_j -n_j | \le \\
2 2^{-k} +2^{-k} \| 2 \nu_{k -1}^D(\n(\nu,\m)) -\nu_k^D(\nu) \|
+| 2^{-k} (\nu_k^D(\nu))_j -2^{-k} n_j | \le \\
2 2^{-k} +c_{30} 2^{-k} +2^{-k} \Gamma^1 = 2^{-k} \Gamma^3
\end{multline*}
с $ \Gamma^3 = 2 +c_{30} +\Gamma^1 > \Gamma^1, $
т.е.
\begin{equation*} \tag{1.4.54}
Q_{k -1, \nu_{k -1}^D(\n(\nu,\m))}^d \subset
\overline Q_{k -1, \nu_{k -1}^D(\n(\nu,\m))}^d
\subset (2^{-k} n +\Gamma^3 2^{-k} B^d).
\end{equation*}

Итак, задавая число $ \Gamma^4 \in \R_+ $ соотношением
$ \Gamma^4 = \max(\Gamma^2, \Gamma^3), $ и обозначая
$$
D_{k,n}^{\prime d,m,D} = 2^{-k} n +\Gamma^4 2^{-k} B^d, \
n \in \Z^d: Q_{k,n}^d \cap G_k^{d,m,D} \ne \emptyset,
$$
в силу (1.4.53), (1.4.54), (1.4.45), (1.4.8), (1.4.21) имеем
\begin{equation*} \tag{1.4.55}
\overline Q_{\mathpzc k^0 +k,\nu^\iota}^d \cup \overline Q_{\mathpzc k^0 +k,\nu^{\iota +1}}^d \subset
D \cap D_{k,n}^{\prime d,m,D}, \ \iota =0,\ldots,\Iota -1.
\end{equation*}
\begin{equation*} \tag{1.4.56}
Q_{k -1, \nu_{k -1}^D(\n(\nu,\m))}^d \subset
D \cap D_{k,n}^{\prime d,m,D},
\end{equation*}
\begin{equation*} \tag{1.4.57}
Q_{k,\nu_k^D(\nu)}^d \subset D \cap D_{k,n}^{\prime d,m,D}, \
\nu \in N_k^{d,m,D}: Q_{k,n}^d \cap \supp g_{k, \nu}^{m,d} \ne
\emptyset, \ \m \in \M^{m,d}(\nu).
\end{equation*}

Из приведенных определений с учётом того, что $ \Gamma^4 > 1, $ видно, что
при $ n \in \Z^d: Q_{k,n}^d \cap G_k^{d,m,D} \ne \emptyset, $
справедливо включение
\begin{equation*} \tag{1.4.58}
Q_{k, n}^d \subset D_{k, n}^{\prime d,m,D},
\end{equation*}
а из (1.4.55) следует, что
\begin{equation*} \tag{1.4.59}
Q_{k,n}^d \subset (2^{-\mathpzc k^0 -k} \nu^{\iota} +(\Gamma^4 +1) 2^{-k} B^d), \
\iota =0,\ldots,\Iota.
\end{equation*}

Учитывая (1.4.58), нетрудно видеть, что существует константа
$ c_{31}(d,m,D) >0 $ такая, что для каждого $ x \in \R^d $ число
\begin{equation*} \tag{1.4.60}
\card \{ n \in \Z^d: Q_{k,n}^d \cap G_k^{d,m,D} \ne \emptyset, \
x \in D_{k,n}^{\prime d,m,D} \} \le c_{31}.
\end{equation*}

Отметим ещё, что вследствие (1.2.1), (1.4.46), (1.4.52) справедлива оценка
\begin{equation*} \tag{1.4.61}
\Iota \le c_{29} \| \nu_k^D(\nu) -
2 \nu_{k -1}^D(\n(\nu,\m)) \| \le c_{29} c_{30} = c_{32}(d,m,D).
\end{equation*}

Принимая во внимание, что при $ \iota =0,\ldots,\Iota -1 $ ввиду (1.2.3),
(1.2.2) для $ Q^\iota = \inter (\overline Q_{\mathpzc k^0 +k,\nu^\iota}^d \cup
\overline Q_{\mathpzc k^0 +k,\nu^{\iota +1}}^d) \subset D $ выполняется равенство
\begin{equation*}
Q^\iota = \begin{cases} 2^{-\mathpzc k^0 -k} \nu^\iota +2^{-\mathpzc k^0 -k} 2^{e_{j^\iota}} I^d,
\text{ при } \epsilon^\iota =1; \\
2^{-\mathpzc k^0 -k} \nu^{\iota +1} +2^{-\mathpzc k^0 -k} 2^{e_{j^\iota}} I^d,

\text{ при } \epsilon^\iota =-1,
\end{cases}
\end{equation*}
определим линейный оператор $ S^\iota: L_1(Q^\iota) \mapsto \mathcal P^{l -1,d}, $
полагая
$$
S^\iota = P_{\delta,x^0,\rho}^{l -1,d,0}
$$
при $ \rho = 2^{e_{j^\iota}}, \ \delta = 2^{-\mathpzc k^0 -k}, $
\begin{equation*}
x^0 = \begin{cases} 2^{-\mathpzc k^0 -k} \nu^\iota, \text{ при } \epsilon^\iota =1; \\
2^{-\mathpzc k^0 -k} \nu^{\iota +1}, \text{ при } \epsilon^\iota =-1 \ ,
\end{cases}
\iota =0,\ldots,\Iota -1.
\end{equation*}

Теперь проведём оценку слагаемых в правой части (1.4.47).
При $ \iota =0,\ldots,\Iota -1, $ благодаря (1.4.59) применяя (1.4.1),
а затем используя (1.4.2), (1.4.3), определение $ Q^\iota, $ и, наконец,
пользуясь (1.4.6) и учитывая (1.4.55), приходим к неравенству
\begin{multline*} \tag{1.4.62}
\| S_{\mathpzc k^0 +k, \nu^\iota}^{l -1,d} f -
S_{\mathpzc k^0 +k, \nu^{\iota +1}}^{l -1,d} f\|_{L_p(Q_{k,n}^d)} = \\
\| S_{\mathpzc k^0 +k, \nu^\iota}^{l -1,d} f -S^\iota f +S^\iota f -
S_{\mathpzc k^0 +k, \nu^{\iota +1}}^{l -1,d} f\|_{L_p(Q_{k,n}^d)} \le \\
\| S_{\mathpzc k^0 +k, \nu^\iota}^{l -1,d} f -S^\iota f \|_{L_p(Q_{k,n}^d)}
+\| S^\iota f -S_{\mathpzc k^0 +k, \nu^{\iota +1}}^{l -1,d} f\|_{L_p(Q_{k,n}^d)} \le \\
c_{33} \| S_{\mathpzc k^0 +k, \nu^\iota}^{l -1,d} f -S^\iota f \|_{L_p(Q_{\mathpzc k^0 +k,\nu^\iota}^d)}
+c_{33} \| S^\iota f -S_{\mathpzc k^0 +k, \nu^{\iota +1}}^{l -1,d} f
\|_{L_p(Q_{\mathpzc k^0 +k,\nu^{\iota +1}}^d)} = \\
c_{33} \| S_{\mathpzc k^0 +k, \nu^\iota}^{l -1,d} (f -S^\iota f)
\|_{L_p(Q_{\mathpzc k^0 +k,\nu^\iota}^d)}
+c_{33} \| S_{\mathpzc k^0 +k, \nu^{\iota +1}}^{l -1,d} (f -S^\iota f)
\|_{L_p(Q_{\mathpzc k^0 +k,\nu^{\iota +1}}^d)} \le \\
c_{34} \| f -S^\iota f \|_{L_p(Q_{\mathpzc k^0 +k,\nu^\iota}^d)}
+c_{34} \| f -S^\iota f\|_{L_p(Q_{\mathpzc k^0 +k,\nu^{\iota +1}}^d)} \le
c_{35} \| f -S^\iota f \|_{L_p(Q^\iota)} \le \\
c_{36} (2^{-\mathpzc k^0 -k})^{-d /p} \biggl(\int_{2^{-\mathpzc k^0 -k} B^d}
\int_{Q_{l \xi}^\iota}| \Delta_\xi^l f(x)|^p dx d\xi\biggr)^{1/p} \le \\
c_{37} 2^{kd /p} \biggl(\int_{2^{-k} B^d}
\int_{Q_{l \xi}^\iota}| \Delta_\xi^l f(x)|^p dx d\xi\biggr)^{1/p} \le \\
c_{37} 2^{kd /p} \biggl(\int_{2^{-k} B^d}
\int_{(D \cap D_{k,n}^{\prime d,m,D})_{l \xi}}
| \Delta_\xi^l f(x)|^p dx d\xi\biggr)^{1/p} \le \\
c_{37} 2^{kd /p} \biggl(\int_{2^{-k +1} B^d}
\int_{D_{l \xi} \cap D_{k,n}^{\prime d,m,D}}
| \Delta_\xi^l f(x)|^p dx d\xi\biggr)^{1/p}.
\end{multline*}

Для оценки последнего слагаемого в правой части (1.4.47), ввиду (1.4.59)
используя (1.4.1), затем пользуясь (1.4.2), (1.4.3), учитывая (1.2.2),
(1.4.46), (1.4.45), применяя (1.4.6) и принимая во внимание (1.4.56), получаем
\begin{multline*} \tag{1.4.63}
\| S_{\mathpzc k^0 +k, \nu^{\Iota}}^{l -1,d} f -
S_{k -1, \nu_{k -1}^D(\n(\nu,\m))}^{l -1,d} f\|_{L_p(Q_{k,n}^d)} \le \\
c_{33} \| S_{\mathpzc k^0 +k, \nu^{\Iota}}^{l -1,d} f -
S_{k -1, \nu_{k -1}^D(\n(\nu,\m))}^{l -1,d} f\|_{L_p(Q_{\mathpzc k^0 +k,\nu^{\Iota}}^d)} = \\
c_{33} \| S_{\mathpzc k^0 +k, \nu^{\Iota}}^{l -1,d} (f -
S_{k -1, \nu_{k -1}^D(\n(\nu,\m))}^{l -1,d} f)\|_{L_p(Q_{\mathpzc k^0 +k,\nu^{\Iota}}^d)} \le \\
c_{38} \| f -S_{k -1, \nu_{k -1}^D(\n(\nu,\m))}^{l -1,d} f
\|_{L_p(Q_{\mathpzc k^0 +k,\nu^{\Iota}}^d)} \le \\
c_{38} \| f -S_{k -1, \nu_{k -1}^D(\n(\nu,\m))}^{l -1,d} f
\|_{L_p(Q_{k -1, \nu_{k -1}^D(\n(\nu,\m))}^d)} \le \\
c_{39} (2^{-(k -1)})^{-d /p} \biggl(\int_{2^{-(k -1)} B^d}
\int_{(Q_{k -1, \nu_{k -1}^D(\n(\nu,\m))}^d)_{l \xi}}
| \Delta_\xi^l f(x)|^p dx d\xi\biggr)^{1/p} \le \\
c_{39} 2^{kd /p} \biggl(\int_{ 2^{-k +1} B^d}
\int_{(D \cap D_{k,n}^{\prime d,m,D})_{l \xi}}
| \Delta_\xi^l f(x)|^p dx d\xi\biggr)^{1/p} \le \\
c_{39} 2^{kd /p} \biggl(\int_{2^{-k +1} B^d}
\int_{D_{l \xi} \cap D_{k,n}^{\prime d,m,D}}
| \Delta_\xi^l f(x)|^p dx d\xi\biggr)^{1/p}.
\end{multline*}

При оценке первого слагаемого в правой части (1.4.47), благодаря (1.4.59)
используя (1.4.1), затем пользуясь (1.4.2), (1.4.3), учитывая (1.2.2),
(1.4.46), применяя (1.4.6) и принимая во внимание (1.4.57), выводим
\begin{multline*} \tag{1.4.64}
\| S_{k, \nu_k^D(\nu)}^{l -1,d} f -
S_{\mathpzc k^0 +k, \nu^0}^{l -1,d} f \|_{L_p(Q_{k,n}^d)} \le \\
c_{33} \| S_{k, \nu_k^D(\nu)}^{l -1,d} f -
S_{\mathpzc k^0 +k, \nu^0}^{l -1,d} f \|_{L_p(Q_{\mathpzc k^0 +k,\nu^0}^d)} = \\
c_{33} \| S_{\mathpzc k^0 +k, \nu^0}^{l -1,d} (f -
S_{k, \nu_k^D(\nu)}^{l -1,d} f)
\|_{L_p(Q_{\mathpzc k^0 +k,\nu^0}^d)} \le \\
c_{40} \| f -S_{k, \nu_k^D(\nu)}^{l -1,d} f
\|_{L_p(Q_{\mathpzc k^0 +k,\nu^0}^d)} \le \\
c_{40} \| f -S_{k,\nu_k^D(\nu)}^{l -1,d} f
\|_{L_p(Q_{k, \nu_k^D(\nu)}^d)} \le \\
c_{41} (2^{-k})^{-d /p} \biggl(\int_{2^{-k} B^d}
\int_{(Q_{k, \nu_k^D(\nu)}^d)_{l \xi}}
| \Delta_\xi^l f(x)|^p dx d\xi\biggr)^{1/p} \le \\
c_{41} 2^{kd /p} \biggl(\int_{ 2^{-k +1} B^d}
\int_{(D \cap D_{k,n}^{\prime d,m,D})_{l \xi}}
| \Delta_\xi^l f(x)|^p dx d\xi\biggr)^{1/p} \le \\
c_{41} 2^{kd /p} \biggl(\int_{2^{-k +1} B^d}
\int_{D_{l \xi} \cap D_{k,n}^{\prime d,m,D}}
| \Delta_\xi^l f(x)|^p dx d\xi\biggr)^{1/p}.
\end{multline*}

Объединяя (1.4.47), (1.4.62), (1.4.63), (1.4.64) и учитывая (1.4.61), при
$ n \in \Z^d: Q_{k,n}^d \cap G_k^{d,m,D} \ne \emptyset, \ \nu \in N_k^{d,m,D}:
Q_{k,n}^d \cap \supp g_{k, \nu}^{m,d} \ne \emptyset, \
\m \in \M^{m,d}(\nu) $ имеем
\begin{multline*} \tag{1.4.65}
\| S_{k, \nu_k^D(\nu)}^{l -1,d} f -
S_{k -1, \nu_{k -1}^D(\n(\nu,\m))}^{l -1,d} f\|_{L_p(Q_{k,n}^d)} \le \\
c_{41} 2^{kd /p} \biggl(\int_{2^{-k +1} B^d}
\int_{D_{l \xi} \cap D_{k,n}^{\prime d,m,D}}
| \Delta_\xi^l f(x)|^p dx d\xi\biggr)^{1/p} +\\
\Iota c_{37} 2^{kd /p} \biggl(\int_{2^{-k +1} B^d}
\int_{D_{l \xi} \cap D_{k,n}^{\prime d,m,D}}
| \Delta_\xi^l f(x)|^p dx d\xi\biggr)^{1/p} +\\
c_{39} 2^{kd /p} \biggl(\int_{2^{-k +1} B^d}
\int_{D_{l \xi} \cap D_{k,n}^{\prime d,m,D}}
| \Delta_\xi^l f(x)|^p dx d\xi\biggr)^{1/p} \le \\
c_{41} 2^{kd /p} \biggl(\int_{2^{-k +1} B^d}
\int_{D_{l \xi} \cap D_{k,n}^{\prime d,m,D}}
| \Delta_\xi^l f(x)|^p dx d\xi\biggr)^{1/p} +\\
c_{32} c_{37} 2^{kd /p} \biggl(\int_{2^{-k +1} B^d}
\int_{D_{l \xi} \cap D_{k,n}^{\prime d,m,D}}
| \Delta_\xi^l f(x)|^p dx d\xi\biggr)^{1/p} +\\
c_{39} 2^{kd /p} \biggl(\int_{2^{-k +1} B^d}
\int_{D_{l \xi} \cap D_{k,n}^{\prime d,m,D}}
| \Delta_\xi^l f(x)|^p dx d\xi\biggr)^{1/p} = \\
c_{42} 2^{kd /p} \biggl(\int_{2^{-k +1} B^d}
\int_{D_{l \xi} \cap D_{k,n}^{\prime d,m,D}}
| \Delta_\xi^l f(x)|^p dx d\xi\biggr)^{1/p}.
\end{multline*}

Из (1.4.44) и (1.4.65) находим, что при $ \mu \in \Z_+^d(\lambda), \ n \in \Z^d:
Q_{k,n}^d \cap G_k^{d,m,D} \ne \emptyset, \ \nu \in N_k^{d,m,D}:
Q_{k,n}^d \cap \supp g_{k, \nu}^{m,d} \ne \emptyset, \
\m \in \M^{m,d}(\nu) $ выполняется неравенство
\begin{multline*}
\| \D^\mu (S_{k, \nu_k^D(\nu)}^{l -1,d} f -
S_{k -1, \nu_{k -1}^D(\n(\nu,\m))}^{l -1,d} f)
\D^{\lambda -\mu} g_{k, \nu}^{m,d} \|_{L_q(Q_{k,n}^d)} \le \\
c_{43} 2^{k (|\lambda| +d /p -d /q)}
2^{kd /p} \biggl(\int_{2^{-k +1} B^d}
\int_{D_{l \xi} \cap D_{k,n}^{\prime d,m,D}}
| \Delta_\xi^l f(x)|^p dx d\xi\biggr)^{1/p}.
\end{multline*}

Подставляя эту оценку в (1.4.43) и учитывая (1.3.17), (1.4.27), а
затем применяя неравенство (1.1.1) при $ a = p/q \le 1, $ и, наконец,
принимая во внимание (1.4.60), при $ \mu \in \Z_+^d(\lambda) $ выводим
\begin{multline*} \tag{1.4.66}
\biggl\| \sum_{\nu \in N_k^{d,m,D}}
\biggl(\sum_{\m \in \M^{m,d}(\nu)} A_{\m}^{m,d}
\D^\mu (S_{k, \nu_k^D(\nu)}^{l -1,d} f -
S_{k -1, \nu_{k -1}^D(\n(\nu,\m))}^{l -1,d} f)\biggr)
\D^{\lambda -\mu} g_{k, \nu}^{m,d} \biggr\|_{L_q(\R^d)}^q \le \\
\sum_{n \in \Z^d: Q_{k,n}^d \cap G_k^{d,m,D} \ne \emptyset}
\biggl(\sum_{\substack{\nu \in N_k^{d,m,D}:\\ Q_{k,n}^d \cap
\supp g_{k, \nu}^{m,d} \ne \emptyset}}
\sum_{\m \in \M^{m,d}(\nu)} A_{\m}^{m,d}
c_{43} 2^{k (|\lambda| +d /p -d /q)} 2^{kd /p} \times \\
\biggl(\int_{2^{-k +1} B^d}
\int_{D_{l \xi} \cap D_{k,n}^{\prime d,m,D}}
| \Delta_\xi^l f(x)|^p dx d\xi\biggr)^{1/p}\biggr)^q \le \\
\sum_{n \in \Z^d: Q_{k,n}^d \cap G_k^{d,m,D} \ne \emptyset}
\biggl(c_{44} c_{43} 2^{k (|\lambda| +d /p -d /q)} 2^{kd /p} \times \\
\biggl(\int_{2^{-k +1} B^d}
\int_{D_{l \xi} \cap D_{k,n}^{\prime d,m,D}}
| \Delta_\xi^l f(x)|^p dx d\xi\biggr)^{1/p}\biggr)^q = \\
(c_{45} 2^{k (|\lambda| +d /p -d /q)} 2^{kd /p})^q
\sum_{n \in \Z^d: Q_{k,n}^d \cap G_k^{d,m,D} \ne \emptyset}
\biggl(\int_{2^{-k +1} B^d}
\int_{D_{l \xi} \cap D_{k,n}^{\prime d,m,D}}
| \Delta_\xi^l f(x)|^p dx d\xi\biggr)^{q/p} \le \\
(c_{45} 2^{k (|\lambda| +d /p -d /q)} 2^{kd /p})^q
\biggl(\sum_{n \in \Z^d: Q_{k,n}^d \cap G_k^{d,m,D} \ne \emptyset}
\int_{2^{-k +1} B^d}
\int_{D_{l \xi} \cap D_{k,n}^{\prime d,m,D}}
| \Delta_\xi^l f(x)|^p dx d\xi\biggr)^{q/p} = \\
(c_{45} 2^{k (|\lambda| +d /p -d /q)} 2^{kd /p})^q
\biggl(\int_{2^{-k +1} B^d}
\sum_{n \in \Z^d: Q_{k,n}^d \cap G_k^{d,m,D} \ne \emptyset}
\int_{D_{l \xi}} \chi_{D_{k,n}^{\prime d,m,D}}(x)
| \Delta_\xi^l f(x)|^p dx d\xi\biggr)^{q/p} = \\
(c_{45} 2^{k (|\lambda| +d /p -d /q)} 2^{kd /p})^q
\biggl(\int_{2^{-k +1} B^d} \int_{D_{l \xi}}
\biggl(\sum_{n \in \Z^d: Q_{k,n}^d \cap G_k^{d,m,D} \ne \emptyset}
\chi_{D_{k,n}^{\prime d,m,D}}(x)\biggr)
| \Delta_\xi^l f(x)|^p dx d\xi\biggr)^{q/p} \le \\
(c_{45} 2^{k (|\lambda| +d /p -d /q)} 2^{kd /p})^q
\biggl(\int_{2^{-k +1} B^d}
\int_{D_{l \xi}} c_{31}
| \Delta_\xi^l f(x)|^p dx d\xi\biggr)^{q/p} = \\
(c_{46} 2^{k (|\lambda| +d /p -d /q)})^q
\biggl(2^{kd /p}
(\int_{2^{-k +1} B^d}
\int_{D_{l \xi}} | \Delta_\xi^l f(x)|^p dx d\xi)^{1/p} \biggr)^q = \\
(c_{47} 2^{k (|\lambda| +d /p -d /q)}
\biggl( (2 2^{-k +1})^{-d}
\int_{2^{-k +1} B^d}
\int_{D_{l \xi}} | \Delta_\xi^l f(x)|^p dx d\xi \biggr)^{1/p})^q = \\
(c_{47} 2^{k (| \lambda| +(d /p -d /q)_+)}
\Omega^{\prime l}(f, 2^{-k +1})_{L_p(D)})^q.
\end{multline*}

Соединяя (1.4.42) с (1.4.66), приходим к (1.4.41) при $ p \le q. $
Для получения неравенства (1.4.41) при $ q < p, $ фиксировав ограниченную
область $ D^\prime $ в $ \R^d, $ для которой $ \supp h \subset D^\prime $
для $ h \in \mathcal P_k^{l -1,d,m,D} $ при $ k \in \Z_+, $ благодаря
неравенству Гёльдера и соблюдению (1.4.41) при $ q = p, $ имеем
\begin{multline*}
\| \D^\lambda \mathcal E_k^{l -1,d,m,D} f \|_{L_q(\R^d)} =
\| \D^\lambda \mathcal E_k^{l -1,d,m,D} f \|_{L_q(D^\prime)} \le \\
c(D^\prime) \| \D^\lambda \mathcal E_k^{l -1,d,m,D} f \|_{L_p(D^\prime)} =
c(D^\prime) \| \D^\lambda \mathcal E_k^{l -1,d,m,D} f \|_{L_p(\R^d)} \le \\
c_{26} 2^{k (|\lambda| +(d /p -d /q)_+)}
\Omega^{\prime l}(f,  2^{-k +1})_{L_p(D)}, \ f \in L_p(D). \square
\end{multline*}

Следствие

В условиях предложения 1.4.9 для $ f \in L_p(D) $ при $ k \in \Z_+: k > K^0(d,D,\e), $
имеет место неравенство
\begin{multline*} \tag{1.4.67}
\| \D^\lambda ((E_k^{l -1,d,m,D} f) \mid_D)
-\D^\lambda ((E_{k -1}^{l -1,d,m,D} f) \mid_D) \|_{L_q(D)} \\
\le c_{26} 2^{k (| \lambda| +(d /p -d /q)_+)}
\Omega^{\prime l}(f, 2^{-k +1})_{L_p(D)}.
\end{multline*}

В самом деле, при соблюдении условий предложения 1.4.9 для $ f \in L_p(D) $
при $ k \in \Z_+: k > K^0, $ согласно (1.3.13), (1.4.40), (1.4.41) имеем
\begin{multline*}
\| \D^\lambda ((E_k^{l -1,d,m,D} f) \mid_D)
-\D^\lambda ((E_{k -1}^{l -1,d,m,D} f) \mid_D) \|_{L_q(D)} = \\
\| \D^\lambda ((E_k^{l -1,d,m,D} f) \mid_D
-(E_{k -1}^{l -1,d,m,D} f) \mid_D) \|_{L_q(D)} = \\
\| \D^\lambda ((E_k^{l -1,d,m,D} f) \mid_D
-(H_{k -1}^{l -1,d,m,D}(E_{k -1}^{l -1,d,m,D} f)) \mid_D) \|_{L_q(D)} = \\
\| \D^\lambda ((E_k^{l -1,d,m,D} f
-H_{k -1}^{l -1,d,m,D}(E_{k -1}^{l -1,d,m,D}
f)) \mid_D) \|_{L_q(D)} = \\
\| \D^\lambda ((\mathcal E_k^{l -1,d,m,D} f) \mid_D) \|_{L_q(D)} =
\| (\D^\lambda (\mathcal E_k^{l -1,d,m,D} f)) \mid_D \|_{L_q(D)} \le \\
\| \D^\lambda \mathcal E_k^{l -1,d,m,D} f \|_{L_q(\R^d)}
\le c_{26} 2^{k (| \lambda| +(d /p -d /q)_+)}
\Omega^{\prime l}(f, 2^{-k +1})_{L_p(D)}.
\end{multline*}

Предложение 1.4.10

Пусть выполнены условия предложения 1.4.9. Тогда существует константа
$ c_{48}(d,l,m,D,\lambda,p,q) >0 $ такая, что
если для функции $ f \in L_p(D) $ функция
\begin{multline*} \tag{1.4.68}
K(f,t) = K_{p,q}^{d,l,m,D,\lambda}(f,t) = \\
t^{-|\lambda| -(d /p -d /q)_+ -1}
\Omega^{\prime l}(f, 2 t)_{L_p(D)} \in L_1(I),
\end{multline*}
то для любого $ k \in \Z_+: k \ge K^0, $ справедливо неравенство
\begin{equation*} \tag{1.4.69}
\|\D^\lambda f -\D^\lambda ((E_k^{l -1,d,m,D} f) \mid_D)\|_{L_q(D)} \le
c_{48} \int_0^{2^{-k}} K(f,t) dt.
\end{equation*}

Доказательство.

Прежде всего заметим, что в условиях предложения 1.4.10 в силу предложения
1.4.6 для $ f \in L_p(D) $ имеет место равенство (1.4.11) с константой $ K^0 $
вместо $ k^0 $ (см. определение 1 при $ \alpha = \e $).
Далее, пусть $ f \in L_p(D) $ и соблюдается условие (1.4.68).
Тогда согласно (1.4.67) при $ k \in \Z_+: k > K^0, $ выполняется
неравенство
\begin{multline*} \tag{1.4.70}
\| \D^\lambda ((E_k^{l -1,d,m,D} f) \mid_D)
-\D^\lambda ((E_{k -1}^{l -1,d,m,D} f) \mid_D) \|_{L_q(D)} \\
\le c_{26} 2^{k (| \lambda| +(d /p -d /q)_+)}
\Omega^{\prime l}(f, 2 2^{-k})_{L_p(D)} \le \\
c_{48} \int_{2^{-k}}^{2^{-k +1}} t^{-| \lambda| -(d /p -d /q)_+ -1}
\Omega^{\prime l}(f, 2 t)_{L_p(D)} dt.
\end{multline*}

Из (1.4.70) и (1.4.68) вытекает, что ряд
$$
\D^\lambda ((E_{K^0}^{l -1,d,m,D} f) \mid_D) +
\sum_{k = K^0 +1}^\infty (\D^\lambda ((E_k^{l -1,d,m,D} f) \mid_D)
-\D^\lambda ((E_{k -1}^{l -1,d,m,D} f) \mid_D))
$$
сходится в $ L_q(D). $ Принимая
во внимание это обстоятельство, а также равенство (1.4.11), для
любой функции $ \phi \in C_0^\infty(D) $ имеем
\begin{multline*}
\langle \D^\lambda f, \phi \rangle =
(-1)^{|\lambda|} \int_{D} f \D^\lambda \phi dx  =
(-1)^{|\lambda|} \int_{D} E_{K^0}^{l -1,d,m,D} f
\D^\lambda \phi dx +\\
\sum_{k = K^0 +1}^\infty (-1)^{|\lambda|} \int_{D}
(E_k^{l -1,d,m,D} f -E_{k -1}^{l -1,d,m,D} f)
\D^\lambda \phi dx = \\
\int_{D} \D^\lambda (E_{K^0}^{l -1,d,m,D} f) \phi dx \\
+\sum_{k = K^0 +1}^\infty \int_{D} (\D^\lambda (E_k^{l -1,d,m,D} f)
-\D^\lambda (E_{k -1}^{l -1,d,m,D} f)) \phi dx \\
= \int_{D} \biggl(\D^\lambda ((E_{K^0}^{l -1,d,m,D} f) \mid_D) + \\
\sum_{k = K^0 +1}^\infty (\D^\lambda ((E_k^{l -1,d,m,D} f) \mid_D) -
\D^\lambda ((E_{k -1}^{l -1,d,m,D} f) \mid_D)) \biggr) \phi dx.
\end{multline*}

А это значит, что для обобщённой производной $ \D^\lambda f $
функции $ f $ в $ L_q(D) $ имеет место равенство
\begin{equation*}
\D^\lambda f = \D^\lambda ((E_{K^0}^{l -1,d,m,D} f) \mid_D)
+\sum_{k = K^0 +1}^\infty \biggl(\D^\lambda ((E_k^{l -1,d,m,D} f) \mid_D) -
\D^\lambda ((E_{k -1}^{l -1,d,m,D} f) \mid_D)\biggr).
\end{equation*}

Исходя из этого равенства, используя  (1.4.70), при $ k \in \Z_+: k \ge K^0 $
выводим оценку
\begin{multline*}
\| \D^\lambda f -\D^\lambda ((E_k^{l -1,d,m,D} f) \mid_D) \|_{L_q(D)} \le \\
\sum_{j = k +1}^\infty \biggl\| \D^\lambda((E_j^{l -1,d,m,D} f) \mid_D)
-\D^\lambda((E_{j -1}^{l -1,d,m,D} f) \mid_D)\biggr\|_{L_q(D)} \\
\le \sum_{j = k +1}^\infty c_{48} \int_{2^{-j}}^{2^{-j +1}} K(f,t) dt =
c_{48} \int_0^{2^{-k}} K(f,t) dt,
\end{multline*}
что совпадает с (1.4.69). $ \square $

Предложение 1.4.11

Пусть $ d \in \N, \ \alpha \in \R_+, \ l = l(\alpha), \ m \in \N, \ \lambda \in
\Z_{+ m}^d, \ D \subset \R^d $ -- ограниченная область $ \e $-типа.
Пусть ещё $ 1 < p < \infty, \ 1 \le q \le \infty $ и соблюдается условие
\begin{equation*} \tag{1.4.71}
\alpha -|\lambda| -(d /p -d /q)_+ >0.
\end{equation*}
Тогда существует константа $ c_{49}(d,\alpha,m,D,\lambda,p,q) >0 $ такая, что
для любой функции $ f \in (\mathcal H_p^\alpha)^\prime(D) $
при $ k \in \Z_+: k \ge K^0, $ имеет место неравенство
\begin{equation*} \tag{1.4.72}
\| \D^\lambda f -\D^\lambda E_k^{l -1,d,m,D} f \|_{L_q(D)}
\le c_{49} 2^{-k(\alpha -|\lambda| -(d /p -d /q)_+)}.
\end{equation*}

Доказательство.

Прежде всего заметим, что для параметров $ d,\alpha,l,m,\lambda,p,q,D, $
указанных в условии предложения, при соблюдении неравенства (1.4.71)
для $ f \in (H_p^\alpha)^\prime(D) $ при $ t \in I $ справедливо соотношение
\begin{multline*} \tag{1.4.73}
K(f,t) = t^{-| \lambda| -(d /p -d /q)_+ -1}
\Omega^{\prime l}(f, 2 t)_{L_p(D)} \\
\le c_{50} t^{\alpha -|\lambda| -(d /p -d /q)_+ -1}
\sup_{\tau >0} \tau^{-\alpha}
\Omega^{\prime l}(f, \tau)_{L_p(D)} \in L_1(I).
\end{multline*}
Принимая во внимание (1.4.68) и (1.4.73), в соответствии с (1.4.69)
заключаем, что для $ f \in (\mathcal H_p^\alpha)^\prime(D) $
при $ k \in \Z_+: k \ge K^0, $ имеет место оценка
\begin{multline*}
\| \D^\lambda f -\D^\lambda ((E_k^{l -1,d,m,D} f) \mid_D)\|_{L_q(D)} \le \\
c_{48} \int_0^{2^{-k}} c_{50} t^{\alpha -|\lambda| -(d /p -d /q)_+ -1} dt = 
c_{49} 2^{-k(\alpha -|\lambda| -(d /p -d /q)_+)}. \square
\end{multline*}

Теорема 1.4.12

Пусть $ d \in \N, \ \alpha \in \R_+, \ l = l(\alpha), \ \lambda \in \Z_+^d. $
Пусть ещё $ 1 < p < \infty, \ 1 \le q \le \infty $ и
соблюдается условие (1.4.71). Тогда существует константа
$ c_{51}(d,\alpha,\lambda,p,q) >0 $ такая, что для любого
$ \delta \in \R_+ $ и любого $ x^0 \in \R^d $ для $ Q = x^0 +\delta I^d $
для любой функции $ f \in (H_p^\alpha)^\prime(Q) $ справедливо неравенство
\begin{multline*} \tag{1.4.74}
\| \D^\lambda f -\D^\lambda P_{\delta,  x^0,\e}^{l -1,d,0} f \|_{L_q(Q)} \le
c_{51} \delta^{-|\lambda| -d /p +d /q} \int_0^1
t^{-| \lambda| -(d /p -d /q)_+ -1} \times \\
\delta^{-d /p} t^{-d /p}
\biggl(\int_{ 2 \delta t B^d}
\int_{ Q_{l \xi}} | \Delta_\xi^l f(x)|^p dx d\xi \biggr)^{1/p} dt.
\end{multline*}

Доказательство.

Сначала установим справедливость (1.4.74) при $ \delta =1, \ x^0 =0. $
Для этого заметим, что куб $ I^d $ является ограниченной областью
$ \e $-типа, причём константа $ K^0 = K^0(d,I^d,\e) $ из определения 1 при
$ \alpha = \e $ равна $ K^0(d,I^d,\e) =2.$
Поэтому в условиях теоремы, фиксировав $ m \in \N, $ для которого $ \lambda \in
\Z_{+ m}^d, $ учитывая, что для $ f \in (H_p^\alpha)^\prime(I^d) $
в силу (1.4.71) имеет место (1.4.73) при $ D = I^d, $ а,
следовательно, соблюдается (1.4.68), используя (1.4.69) при $ k =2, \ D = I^d, $ получаем
\begin{multline*} \tag{1.4.75}
\| \D^\lambda f -\D^\lambda P_{1,0,\e}^{l -1,d,0} f \|_{L_q(I^d)} \le
\| \D^\lambda f -\D^\lambda ((E_2^{l -1,d,m,I^d} f) \mid_{I^d})\|_{L_q(I^d)} +\\
\| \D^\lambda ((E_2^{l -1,d,m,I^d} f) \mid_{I^d}) -
\D^\lambda ((P_{1,0,\e}^{l -1,d,0} f) \mid_{I^d}) \|_{L_q(I^d)} \le \\
c_{48} \int_0^{2^{-2}} t^{-|\lambda| -(d /p -d /q)_+ -1}
\Omega^{\prime l}(f, 2 t)_{L_p(I^d)} dt + \\
\| \D^\lambda ((E_2^{l -1,d,m,I^d} f) \mid_{I^d}) -
\D^\lambda ((P_{1,0,\e}^{l -1,d,0} f) \mid_{I^d}) \|_{L_q(I^d)} = \\
c_{52} \int_0^{2^{-2}} t^{-|\lambda| -(d /p -d /q)_+ -1}
\biggl( t^{-d} \int_{ 2 t B^d} \int_{I^d_{l \xi}}
|\Delta_{\xi}^l f(x)|^p dx d\xi\biggr)^{1 /p} dt + \\
\| \D^\lambda ((E_2^{l -1,d,m,I^d} f) \mid_{I^d}) -
\D^\lambda ((P_{1,0,\e}^{l -1,d,0} f) \mid_{I^d}) \|_{L_q(I^d)} = \\
c_{53} \int_0^{1} t^{-|\lambda| -(d /p -d /q)_+ -1}
t^{-d /p} \biggl( \int_{ 2^{-1} t B^d} \int_{I^d_{l \xi}}
|\Delta_{\xi}^l f(x)|^p dx d\xi\biggr)^{1 /p} dt +\\
\| \D^\lambda ((E_2^{l -1,d,m,I^d} f) \mid_{I^d}) -
\D^\lambda ((P_{1,0,\e}^{l -1,d,0} f) \mid_{I^d}) \|_{L_q(I^d)}.
\end{multline*}

Оценим второе слагаемое в правой части (1.4.75). Принимая во внимание (1.3.9),
(1.4.33), имеем
\begin{multline*} \tag{1.4.76}
\biggl\| \D^\lambda ((E_2^{l -1,d,m,I^d} f) \mid_{I^d}) -
\D^\lambda ((P_{1,0,\e}^{l -1,d,0} f) \mid_{I^d}) \biggr\|_{L_q(I^d)} = \\
\biggl\| \D^\lambda \biggl((\sum_{\nu \in N_2^{d,m,I^d}}
(S_{2,\nu_2^{I^d}(\nu)}^{l -1,d} f ) g_{2, \nu}^{m,d}) \mid_{I^d}\biggr) -
\D^\lambda \biggl((P_{1,0,\e}^{l -1,d,0} f) \mid_{I^d}
(\sum_{\nu \in N_2^{d,m,I^d}} g_{2, \nu}^{m,d}) \mid_{I^d}\biggr) \biggr\|_{L_q(I^d)} = \\
\biggl\| \D^\lambda \biggl((\sum_{\nu \in N_2^{d,m,I^d}}
(S_{2,\nu_2^{I^d}(\nu)}^{l -1,d} f ) g_{2, \nu}^{m,d}) \mid_{I^d}\biggr) -
\D^\lambda \biggl((\sum_{\nu \in N_2^{d,m,I^d}} (S_{0,0}^{l -1,d} f)
g_{2, \nu}^{m,d}) \mid_{I^d}\biggr) \biggr\|_{L_q(I^d)} = \\
\biggl\| \D^\lambda \biggl((\sum_{\nu \in N_2^{d,m,I^d}}
(S_{2,\nu_2^{I^d}(\nu)}^{l -1,d} f ) g_{2, \nu}^{m,d} -
\sum_{\nu \in N_2^{d,m,I^d}} (S_{0,0}^{l -1,d} f)
g_{2, \nu}^{m,d}) \mid_{I^d}\biggr) \biggr\|_{L_q(I^d)} = \\
\biggl\| \D^\lambda \biggl((\sum_{\nu \in N_2^{d,m,I^d}}
(S_{2,\nu_2^{I^d}(\nu)}^{l -1,d} f -
S_{0,0}^{l -1,d} f) g_{2, \nu}^{m,d}) \mid_{I^d}\biggr) \biggr\|_{L_q(I^d)} = \\
\biggl\| \biggl(\D^\lambda (\sum_{\nu \in N_2^{d,m,I^d}}
(S_{2,\nu_2^{I^d}(\nu)}^{l -1,d} f -
S_{0,0}^{l -1,d} f) g_{2, \nu}^{m,d})\biggr) \mid_{I^d} \biggr\|_{L_q(I^d)} = \\
\biggl\| \sum_{\nu \in N_2^{d,m,I^d}}
\D^\lambda \biggl((S_{2,\nu_2^{I^d}(\nu)}^{l -1,d} f -
S_{0,0}^{l -1,d} f) g_{2, \nu}^{m,d}\biggr) \biggr\|_{L_q(I^d)} = \\
\biggl\| \sum_{\nu \in N_2^{d,m,I^d}} \sum_{\mu \in \Z_+^d(\lambda)} C_\lambda^\mu
\D^\mu (S_{2,\nu_2^{I^d}(\nu)}^{l -1,d} f -
S_{0,0}^{l -1,d} f) \D^{\lambda -\mu} g_{2, \nu}^{m,d} \biggr\|_{L_q(I^d)} = \\
\biggl\| \sum_{\mu \in \Z_+^d(\lambda)} C_\lambda^\mu \sum_{\nu \in N_2^{d,m,I^d}}
\D^\mu (S_{2,\nu_2^{I^d}(\nu)}^{l -1,d} f -
S_{0,0}^{l -1,d} f) \D^{\lambda -\mu} g_{2, \nu}^{m,d} \biggr\|_{L_q(I^d)} \le \\
\sum_{\mu \in \Z_+^d(\lambda)} C_\lambda^\mu
\biggl\| \sum_{\nu \in N_2^{d,m,I^d}}
\D^\mu (S_{2,\nu_2^{I^d}(\nu)}^{l -1,d} f -
S_{0,0}^{l -1,d} f) \D^{\lambda -\mu} g_{2, \nu}^{m,d} \biggr\|_{L_q(I^d)}.
\end{multline*}

Оценивая правую часть (1.4.76), при  $ \mu \in \Z_+^d(\lambda) $ имеем
\begin{multline*} \tag{1.4.77}
\biggl\| \sum_{\nu \in N_2^{d,m,I^d}}
 \D^\mu (S_{2,\nu_2^{I^d}(\nu)}^{l -1,d} f -
S_{0,0}^{l -1,d} f) \D^{\lambda -\mu} g_{2, \nu}^{m,d} \biggr\|_{L_q(I^d)}^q = \\
\int_{I^d} \biggl| \sum_{\nu \in N_2^{d,m,I^d}}
\D^\mu (S_{2,\nu_2^{I^d}(\nu)}^{l -1,d} f -
S_{0,0}^{l -1,d} f) \D^{\lambda -\mu} g_{2, \nu}^{m,d} \biggr|^q dx = \\
\int_{ \cup_{n \in \Nu_{0, 2^2 -1}^d} Q_{2,n}^d}
\biggl| \sum_{\nu \in N_2^{d,m,I^d}}
\D^\mu (S_{2,\nu_2^{I^d}(\nu)}^{l -1,d} f -
S_{0,0}^{l -1,d} f) \D^{\lambda -\mu} g_{2, \nu}^{m,d} \biggr|^q dx = \\
\sum_{n \in \Nu_{0, 3}^d} \int_{Q_{2,n}^d}
\biggl| \sum_{\nu \in N_2^{d,m,I^d}} \D^\mu (S_{2,\nu_2^{I^d}(\nu)}^{l -1,d} f -
S_{0,0}^{l -1,d} f) \D^{\lambda -\mu} g_{2, \nu}^{m,d} \biggr|^q dx = \\
\sum_{n \in \Nu_{0, 3}^d} \int_{Q_{2,n}^d}
\biggl| \sum_{\nu \in N_2^{d,m,I^d}: Q_{2,n}^d \cap \supp g_{2,\nu}^{m,d} \ne
\emptyset}
\D^\mu (S_{2,\nu_2^{I^d}(\nu)}^{l -1,d} f -
S_{0,0}^{l -1,d} f) \D^{\lambda -\mu} g_{2, \nu}^{m,d} \biggr|^q dx = \\
\sum_{n \in \Nu_{0, 3}^d}
\biggl\| \sum_{\nu \in N_2^{d,m,I^d}: Q_{2,n}^d \cap \supp g_{2,\nu}^{m,d} \ne
\emptyset}
\D^\mu (S_{2,\nu_2^{I^d}(\nu)}^{l -1,d} f -
S_{0,0}^{l -1,d} f) \D^{\lambda -\mu} g_{2, \nu}^{m,d} \biggr\|_{L_q(Q_{2,n}^d)}^q \le \\
\sum_{n \in \Nu_{0, 3}^d}
\biggl( \sum_{\nu \in N_2^{d,m,I^d}: Q_{2,n}^d \cap \supp g_{2,\nu}^{m,d} \ne
\emptyset}
\|\D^\mu (S_{2,\nu_2^{I^d}(\nu)}^{l -1,d} f -
S_{0,0}^{l -1,d} f) \D^{\lambda -\mu} g_{2, \nu}^{m,d} \|_{L_q(Q_{2,n}^d)}\biggr)^q.
\end{multline*}

Далее, при $ \mu \in \Z_+^d(\lambda) $ для $ n \in \Nu_{0, 3}^d, \
\nu \in N_2^{d,m,I^d}: Q_{2,n}^d \cap \supp g_{2,\nu}^{m,d} \ne \emptyset, $
ввиду (1.3.10) выполняется неравенство
\begin{multline*} \tag{1.4.78}
\|\D^\mu (S_{2,\nu_2^{I^d}(\nu)}^{l -1,d} f -S_{0,0}^{l -1,d} f)
\D^{\lambda -\mu} g_{2, \nu}^{m,d} \|_{L_q(Q_{2,n}^d)} \le \\
\| \D^{\lambda -\mu} g_{2, \nu}^{m,d} \|_{L_\infty(\R^d)}
\|\D^\mu (S_{2,\nu_2^{I^d}(\nu)}^{l -1,d} f -
S_{0,0}^{l -1,d} f) \|_{L_q(Q_{2,n}^d)} = \\
c_{54} 2^{2 |\lambda -\mu|}
\|\D^\mu (S_{2,\nu_2^{I^d}(\nu)}^{l -1,d} f -
S_{0,0}^{l -1,d} f) \|_{L_q(Q_{2,n}^d)},
\end{multline*}
а с учётом того, что $ \nu_2^{I^d}(\nu) \in \Nu_{0,3}^d, $ в силу (1.4.1),
(1.4.2), (1.4.3), (1.4.6) соблюдается неравенство
\begin{multline*} \tag{1.4.79}
\|\D^\mu (S_{2,\nu_2^{I^d}(\nu)}^{l -1,d} f -
S_{0,0}^{l -1,d} f) \|_{L_q(Q_{2,n}^d)} \le \\
c_{55} (2^{-2})^{-|\mu| -d /p +d /q}
\| S_{2,\nu_2^{I^d}(\nu)}^{l -1,d} f -
S_{0,0}^{l -1,d} f \|_{L_p(Q_{2,\nu_2^{I^d}(\nu)}^d)} = \\
c_{55} (2^{-2})^{-|\mu| -d /p +d /q}
\| S_{2,\nu_2^{I^d}(\nu)}^{l -1,d} (f -
S_{0,0}^{l -1,d} f) \|_{L_p(Q_{2,\nu_2^{I^d}(\nu)}^d)} \le \\
c_{56} (2^{-2})^{-|\mu| -d /p +d /q}
\| f -S_{0,0}^{l -1,d} f \|_{L_p(Q_{2,\nu_2^{I^d}(\nu)}^d)} \le \\
c_{56} (2^{-2})^{-|\mu| -d /p +d /q}
\| f -S_{0,0}^{l -1,d} f \|_{L_p(I^d)} \le \\
c_{56} (2^{-2})^{-|\mu| -d /p +d /q}
c_5 \left(\int_{B^d} \int_{I^d_{l \xi}}
| \Delta_\xi^l f(x)|^p \,dx \,d\xi \right)^{1/p} = \\
c_{57} (2^{-2})^{-|\mu| -d /p +d /q}
\left(\int_{B^d} \int_{I^d_{l \xi}}
| \Delta_\xi^l f(x)|^p \,dx \,d\xi \right)^{1/p}.
\end{multline*}

Объединяя (1.4.78), (1.4.79), при $ \mu \in \Z_+^d(\lambda) $ для $ n \in
\Nu_{0, 3}^d, \ \nu \in N_2^{d,m,I^d}: Q_{2,n}^d \cap \supp g_{2,\nu}^{m,d} \ne
\emptyset, $ получаем оценку
\begin{multline*}
\|\D^\mu (S_{2,\nu_2^{I^d}(\nu)}^{l -1,d} f -S_{0,0}^{l -1,d} f)
\D^{\lambda -\mu} g_{2, \nu}^{m,d} \|_{L_q(Q_{2,n}^d)} \le \\
c_{54} 2^{2 |\lambda -\mu|}
c_{57} (2^{-2})^{-|\mu| -d /p +d /q}
\left(\int_{B^d} \int_{I^d_{l \xi}}
| \Delta_\xi^l f(x)|^p \,dx \,d\xi \right)^{1/p} = \\
c_{58} \left(\int_{B^d} \int_{I^d_{l \xi}}
| \Delta_\xi^l f(x)|^p \,dx \,d\xi \right)^{1/p}.
\end{multline*}

Подставляя эту оценку в (1.4.77) и учитывая (1.3.8), находим, что при
$ \mu \in \Z_+^d(\lambda) $ справедливо неравенство
\begin{multline*}
\biggl\| \sum_{\nu \in N_2^{d,m,I^d}}
\D^\mu (S_{2,\nu_2^{I^d}(\nu)}^{l -1,d} f -
S_{0,0}^{l -1,d} f) \D^{\lambda -\mu} g_{2, \nu}^{m,d} \biggr\|_{L_q(I^d)}^q \le \\
\sum_{n \in \Nu_{0, 3}^d}
\biggl( \sum_{\nu \in N_2^{d,m,I^d}: Q_{2,n}^d \cap \supp g_{2,\nu}^{m,d} \ne
\emptyset} c_{58} \biggl(\int_{B^d} \int_{I^d_{l \xi}}
| \Delta_\xi^l f(x)|^p \,dx \,d\xi \biggr)^{1/p}\biggr)^q \le \\
\sum_{n \in \Nu_{0, 3}^d}
\biggl( \card \{ \nu \in \Z^d: Q_{2,n}^d \cap \supp g_{2,\nu}^{m,d} \ne \emptyset\}
c_{58} \biggl(\int_{B^d} \int_{I^d_{l \xi}}
| \Delta_\xi^l f(x)|^p \,dx \,d\xi \biggr)^{1/p}\biggr)^q = \\
\sum_{n \in \Nu_{0, 3}^d}
\biggl( \card \Nu_{-m,0}^d c_{58} \biggl(\int_{B^d} \int_{I^d_{l \xi}}
| \Delta_\xi^l f(x)|^p \,dx \,d\xi \biggr)^{1/p}\biggr)^q = \\
\biggl(c_{59} \biggl(\int_{B^d} \int_{I^d_{l \xi}}
| \Delta_\xi^l f(x)|^p \,dx \,d\xi \biggr)^{1/p}\biggr)^q.
\end{multline*}

Соединяя эту оценку с (1.4.76), выводим
\begin{multline*}
\| \D^\lambda ((E_2^{l -1,d,m,I^d} f) \mid_{I^d}) -
\D^\lambda ((P_{1,0,\e}^{l -1,d,0} f) \mid_{I^d}) \|_{L_q(I^d)} \le \\
\sum_{\mu \in \Z_+^d(\lambda)} C_\lambda^\mu
c_{59} \left(\int_{B^d} \int_{I^d_{l \xi}}
| \Delta_\xi^l f(x)|^p \,dx \,d\xi \right)^{1/p} =
c_{60} \left(\int_{B^d} \int_{I^d_{l \xi}}
| \Delta_\xi^l f(x)|^p \,dx \,d\xi \right)^{1/p} = \\
c_{61} \int_{2^{-1}}^1 \left(\int_{B^d} \int_{I^d_{l \xi}}
| \Delta_\xi^l f(x)|^p \,dx \,d\xi \right)^{1/p} dt \le
c_{61} \int_{2^{-1}}^1 \left(\int_{ 2 t B^d} \int_{I^d_{l \xi}}
| \Delta_\xi^l f(x)|^p \,dx \,d\xi \right)^{1/p} dt \le \\
c_{61} \int_0^1 \left(\int_{ 2 t B^d} \int_{I^d_{l \xi}}
| \Delta_\xi^l f(x)|^p \,dx \,d\xi \right)^{1/p} dt \le \\
c_{61} \int_0^1 t^{-|\lambda| -(d /p -d /q)_+ -1}
t^{-d /p} \left(\int_{ 2 t B^d} \int_{I^d_{l \xi}}
| \Delta_\xi^l f(x)|^p \,dx \,d\xi \right)^{1/p} dt.
\end{multline*}

Подставляя эту оценку в (1.4.75), приходим к неравенству
\begin{multline*}
\| \D^\lambda f -\D^\lambda P_{1,0,\e}^{l -1,d,0} f \|_{L_q(I^d)} \le \\
c_{53} \int_0^{1} t^{-|\lambda| -(d /p -d /q)_+ -1}
t^{-d /p} \biggl( \int_{ 2^{-1} t B^d} \int_{I^d_{l \xi}}
|\Delta_{\xi}^l f(x)|^p dx d\xi\biggr)^{1 /p} dt +\\
c_{61} \int_0^1 t^{-|\lambda| -(d /p -d /q)_+ -1}
t^{-d /p} \left(\int_{ 2 t B^d} \int_{I^d_{l \xi}}
| \Delta_\xi^l f(x)|^p \,dx \,d\xi \right)^{1/p} dt \le \\
c_{51} \int_0^{1} t^{-|\lambda| -(d /p -d /q)_+ -1}
t^{-d /p} \biggl( \int_{ 2 t B^d} \int_{I^d_{l \xi}}
|\Delta_{\xi}^l f(x)|^p dx d\xi\biggr)^{1 /p} dt,
\end{multline*}
которое совпадает с (1.4.74) при $ \delta = 1, \ x^0 =0. $ А отсюда выводим неравенство
(1.4.74) при произвольных $ \delta \in \R_+ $ и $ x^0 \in \R^d $
путём перехода от функции $ f \in (H_p^\alpha)^\prime(Q) $ к функции
$ (h_{\delta,x^0} f)(\cdot) = f(x^0 +\delta \cdot) \in
(H_p^\alpha)^\prime(I^d), $ и, учитывая тот факт, что
$ P_{\delta,x^0,\e}^{l -1,d,0} f = h_{\delta, x^0}^{-1} P_{1, 0,\e}^{l -1,d,0}
h_{\delta, x^0} f, \ f \in L_1(Q). \square $

Теорема 1.4.13

Пусть $ d \in \N, \ \alpha \in \R_+, \ 1 < p < \infty, \ 1 \le q \le \infty, \
\lambda \in \Z_+^d, \ D \subset \R^d $ -- ограниченная область $ \e $-типа.
Тогда если соблюдается условие (1.4.71), то существует константа
$ c_{62}(d,\alpha,p,q,\lambda,D) >0 $ такая, что для любой функции $ f \in
(H_p^\alpha)^\prime(D) $ выполняется неравенство
\begin{equation*} \tag{1.4.80}
\| \D^\lambda f \|_{L_q(D)} \le c_{62} \| f \|_{(H_p^\alpha)^\prime(D)},
\end{equation*}
и если
\begin{equation*} \tag{1.4.81}
\alpha -|\lambda| -d /p >0,
\end{equation*}
то для любой функции $ f \in (H_p^\alpha)^\prime(D) $ имеет место включение
\begin{equation*} \tag{1.4.82}
\D^\lambda f \in C(D).
\end{equation*}

Доказательство.

Неравенство (1.4.80) без труда выводится, исходя из (1.4.72) и (1.4.34)
при соответствующем выборе значений  параметров $ l,m,k,s. $

Для доказательства (1.4.82) достаточно заметить, что при выполнении (1.4.81)
и $ m \in \N: \lambda \in \Z_{+ m -1}^d, $ в силу (1.4.72) каждая
функция $ \D^\lambda f $ при $ f \in (\mathcal H_p^\alpha)^\prime(D), $
а, следовательно, и при $ f \in (H_p^\alpha)^\prime(D), $
является пределом в $ L_\infty(D) $ равномерно сходящейся на $ \overline D $
последовательности непрерывных в $ \R^d $ функций $ \{\D^\lambda E_k^{l -1,d,m,D} f, \
k \in \Z_+: k \ge K^0\}. \square $
\bigskip

\centerline{\S 2. Продолжение функций из пространств
$ (H_p^\alpha)^\prime(D) $ и $ (B_{p,\theta}^\alpha)^\prime(D) $}
\bigskip

2.1. В этом пункте доказывается один из важнейших результатов работы

Теорема 2.1.1

Пусть $ d \in \N, \ \alpha \in \R_+, \ 1 < p < \infty, \ 1 \le \theta \le \infty $
и $ D \subset \R^d $ -- ограниченная область $ \e $-типа.
Тогда существует непрерывное линейное отображение
$ E^{d, \alpha, p, \theta,D}: (B_{p, \theta}^\alpha)^\prime(D) \mapsto
L_p(\R^d) $ такое, что для любой функции $ f \in (B_{p, \theta}^\alpha)^\prime(D) $
соблюдается равенство
\begin{equation*} \tag{2.1.1}
(E^{d, \alpha, p, \theta,D} f) \mid_{D} = f,
\end{equation*}
а при $ \l \in \Z_+: \l < \alpha, $ существует константа
$ c_1(d,\alpha,p,\theta,D,\l) > 0, $ для которой выполняется неравенство
\begin{equation*} \tag{2.1.2}
\| E^{d, \alpha, p, \theta,D} f\|_{(B_{p, \theta}^\alpha)^{\l}(\R^d)} \le
c_1 \| f\|_{(B_{p, \theta}^\alpha)^\prime(D)}.
\end{equation*}

Доказательство.

В условиях теоремы положим $ l = l(\alpha) $ и фиксируем $ m \in \N $ так,
чтобы $ l \le m. $ Отметим ещё, что в качестве $ K^0 $ рассматривается
константа из определения 1 при $ \alpha = \e. $ Сначала рассмотрим случай,
когда $ \theta \ne \infty. $

Прежде всего заметим, что для $ f \in (B_{p, \theta}^\alpha)^\prime(D) $ ряд
\begin{equation*}
\sum_{k = K^0 +1}^\infty \| \mathcal E_k^{l -1,d,m,D} f\|_{L_p(\R^d)}
\end{equation*}
сходится, ибо в силу (1.4.41) при $ q = p, \ \lambda =0 $ и неравенства
Гёльдера имеет место оценка
\begin{multline*} \tag{2.1.3}
\sum_{ k \in \N: n \le k \le n +r}
\| \mathcal E_k^{l -1,d,m,D} f \|_{L_p(\R^d)} \le
\sum_{ k \in \N: n \le k \le n +r} c_{2}
\Omega^{\prime l}(f, 2 2^{-k})_{L_p(D)} = \\
c_{2} \sum_{ k \in \n: n \le k \le n +r}
\Omega^{\prime l}(f, 2 2^{-k})_{L_p(D)} = 
c_{2} \sum_{ k \in \N: n \le k \le n +r}
2^{-k \alpha} 2^{k \alpha} \Omega^{\prime l}(f, 2 2^{-k})_{L_p(D)} \le \\
c_{2} \biggl( \sum_{ k \in \N: n \le k \le n +r}
2^{-\theta^\prime k \alpha} \biggr)^{1 / \theta^\prime}
\biggl( \sum_{ k \in \N: n \le k \le n +r} (2^{k \alpha} \Omega^{\prime l}(f, 2
2^{-k})_{L_p(D)})^{\theta} \biggr)^{1 / \theta} \le \\
c_{3} \biggl( \sum_{ k \in \N: n \le k \le n +r}
2^{-\theta^\prime k \alpha} \biggr)^{1 / \theta^\prime} \biggl( \sum_{k \in \N}
((2^{-(k -1) })^{-\alpha}
\Omega^{\prime l}(f, 2 2^{-k})_{L_p(D)})^{\theta} \biggr)^{1 / \theta} \le \\
c_{4} \biggl( \sum_{ k \in \N: n \le k \le n +r}
2^{-\theta^\prime k \alpha} \biggr)^{1 / \theta^\prime} \biggl( \sum_{k \in \N}
\int_{2^{-k}}^{2^{-(k -1)}} t^{-1 -\theta \alpha}
(\Omega^{\prime l}(f, 2 t)_{L_p(D)})^{\theta} dt \biggr)^{1 / \theta} \le \\
c_{4} \biggl( \sum_{ k \in \N: n \le k \le n +r}
2^{-\theta^\prime k \alpha} \biggr)^{1 / \theta^\prime}
\biggl(\int_{\R_+} t^{-1 -\theta \alpha}
(\Omega^{\prime l}(f, 2 t)_{L_p(D)})^{\theta} dt \biggr)^{1 / \theta} \le \\
c_{5} \biggl( \sum_{ k \in \N: n \le k \le n +r}
2^{-\theta^\prime k \alpha} \biggr)^{1 / \theta^\prime}
\biggl(\int_{\R_+} t^{-1 -\theta \alpha}
(\Omega^{\prime l}(f, t)_{L_p(D)})^{\theta} dt \biggr)^{1 / \theta},\\
n \in \N: n > K^0, \ r \in \Z_+,
\end{multline*}
где $ \theta^\prime = \theta/(\theta -1). $

Поэтому для $ f \in (B_{p,\theta}^\alpha)^\prime(D) $ ряд
$ \sum_{ k = K^0 +1}^\infty \mathcal E_k^{l -1,d,m,D} f $
сходится в $ L_p(\R^d). $ Определим отображение $ E^{d,\alpha,p,\theta,D}:
(B_{p,\theta}^\alpha)^\prime(D) \mapsto L_p(\R^d), $ полагая для
$ f \in (B_{p,\theta}^\alpha)^\prime(D) $ значение
\begin{equation*} \tag{2.1.4}
E^{d,\alpha,p,\theta,D} f = E_{K^0}^{l -1,d,m,D} f +
\sum_{ k = K^0 +1}^\infty \mathcal E_k^{l -1,d,m,D} f =
E_{K^0}^{l -1,d,m,D} f +E^{\prime d,\alpha,p,\theta,D} f,
\end{equation*}
где
\begin{equation*} \tag{2.1.5}
E^{\prime d,\alpha,p,\theta,D} f = \sum_{ k = K^0 +1}^\infty
\mathcal E_k^{l -1,d,m,D} f.
\end{equation*}
Ясно, что отображение $ E^{d,\alpha,p,\theta,D} $ линейно, а, вследствие
(2.1.4), (1.4.40), (1.3.13), (1.4.11), в $ L_p(D) $ выполняется равенство
\begin{multline*}
(E^{d,\alpha,p,\theta,D} f) \mid_{D} =
(E_{K^0}^{l -1,d,m,D} f) \mid_D +\sum_{ k = K^0 +1}^\infty
(\mathcal E_k^{l -1,d,m,D} f) \mid_{D} = \\
(E_{K^0}^{l -1,d,m,D} f) \mid_D +\sum_{ k = K^0 +1}^\infty
((E_k^{l -1,d,m,D} f) \mid_D -
(H_{k -1}^{l -1,d,m,D} (E_{k -1}^{l -1,d,m,D} f)) \mid_D) = \\
(E_{K^0}^{l -1,d,m,D} f) \mid_D +\sum_{ k = K^0 +1}^\infty
((E_k^{l -1,d,m,D} f) \mid_D -
(E_{k -1}^{l -1,d,m,D} f) \mid_D) = \\
(E_{K^0}^{l -1,d,m,D} f) \mid_D +\sum_{ k = K^0 +1}^\infty
(E_k^{l -1,d,m,D} f -E_{k -1}^{l -1,d,m,D} f) \mid_D = f,
\end{multline*}
т.е. соблюдается (2.1.1).

Покажем, что при $ \l \in \Z_+: \l < \alpha, $ отображение, которое каждой
функции $ f \in (B_{p,\theta}^\alpha)^\prime(D) $ ставит в соответствие функцию
$ E_{K^0}^{l -1,d,m,D} f $
является непрерывным отображением пространства $ (B_{p,\theta}^\alpha)^\prime(D) $
в $ (B_{p,\theta}^\alpha)^{\l}(\R^d). $ Действительно, для $ \mathcal F =
E_{K^0}^{l -1,d,m,D} f, \ f \in (B_{p,\theta}^\alpha)^\prime(D), $
при $ \lambda \in \Z_{+ \l}^{\prime d} $ и $ h \in \R^d, $ согласно (1.1.3),
(1.4.34), с одной стороны имеет место неравенство
\begin{multline*}
\| \Delta_h^{l -\l} (\D^\lambda \mathcal F)\|_{L_p(\R^d)} =
\| \Delta_h^{l -\l} (\D^\lambda  E_{K^0}^{l -1,d,m,D} f) \|_{L_p(\R^d)} \le \\
(l -\l)! \cdot \|h\|^{l -\l} \sum_{\mu \in \Z_{+ l -\l}^{\prime d}}
\| \D^\mu (\D^\lambda E_{K^0}^{l -1,d,m,D} f) \|_{L_p(\R^d)} \le \\
(l -\l)! \cdot \|h\|^{l -\l} \sum_{\mu \in \Z_{+ l -\l}^{\prime d}}
c_6 2^{K^0 |\lambda +\mu|} \| f\|_{L_p(D)} =
c_7 \|h\|^{l -\l} \| f \|_{L_p(D)},
\end{multline*}
а с другой стороны --
\begin{equation*}
\| \Delta_h^{l -\l} (\D^\lambda  \mathcal F) \|_{L_p(\R^d)} \le
c_8 \| \D^\lambda \mathcal F \|_{L_p(\R^d)} =
c_8 \| \D^\lambda E_{K^0}^{l -1,d,m,D} f \|_{L_p(\R^d)} \le
c_{9} \| f \|_{L_p(D)}.
\end{equation*}
Отсюда видим, что при $ \lambda \in \Z_{+ \l}^{\prime d}, \ t \in \R_+ $
выполняются неравенства
\begin{equation*}
\Omega^{l -\l}(\D^\lambda \mathcal F, t)_{L_p(\R^d)} \le
c_7 t^{l -\l} \| f \|_{L_p(D)}; \\
\Omega^{l -\l}(\D^\lambda \mathcal F, t)_{L_p(\R^d)} \le
c_{9} \| f \|_{L_p(D)}.
\end{equation*}
Используя эти неравенства, получаем, что
\begin{multline*}
\biggl(\int_{\R_+} t^{-1 -\theta (\alpha -\l)}
(\Omega^{l -\l}(\D^\lambda \mathcal F, t)_{L_p(\R^d)})^{\theta} dt \biggr)^{1 / \theta} = \\
\biggl(\int_{I} t^{-1 -\theta (\alpha -\l)}
(\Omega^{l -\l}(\D^\lambda \mathcal F, t)_{L_p(\R^d)})^{\theta} dt +
\int_1^\infty t^{-1 -\theta (\alpha -\l)}
(\Omega^{l -\l}(\D^\lambda \mathcal F, t)_{L_p(\R^d)})^{\theta} dt \biggr)^{1 / \theta} \le \\
\biggl(\int_{I} t^{-1 -\theta (\alpha -\l)}
(c_7 t^{l -\l} \| f \|_{L_p(D)})^{\theta} dt +
\int_1^\infty t^{-1 -\theta (\alpha -\l)}
(c_{9} \| f \|_{L_p(D)})^{\theta} dt \biggr)^{1 / \theta} = \\
\biggl(c_7^\theta \int_{I} t^{-1 +\theta (l -\alpha)} dt +c_{9}^\theta
\int_1^\infty t^{-1 -\theta (\alpha -\l)} dt \biggr)^{1 / \theta} \| f \|_{L_p(D)} = \\
c_{10} \| f \|_{L_p(D)} \le c_{10} \| f \|_{(B_{p,\theta}^\alpha)^\prime(D)},
\ f \in (B_{p,\theta}^\alpha)^\prime(D), \ \lambda \in \Z_{+ \l}^{\prime d}.
\end{multline*}
Из сказанного с учётом (1.4.34) вытекает, что для $ f \in (B_{p,\theta}^\alpha)^\prime(D) $
справедливо неравенство
\begin{equation*} \tag{2.1.6}
\| E_{K^0}^{l -1,d,m,D} f \|_{(B_{p,\theta}^\alpha)^{\l}(\R^d)}
\le c_{11} \| f \|_{(B_{p,\theta}^\alpha)^\prime(D)}.
\end{equation*}

Теперь проверим, что существует константа
$ c_{12}(d,\alpha,p,\theta,D) > 0 $ такая, что для
$ f \in (B_{p, \theta}^\alpha)^\prime(D) $ выполняется неравенство
\begin{equation*} \tag{2.1.7}
\| E^{\prime d,\alpha,p,\theta,D} f\|_{(B_{p, \theta}^\alpha)^{\l}(\R^d)} \le
c_{12} \| f\|_{(B_{p, \theta}^\alpha)^\prime(D)}.
\end{equation*}

В самом деле, в силу (2.1.5), (2.1.3) справедливо неравенство
\begin{multline*} \tag{2.1.8}
\| E^{\prime d,\alpha,p,\theta,D} f \|_{L_p(\R^d)} \le \sum_{ k = K^0 +1}^\infty
\| \mathcal E_k^{l -1,d,m,D} f\|_{L_p(\R^d)} \le \\
c_5 \biggl(\sum_{ k = K^0 +1}^\infty 2^{-\theta^\prime k \alpha}\biggr)^{1 / \theta^\prime}
\biggl(\int_{\R_+} t^{-1 -\theta \alpha}
(\Omega^{\prime l}(f, t)_{L_p(D)})^{\theta} dt \biggr)^{1 / \theta} \le \\
c_{13} \| f \|_{(B_{p, \theta}^\alpha)^\prime(D)}.
\end{multline*}

Далее, пользуясь тем, что для $ f \in (B_{p,\theta}^\alpha)^\prime(D) $
при $ \lambda \in \Z_{+ \l}^{\prime d}, \ k \in \Z_+: k > K^0, $ ввиду (1.4.41) справедлива оценка
\begin{multline*}
\| \D^\lambda \mathcal E_k^{l -1,d,m,D} f \|_{L_p(\R^d)}
\le c_{14} 2^{k |\lambda|} \Omega^{\prime l}(f, 2 2^{-k})_{L_p(D)} \le \\
c_{15} 2^{k \l} 2^{-k \alpha} (2 2^{-k})^{-\alpha}
\Omega^{\prime l}(f, 2 2^{-k})_{L_p(D)} \le
c_{15} 2^{-k(\alpha -\l)} \sup_{t \in \R_+}
t^{-\alpha} \Omega^{\prime l}(f, t)_{L_p(D)},
\end{multline*}
с учётом (1.1.4) заключаем, что ряд
\begin{equation*}
\sum_{k = K^0 +1}^\infty \D^\lambda  \mathcal E_k^{l -1,d,m,D} f
\end{equation*}
сходится в $ L_p(\R^d). $ Отсюда и из (2.1.5) приходим к выводу, что в $ L_p(\R^d) $
имеет место равенство
\begin{equation*} \tag{2.1.9}
\D^\lambda E^{\prime d,\alpha,p,\theta,D} f = \sum_{ k = K^0 +1}^\infty
\D^\lambda \mathcal E_k^{l -1,d,m,D} f, \ f \in (B_{p,\theta}^\alpha)^\prime(D), \
\lambda \in \Z_{+ \l}^{\prime d}.
\end{equation*}
Принимая во внимание (2.1.9), (1.4.41), а также неравенство Гёльдера и (2.1.3),
для $ f \in (B_{p,\theta}^\alpha)^\prime(D) $ при
$ \lambda \in \Z_{+ \l}^{\prime d} $ имеем
\begin{multline*} \tag{2.1.10}
\| \D^\lambda E^{\prime d,\alpha,p,\theta,D} f \|_{L_p(\R^d)}
= \biggl\| \sum_{ k = K^0 +1}^\infty \D^\lambda \mathcal E_k^{l -1,d,m,D} f \biggr\|_{L_p(\R^d)} \le \\
\sum_{ k = K^0 +1}^\infty \| \D^\lambda \mathcal E_k^{l -1,d,m,D} f \|_{L_p(\R^d)} \le
\sum_{ k = K^0 +1}^\infty c_{14} 2^{k |\lambda|} \Omega^{\prime l}(f, 2 2^{-k})_{L_p(D)} = \\
c_{14} \sum_{ k = K^0 +1}^\infty 2^{k \l}
\Omega^{\prime l}(f, 2 2^{-k})_{L_p(D)} =
c_{14} \sum_{ k = K^0 +1}^\infty 2^{-k(\alpha -\l)} 2^{k \alpha}
\Omega^{\prime l}(f, 2  2^{-k})_{L_p(D)} \le \\
c_{14} \biggl( \sum_{ k = K^0 +1}^\infty
2^{-k \theta^\prime (\alpha -\l)} \biggr)^{1 / \theta^\prime}
\biggl( \sum_{ k = K^0 +1}^\infty (2^{k \alpha} \Omega^{\prime l}(f, 2
2^{-k})_{L_p(D)})^{\theta} \biggr)^{1 / \theta} \le \\
c_{16} \biggl( \sum_{ k = K^0 +1}^\infty
2^{-k \theta^\prime (\alpha -\l)} \biggr)^{1 / \theta^\prime}
\biggl(\int_{\R_+} t^{-1 -\theta \alpha}
(\Omega^{\prime l}(f, t)_{L_p(D)})^{\theta} dt \biggr)^{1 / \theta} \le \\
c_{17} \| f\|_{(B_{p,\theta}^\alpha)^\prime(D)}.
\end{multline*}

Наконец, для функции $ F = E^{\prime d,\alpha,p,\theta,D} f, \
f \in (B_{p,\theta}^\alpha)^\prime(D), $
при $ \lambda \in \Z_{+ \l}^{\prime d} $ оценим \\
\begin{multline*} \tag{2.1.11}
\biggl(\int_{\R_+} t^{-1 -\theta (\alpha -\l)}
(\Omega^{l -\l}(\D^\lambda F, t)_{L_p(\R^d)})^{\theta} dt\biggr)^{1 /\theta} =\\
\biggl(\int_{I} t^{-1 -\theta (\alpha -\l)} (\Omega^{l -\l}(\D^\lambda F, t)_{L_p(\R^d)})^{\theta} dt +
\int_1^\infty t^{-1 -\theta (\alpha -\l)}
(\Omega^{l -\l}(\D^\lambda F, t)_{L_p(\R^d)})^{\theta} dt\biggr)^{1 /\theta}.
\end{multline*}

Учитывая, что при $ \lambda \in \Z_{+ \l}^{\prime d} $ для $ t \in \R_+ $
справедливо неравенство
\begin{multline*}
\Omega^{l -\l}(\D^\lambda F, t)_{L_p(\R^d)} =
\sup_{ \{ h \in \R^d: h \in t B^d \}} \| \Delta_h^{l -\l} \D^\lambda F \|_{L_p(\R^d)} \le \\
\sup_{ \{ h \in \R^d: h \in t B^d \}} c_{8} \| \D^\lambda F \|_{L_p(\R^d)} =
c_{8} \| \D^\lambda F \|_{L_p(\R^d)},
\end{multline*}
находим, что
\begin{multline*}
\int_1^\infty t^{-1 -\theta (\alpha -\l)}
(\Omega^{l -\l}(\D^\lambda F, t)_{L_p(\R^d)})^{\theta} dt \le
\int_1^\infty t^{-1 -\theta (\alpha -\l)}
(c_{8} \| \D^\lambda F \|_{L_p(\R^d)})^{\theta} dt = \\
(c_{8} \| \D^\lambda F \|_{L_p(\R^d)})^{\theta}
\int_1^\infty t^{-1 -\theta (\alpha -\l)} dt =
(c_{18} \| \D^\lambda F \|_{L_p(\R^d)})^{\theta}.
\end{multline*}

Подставляя последнюю оценку в (2.1.11) и применяя неравенство (1.1.1),
получаем, что

\begin{multline*} \tag{2.1.12}
\biggl(\int_{\R_+} t^{-1 -\theta (\alpha -\l)}
(\Omega^{l -\l}(\D^\lambda F, t)_{L_p(\R^d)})^{\theta} dt\biggr)^{1 /\theta} \le \\
\biggl(\int_{I} t^{-1 -\theta (\alpha -\l)} (\Omega^{l -\l}(\D^\lambda F, t)_{L_p(\R^d)})^{\theta} dt +
(c_{18} \| \D^\lambda F \|_{L_p(\R^d)})^{\theta}\biggr)^{1 /\theta} \le \\
(\int_{I} t^{-1 -\theta (\alpha -\l)}
(\Omega^{l -\l}(\D^\lambda F, t)_{L_p(\R^d)})^{\theta} dt)^{1 /\theta} +
c_{18} \| \D^\lambda F \|_{L_p(\R^d)}, \ \lambda \in \Z_{+ \l}^{\prime d}.
\end{multline*}

Таким образом, приходим к необходимости оценки
$$
\biggl(\int_{I} t^{-1 -\theta (\alpha -\l)}
(\Omega^{l -\l}(\D^\lambda F, t)_{L_p(\R^d)})^{\theta} dt\biggr)^{1 /\theta}
$$
для $ \lambda \in \Z_{+ \l}^{\prime d}.$ Проведём эту оценку.

Пользуясь тем, что имеет место представление
$ I = \cup_{k \in \N} [2^{-k}, 2^{-k +1}), $
получаем, что при $ \lambda \in \Z_{+ \l}^{\prime d} $ соблюдается неравенство
\begin{multline*} \tag{2.1.13}
\biggl(\int_I t^{-1 -\theta (\alpha -\l)}
(\Omega^{l -\l}(\D^\lambda F, t)_{L_p(\R^d)})^{\theta} dt\biggr)^{1/\theta}  = \\
\biggl(\sum_{ k \in \N} \int_{[2^{-k}, 2^{-k +1})}
t^{-1 -\theta (\alpha -\l)}
(\Omega^{l -\l}(\D^\lambda F, t)_{L_p(\R^d)})^{\theta} dt\biggr)^{1 /\theta} \le \\
\biggl(\sum_{ k \in \N} \!\int_{[2^{-k}, 2^{-k +1})}
(2^{-k}\!)^{-1 -\theta(\alpha -\l)}
(\Omega^{l -\l}(\D^\lambda F, 2^{-k +1}\!)_{L_p(\R^d)})^{\theta} dt\biggr)^{1 /\theta} = \\
\biggl(\sum_{ k \in \N} 2^{k \theta(\alpha -\l)}
(\Omega^{l -\l}(\D^\lambda F, 2^{-k +1})_{L_p(\R^d)})^{\theta}\biggr)^{1 /\theta}.
\end{multline*}

При $ \lambda \in \Z_{+ \l}^{\prime d} , \ k \in \N, $ чтобы оценить
$ \Omega^{l -\l}(\D^\lambda F, 2^{-k +1})_{L_p(\R^d)}, $
для $ h \in (2^{-k +1} B^d) $ с учётом (2.1.9), (2.1.10) имеем
\begin{multline*} \tag{2.1.14}
\| \Delta_h^{l -\l} \D^\lambda F \|_{L_p(\R^d)} =
\biggl\| \Delta_h^{l -\l} \biggl(\sum_{ \kappa = K^0 +1}^\infty
\D^\lambda \mathcal E_\kappa^{l -1,d,m,D} f\biggr) \biggr\|_{L_p(\R^d)} = \\
\biggl\| \sum_{ \kappa = K^0 +1}^\infty \Delta_h^{l -\l}
(\D^\lambda \mathcal E_\kappa^{l -1,d,m,D} f) \biggr\|_{L_p(\R^d)} \le
\sum_{ \kappa = K^0 +1}^\infty \| \Delta_h^{l -\l}
(\D^\lambda \mathcal E_\kappa^{l -1,d,m,D} f)\|_{L_p(\R^d)} = \\
\sum_{ \kappa: K^0 +1 \le \kappa \le k} \| \Delta_h^{l -\l}
(\D^\lambda \mathcal E_\kappa^{l -1,d,m,D} f)\|_{L_p(\R^d)} +
\sum_{ \kappa \ge \max(K^0 +1, k +1)} \| \Delta_h^{l -\l}
(\D^\lambda \mathcal E_\kappa^{l -1,d,m,D} f)\|_{L_p(\R^d)}.
\end{multline*}

При оценке первой суммы в правой части (2.1.14), пользуясь (1.1.3), (1.4.41),
выводим
\begin{multline*} \tag{2.1.15}
\sum_{ \kappa: K^0 +1 \le \kappa \le k} \| \Delta_h^{l -\l}
(\D^\lambda \mathcal E_\kappa^{l -1,d,m,D} f)\|_{L_p(\R^d)} \le \\
\sum_{ \kappa: K^0 +1 \le \kappa \le k} (l -\l)! \cdot \|h\|^{l -\l} \sum_{\mu \in \Z_{+ l -\l}^{\prime d}}
\| \D^\mu (\D^\lambda \mathcal E_{\kappa}^{l -1,d,m,D} f) \|_{L_p(\R^d)} = \\
\sum_{ \kappa: K^0 +1 \le \kappa \le k} (l -\l)! \cdot \|h\|^{l -\l}
\sum_{\mu \in \Z_{+ l -\l}^{\prime d}}
\| \D^{\lambda +\mu} \mathcal E_{\kappa}^{l -1,d,m,D} f) \|_{L_p(\R^d)} \le \\
\sum_{ \kappa: K^0 +1 \le \kappa \le k} c_{19} 2^{-k(l -\l)}
\sum_{\mu \in \Z_{+ l -\l}^{\prime d}}
c_{20} 2^{\kappa |\lambda +\mu|}
\Omega^{\prime l}(f, 2 2^{-\kappa})_{L_p(D)} = \\
c_{21} 2^{-k(l -\l)} \sum_{ \kappa: K^0 +1 \le \kappa \le k}
\sum_{\mu \in \Z_{+ l -\l}^{\prime d}}
2^{\kappa l} \Omega^{\prime l}(f, 2 2^{-\kappa})_{L_p(D)} \le \\
c_{22} 2^{-k(l -\l)} \sum_{ \kappa =1}^k 2^{\kappa l}
\Omega^{\prime l}(f, 2 2^{-\kappa})_{L_p(D)},\\
h \in (2^{-k +1} B^d), \ k \in \N, \ \lambda \in \Z_{+ \l}^{\prime d}.
\end{multline*}

Оценивая вторую сумму в правой части (2.1.14), в силу (1.4.41) получаем
\begin{multline*} \tag{2.1.16}
\sum_{ \kappa \ge \max(K^0 +1, k +1)} \| \Delta_h^{l -\l}
(\D^\lambda \mathcal E_\kappa^{l -1,d,m,D} f)\|_{L_p(\R^d)} \le \\
\sum_{ \kappa \ge \max(K^0 +1, k +1)} c_{8}
\| \D^\lambda \mathcal E_\kappa^{l -1,d,m,D} f\|_{L_p(\R^d)} \le \\
c_{8} \sum_{ \kappa \ge \max(K^0 +1, k +1)}
c_{14} 2^{\kappa |\lambda|} \Omega^{\prime l}(f, 2 2^{-\kappa})_{L_p(D)} \le
c_{23} \sum_{ \kappa = k +1}^\infty
2^{\kappa \l} \Omega^{\prime l}(f, 2 2^{-\kappa})_{L_p(D)},\\
h \in (2^{-k +1} B^d), \ k \in \N, \ \lambda \in \Z_{+ \l}^{\prime d}.
\end{multline*}

Соединяя (2.1.14), (2.1.15), (2.1.16), заключаем, что для $ k \in \N, \
\lambda \in \Z_{+ \l}^{\prime d} $ соблюдается неравенство
\begin{multline*}
\Omega^{l -\l}(\D^\lambda F, 2^{-k +1})_{L_p(\R^d)} =
\sup_{ h \in (2^{-k +1} B^d)}
\| \Delta_h^{l -\l} \D^\lambda F \|_{L_p(\R^d)} \le \\
c_{22} 2^{-k(l -\l)} \sum_{ \kappa =1}^k 2^{\kappa l}
\Omega^{\prime l}(f, 2 2^{-\kappa})_{L_p(D)} +
c_{23} \sum_{ \kappa = k +1}^\infty 2^{\kappa \l}
\Omega^{\prime l}(f, 2 2^{-\kappa})_{L_p(D)} \le \\
c_{24} \biggl(2^{-k(l -\l)}
\sum_{ \kappa =1}^k 2^{\kappa l}
\Omega^{\prime l}(f, 2 2^{-\kappa})_{L_p(D)} +
\sum_{ \kappa = k +1}^\infty 2^{\kappa \l}
\Omega^{\prime l}(f, 2 2^{-\kappa})_{L_p(D)}\biggr).
\end{multline*}

При $ \lambda \in \Z_{+ \l}^{\prime d},$ подставляя эту оценку в (2.1.13) и
применяя неравенство Гёльдера, приходим к соотношению
\begin{multline*} \tag{2.1.17}
\biggl(\int_I t^{-1 -\theta (\alpha -\l)}
(\Omega^{l -\l}(\D^\lambda F, t)_{L_p(\R^d)})^{\theta} dt\biggr)^{1/\theta} \le \\
\biggl(\sum_{ k \in \N} 2^{k \theta(\alpha -\l)}
(c_{24} \biggl(2^{-k(l -\l)}
\sum_{ \kappa =1}^k 2^{\kappa l}
\Omega^{\prime l}(f, 2 2^{-\kappa})_{L_p(D)} + \\
\sum_{ \kappa = k +1}^\infty 2^{\kappa \l}
\Omega^{\prime l}(f, 2 2^{-\kappa})_{L_p(D)} \biggr))^\theta\biggr)^{1 /\theta} = \\
c_{24} \biggl(\sum_{ k \in \N} 2^{k \theta(\alpha -\l)}
\biggl(2^{-k(l -\l)} \sum_{ \kappa =1}^k 2^{\kappa l}
\Omega^{\prime l}(f, 2 2^{-\kappa})_{L_p(D)} + \\
\sum_{ \kappa = k +1}^\infty 2^{\kappa \l}
\Omega^{\prime l}(f, 2 2^{-\kappa})_{L_p(D)}\biggr)^\theta\biggr)^{1 /\theta} \le \\
c_{24} \biggl(\sum_{ k \in \N} 2^{k \theta(\alpha -\l)}
c_{25} \biggl((2^{-k(l -\l)} \sum_{ \kappa =1}^k 2^{\kappa l}
\Omega^{\prime l}(f, 2 2^{-\kappa})_{L_p(D)})^\theta + \\
(\sum_{ \kappa = k +1}^\infty 2^{\kappa \l}
\Omega^{\prime l}(f, 2 2^{-\kappa})_{L_p(D)})^\theta\biggr)\biggr)^{1 /\theta} = \\
c_{26} \biggl(\sum_{ k \in \N}
\biggl(2^{-k(l -\alpha) \theta} (\sum_{ \kappa =1}^k 2^{\kappa l}
\Omega^{\prime l}(f, 2 2^{-\kappa})_{L_p(D)})^\theta + \\
2^{k \theta(\alpha -\l)} (\sum_{ \kappa = k +1}^\infty 2^{\kappa \l}
\Omega^{\prime l}(f, 2 2^{-\kappa})_{L_p(D)})^\theta\biggr)\biggr)^{1 /\theta} = \\
c_{26} \biggl(\sum_{ k \in \N}
2^{-k(l -\alpha) \theta} \biggl(\sum_{ \kappa =1}^k 2^{\kappa l}
\Omega^{\prime l}(f, 2 2^{-\kappa})_{L_p(D)}\biggr)^\theta + \\
\sum_{ k \in \N} 2^{k \theta(\alpha -\l)}
\biggl(\sum_{ \kappa = k +1}^\infty 2^{\kappa \l}
\Omega^{\prime l}(f, 2 2^{-\kappa})_{L_p(D)}\biggr)^\theta\biggr)^{1 /\theta}.
\end{multline*}

Далее, проводя оценку правой части (2.1.17), фиксируем $ \epsilon \in I, $
для которого выполняются неравенства
$ l -\alpha -\epsilon >0, \ \alpha -\l -\epsilon >0. $
Тогда при $ k \in \N, $ применяя неравенство Гёльдера, имеем
\begin{multline*} \tag{2.1.18}
\biggl(\sum_{ \kappa =1}^k 2^{\kappa l}
\Omega^{\prime l}(f, 2 2^{-\kappa})_{L_p(D)}\biggr)^\theta =
\biggl(\sum_{ \kappa =1}^k 2^{\kappa l -\alpha \kappa -\epsilon \kappa}
2^{\epsilon \kappa} 2^{\alpha \kappa}
\Omega^{\prime l}(f, 2 2^{-\kappa})_{L_p(D)}\biggr)^\theta \le \\
\biggl(\biggl(\sum_{ \kappa =1}^k 2^{\kappa (l -\alpha -\epsilon) \theta^\prime}\biggr)^{1 /\theta^\prime}
\biggl(\sum_{ \kappa =1}^k 2^{\kappa \epsilon \theta} 2^{\kappa \alpha \theta}
(\Omega^{\prime l}(f, 2 2^{-\kappa})_{L_p(D)})^\theta\biggr)^{1 /\theta}\biggr)^\theta \le \\
(c_{27} 2^{k(l -\alpha -\epsilon) \theta^\prime})^{\theta /\theta^\prime}
\sum_{ \kappa =1}^k 2^{\kappa \epsilon \theta} 2^{\kappa \alpha \theta}
(\Omega^{\prime l}(f, 2 2^{-\kappa})_{L_p(D)})^\theta = \\
c_{28} 2^{k(l -\alpha -\epsilon) \theta}
\sum_{ \kappa =1}^k 2^{\kappa \epsilon \theta} 2^{\kappa \alpha \theta}
(\Omega^{\prime l}(f, 2 2^{-\kappa})_{L_p(D)})^\theta.
\end{multline*}
Используя (2.1.18) и принимая во внимание вывод неравенства (2.1.3),
получаем
\begin{multline*} \tag{2.1.19}
\sum_{ k \in \N}
2^{-k(l -\alpha) \theta} \biggl(\sum_{ \kappa =1}^k 2^{\kappa l}
\Omega^{\prime l}(f, 2 2^{-\kappa})_{L_p(D)}\biggr)^\theta \le \\
\sum_{ k \in \N}
2^{-k(l -\alpha) \theta} c_{28} 2^{k(l -\alpha -\epsilon) \theta}
\sum_{ \kappa =1}^k 2^{\kappa \epsilon \theta} 2^{\kappa \alpha \theta}
(\Omega^{\prime l}(f, 2 2^{-\kappa})_{L_p(D)})^\theta = \\
c_{28} \sum_{ k \in \N} 2^{-k \epsilon \theta}
\sum_{ \kappa =1}^k 2^{\kappa \epsilon \theta} 2^{\kappa \alpha \theta}
(\Omega^{\prime l}(f, 2 2^{-\kappa})_{L_p(D)})^\theta = \\
c_{28} \sum_{ k \in \N} \sum_{ \kappa =1}^k
2^{-k \epsilon \theta} 2^{\kappa \epsilon \theta} 2^{\kappa \alpha \theta}
(\Omega^{\prime l}(f, 2 2^{-\kappa})_{L_p(D)})^\theta = \\
c_{28} \sum_{ k \in \N, \kappa \in \N: \kappa \le k}
2^{-k \epsilon \theta} 2^{\kappa \epsilon \theta} 2^{\kappa \alpha \theta}
(\Omega^{\prime l}(f, 2 2^{-\kappa})_{L_p(D)})^\theta = \\
c_{28} \sum_{\kappa \in \N} \biggl(\sum_{k = \kappa}^\infty
2^{-k \epsilon \theta}\biggr) 2^{\kappa \epsilon \theta} 2^{\kappa \alpha \theta}
(\Omega^{\prime l}(f, 2 2^{-\kappa})_{L_p(D)})^\theta \le \\
c_{28} \sum_{\kappa \in \N} c_{29} 2^{-\kappa \epsilon \theta}
2^{\kappa \epsilon \theta} 2^{\kappa \alpha \theta}
(\Omega^{\prime l}(f, 2 2^{-\kappa})_{L_p(D)})^\theta = \\
c_{30} \sum_{\kappa \in \N} (2^{\kappa \alpha}
\Omega^{\prime l}(f, 2 2^{-\kappa})_{L_p(D)})^\theta \le \\
c_{31} \int_{\R_+} t^{-1 -\theta \alpha}
(\Omega^{\prime l}(f, t)_{L_p(D)})^\theta dt.
\end{multline*}

Продолжая оценку правой части (2.1.17), при $ k \in \N, $ пользуясь с учётом
(2.1.3) неравенством Гёльдера. выводим
\begin{multline*} \tag{2.1.20}
\biggl(\sum_{ \kappa = k +1}^\infty 2^{\kappa \l}
\Omega^{\prime l}(f, 2 2^{-\kappa})_{L_p(D)}\biggr)^\theta =
\biggl(\sum_{ \kappa = k +1}^\infty 2^{\kappa \l} 2^{\kappa \epsilon}
2^{-\kappa \alpha} 2^{-\kappa \epsilon} 2^{\kappa \alpha}
\Omega^{\prime l}(f, 2 2^{-\kappa})_{L_p(D)}\biggr)^\theta \le \\
\biggl(\biggl(\sum_{ \kappa = k +1}^\infty 2^{-\kappa(\alpha -\l -\epsilon) \theta^\prime}\biggr)^{1 /\theta^\prime}
\biggl(\sum_{ \kappa = k +1}^\infty 2^{-\kappa \epsilon \theta} 2^{\kappa \alpha \theta}
(\Omega^{\prime l}(f, 2 2^{-\kappa})_{L_p(D)})^\theta\biggr)^{1 /\theta}\biggr)^\theta \le \\
(c_{32} 2^{-k (\alpha -\l -\epsilon) \theta^\prime})^{\theta /\theta^\prime}
\sum_{ \kappa = k +1}^\infty 2^{-\kappa \epsilon \theta} 2^{\kappa \alpha \theta}
(\Omega^{\prime l}(f, 2 2^{-\kappa})_{L_p(D)})^\theta = \\
c_{33} 2^{-k(\alpha -\l -\epsilon) \theta}
\sum_{ \kappa = k +1}^\infty 2^{-\kappa \epsilon \theta} 2^{\kappa \alpha \theta}
(\Omega^{\prime l}(f, 2 2^{-\kappa})_{L_p(D)})^\theta.
\end{multline*}

Пользуясь (2.1.20), находим, что
\begin{multline*} \tag{2.1.21}
\sum_{ k \in \N} 2^{k \theta(\alpha -\l)}
\biggl(\sum_{ \kappa = k +1}^\infty 2^{\kappa \l}
\Omega^{\prime l}(f, 2 2^{-\kappa})_{L_p(D)}\biggr)^\theta \le \\
\sum_{ k \in \N} 2^{k \theta(\alpha -\l)}
c_{33} 2^{-k(\alpha -\l -\epsilon) \theta}
\sum_{ \kappa = k +1}^\infty 2^{-\kappa \epsilon \theta} 2^{\kappa \alpha \theta}
(\Omega^{\prime l}(f, 2 2^{-\kappa})_{L_p(D)})^\theta = \\
c_{33} \sum_{ k \in \N} 2^{k \epsilon \theta}
\sum_{ \kappa = k +1}^\infty 2^{-\kappa \epsilon \theta} 2^{\kappa \alpha \theta}
(\Omega^{\prime l}(f, 2 2^{-\kappa})_{L_p(D)})^\theta \le \\
c_{33} \sum_{ k \in \N} 2^{k \epsilon \theta}
\sum_{ \kappa = k}^\infty 2^{-\kappa \epsilon \theta}
(2^{\kappa \alpha} \Omega^{\prime l}(f, 2 2^{-\kappa})_{L_p(D)})^\theta = \\
c_{33} \sum_{ k \in \N, \kappa \in \N: \kappa \ge k} 2^{k \epsilon \theta}
2^{-\kappa \epsilon \theta}
(2^{\kappa \alpha} \Omega^{\prime l}(f, 2 2^{-\kappa})_{L_p(D)})^\theta = \\
c_{33} \sum_{\kappa \in \N} \biggl(\sum_{k =1}^\kappa 2^{k \epsilon \theta}\biggr)
2^{-\kappa \epsilon \theta}
(2^{\kappa \alpha} \Omega^{\prime l}(f, 2 2^{-\kappa})_{L_p(D)})^\theta \le \\
c_{33} \sum_{\kappa \in \N} c_{34} 2^{\kappa \epsilon \theta}
2^{-\kappa \epsilon \theta}
(2^{\kappa \alpha} \Omega^{\prime l}(f, 2 2^{-\kappa})_{L_p(D)})^\theta = \\
c_{35} \sum_{\kappa \in \N} (2^{\kappa \alpha}
\Omega^{\prime l}(f, 2 2^{-\kappa})_{L_p(D)})^\theta \le \\
c_{36} \int_{\R_+} t^{-1 -\theta \alpha}
(\Omega^{\prime l}(f, t)_{L_p(D)})^\theta dt (\text{ см. вывод } (2.1.3)).
\end{multline*}

Соединяя (2.1.17), (2.1.19), (2.1.21), приходим к неравенству
\begin{multline*}
\biggl(\int_I t^{-1 -\theta (\alpha -\l)}
(\Omega^{l -\l}(\D^\lambda F, t)_{L_p(\R^d)})^{\theta} dt\biggr)^{1/\theta} \le \\
c_{26} \biggl(c_{37} \int_{\R_+} t^{-1 -\theta \alpha}
(\Omega^{\prime l}(f, t)_{L_p(D)})^\theta dt\biggr)^{1 /\theta} \le \\
c_{38} \biggl(\int_{\R_+} t^{-1 -\theta \alpha}
(\Omega^{\prime l}(f, t)_{L_p(D)})^\theta dt\biggr)^{1 /\theta},
\end{multline*}
которое в сочетании с (2.1.12) влечёт неравенство
\begin{multline*} \tag{2.1.22}
\biggl(\int_{\R_+} t^{-1 -\theta (\alpha -\l)}
(\Omega^{l -\l}(\D^\lambda F, t)_{L_p(\R^d)})^{\theta} dt\biggl)^{1 /\theta} \le \\
c_{39} \biggl(\biggl(\int_{\R_+} t^{-1 -\theta \alpha}
(\Omega^{\prime l}(f, t)_{L_p(D)})^\theta dt\biggr)^{1 /\theta} +
\| \D^\lambda F \|_{L_p(\R^d)}\biggr) = \\
c_{39} \biggl(\biggl(\int_{\R_+} t^{-1 -\theta \alpha}
(\Omega^{\prime l}(f, t)_{L_p(D)})^\theta dt\biggr)^{1 /\theta} +
\| \D^\lambda (E^{\prime d,\alpha,p,\theta,D} f) \|_{L_p(\R^d)}\biggr),\\
\lambda \in \Z_{+ \l}^{\prime d}.
\end{multline*}
Из (2.1.8), (2.1.10), (2.1.22) следует (2.1.7). Соединяя (2.1.6), (2.1.7) и
учитывая (2.1.4), приходим к (2.1.2) при $ \theta \ne \infty. $
При $ \theta = \infty $ доказательство теоремы проводится по той же схеме с
заменой в соответствующих местах операции суммирования на операцию
взятия супремума или максимума. $ \square $

Следствие

В условиях теоремы 2.1.1 имеет место включение
\begin{equation*} \tag{2.1.23}
(B_{p, \theta}^\alpha)^\prime(D) \subset (B_{p, \theta}^\alpha)^{\l}(D),
\end{equation*}
и для любой функции $ f \in (B_{p, \theta}^\alpha)^\prime(D) $
выполняется неравенство
\begin{equation*} \tag{2.1.24}
\| f\|_{(B_{p, \theta}^\alpha)^{\l}(D)} \le c_1
\| f\|_{(B_{p, \theta}^\alpha)^\prime(D)}.
\end{equation*}

Для получения (2.1.23), (2.1.24) достаточно применить теорему 2.1.1
и неравенство
\begin{equation*}
\| (E^{d,\alpha,p,\theta,D} f) \mid_{D} \|_{(B_{p, \theta}^\alpha)^{\l}(D)} \le
\| E^{d,\alpha,p,\theta,D} f\|_{(B_{p, \theta}^\alpha)^{\l}(\R^d)}.
\end{equation*}

Из (1.1.5), (1.1.6) и (2.1.23), (2.1.24) вытекает, что в условиях теоремы 2.1.1

имеет место равенство $ (B_{p, \theta}^\alpha)^\prime(D) =
(B_{p, \theta}^\alpha)^{\l}(D) $ и нормы
$ \| \cdot \|_{(B_{p, \theta}^\alpha)^{\l}(D)},
\| \cdot \|_{(B_{p, \theta}^\alpha)^\prime(D)} $ эквивалентны.

Отметим ещё, что с помощью теоремы 2.1.1 и теоремы 2.1.1 из [9]
(см. также [10]), опираясь на п. 5.6.2 из [8], легко проверить, что в условиях
теоремы 2.1.1 справедливо равенство $ (B_{p, \theta}^\alpha)^\prime(D) =
(B_{p, \theta}^{(\alpha, \ldots, \alpha)})^{\prime}(D) $ и нормы
$ \| \cdot \|_{(B_{p, \theta}^{(\alpha, \ldots, \alpha)})^{\prime}(D)},
\| \cdot \|_{(B_{p, \theta}^\alpha)^\prime(D)} $ эквивалентны. Заметим
однако, что сами по себе эти факты не позволяют получить результаты,
установленные в \S 3, \S 4 и \S 5, как следствия соответствующих утверждений из [9].
\bigskip

2.2. В этом пункте приведём используемые в следующих параграфах
соотношения между нормами образов и прообразов при некоторых отображениях
рассматриваемых пространств.

При $ d \in \N $ для $ \delta \in \R_+ $ и $ x^0 \in \R^d $ обозначим через
$ h_{\delta, x^0} $ отображение, которое каждой функции $ f, $ заданной на
некотором множестве $ S \subset \R^d, $  ставит в соответствие функцию
$ h_{\delta, x^0} f, $ определяемую на множестве $ \{ x \in \R^d: x^0 +\delta
x \in S\} = \delta^{-1} (S -x^0) $ равенством $ (h_{\delta, x^0} f)(x) =
f(x^0 +\delta x). $ Так как для $ \delta \in \R_+, \ x^0 \in \R^d $
отображение $ \R^d \ni x \mapsto x^0 +\delta x \in \R^d $ ---
взаимно однозначно, то отображение $ h_{\delta, x^0} $ является
биекцией на себя  множества всех функций с областью определения
в $ \R^d. $ При этом обратное  отображение $ h_{\delta, x^0}^{-1} $
задаётся равенством
\begin{equation*} \tag{2.2.1}
(h_{\delta, x^0}^{-1} f)(x) = f(\delta^{-1} (x -x^0)) =
(h_{\delta^\prime, x^{\prime 0}} f)(x)  \text{ с } \delta^\prime = \delta^{-1}, \
x^{\prime 0} =-\delta^{-1} x^0.
\end{equation*}

Отметим, что при $ 1 \le p \le \infty, \
l \in \Z_+ $ для $ f \in W_p^l(x^0 +\delta D), \ \lambda \in \Z_{+ l}^{\prime d}, $
где $ D $ -- область в $ \R^d, \ \delta \in \R_+, \ x^0 \in \R^d, $ имеет место
равенство
\begin{multline*} \tag{2.2.2}
\| \D^\lambda (h_{\delta,x^0} f)\|_{L_p(D)} =
\biggl(\int_D |\D^\lambda (f(x^0 +\delta y))|^p dy\biggr)^{1/p} = \\
\biggl(\int_D | \delta^{|\lambda|} (\D^\lambda f)(x^0 +\delta y)|^p dy\biggr)^{1/p}
= \delta^{|\lambda|} \biggl(\int_{ x^0 +\delta D} |\D^\lambda f(x)|^p \delta^{-d} dx\biggr)^{1/p} = \\
\delta^{|\lambda| -d /p} \biggl(\int_{x^0 +\delta D} |\D^\lambda f(x)|^p dx\biggr)^{1/p} =
\delta^{|\lambda| -d /p} \|\D^\lambda f\|_{L_p(x^0 +\delta D)},
\end{multline*}
а, следовательно, для $ f \in L_p(D) $ выполняется равенство
\begin{equation*} \tag{2.2.3}
\| h_{\delta,x^0}^{-1} f\|_{L_p(x^0 +\delta D)} = \delta^{d /p}
\| h_{\delta,x^0} h_{\delta,x^0}^{-1} f \|_{L_p(D)} =
\delta^{d /p} \|f\|_{L_p(D)}.
\end{equation*}

Лемма  2.2.1

Пусть $ d \in \N, \ l \in \Z_+, \ D $ -- область в $ \R^d, \ 1 \le p < \infty, \
\delta \in \R_+, \ x^0 \in \R^d. $ Тогда при $ t \in \R_+ $ для $ f \in
L_p(x^0 +\delta D) $ справедливо равенство
\begin{equation*} \tag{2.2.4}
\Omega^{\prime l}((h_{\delta, x^0} f),t)_{L_p(D)} =
\delta^{-d /p} \Omega^{\prime l}(f, \delta t)_{L_p(x^0 +\delta D)}.
\end{equation*}

Доказательство.

В условиях леммы для $ f \in L_p(x^0 +\delta D), $
делая замену переменных $ \eta = \delta^{-1} \xi, \ y = \delta^{-1}(x -x^0), $
находим, что соблюдается равенство
\begin{multline*}
\int_{ t B^d} \int_{ D_{l \eta}} | (\Delta_{\eta}^{l}
(h_{\delta, x^0} f))(y)|^p dy d\eta = \\
\int_{\{(\eta,y): \eta \in t B^d, \ y \in D_{l \eta}\}}
\biggl|\sum_{i =0}^{l} C_l^i (-1)^{l -i}
(h_{\delta, x^0} f)(y +i \eta)\biggr|^p d\eta dy = \\
\int_{\{(\eta,y): \eta \in t B^d, \ y \in D_{l \eta}\}}
\biggl|\sum_{i =0}^{l} C_{l}^i (-1)^{l -i}
f(x^0 +\delta y +i \delta \eta)\biggr|^p d\eta dy = \\
\int_{\{ (\xi,x): \xi \in \delta t B^d, \ x \in (x^0 +\delta D)_{l \xi} \}}
\biggl|\sum_{i =0}^{l} C_{l}^i (-1)^{l -i}
f(x +i \xi)\biggr|^p \delta^{-d} \ \delta^{-d} d\xi dx = \\
\delta^{-d} \delta^{-d} \int_{\delta t B^d} \int_{ (x^0 +\delta D)_{l \xi}}
|\Delta_{\xi}^{l} f(x)|^p dx d\xi.
\end{multline*}
Откуда выводим
\begin{multline*}
\Omega^{\prime l}((h_{\delta, x^0} f),t)_{L_p(D)} =
((2 t)^{-d} \int_{ t B^d} \int_{D_{l \eta}}
| (\Delta_{\eta}^{l} (h_{\delta, x^0} f))(y)|^p dy d\eta)^{1 /p} = \\
((2 t)^{-d} \delta^{-d} \delta^{-d} \int_{\delta t B^d}
\int_{ (x^0 +\delta D)_{l \xi}}
|\Delta_{\xi}^{l} f(x)|^p dx d\xi)^{1 /p} = \\
\delta^{-d /p} ((2 \delta t)^{-d} \int_{\delta t B^d}
\int_{ (x^0 +\delta D)_{l \xi}}
|\Delta_{\xi}^{l} f(x)|^p dx d\xi)^{1 /p} = \\
\delta^{-d /p} \Omega^{\prime l}(f, \delta t)_{L_p(x^0 +\delta D)}, \
t \in \R_+, \ f \in L_p(x^0 +\delta D),
\end{multline*}
что совпадает с (2.2.4). $ \square $

Лемма 2.2.2

Пусть $ d \in \N, \ D $ -- область в $ \R^d, \ \alpha \in \R_+, \ 1 \le p < \infty, \
1 \le \theta \le \infty, \ \delta \in \R_+, \ x^0 \in \R^d. $ Тогда существуют
константы $ c_1(d,\alpha,p,\delta) > 0, \ c_2(d,\alpha,p,\delta) > 0 $
такие, что для любой функции $ f \in (B_{p,\theta}^\alpha)^\prime(x^0 +\delta D) $
соблюдается неравенство
\begin{equation*} \tag{2.2.5}
\| h_{\delta, x^0} f \|_{(B_{p,\theta}^\alpha)^\prime(D)} \le
c_1 \| f \|_{(B_{p,\theta}^\alpha)^\prime(x^0 +\delta D)},
\end{equation*}
а для $ f \in (B_{p,\theta}^\alpha)^\prime(D) $ выполняется неравенство
\begin{equation*} \tag{2.2.6}
\| h_{\delta, x^0}^{-1} f \|_{(B_{p,\theta}^\alpha)^\prime(x^0 +\delta D)} \le
c_2 \| f \|_{(B_{p,\theta}^\alpha)^\prime(D)}.
\end{equation*}

Доказательство.

В условиях леммы, полагая $ l = l(\alpha), $ для
$ f \in (B_{p,\theta}^\alpha)^\prime(x^0 +\delta D) $ в силу (2.2.4) функция
$ t^{-1 -\theta \alpha} (\Omega^{\prime l}((h_{\delta, x^0}
f),t)_{L_p(D)})^\theta $ суммируема на $ \R_+, $ и справедливо соотношение
\begin{multline*} \tag{2.2.7}
\biggl(\int_{\R_+} t^{-1 -\theta \alpha} (\Omega^{\prime l}((h_{\delta, x^0}
f),t)_{L_p(D)})^\theta dt\biggr)^{1 /\theta} = \\
\biggl(\int_{\R_+} t^{-1 -\theta \alpha} (\delta^{-d /p}
\Omega^{\prime l}(f, \delta t)_{L_p(x^0 +\delta D)})^\theta dt\biggr)^{1 /\theta} = \\
\delta^{-d /p} \biggl(\int_{\R_+} t^{-1 -\theta \alpha}
(\Omega^{\prime l}(f, \delta t)_{L_p(x^0 +\delta D)})^\theta dt\biggr)^{1 /\theta} = \\
\delta^{-d /p} \biggl(\int_{\R_+} (\delta^{-1} \tau)^{-1 -\theta \alpha}
(\Omega^{\prime l}(f, \delta \delta^{-1} \tau)_{L_p(x^0 +\delta D)})^\theta
\delta^{-1} d\tau\biggr)^{1 /\theta} = \\
\delta^{-d /p} \biggl(\int_{\R_+} \delta^{\theta \alpha} \tau^{-1 -\theta \alpha}
(\Omega^{\prime l}(f, \tau)_{L_p(x^0 +\delta D)})^\theta d\tau\biggr)^{1 /\theta} = \\
\delta^{-d /p} \delta^{\alpha} \biggl(\int_{\R_+} \tau^{-1 -\theta \alpha}
(\Omega^{\prime l}(f, \tau)_{L_p(x^0 +\delta D)})^\theta d\tau\biggr)^{1 /\theta}.
\end{multline*}
Объединяя равенства (2.2.7) и (2.2.2) при $ \lambda =0, $ приходим к (2.2.5).
Наконец, для $ f \in (B_{p,\theta}^\alpha)^\prime(D), $ применяя (2.2.1),
(2.2.5), получаем, что $ h_{\delta,x^0}^{-1} f \in
(B_{p,\theta}^\alpha)^\prime(x^0 +\delta D) $ и
соблюдается (2.2.6). $ \square $
\bigskip

\centerline{\S 3. Восстановление функций вместе с их производными}
\centerline{из изотропных классов Никольского и Бесова}
\centerline{по значениям функций в заданном числе точек}
\bigskip

3.1. В этом пункте будет установлена оценка сверху величины
наилучшей точности линейного восстановления в пространстве
$ W_q^{\mathpzc m}(D) $ по значениям в $ n $ точках функций $ f $
из классов $ (\mathcal H_p^\alpha)^\prime(D) $ и
$ (\mathcal B_{p,\theta}^\alpha)^\prime(D), $ определённых в
ограниченной области $ D \ \e$-типа.

Но сначала опишем постановку задачи.

Пусть $ C(T) $ --- пространство непрерывных вещественных
функций на топологическом пространстве $ T, $ и $ X $ --- некоторое
линейное нормированное пространство вещественных функций,
определённых на $ T. $
Пусть ещё $ \mathcal K \subset C(T) \cap X $ -- некоторый класс функций.

Для $ n \in \N $ через $ \Phi_n(C(T)) $ обозначим совокупность всех отображений
$ \phi: C(T) \mapsto \R^n, $ для каждого из которых существует набор точек
$ \{t^j \in T, \ j =1,\ldots,n\} $ такой, что значение
$ \phi(f) = (f(t^1), \ldots, f(t^n)), \ f \in C(T), $ а также через
$ \mathcal A^n(X) (\overline{\mathcal A}^n(X)) $ обозначим
множество всех отображений (всех линейных отображений)
$ A: \R^n \mapsto X. $

Тогда при $ n \in \N $ положим
\begin{equation*}
\sigma_n(\mathcal K,C(T),X) = \inf_{A \in \mathcal A^n(X), \ \phi \in \Phi_n(C(T))}
\sup_{f \in \mathcal K} \| f -A \circ \phi(f) \|_X,
\end{equation*}
а
\begin{equation*}
\overline \sigma_n(\mathcal K,C(T),X) = \inf_{ A \in \overline{\mathcal A}^n(X), \
\phi \in \Phi_n(C(T))} \sup_{f \in \mathcal K} \| f -A \circ \phi(f) \|_X.
\end{equation*}

Для доказательства основного утверждения этого пункта, теоремы
3.1.2, потребуется взятая из [12] лемма 3.1.1.

Лемма 3.1.1

Пусть $ d \in \N $ и $ l \in \Z_+. $ Тогда справедливы следующие
утверждения:

1) отображение $ \J^{l,d}: \mathcal P^{l,d} \mapsto \R^{\card \Z_{+ l}^d}, $
которое каждому полиному $ f \in \mathcal P^{l,d} $ ставит в соответствие
набор вещественных чисел $ y = \{y_\lambda = f(\lambda), \ \lambda \in \Z_{+ l}^d\}, $
является изоморфизмом;

2) для $ \rho, \ \sigma, \ \tau \in \R_+ $ существует константа
$ c_1(d,l,\rho,\sigma,\tau) >0 $ такая, что для любого $ \delta \in \R_+, $
для любых точек $ x^0, \ y^0 \in \R^d $ таких, что $ y^0 \in
(x^0 +\sigma \delta B^d),$ и любого множества $ Q \subset
(y^0 +\tau \delta B^d), $ для любого полинома $ f \in \mathcal P^{l,d} $
справедлива оценка
\begin{equation*} \tag{3.1.1}
\sup_{x \in Q} | f(x)| \le c_1 \max_{\lambda \in \Z_{+ l}^d} | f(x^0
+\delta \rho \lambda)|.
\end{equation*}

Теорема 3.1.2

Пусть $ d \in \N, \ \alpha \in \R_+, \ D \subset \R^d $ -- ограниченная область
$ \e $-типа, $ \mathpzc m \in \Z_+, \ 1 < p < \infty, \ 1 \le q \le \infty $
таковы, что выполняются условия (1.4.81) при $ \lambda =0 $ и
\begin{equation*} \tag{3.1.2}
\alpha -\mathpzc m -(d /p -d /q)_+ >0.
\end{equation*}
Пусть ещё $ T = D, \ X = W_q^{\mathpzc m}(D),
\ \mathcal K = (\mathcal H_p^\alpha)^\prime(D), \
(\mathcal B_{p,\theta}^\alpha)^\prime(D), $
где $ \theta \in \R: \theta \ge 1. $ Тогда существуют константы
$ c_2(d,\alpha,D,p,q,\mathpzc m) >0, \ n_0 \in \N $ такие, что для любого
$ n \in \N: n \ge n_0, $ соблюдается неравенство
\begin{equation*} \tag{3.1.3}
\overline \sigma_n(\mathcal K,C(T),X) \le c_2 n^{-(\alpha -\mathpzc m) /d +(1 /p -1 /q)_+}.
\end{equation*}

Доказательство.

Учитывая включение (1.1.4), доказательство достаточно провести
лишь в случае $ \mathcal K = (\mathcal H_p^\alpha)^\prime(D). $ Так что
рассмотрим этот случай.

Прежде всего отметим, что соблюдение условий (1.4.81) при $ \lambda =0 $ и
(3.1.2) обеспечивает включение $ \mathcal K \subset C(D) \cap W_q^{\mathpzc m}(D) $
(см. теорему 1.4.13), и, следовательно,
величина $ \overline \sigma_n(\mathcal K,C(T),X) $ определена.

В условиях теоремы фиксируем $ l = l(\alpha) \in \N, \ m \in \N: \mathpzc m \le m, $ и
константу $ K^0 = K^0(d,D,\e) \in \Z_+ $ (см. п. 1.2., определение 1 при
$ \alpha = \e $). При $ k \in \Z_+: k \ge K^0, $ положим
$$
n(k,D) = n^{l -1,d}(k,D) = (\card \Z_{+ l -1}^d) \cdot \card \Nu_k(D),
$$
где $ \Nu_k(D) = \{\nu \in \Z^d: \overline Q_{k,\nu}^d \subset D\}. $

Пусть теперь $ n \in \N: n \ge n(K^0,D). $ Тогда выберем $ k \in \Z_+: k \ge
K^0, $ так, чтобы имело место соотношение
\begin{equation*} \tag{3.1.4}
n(k,D) \le n \le n(k +1,D).
\end{equation*}

Далее, для этого $ k $ построим систему точек
\begin{equation*}
x_{k,\nu}^\mu = 2^{-k} (\nu +\frac{1} {4} \e
+\frac {1} {2} l^{-1} \mu) \in Q_{k,\nu}^d, \ \nu \in \Nu_k(D), \
\mu \in \Z_{+ l -1}^d,
\end{equation*}
и определим отображение $ \phi \in \Phi_{n(k,D)}(C(D)), $ полагая для $ f \in C(D) $ значение
$$
\phi(f) = \{ f(x_{k,\nu}^\mu), \ \nu \in \Nu_k(D), \ \mu \in \Z_{+ l -1}^d\}.
$$

Пользуясь леммой 3.1.1, возьмём систему полиномов $ \{\pi_\mu^{l -1,d} \in
\mathcal P^{l -1,d}, \ \mu \in \Z_{+ l -1}^d\} $ такую, что для
$ \mu, \mu^\prime \in \Z_{+ l -1}^d $ значение
\begin{equation*}
\pi_\mu^{l -1,d}(\mu^\prime) = \begin{cases} 1, \text{ при } \mu^\prime = \mu; \\
0, \text{ при } \mu^\prime \ne \mu, \end{cases}
\end{equation*}
и обозначим через $ R_{k,\nu,\mu}^{l -1,d} $ полином
из $ \mathcal P^{l -1,d} $ равный $ R_{k,\nu,\mu}^{l -1,d}(x) =
\pi_\mu^{l -1,d}(2l 2^k (x -x_{k,\nu}^0)), \ \nu \in \Nu_k(D), \
\mu \in \Z_{+ l -1}^d. $

Нетрудно видеть, что при $ \nu \in \Nu_k(D), \ \mu, \mu^\prime \in \Z_{+ l -1}^d $
значение $ R_{k,\nu,\mu}^{l -1,d} (x_{k,\nu}^{\mu^\prime}) =
\pi_\mu^{l -1,d}(\mu^\prime). $

Определим для $ \nu \in \Nu_k(D) $ линейный оператор
$ R_{k,\nu}^{l -1,d}: \R^{n(k,D)} \mapsto \mathcal P^{l -1,d}, $
полагая для $ t = \{t_{\nu^\prime,\mu} \in \R, \ \nu^\prime \in \Nu_k(D), \
\mu \in \Z_{+ l -1}^d\} $ значение
\begin{equation*}
R_{k,\nu}^{l -1,d} t = \sum_{\mu \in \Z_{+ l -1}^d}
t_{\nu,\mu} R_{k,\nu,\mu}^{l -1,d}.
\end{equation*}

Принимая во внимание сказанное выше,  получаем, что при
$ \nu \in \Nu_k(D) $ для $ f \in C(D) $ и $ \mu \in \Z_{+ l -1}^d $ имеет
место равенство
\begin{equation*} \tag{3.1.5}
(R_{k,\nu}^{l -1,d} \circ \phi(f))  (x_{k,\nu}^\mu)
= f(x_{k,\nu}^\mu).
\end{equation*}

Теперь определим оператор $ A \in \overline{\mathcal A}^{n(k,D)}(X), $
полагая для $ t = \{t_{\nu,\mu} \in \R, \ \nu \in \Nu_k(D), \ \mu \in \Z_{+ l -1}^d\} $
значение
\begin{equation*}
A t = \biggl(\sum_{\nu \in N_k^{d,m,D}}
(R_{k, \nu_k^D(\nu)}^{l -1,d} t) g_{k, \nu}^{m,d}\biggr)
\text{ (см. п. 1.4.).}
\end{equation*}

Тогда для $ f \in \mathcal K $ при $ \lambda \in \Z_{+ \mathpzc m}^d $ имеем
\begin{multline*} \tag{3.1.6}
\| \D^\lambda f -\D^\lambda A \circ \phi(f) \|_{L_q(D)} \le
\| \D^\lambda f -\D^\lambda E_k^{l -1,d,m,D}(f) \|_{L_q(D)} \\
+\| \D^\lambda E_k^{l -1,d,m,D} f -\D^\lambda A \circ \phi(f)\|_{L_q(D)}.
\end{multline*}

Оценку первого слагаемого в правой части (3.1.6) даёт (1.4.72).
При оценке второго слагаемого в правой части (3.1.6) будем действовать так же,
как при доказательстве предложения 1.4.9, опираясь на те же вспомогательные
объекты, что и там. Используя (1.4.33), имеем
\begin{multline*} \tag{3.1.7}
\| \D^\lambda E_k^{l -1,d,m,D} f -\D^\lambda A \circ \phi(f)\|_{L_q(D)} = \\
\biggl\|\D^\lambda E_k^{l -1,d,m,D} f -\D^\lambda \biggl(\sum_{\nu \in N_k^{d,m,D}}
(R_{k, \nu_k^D(\nu)}^{l -1,d} \phi(f)) g_{k, \nu}^{m,d}\biggr) \biggr\|_{L_q(D)} \le \\
\biggl\|\D^\lambda E_k^{l -1,d,m,D} f -\D^\lambda \biggl(\sum_{\nu \in N_k^{d,m,D}}
(R_{k, \nu_k^D(\nu)}^{l -1,d} \phi(f)) g_{k, \nu}^{m,d}\biggr) \biggr\|_{L_q(\R^d)} = \\
\biggl\|\D^\lambda ((E_k^{l -1,d,m,D} f) -\sum_{\nu \in N_k^{d,m,D}}
(R_{k, \nu_k^D(\nu)}^{l -1,d} \phi(f)) g_{k, \nu}^{m,d}) \biggr\|_{L_q(\R^d)} = \\
\biggl\|\D^\lambda \biggl(\sum_{\nu \in N_k^{d,m,D}}
(S_{k, \nu_k^D(\nu)}^{l -1,d} f) g_{k, \nu}^{m,d}
-\sum_{\nu \in N_k^{d,m,D}}
(R_{k, \nu_k^D(\nu)}^{l -1,d} \phi(f)) g_{k, \nu}^{m,d}\biggr) \biggr\|_{L_q(\R^d)} = \\
\biggl\|\D^\lambda \biggl(\sum_{\nu \in N_k^{d,m,D}}
((S_{k, \nu_k^D(\nu)}^{l -1,d} f)
-(R_{k, \nu_k^D(\nu)}^{l -1,d} \phi(f))) g_{k, \nu}^{m,d}\biggr) \biggr\|_{L_q(\R^d)} = \\
\biggl\| \sum_{\nu \in N_k^{d,m,D}}
\D^\lambda (((S_{k, \nu_k^D(\nu)}^{l -1,d} f)
-(R_{k, \nu_k^D(\nu)}^{l -1,d} \phi(f))) g_{k, \nu}^{m,d}) \biggr\|_{L_q(\R^d)} = \\
\biggl\| \sum_{\nu \in N_k^{d,m,D}} \sum_{\mu \in \Z_+^d(\lambda)} C_\lambda^\mu
\D^\mu ((S_{k, \nu_k^D(\nu)}^{l -1,d} f)
-(R_{k, \nu_k^D(\nu)}^{l -1,d} \phi(f)))
\D^{\lambda -\mu} g_{k, \nu}^{m,d} \biggr\|_{L_q(\R^d)} = \\
\biggl\| \sum_{\mu \in \Z_+^d(\lambda)} C_\lambda^\mu
\sum_{\nu \in N_k^{d,m,D}}
\D^\mu ((S_{k, \nu_k^D(\nu)}^{l -1,d} f)
-(R_{k, \nu_k^D(\nu)}^{l -1,d} \phi(f)))
\D^{\lambda -\mu} g_{k, \nu}^{m,d} \biggr\|_{L_q(\R^d)} \le \\
\sum_{\mu \in \Z_+^d(\lambda)} C_\lambda^\mu
\biggl\| \sum_{\nu \in N_k^{d,m,D}}
\D^\mu ((S_{k, \nu_k^D(\nu)}^{l -1,d} f)
-(R_{k, \nu_k^D(\nu)}^{l -1,d} \phi(f)))
\D^{\lambda -\mu} g_{k, \nu}^{m,d} \biggr\|_{L_q(\R^d)}.
\end{multline*}

Оценивая правую часть (3.1.7), при $ \mu \in \Z_+^d(\lambda) $ с учётом
(1.4.14), (1.4.12) получаем
\begin{multline*} \tag{3.1.8}
\biggl\| \sum_{\nu \in N_k^{d,m,D}}
\D^\mu ((S_{k, \nu_k^D(\nu)}^{l -1,d} f) -
(R_{k, \nu_k^D(\nu)}^{l -1,d} \phi(f)))
\D^{\lambda -\mu} g_{k, \nu}^{m,d} \biggr\|_{L_q(\R^d)}^q = \\
\int_{\R^d} \biggl| \sum_{\nu \in N_k^{d,m,D}}
\D^\mu ((S_{k, \nu_k^D(\nu)}^{l -1,d} f) -
(R_{k, \nu_k^D(\nu)}^{l -1,d} \phi(f)))
\D^{\lambda -\mu} g_{k, \nu}^{m,d} \biggr|^q dx = \\
\int_{G_k^{d,m,D}} \biggl| \sum_{\nu \in N_k^{d,m,D}}
\D^\mu ((S_{k, \nu_k^D(\nu)}^{l -1,d} f) -
(R_{k, \nu_k^D(\nu)}^{l -1,d} \phi(f)))
\D^{\lambda -\mu} g_{k, \nu}^{m,d} \biggr|^q dx = \\
\sum_{\substack{n \in \Z^d:\\ Q_{k,n}^d \cap G_k^{d,m,D} \ne \emptyset}}
\int_{Q_{k,n}^d} \biggl| \sum_{\nu \in N_k^{d,m,D}}
\D^\mu ((S_{k, \nu_k^D(\nu)}^{l -1,d} f) -
(R_{k, \nu_k^D(\nu)}^{l -1,d} \phi(f)))
\D^{\lambda -\mu} g_{k, \nu}^{m,d} \biggr|^q dx = \\
\sum_{\substack{n \in \Z^d:\\ Q_{k,n}^d \cap G_k^{d,m,D} \ne \emptyset}}
\int_{Q_{k,n}^d} \biggl| \sum_{\nu \in N_k^{d,m,D}: Q_{k,n}^d \cap
\supp g_{k, \nu}^{m,d} \ne \emptyset}
\D^\mu ((S_{k, \nu_k^D(\nu)}^{l -1,d} f) -\\
(R_{k, \nu_k^D(\nu)}^{l -1,d} \phi(f)))
\D^{\lambda -\mu} g_{k, \nu}^{m,d} \biggr|^q dx = \\
\sum_{n \in \Z^d: Q_{k,n}^d \cap G_k^{d,m,D} \ne \emptyset}
\biggl\| \sum_{\nu \in N_k^{d,m,D}: Q_{k,n}^d \cap
\supp g_{k, \nu}^{m,d} \ne \emptyset}
\D^\mu ((S_{k, \nu_k^D(\nu)}^{l -1,d} f) -\\
(R_{k, \nu_k^D(\nu)}^{l -1,d} \phi(f)))
\D^{\lambda -\mu} g_{k, \nu}^{m,d} \biggr\|_{L_q(Q_{k,n}^d)}^q \le \\
\sum_{n \in \Z^d: Q_{k,n}^d \cap G_k^{d,m,D} \ne \emptyset}
\biggl(\sum_{\nu \in N_k^{d,m,D}: Q_{k,n}^d \cap
\supp g_{k, \nu}^{m,d} \ne \emptyset}
\| \D^\mu (S_{k, \nu_k^D(\nu)}^{l -1,d} f -\\
R_{k, \nu_k^D(\nu)}^{l -1,d} \phi(f))
\D^{\lambda -\mu} g_{k, \nu}^{m,d} \|_{L_q(Q_{k,n}^d)}\biggr)^q.
\end{multline*}

Для оценки правой части (3.1.8), используя сначала (1.3.10), а
затем применяя (1.4.1), при $ \lambda \in \Z_{+ \mathpzc m}^d, \ \mu \in \Z_+^d(\lambda), \
n \in \Z^d: Q_{k,n}^d \cap G_k^{d,m,D} \ne \emptyset, \
\nu \in N_k^{d,m,D}: Q_{k,n}^d \cap \supp g_{k, \nu}^{m,d} \ne \emptyset, $
выводим
\begin{multline*} \tag{3.1.9}
\| \D^\mu (S_{k, \nu_k^D(\nu)}^{l -1,d} f -
R_{k, \nu_k^D(\nu)}^{l -1,d} \phi(f))
\D^{\lambda -\mu} g_{k, \nu}^{m,d} \|_{L_q(Q_{k,n}^d)} \le \\
\| \D^{\lambda -\mu} g_{k, \nu}^{m,d} \|_{L_\infty(\R^d)}
\| \D^\mu (S_{k, \nu_k^D(\nu)}^{l -1,d} f -
R_{k, \nu_k^D(\nu)}^{l -1,d} \phi(f))\|_{L_q(Q_{k,n}^d)} = \\
c_3 2^{k |\lambda -\mu|}
\| \D^\mu (S_{k, \nu_k^D(\nu)}^{l -1,d} f -
R_{k, \nu_k^D(\nu)}^{l -1,d} \phi(f))\|_{L_q(Q_{k,n}^d)} \le \\
c_4 2^{k |\lambda -\mu|} 2^{k(|\mu| -d /q)}
\| S_{k, \nu_k^D(\nu)}^{l -1,d} f -
R_{k, \nu_k^D(\nu)}^{l -1,d} \phi(f)\|_{L_\infty(Q_{k,n}^d)} = \\
c_4 2^{k(|\lambda| -d /q)}
\| S_{k, \nu_k^D(\nu)}^{l -1,d} f -
R_{k, \nu_k^D(\nu)}^{l -1,d} \phi(f)\|_{L_\infty(Q_{k,n}^d)}.
\end{multline*}

Далее, с учётом того, что в силу (1.4.21) имеет место включение
\begin{multline*}
Q_{k,n}^d \subset (2^{-k} \nu_k^D(\nu) +
(\gamma^1 +1) 2^{-k} B^d),\\
n \in \Z^d: Q_{k,n}^d \cap G_k^{d,m,D} \ne \emptyset, \
\nu \in N_k^{d,m,D}: \supp g_{k, \nu}^{m,d} \cap Q_{k,n}^d
\ne \emptyset,
\end{multline*}
применяя (3.1.1), затем используя (3.1.5), (1.4.74), наконец, делая
замену переменной, получаем, что при $ n \in \Z^d: Q_{k,n}^d \cap
G_k^{d,m,D} \ne \emptyset, \
\nu \in N_k^{d,m,D}: \supp g_{k, \nu}^{m,d} \cap Q_{k,n}^d
\ne \emptyset,, $ выполняется неравенство
\begin{multline*} \tag{3.1.10}
\| S_{k, \nu_k^D(\nu)}^{l -1,d} f -
R_{k, \nu_k^D(\nu)}^{l -1,d} \phi(f)\|_{L_\infty(Q_{k,n}^d)} = \\
\sup_{x \in Q_{k,n}^d}
| (S_{k, \nu_k^D(\nu)}^{l -1,d} f)(x) -
(R_{k, \nu_k^D(\nu)}^{l -1,d} \phi(f))(x)| \le \\
c_5 \max_{\mu \in \Z_{+ l -1}^d}
\biggl| (S_{k, \nu_k^D(\nu)}^{l -1,d} f)
(x_{k,\nu_k^D(\nu)}^0 +\frac {1} {2} l^{-1} 2^{-k} \mu) -
(R_{k, \nu_k^D(\nu)}^{l -1,d} \phi(f))
(x_{k,\nu_k^D(\nu)}^0 +\frac {1} {2} l^{-1} 2^{-k} \mu) \biggr| = \\
c_5 \max_{\mu \in \Z_{+ l -1}^d}
| (S_{k, \nu_k^D(\nu)}^{l -1,d} f)
(x_{k,\nu_k^D(\nu)}^\mu) -
(R_{k, \nu_k^D(\nu)}^{l -1,d} \phi(f))
(x_{k,\nu_k^D(\nu)}^\mu)| = \\
c_5 \max_{\mu \in \Z_{+ l -1}^d}
| (S_{k, \nu_k^D(\nu)}^{l -1,d} f)
(x_{k,\nu_k^D(\nu)}^\mu) -f(x_{k,\nu_k^D(\nu)}^\mu)| \le \\
c_5 \sup_{x \in Q_{k,\nu_k^D(\nu)}^d}
| (S_{k, \nu_k^D(\nu)}^{l -1,d} f)(x) -f(x)| =
c_5 \| f -S_{k, \nu_k^D(\nu)}^{l -1,d} f\|_{L_\infty(Q_{k,\nu_k^D(\nu)}^d)} \le \\
c_6 2^{k d / p} \int_0^1 t^{-d /p -1}
2^{k d /p} t^{-d /p} \biggl(\int_{ 2 2^{-k} t B^d}
\int_{ (Q_{k,\nu_k^D(\nu)}^d)_{l \xi}}
| \Delta_{\xi}^{l} f(x)|^p dx d\xi \biggr)^{1/p} dt = \\
c_6 2^{k d / p} \int_0^1 t^{-d /p -1}
(2^{-k} t)^{-d /p} \biggl(\int_{ 2 (2^{-k} t) B^d}
\int_{(Q_{k,\nu_k^D(\nu)}^d)_{l \xi}}
| \Delta_{\xi}^{l} f(x)|^p dx d\xi\biggr)^{1/p} dt = \\
c_6 \int_0^{2^{-k}} t^{-d /p -1}
 t^{-d /p} \biggl(\int_{2 t B^d}
\int_{(Q_{k,\nu_k^D(\nu)}^d)_{l \xi}}
| \Delta_{\xi}^{l} f(x)|^p dx d\xi\biggr)^{1/p} dt.
\end{multline*}

Объединяя (3.1.9), (3.1.10) и учитывая (1.4.57), выводим
\begin{multline*} \tag{3.1.11}
\| \D^\mu (S_{k, \nu_k^D(\nu)}^{l -1,d} f -
R_{k, \nu_k^D(\nu)}^{l -1,d} \phi(f))
\D^{\lambda -\mu} g_{k, \nu}^{m,d} \|_{L_q(Q_{k,n}^d)} \le \\
c_4 2^{k(|\lambda| -d /q)}
c_6 \int_0^{2^{-k}} t^{-d /p -1}  t^{-d /p} \biggl(\int_{2 t B^d}
\int_{(Q_{k,\nu_k^D(\nu)}^d)_{l \xi}}
| \Delta_{\xi}^{l} f(x)|^p dx d\xi\biggr)^{1/p} dt \le \\
c_7 2^{k(| \lambda| -d /q)} \int_0^{2^{-k}} t^{-d /p -1}
 t^{-d /p} \biggl(\int_{2 t B^d}
\int_{(D \cap D_{k,n}^{\prime d,m,D})_{l \xi}}
| \Delta_{\xi}^{l} f(x)|^p dx d\xi\biggr)^{1/p} dt \le \\
c_7 2^{k(| \lambda| -d /q)} \int_0^{2^{-k}} t^{-d /p -1}
 t^{-d /p} \biggl(\int_{2 t B^d}
\int_{D_{l \xi} \cap D_{k,n}^{\prime d,m,D}}
| \Delta_{\xi}^{l} f(x)|^p dx d\xi\biggr)^{1/p} dt,\\
n \in \Z^d: Q_{k,n}^d \cap G_k^{d,m,D} \ne \emptyset, \\
\nu \in N_k^{d,m,D}: \supp g_{k, \nu}^{m,d} \cap Q_{k,n}^d
\ne \emptyset, \ \lambda \in \Z_{+ \mathpzc m}^d, \ \mu \in \Z_+^d(\lambda).
\end{multline*}

Выбирая $ \epsilon >0 $ так, чтобы соблюдалось неравенство
$ \alpha -d /p -\epsilon >0, $ из (3.1.11), пользуясь
неравенством Гёльдера, выводим оценку
\begin{multline*}
\| \D^\mu (S_{k, \nu_k^D(\nu)}^{l -1,d} f -
R_{k, \nu_k^D(\nu)}^{l -1,d} \phi(f))
\D^{\lambda -\mu} g_{k, \nu}^{m,d} \|_{L_q(Q_{k,n}^d)} \le \\
c_7 2^{k(| \lambda| -d /q)} \int_0^{2^{-k}} t^{\alpha -d /p -\epsilon -1 /q^\prime}
t^{-(\alpha -\epsilon) -1 /q} \\
\times t^{-d /p} \biggl(\int_{2 t B^d}
\int_{D_{l \xi} \cap D_{k,n}^{\prime d,m,D}}
| \Delta_{\xi}^{l} f(x)|^p dx d\xi\biggr)^{1/p} dt \le \\
c_{8} 2^{-k(\alpha -|\lambda| -d /p +d /q -\epsilon)}
\biggl(\int_0^{2^{-k}} t^{-q(\alpha -\epsilon) -1} t^{-q d /p} \\
\times \biggl(\int_{ 2 t B^d} \int_{D_{l \xi} \cap D_{k,n}^{\prime d,m,D}}
| \Delta_{\xi}^{l} f(x)|^p dx d\xi \biggr)^{q /p} dt \biggr)^{1/q},\\
n \in \Z^d: Q_{k,n}^d \cap G_k^{d,m,D} \ne \emptyset, \\
\nu \in N_k^{d,m,D}: \supp g_{k, \nu}^{m,d} \cap Q_{k,n}^d
\ne \emptyset, \ \lambda \in \Z_{+ \mathpzc m}^d, \ \mu \in \Z_+^d(\lambda).
\end{multline*}

Подставляя эту оценку в (3.1.8) и  проводя соответствующие выкладки,
в частности, используя (1.4.27), применяя неравенство (1.1.1)
при $ a = p /q \le 1 $ и (1.4.60), для $ f \in \mathcal K $ при $ p \le q, \
\lambda \in \Z_{+ \mathpzc m}^d, \ \mu \in \Z_+^d(\lambda) $ получаем
\begin{multline*} \tag{3.1.12}
\biggl\| \sum_{\nu \in N_k^{d,m,D}}
\D^\mu ((S_{k, \nu_k^D(\nu)}^{l -1,d} f) -
(R_{k, \nu_k^D(\nu)}^{l -1,d} \phi(f)))
\D^{\lambda -\mu} g_{k, \nu}^{m,d} \biggr\|_{L_q(\R^d)}^q \le \\
\sum_{n \in \Z^d: Q_{k,n}^d \cap G_k^{d,m,D} \ne \emptyset}
\biggl(\sum_{\nu \in N_k^{d,m,D}: Q_{k,n}^d \cap
\supp g_{k, \nu}^{m,d} \ne \emptyset}
c_{8} 2^{-k(\alpha -|\lambda| -d /p +d /q -\epsilon)} \\
\times \biggl(\int_0^{2^{-k}} t^{-q(\alpha -\epsilon) -1} t^{-q d /p} \\
\times \biggl(\int_{ 2 t B^d}
\int_{D_{l \xi} \cap D_{k,n}^{\prime d,m,D}}
| \Delta_{\xi}^{l} f(x)|^p dx d\xi \biggr)^{q /p} dt \biggr)^{1/q}\biggr)^q \le \\
\sum_{n \in \Z^d: Q_{k,n}^d \cap G_k^{d,m,D} \ne \emptyset}
\biggl(c_{9} 2^{-k(\alpha -|\lambda| -d /p +d /q -\epsilon)}
\biggl(\int_0^{2^{-k}} t^{-q(\alpha -\epsilon) -1} t^{-q d /p} \\
\times \biggl(\int_{ 2 t B^d}
\int_{D_{l \xi} \cap D_{k,n}^{\prime d,m,D}}
| \Delta_{\xi}^{l} f(x)|^p dx d\xi \biggr)^{q /p} dt \biggr)^{1/q} \biggr)^q = \\
(c_{9} 2^{-k(\alpha -|\lambda| -d /p +d /q -\epsilon)})^q
\sum_{n \in \Z^d: Q_{k,n}^d \cap G_k^{d,m,D} \ne \emptyset}
\int_0^{2^{-k}} t^{-q(\alpha -\epsilon) -1} t^{-q d /p} \\
\times \biggl(\int_{ 2 t B^d}
\int_{D_{l \xi} \cap D_{k,n}^{\prime d,m,D}}
| \Delta_{\xi}^{l} f(x)|^p dx d\xi \biggr)^{q /p} dt = \\
(c_{9} 2^{-k(\alpha -|\lambda| -d /p +d /q -\epsilon)})^q
\int_0^{2^{-k}} t^{-q(\alpha -\epsilon) -1} t^{-q d /p}  \times \\
\sum_{n \in \Z^d: Q_{k,n}^d \cap G_k^{d,m,D} \ne \emptyset}
\biggl(\int_{ 2 t B^d} \int_{D_{l \xi} \cap D_{k,n}^{\prime d,m,D}}
| \Delta_{\xi}^{l} f(x)|^p dx d\xi \biggr)^{q /p} dt \le \\
(c_{9} 2^{-k(\alpha -|\lambda| -d /p +d /q -\epsilon)})^q
\int_0^{2^{-k}} t^{-q(\alpha -\epsilon) -1} t^{-q d /p}  \\
\times \biggl(\sum_{n \in \Z^d: Q_{k,n}^d \cap G_k^{d,m,D} \ne \emptyset}
\int_{ 2 t B^d} \int_{D_{l \xi} \cap D_{k,n}^{\prime d,m,D}}
| \Delta_{\xi}^{l} f(x)|^p dx d\xi \biggr)^{q /p} dt = \\
(c_{9} 2^{-k(\alpha -|\lambda| -d /p +d /q -\epsilon)})^q
\int_0^{2^{-k}} t^{-q(\alpha -\epsilon) -1} t^{-q d /p} \\
\times \biggl(\int_{ 2 t B^d}
\sum_{n \in \Z^d: Q_{k,n}^d \cap G_k^{d,m,D} \ne \emptyset}
\int_{D_{l \xi}} \chi_{ D_{k,n}^{\prime d,m,D}}(x)
| \Delta_{\xi}^{l} f(x)|^p dx d\xi \biggr)^{q /p} dt = \\
(c_{9} 2^{-k(\alpha -|\lambda| -d /p +d /q -\epsilon)})^q
\int_0^{2^{-k}} t^{-q(\alpha -\epsilon) -1} t^{-q d /p}  \\
\times \biggl(\int_{ 2 t B^d} \int_{D_{l \xi}}
(\sum_{n \in \Z^d: Q_{k,n}^d \cap
G_k^{d,m,D} \ne \emptyset}
\chi_{ D_{k,n}^{\prime d,m,D}}(x))
| \Delta_{\xi}^{l} f(x)|^p dx d\xi \biggr)^{q /p} dt \le \\
(c_{9} 2^{-k(\alpha -|\lambda| -d /p +d /q -\epsilon)})^q
\int_0^{2^{-k}} t^{-q(\alpha -\epsilon) -1} t^{-q d /p}  \\
\times \biggl(\int_{ 2 t B^d} \int_{D_{l \xi}} c_{10}
| \Delta_{\xi}^{l} f(x)|^p dx d\xi \biggr)^{q /p} dt \le \\
(c_{11} 2^{-k(\alpha -|\lambda| -d /p +d /q -\epsilon)})^q
\int_0^{2^{-k}} t^{-q(\alpha -\epsilon) -1}
\biggl((2 2 t)^{-d /p} \\
\times \biggl(\int_{ 2 t B^d} \int_{D_{l \xi}}
| \Delta_{\xi}^{l} f(x)|^p dx d\xi\biggr)^{1 /p} \biggr)^q dt = \\
(c_{11} 2^{-k(\alpha -|\lambda| -d /p +d /q -\epsilon)})^q
\int_0^{2^{-k}} t^{-q(\alpha -\epsilon) -1}
\biggl(\Omega^{\prime l}(f, 2 t)_{L_p(D)} \biggr)^q dt \le \\
(c_{11} 2^{-k(\alpha -|\lambda| -d /p +d /q -\epsilon)})^q
\int_0^{2^{-k}} t^{-q(\alpha -\epsilon) -1} c_{12} t^{\alpha q} dt = \\
(c_{13} 2^{-k(\alpha -|\lambda| -d /p +d /q -\epsilon)})^q
\int_0^{2^{-k}} t^{q \epsilon -1} dt =
(c_{13} 2^{-k(\alpha -|\lambda| -d /p +d /q -\epsilon)})^q
c_{14} 2^{-k q \epsilon} = \\
(c_{15} 2^{-k(\alpha -|\lambda| -d /p +d /q)})^q.
\end{multline*}

Соединяя (3.1.7) и (3.1.12), приходим к оценке
\begin{equation*} \tag{3.1.13}
\| \D^\lambda E_k^{l -1,d,m,D} f -\D^\lambda A \circ \phi(f)\|_{L_q(D)} \le
c_{16} 2^{-k(\alpha -|\lambda| -d /p +d /q)}.
\end{equation*}

Подстановка оценок (1.4.72) и (3.1.13) в (3.1.6) приводит к неравенству
\begin{equation*}
\| \D^\lambda f -\D^\lambda A \circ \phi(f) \|_{L_q(D)} \le
c_{17} 2^{-k(\alpha -|\lambda| -(d /p -d /q)_+)}, \
f \in \mathcal K, \ \lambda \in \Z_{+ \mathpzc m}^d, \ p \le q,
\end{equation*}
и, следовательно,
\begin{equation*} \tag{3.1.14}
\| f -A \circ \phi(f) \|_{W_q^{\mathpzc m}(D)} \le
c_{17} 2^{-k(\alpha -\mathpzc m -(d /p -d /q)_+)}, \
f \in \mathcal K, \ p \le q.
\end{equation*}

Из (3.1.14) с учётом (3.1.4) следует, что
\begin{equation*}
\overline \sigma_n(\mathcal K,C(T),X) \le \overline \sigma_{n(k,D)}(\mathcal K,C(T),X) \le
c_{17} 2^{-k(\alpha -\mathpzc m -(d /p -d /q)_+)}.
\end{equation*}
Отсюда, пользуясь тем, что для $ n, $ удовлетворяющего (3.1.4),
справедливо соотношение $ n \le c_{18}(D) 2^{kd}, $
приходим к (3.1.3) при $ p \le q. $

Если же $ q < p, $ то ввиду ограниченности области $ D $
пространство $ W_p^{\mathpzc m}(D) \subset W_q^{\mathpzc m}(D) $ и для
$ g \in W_p^{\mathpzc m}(D) $ норма
$ \| g \|_{W_q^{\mathpzc m}(D)} \le c_{19}(d,D,p,q,\mathpzc m)
\| g \|_{W_p^{\mathpzc m}(D)}. $
Поэтому $ \overline \sigma_n(\mathcal K,C(D),W_q^{\mathpzc m}(D)) \le c_{19}
\overline \sigma_n(\mathcal K,C(D),W_p^{\mathpzc m}(D)), $
и, следовательно, (3.1.3) справедливо и при $ q < p. \square $
\bigskip

3.2. В этом пункте выводится оценка  снизу величины наилучшей точности
восстановления в $ W_q^{\mathpzc m}(D) $ по значениям в $ n $ точках функций $ f $
из $ (\mathcal H_p^\alpha)^\prime(D),
(\mathcal B_{p,\theta}^\alpha)^\prime(D). $

Теорема 3.2.1

В условиях теоремы 3.1.2 существует константа $ c_1(d,\alpha,D,p,\theta,q,\mathpzc m) >0 $
такая, что при $ n \in \N $ выполняется неравенство
\begin{equation*} \tag{3.2.1}
\sigma_n(\mathcal K,C(T),X) \ge c_1 n^{-(\alpha -\mathpzc m) /d +(1 /p -1 /q)_+}.
\end{equation*}

В основе доказательства теоремы 3.2.1 лежит лемма 3.2.2, установленная в [5].

Лемма 3.2.2

Пусть $ d \in \N, \ \alpha \in \R_+, \ \theta \in \R: \theta \ge 1, \
1 \le p < \infty, \ 1 \le q \le \infty, \ \lambda \in \Z_+^d $ и
$ \l = \l(\alpha) \in \Z_+. $
Тогда существует константа $ c_2(d,\alpha,p,\theta,q,\lambda) >0 $ такая,
что для любого $ n \in \N $ для любого набора точек $ x^1,\ldots,x^n \in I^d $
можно построить функцию
$ f \in C_0^\infty(I^d) \cap (\mathcal B_{p,\theta}^\alpha)^{\l}(\R^d), $
для которой
\begin{equation*} \tag{3.2.2}
f(x^i) =0, \ i =1,\ldots,n,
\end{equation*}
и
\begin{equation*} \tag{3.2.3}
\| \D^\lambda f\|_{L_q(I^d)} \ge
c_2 n^{-(\alpha -|\lambda|) /d +(1 /p -1 /q)_+}.
\end{equation*}

Отметим, что функция $ f $ в формулировке леммы 3.2.2, приведенной в
[5], принадлежит $ C_0^\infty(\R^d) \cap
(\mathcal B_{p,\theta}^\alpha)^{\l}(I^d). $ Однако, как видно из
доказательства леммы, на самом деле $ f \in C_0^\infty(I^d) \cap
(\mathcal B_{p,\theta}^\alpha)^{\l}(\R^d). $

Доказательство теоремы 3.2.1.

Принимая во внимание включение (1.1.4), доказательство проводим в случае
$ \mathcal K = (\mathcal B_{p,\theta}^\alpha)^\prime(D). $

В условиях теоремы 3.1.2 фиксируем точку $ x^0 \in \R^d $ и число $ \delta \in
\R_+ $ такие, что $ Q = (x^0 +\delta I^d) \subset D, $ и пусть $ n \in \N, \
A \in \mathcal A^n(X) $ и $ \phi \in \Phi_n(C(D)) $ определяется системой
точек $ y^1,\ldots,y^n \in D. $ Для некоторого фиксированного
$ \lambda \in \Z_{+ \mathpzc m}^{\prime d} $ возьмём функцию $ f \in
(\mathcal B_{p,\theta}^\alpha)^{\l}(\R^d) \cap
C_0^\infty(I^d), $ обладающую свойствами (3.2.2) и (3.2.3), где
$ x^i = \delta^{-1}(y^i -x^0), \ i =1,\ldots,n. $ Рассмотрим функцию
$ F = h_{\delta,x^0}^{-1} f $ (см. п. 2.2.). Тогда, поскольку ввиду (1.1.5)
функция $ f \in C_0 (\mathcal B_{p,\theta}^\alpha)^\prime(\R^d) \cap
C_0^\infty(I^d), $ то в силу (2.2.1), (2.2.7) функция
$ F \in (c_3 (\mathcal B_{p,\theta}^\alpha)^\prime(\R^d)) \cap C_0^\infty(Q). $
Полагая $ \mathcal F = (1 / c_3) F \mid_D, $ получаем, что
$ \supp \mathcal F \subset Q, \ \mathcal F \in (\mathcal B_{p,\theta}^\alpha)^\prime(D), \
\mathcal F(y^i) = (1 / c_3) F(y^i) = (1 / c_3) f(\delta^{-1}(y^i -x^0)) =
(1 / c_3) f(x^i) =0, \ i =1,\ldots,n. $

При этом ввиду (3.2.3), (2.2.2) соблюдается неравенство
\begin{multline*}
c_2 n^{-(\alpha -\mathpzc m) /d +(1 /p -1 /q)_+} =
c_2 n^{-(\alpha -|\lambda|) /d +(1 /p -1 /q)_+}
\le \|\D^\lambda f\|_{L_q(I^d)} \\
= \|\D^\lambda (h_{\delta,x^0} F)\|_{L_q(I^d)} =
c_4(\delta,\lambda,q) \| \D^\lambda F\|_{L_q(Q)} =
c_5 \| \D^\lambda \mathcal F\|_{L_q(Q)} = c_5 \| \D^\lambda \mathcal F\|_{L_q(D)} \le \\
c_5 \| \mathcal F\|_{W_q^{\mathpzc m}(D)} = c_5 \| \mathcal F -A \circ \phi(\mathcal F)
+A \circ \phi(0) -0\|_{W_q^{\mathpzc m}(D)}
\le 2 c_5 \sup_{g \in \mathcal K} \| g -A \circ \phi(g)\|_X.
\end{multline*}
Отсюда, в силу произвольности $ A \in \mathcal A^n(X) $ и $ \phi \in
\Phi_n(C(D)), $ заключаем, что имеет место (3.2.1). $\square$
\bigskip

\centerline{\S 4. Задача С.Б. Стечкина}
\centerline{для оператора частного дифференцирования}
\centerline{на изотропных классах функций Никольского и Бесова}
\bigskip

4.1. В этом пункте будет получена оценка сверху наилучшей точности
приближения в $ L_q(D) $ оператора $ \D^\lambda $ ограниченными
операторами, действующими из $ L_s(D) $ в $ L_q(D), $ на
классах $ (\mathcal H_p^\alpha)^\prime(D) $ и
$ (\mathcal B_{p,\theta}^\alpha)^\prime(D), $ где $ D $ -- ограниченная
область $ \e $-типа в $ \R^d. $

Напомним постановку общей задачи,  частным случаем которой
является задача, рассматриваемая в этом параграфе.

Пусть $ X,Y $ -- банаховы пространства, $ \mathcal B(X, Y) $ --
банахово пространство непрерывных линейных операторов $ V: X \mapsto Y $
с обычной нормой, $ A: D(A) \mapsto Y $ -- линейный оператор с областью
определения $ D(A) \subset X. $ Пусть ещё множество $ \mathcal K \subset D(A) $ и
$ \rho >0. $

Требуется описать поведение в зависимости от $ \rho $ величины
\begin{equation*}
E(A,X,Y,\mathcal K,\rho) = \inf_{\{ V \in \mathcal B(X,Y): \|V\|_{\mathcal B(X,Y)}
\le \rho\}} \sup_{x \in \mathcal K} \|Ax -Vx\|_Y.
\end{equation*}

Теорема 4.1.1

Пусть $ d \in \N, \ \alpha \in \R_+, \ 1 < p < \infty, \ 1 \le q,s \le \infty $
и $ \lambda \in \Z_+^d $ удовлетворяют  условию (1.4.71), а также при $ p < s $
соблюдается неравенство
\begin{equation*}
\alpha -(d /p -d /s) >0,
\end{equation*}
и, кроме того, выполняется неравенство
$$
|\lambda| +(d /s -d /q)_+ >0.
$$
Пусть ещё $ D \subset \R^d $ -- ограниченная область $ \e $-типа,
$ X = L_s(D), \ Y = L_q(D), \ A = \D^\lambda, \ D(A) = \{f \in L_s(D):
\D^\lambda f \in L_q(D)\}, \ \mathcal K =
(\mathcal H_p^\alpha)^\prime(D), \
(\mathcal B_{p,\theta}^\alpha)^\prime(D), $ где $ \theta \in \R:
\theta \ge 1, \ \gamma = \alpha -|\lambda| -(d /p -d /q)_+, \
\tau = |\lambda| +(d /s -d /q)_+. $ Тогда существуют
константы $ \rho_1(d,\alpha,D,q,s,\lambda) >0 $ и
$ c_1(d,\alpha,D,p,q,s,\lambda) >0 $ такие, что для $ \rho \ge \rho_1 $
справедливо неравенство
\begin{equation*} \tag{4.1.1}
E(A,X,Y,\mathcal K,\rho) \le c_1 \rho^{-\gamma /\tau}.
\end{equation*}

Для доказательства (4.1.1) достаточно заметить, что в условиях теоремы
вследствие теоремы 1.4.13 и (1.1.4) множество $ \mathcal K \subset D(A), $
и если взять $ l = l(\alpha), \ m \in \N: \lambda \in \Z_{+ m}^d, $ то для каждого
достаточно большого $ k \in \Z_+ $
для оператора $ V: L_s(D) \mapsto L_q(D), $ определяемого равенством
$$
V f = (\D^\lambda (E_k^{l -1,d,m,D} f)) \mid_D,
$$
в силу (1.4.72), (1.4.34) (см. также (1.1.4)) выполняются соотношения
\begin{equation*}
\sup_{f \in \mathcal K} \| \D^\lambda f -V f \|_{L_q(D)} \le
c_2 2^{-k(\alpha -|\lambda| -(d /p -d /q)_+)},
\end{equation*}
и
\begin{equation*}
\| V \|_{\mathcal B(L_s(D), L_q(D))} \le
c_3 2^{k (|\lambda| +(d /s -d /q)_+)},
\end{equation*}
а, следовательно,
\begin{multline*}
\inf_{\{\mathcal V \in \mathcal B(X,Y): \| \mathcal V \|_{\mathcal B(X,Y)} \le
c_3 2^{k (|\lambda| +(d /s -d /q)_+)}\}}
\sup_{f \in \mathcal K} \| \D^\lambda f -\mathcal V f \|_{Y} \le \\
c_2 2^{-k(\alpha -|\lambda| -(d /p -d /q)_+)}. \square
\end{multline*}
\bigskip

4.2. В этом пункте проводится оценка снизу величины, которая
рассматривалась в п. 4.1. Для доказательства этой оценки используем
теорему 4.2.1, справедливость которой по существу установлена в [6]
при доказательстве теоремы 2.2.4, опирающейся на теорему 2.2.1 и
леммы 2.2.2, 2.2.3.
При этом будем пользоваться следующим обозначением. Для множества $ S, $
состоящего из функций $ f, $ область определения которых содержит множество
$ D \subset \R^d, $ через $ S \mid_D $ обозначим множество $ S \mid_D =
\{ f \mid_D: f \in S\}. $

Теорема 4.2.1

Пусть $ d, \alpha, p, q,s, \lambda $  удовлетворяют условиям теоремы 4.1.1 и
$ \l = \l(\alpha) $ (см. п. 1.1.). Пусть ещё $ A = \D^\lambda, \
D(A) = \{f \in L_s(I^d): \D^\lambda f \in L_q(I^d)\}, \
X = L_s(I^d), \ Y = L_q(I^d), \ \mathcal K =
((\mathcal B_{p,\theta}^\alpha)^{\l}(\R^d) \cap C_0^\infty(I^d)) \mid_{I^d}, $
где $ \theta \in \R: \theta \ge 1, \ \gamma = \alpha -|\lambda| -(d /p -d /q)_+), \
\tau = |\lambda| +(d /s -d /q)_+. $
Тогда существуют константы $ c_1(d,\alpha,p,\theta,q,s,\lambda) >0 $ и
$ \rho_1(d,\alpha,q,s,\lambda) >0 $ такие, что для $ \rho > \rho_1 $
выполняется неравенство
\begin{equation*} \tag{4.2.1}
E(A,X,Y,\mathcal K,\rho) \ge c_1 \rho^{-\gamma /\tau}.
\end{equation*}

Теорема 4.2.2

В условиях и обозначениях теоремы 4.1.1 существуют константы
$ \rho_2(d,\alpha,D,q,s,\lambda) >0, \ c_2(d,\alpha,D,p,\theta,q,s,\lambda) >0 $
такие, что при $ \rho > \rho_2 $ соблюдается неравенство (4.2.1) с
константой $ c_2 $ вместо $ c_1. $

Доказательство.

Ввиду включения (1.1.4) рассмотрим лишь случай, когда $ \mathcal K =
(\mathcal B_{p,\theta}^\alpha)^\prime(D), \ \theta \ne \infty. $
Фиксируем $ x^0 \in \R^d $ и $ \delta \in \R_+ $ такие, что $ Q =
(x^0 +\delta I^d) \subset D. $
Прежде всего, отметим, что в силу (1.1.5), (4.2.1) при $ \rho > \rho_1$
верно неравенство
\begin{multline*} \tag{4.2.2}
c_1 \rho^{-\gamma /\tau} \le E(\D^\lambda,L_s(I^d),L_Q(I^d),
C_0 ((\mathcal B_{p,\theta}^\alpha)^\prime(\R^d) \cap C_0^\infty(I^d)) \mid_{I^d}, \rho) = \\
C_0 E(\D^\lambda,L_s(I^d),L_Q(I^d),
((\mathcal B_{p,\theta}^\alpha)^\prime(\R^d) \cap C_0^\infty(I^d)) \mid_{I^d}, \rho).
\end{multline*}

Далее, заметим, что для $ f \in L_s(I^d): \D^\lambda f \in L_q(I^d), \
V \in \mathcal B(L_s(Q),L_q(Q)) $ имеет место соотношение
\begin{multline*}
\| \D^\lambda f -(h_{\delta,x^0} V h_{\delta,x^0}^{-1}) f\|_{L_q(I^d)} =
\biggl(\int_{I^d}| \D^\lambda f -(h_{\delta,x^0} V h_{\delta,x^0}^{-1}) f|^q dx\biggr)^{1 /q} = \\
\biggl(\int_{Q}| h_{\delta,x^0}^{-1}(\D^\lambda f) -V (h_{\delta,x^0}^{-1} f)|^q \delta^{-d} dy\biggr)^{1 /q} =\\
\delta^{-d /q} \biggl(\int_{Q}| \delta^{|\lambda|} \D^\lambda (h_{\delta,x^0}^{-1} f) -
V (h_{\delta,x^0}^{-1} f)|^q dy\biggr)^{1 /q} = \\
\delta^{|\lambda| -d /q} \| \D^\lambda (h_{\delta,x^0}^{-1} f) -
\delta^{-|\lambda|} V (h_{\delta,x^0}^{-1} f)\|_{L_q(Q)}.
\end{multline*}
Отметим ещё, что для $ V \in \mathcal B(L_s(Q),L_q(Q)), $ вследствие
(2.2.2) при $ \lambda =0, $ (2.2.3), норма
\begin{multline*}
\|h_{\delta,x^0} V h_{\delta,x^0}^{-1}\|_{\mathcal B(L_s(I^d),L_q(I^d))} \le
\|h_{\delta,x^0} \|_{\mathcal B(L_q(Q),L_q(I^d))} \\
\| V\|_{\mathcal B(L_s(Q),L_q(Q))} \| h_{\delta,x^0}^{-1}\|_{\mathcal B(L_s(I^d),L_s(Q))} \le \\
\delta^{-d /q} \| V\|_{\mathcal B(L_s(Q),L_q(Q))} \delta^{d /s} =
\delta^{d /s -d /q} \| V\|_{\mathcal B(L_s(Q),L_q(Q))}.
\end{multline*}
Учитывая эти обстоятельства, а также (2.2.1), (2.2.7), при $ \rho \in \R_+ $
для $ V \in \mathcal B(L_s(Q),L_q(Q)): \|V\|_{\mathcal B(L_s(Q),L_q(Q))} \le
\delta^{|\lambda|} \rho, $ получаем, что
\begin{multline*}
\inf_{\substack{\{ U \in \mathcal B(L_s(I^d),L_q(I^d)):\\
\|U\|_{\mathcal B(L_s(I^d),L_q(I^d))} \le \delta^{|\lambda| +d /s -d /q} \rho\}}}
\sup_{f \in ((\mathcal B_{p,\theta}^\alpha)^\prime(\R^d) \cap C_0^\infty(I^d)) \mid_{I^d}}
\| \D^\lambda f -U f\|_{L_q(I^d)} \le \\
\sup_{f \in ((\mathcal B_{p,\theta}^\alpha)^\prime(\R^d) \cap C_0^\infty(I^d)) \mid_{I^d}}
\| \D^\lambda f -(h_{\delta,x^0} V h_{\delta,x^0}^{-1}) f\|_{L_q(I^d)} = \\
\sup_{f \in ((\mathcal B_{p,\theta}^\alpha)^\prime(\R^d) \cap C_0^\infty(I^d)) \mid_{I^d}}
\delta^{|\lambda| -d /q} \| \D^\lambda (h_{\delta,x^0}^{-1} f) -
\delta^{-|\lambda|} V (h_{\delta,x^0}^{-1} f)\|_{L_q(Q)} = \\
\delta^{|\lambda| -d /q}
\sup_{f \in ((\mathcal B_{p,\theta}^\alpha)^\prime(\R^d) \cap C_0^\infty(I^d)) \mid_{I^d}}
\| \D^\lambda (h_{\delta,x^0}^{-1} f) -
\delta^{-|\lambda|} V (h_{\delta,x^0}^{-1} f)\|_{L_q(Q)} \le \\
\delta^{|\lambda| -d /q}
\sup_{F \in (c_3 (\mathcal B_{p,\theta}^\alpha)^\prime(\R^d)) \cap C_0^\infty(Q)}
\| \D^\lambda (F \mid_Q) -\delta^{-|\lambda|} V (F \mid_Q)\|_{L_q(Q)},
\end{multline*}
и, значит,
\begin{multline*} \tag{4.2.3}
c_3 \inf_{\substack{ \mathcal V \in \mathcal B(L_s(Q),L_q(Q)):\\
\|\mathcal V\|_{\mathcal B(L_s(Q),L_q(Q))} \le \rho}}
\sup_{F \in (\mathcal B_{p,\theta}^\alpha)^\prime(\R^d) \cap C_0^\infty(Q)}
\| \D^\lambda (F \mid_Q) -\mathcal V (F \mid_Q)\|_{L_q(Q)} \ge\\
\delta^{-|\lambda| +d /q}
\inf_{\substack{\{ U \in \mathcal B(L_s(I^d),L_q(I^d)):\\
\|U\|_{\mathcal B(L_s(I^d),L_q(I^d))} \le \delta^{|\lambda| +d /s -d /q} \rho\}}}
\sup_{f \in ((\mathcal B_{p,\theta}^\alpha)^\prime(\R^d) \cap C_0^\infty(I^d)) \mid_{I^d}}
\| \D^\lambda f -U f\|_{L_q(I^d)}.
\end{multline*}

Пользуясь тем, что для $ F \in C_0^\infty(Q) $ справедливо равенство
$ F \mid_D = (\mathcal I^Q(F \mid_Q)) \mid_D \ (\mathcal I^Q $ см. в п. 1.4.),
а для $ \mathcal V \in \mathcal B(L_s(D),L_q(D)) $ имеет место неравенство
\begin{multline*}
\| (\mathcal V((\mathcal I^Q f) \mid_D)) \mid_Q \|_{L_q(Q)} \le
\| \mathcal V((\mathcal I^Q f) \mid_D) \|_{L_q(D)} \le\\
\| \mathcal V\|_{\mathcal B(L_s(D),L_q(D))}
\| (\mathcal I^Q f) \mid_D \|_{L_s(D)} =
\| \mathcal V\|_{\mathcal B(L_s(D),L_q(D))}
\| ((\mathcal I^Q f) \mid_D) \mid_Q \|_{L_s(Q)} = \\
\| \mathcal V\|_{\mathcal B(L_s(D),L_q(D))} \| f\|_{L_s(Q)}, \
f \in L_s(Q),
\end{multline*}
для $ \mathcal V \in \mathcal B(L_s(D),L_q(D)):
\| \mathcal V\|_{\mathcal B(L_s(D),L_q(D))} \le \rho, $ выводим
\begin{multline*}
\sup_{f \in (\mathcal B_{p,\theta}^\alpha)^\prime(D)}
\| \D^\lambda f -\mathcal V f\|_{L_q(D)} \ge \\
\sup_{F \in (\mathcal B_{p,\theta}^\alpha)^\prime(\R^d) \cap C_0^\infty(Q)}
\| \D^\lambda (F \mid_D) -\mathcal V (F \mid_D)\|_{L_q(D)} \ge \\
\sup_{F \in (\mathcal B_{p,\theta}^\alpha)^\prime(\R^d) \cap C_0^\infty(Q)}
\| (\D^\lambda (F \mid_D)) \mid_Q -(\mathcal V (F \mid_D)) \mid_Q\|_{L_q(Q)} = \\
\sup_{F \in (\mathcal B_{p,\theta}^\alpha)^\prime(\R^d) \cap C_0^\infty(Q)}
\| \D^\lambda (F \mid_Q) -(\mathcal V ((\mathcal I^Q(F \mid_Q)) \mid_D)) \mid_Q\|_{L_q(Q)} \ge \\
\inf_{ V \in \mathcal B(L_s(Q),L_q(Q)): \|V\|_{\mathcal B(L_s(Q),L_q(Q))} \le \rho}
\sup_{F \in (\mathcal B_{p,\theta}^\alpha)^\prime(\R^d) \cap C_0^\infty(Q)}
\| \D^\lambda (F \mid_Q) -V (F \mid_Q)\|_{L_q(Q)},
\end{multline*}
а, следовательно,
\begin{multline*} \tag{4.2.4}
\inf_{ \mathcal V \in \mathcal B(L_s(D),L_q(D)):
\| \mathcal V\|_{\mathcal B(L_s(D),L_q(D))} \le \rho}
\sup_{f \in (\mathcal B_{p,\theta}^\alpha)^\prime(D)}
\| \D^\lambda f -\mathcal V f\|_{L_q(D)} \ge \\
\inf_{ V \in \mathcal B(L_s(Q),L_q(Q)): \|V\|_{\mathcal B(L_s(Q),L_q(Q))} \le \rho}
\sup_{F \in (\mathcal B_{p,\theta}^\alpha)^\prime(\R^d) \cap C_0^\infty(Q)}
\| \D^\lambda (F \mid_Q) -V (F \mid_Q)\|_{L_q(Q)}.
\end{multline*}

Соединяя (4.2.4), (4.2.3), (4.2.2), заключаем, что при $ \rho > \rho_2 =
\delta^{-|\lambda| -d /s +d /q} \rho_1 $ соблюдается неравенство (4.2.1) с
константой $ c_2 $ вместо $ c_1, $ где $ A,X,Y,\mathcal K $ имеют тот же смысл,
что в формулировке теоремы 4.1.1. $ \square $
\bigskip

\centerline{\S 5. Поперечники классов $ B((H_p^\alpha)^\prime(D)) $ и
$B((B_{p,\theta}^\alpha)^\prime(D))$}
\centerline{в пространстве $ W_q^{\m}(D)$}
\bigskip

5.1. В этом пункте введём в рассмотрение подходящие нам подпространства
кусочно-полиномиальных функций и установим некоторые вспомогательные
утверждения, касающиеся самих подпространств и линейных отображений в эти подпространства,
с помощью которых проводятся верхние оценки изучаемых поперечников.

Сначала заметим, что при $ m \in \N, \ d \in \N $ для каждого $ \nu \in \Z^d $
существует единственная пара $ (\rho,\sigma): \rho \in \Z^d, \sigma \in
\Nu_{0,m}^d, $ такая, что $ \nu = (m +1) \rho +\sigma. $
А это значит, что множество $ \Z^d $ представляется в виде
\begin{multline*} \tag{5.1.1}
\Z^d = \bigcup_{\sigma \in \Nu_{0,m}^d} (\sigma +(m +1) \Z^d),\text{ причём }
(\sigma +(m +1) \Z^d) \cap (\sigma^{\prime} +(m +1) \Z^d) = \emptyset, \\
\text{ для } \sigma, \sigma^{\prime} \in \Nu_{0,m}^d: \sigma \ne
\sigma^{\prime}.
\end{multline*}

Для $ d, m \in \N, \ k \in \Z_+, \ \sigma \in \Nu_{0,m}^d $ и области
$ D \subset \R^d $ обозначим
$$
N_{k,\sigma}^{d,m,D} = (\sigma +(m +1) \Z^d) \cap N_k^{d,m,D}
$$
и
$$
r_{k,\sigma}^{d,m,D} = \card N_{k,\sigma}^{d,m,D}.
$$

Принимая во внимание, что
\begin{multline*}
r_{k,\sigma}^{d,m,D} = \card \{\nu = (m +1) \rho +\sigma:
\rho \in \Z^d, \ (2^{-k} ((m +1) \rho +\sigma) +2^{-k} (m +1) \overline I^d)
\cap D \ne \emptyset \} = \\
\card \{\rho \in \Z^d: (\rho +I^d) \cap (-(m +1)^{-1} \sigma +2^k (m +1)^{-1} D)
\ne \emptyset \},
\end{multline*}
нетрудно видеть, что существуют константы $ c_1(d,m,D) >0 $ и $ c_2(d,m,D) >0 $
такие, что для $ d,m \in \N, \ k \in \Z_+, \ \sigma \in \Nu_{0,m}^d, $
ограниченной области $ D \subset \R^d $
выполняется неравенство
\begin{equation*} \tag{5.1.2}
c_1 2^{kd} \le r_{k,\sigma}^{d,m,D} \le c_2 2^{kd}.
\end{equation*}

Для $ d, m \in \N, \ k, l \in \Z_+, \ \sigma \in \Nu_{0,m}^d $ и ограниченной
области $ D \subset \R^d $ обозначим через
$ \mathcal P_{k,\sigma}^{l,d,m,D} $ линейное пространство, состоящее из
функций $ f: \R^d \mapsto \R, $ для  каждой из которых существует набор

полиномов $ \{f_{\nu} \in \mathcal  P^{l,d}, \ \nu \in N_{k,\sigma}^{d,m,D}\} $
такой, что для всех $ x \in \R^d $
соблюдается равенство
\begin{equation*} \tag{5.1.3}
f(x) = \sum_{\nu \in N_{k,\sigma}^{d,m,D}} f_{\nu}(x) g_{k,\nu}^{m,d}(x).
\end{equation*}

В силу (5.1.1) для $ d,m \in \N, \ k,l \in \Z_+, $ ограниченной области
$ D \subset \R^d $ имеет место равенство
\begin{equation*} \tag{5.1.4}
\mathcal P_k^{l,d,m,D} = \sum_{\sigma \in \Nu_{0,m}^d}
\mathcal P_{k,\sigma}^{l,d,m,D}.
\end{equation*}

Для доказательства предложения 5.1.2 нам понадобится лемма 5.1.1, являющаяся
очевидным следствием того факта, что  все нормы на конечномерном линейном
пространстве эквивалентны.

Лемма 5.1.1

Пусть $ d,m \in \N, \ l \in \Z_+, \ 1 \le  p \le \infty. $
Тогда существуют константы $ c_3(l,d,m,p) >0 $ и $ c_4(l,d,m,p) >0 $
такие, что для $ f \in \mathcal P^{l,d} $ соблюдается соотношение
\begin{equation*} \tag{5.1.5}
c_3 \| f \psi^{m,d} \|_{L_p((m +1) I^d)}\le \|f\|_{L_p((m+1) I^d)}
\le c_4 \| f \psi^{m,d} \|_{L_p((m +1) I^d)}.
\end{equation*}

Предложение 5.1.2

Пусть $ d,m \in \N, \ l \in \Z_+, \ D $ -- ограниченная область в $ \R^d $ и
$ 1 \le p \le \infty. $ Тогда существуют константы $ c_5(l,d,m,p) >0 $ и
$ c_6(l,d,m,p) >0 $ такие, что для $ k \in \Z_+ $ и $ \sigma \in \Nu_{0,m}^d, $
для любого набора полиномов $ \{f_{\nu} \in \mathcal P^{l,d}, \ \nu \in N_{k,\sigma}^{d,m,D}\} $
для функции $ f$, определяемой для всех
$ x \in \R^d $ равенством (5.1.3), выполняется неравенство
\begin{equation*} \tag{5.1.6}
c_5 \|f\|_{L_p(\R^d)} \le \biggl(\sum_{\nu \in N_{k,\sigma}^{d,m,D}}
\|f_{\nu}\|_{L_p(Q_{k,\nu}^{m,d})}^p\biggr)^{1/p} \le
c_6 \|f\|_{L_p(\R^d)}, \ p \ne \infty,
\end{equation*}

\begin{equation*} \tag{5.1.6'}
c_5 \|f\|_{L_{\infty}(\R^d)} \le \max_{\nu \in N_{k,\sigma}^{d,m,D}}
\|f_{\nu}\|_{L_{\infty}(Q_{k,\nu}^{m,d})} \le
c_6 \|f\|_{L_{\infty}(\R^d)},
\end{equation*}

где $ Q_{k,\nu}^{m,d} = 2^{-k} \nu +2^{-k} (m +1) I^d. $

Доказательство.

Проверим (5.1.6) (аналогично проверяется (5.1.6')).
В условиях предложения, полагая
$ G_{k,\sigma}^{d,m,D} = \cup_{\nu \in N_{k,\sigma}^{d,m,D}} \supp g_{k,\nu}^{m,d}, $
и ввиду (1.3.7) учитывая, что
$$
\supp g_{k,\nu}^{m,d} = Q_{k,\nu}^{m,d} \cup (\overline Q_{k,\nu}^{m,d}
\setminus Q_{k,\nu}^{m,d}),
$$
причём
$$
\mes (\overline Q_{k,\nu}^{m,d} \setminus Q_{k,\nu}^{m,d}) =0,
$$
имеем
$$
G_{k,\sigma}^{d,m,D} = (\bigcup_{\nu \in N_{k,\sigma}^{d,m,D}} Q_{k,\nu}^{m,d})
\cup A_{k,\sigma}^{d,m,D},
$$
причём $ \mes A_{k,\sigma}^{d,m,D} =0, $ а $ Q_{k,\nu}^{m,d} \cap
Q_{k,\nu^{\prime}}^{m,d} =
Q_{k,\nu}^{m,d} \cap \supp g_{k,\nu^{\prime}}^{m,d} = \emptyset $
для $ \nu, \nu^{\prime} \in (\sigma +(m +1) \Z^d): \nu \ne
\nu^{\prime}. $

Принимая во внимание сказанное, получаем, что в рассматриваемой ситуации
\begin{multline*} \tag{5.1.7}
\|f\|_{L_p(\R^d)}^p = \int_{\R^d} \biggl|\sum_{\nu^{\prime} \in
N_{k,\sigma}^{d,m,D}} f_{\nu^{\prime}}
g_{k,\nu^{\prime}}^{m,d}\biggr|^p dx = \\
\int_{G_{k,\sigma}^{d,m,D}} \biggl|\sum_{\nu^{\prime} \in
N_{k,\sigma}^{d,m,D}} f_{\nu^{\prime}}
g_{k,\nu^{\prime}}^{m,d}\biggr|^p dx =
\sum_{\nu \in N_{k,\sigma}^{d,m,D}} \int_{ Q_{k,\nu}^{m,d}} |f_{\nu}
g_{k,\nu}^{m,d}|^p dx.
\end{multline*}

Теперь для каждого $ \nu \in N_{k,\sigma}^{d,m,D}, $ учитывая (5.1.5),
имеем
\begin{multline*}
\int_{ Q_{k,\nu}^{m,d} } |f_{\nu}
g_{k,\nu}^{m,d}|^p dx = \int_{(2^{-k} \nu +2^{-k} (m +1) I^d)}
|f_{\nu}(x) \psi^{m,d}(2^k
x -\nu)|^p dx \\
= 2^{-kd} \int_{(m +1) I^d} |f_{\nu}(2^{-k} \nu +2^{-k} y) \psi^{m,d}(y)|^p dy \ge
2^{-kd} c_7 \int_{(m +1) I^d} |f_{\nu}(2^{-k} \nu +2^{-k} y)|^p dy \\
= c_7 \int_{ Q_{k,\nu}^{m,d}} |f_{\nu}(x)|^p dx = c_7
\|f_{\nu}\|_{L_p( Q_{k,\nu}^{m,d})}^p.
\end{multline*}

Подставляя эту оценку в (5.1.7), получаем второе неравенство в
(5.1.6).

Далее, для каждого $ \nu \in N_{k,\sigma}^{d,m,D} $ справедливо
неравенство
$$
\int_{ Q_{k,\nu}^{m,d} } |f_{\nu} g_{k,\nu}^{m,d}|^p dx
= \int_{ Q_{k,\nu}^{m,d}} |f_{\nu}(x) \psi^{m,d}(2^k x -\nu)|^p dx \le
\|\psi^{m,d}\|_{L_{\infty}(\R^d)}^p \cdot
\|f_{\nu}\|_{L_p(Q_{k,\nu}^{m,d})}^p.
$$

Соединение этой оценки с (5.1.7)  приводит к первому неравенству
в (5.1.6). $\square$

Из (5.1.6) вытекает, что в условиях предложения 5.1.2 линейный оператор,
определяемый равенством (5.1.3), является изоморфизмом прямого
произведения $ r_{k,\sigma}^{d,m,D} $   экземпляров пространства
$ \mathcal P^{l,d} $ на пространство $ \mathcal P_{k,\sigma}^{l,d,m,D}. $
Поэтому, полагая
$$
R_k^{l,d,m,D} = \dim \mathcal  P_k^{l,d,m,D}, \
R_{k,\sigma}^{l,d,m,D} = \dim \mathcal  P_{k,\sigma}^{l,d,m,D}, \
c(l,d) = \dim \mathcal P^{l,d} = \card \Z_{+ l}^d,
$$
заключаем, что для $ d,m \in \N, \ l,k \in \Z_+, \ \sigma \in \Nu_{0,m}^d $
и ограниченной области $ D \subset \R^d $ справедливы соотношения
$$
R_{k,\sigma}^{l,d,m,D} = c(l,d) r_{k,\sigma}^{d,m,D},
$$
и, значит, в силу (5.1.2) --
\begin{equation*} \tag{5.1.8}
c_8 2^{kd} < R_{k,\sigma}^{l,d,m,D} < c_9 2^{kd},
\end{equation*}
где $ c_8 > 0 $ и $ c_9 > 0 $ зависят только от $ l,d,m,D. $

Из (5.1.4) следует, что при $ d,m \in \N, \ l,k \in \Z_+, \ \sigma \in \Nu_{0,m}^d $
и ограниченной области $ D \subset \R^d $ выполняется соотношение
$$
R_{k,\sigma}^{l,d,m,D} \le R_k^{l,d,m,D} \le
\sum_{\tau \in \Nu_{0,m}^d} R_{k,\tau}^{l,d,m,D},
$$
которое в соединении с (5.1.8) даёт оценку
\begin{equation*} \tag{5.1.9}
c_8 2^{kd} < R_k^{l,d,m,D} < c_{10} 2^{kd},
\end{equation*}
с константой $ c_{10} > 0, $ зависящей
только от $ l,d,m,D.$

Ещё одним простым следствием эквивалентности норм на конечномерных
линейных пространствах является следующая лемма.

Лемма 5.1.3

Пусть $ l,m \in \Z_+, \ d \in \N, \ 1 \le p \le \infty. $ Тогда существуют
константы $ c_{11}(l,m,d,p) >0 $ и $ c_{12}(l,m,d,p) >0 $ такие, что для
$ f \in \mathcal P^{l,d} $ имеют место
неравенства
\begin{equation*} \tag{5.1.10}
c_{11} \|f\|_{L_p((m +1) I^d)} \le \|\{f(\lambda), \ \lambda \in
\Z_{+ l}^d \}\|_{l_p^{c(l,d)}} \le c_{12} \|f\|_{L_p((m +1) I^d)}.
\end{equation*}

Для $ l,k \in \Z_+, \ d,m \in \N, \ \sigma \in \Nu_{0,m}^d, $ ограниченной
области $ D \subset \R^d $ обозначим через $ \I_{k,\sigma}^{l,d,m,D}:
\mathcal P_{k,\sigma}^{l,d,m,D} \mapsto \R^{R_{k,\sigma}^{l,d,m,D}}, $
линейный оператор, который каждой функции
$ f \in \mathcal P_{k,\sigma}^{l,d,m,D} $ ставит в соответствие набор
значений $ \{f_{\nu}(2^{-k} \nu +2^{-k} \lambda), \ \nu \in
N_{k,\sigma}^{d,m,D}, \ \lambda \in \Z_{+ l}^d\}, $ где
$ \{f_{\nu} \in \mathcal P^{l,d}, \ \nu \in N_{k,\sigma}^{d,m,D}\} $ -- система
полиномов, удовлетворяющих (5.1.3).

Предложение 5.1.4

Пусть $ l \in \Z_+, \ d,m \in \N, \ 1 \le p \le \infty, \ D $ -- ограниченная
область в $ \R^d. $ Тогда существуют константы $ c_{13}(l,d,m,p) >0 $ и
$ c_{14}(l,d,m,p) >0 $ такие, что при $ k \in \Z_+, \ \sigma \in \Nu_{0,m}^d $
для $ f \in \mathcal P_{k,\sigma}^{l,d,m,D} $ соблюдается
соотношение
\begin{equation*} \tag{5.1.11}
c_{13} \|f\|_{L_p(\R^d)} \le 2^{-(kd) /p}
\| \I_{k,\sigma}^{l,d,m,D} f \|_{l_p^{R_{k,\sigma}^{l,d,m,D}}} \le
c_{14} \|f\|_{L_p(\R^d)}.
\end{equation*}

Доказательство.

Для доказательства (5.1.11) достаточно воспользоваться
соотношениями (5.1.6) (или (5.1.6') при $ p = \infty $),
а затем, пользуясь тем, что
$$
\|f_{\nu}\|_{L_p(Q_{k,\nu}^{m,d})} = 2^{-kd /p}
\|f_{\nu}(2^{-k} \nu +2^{-k} \cdot)\|_{L_p((m +1) I^d)},
$$
и тем, что $ f_{\nu}(2^{-k} \nu +2^{-k} \cdot) \in \mathcal P^{l,d} $ при
$ \nu \in N_{k,\sigma}^{d,m,D}, \ k \in \Z_+, $ применить (5.1.10).$\square$

Лемма 5.1.5

Пусть $ l \in \Z_+, \ d,m \in \N, \ 1 \le q \le \infty, \ D $ -- ограниченная
область в $ \R^d. $ Тогда при $ \mathpzc m \in \Z_+: \mathpzc m \le m, \ k \in \Z_+, \
\sigma \in \Nu_{0,m}^d $ отображение
\begin{equation*}
U = U_{k,\sigma}^{l,d,m,D,\mathpzc m,q}: \mathcal P_{k,\sigma}^{l,d,m,D} \cap
W_q^{\mathpzc m}(\R^d) \mapsto \{f \mid_D: f \in \mathcal P_{k,\sigma}^{l,d,m,D}\} \cap W_q^{\mathpzc m}(D),
\end{equation*}
задаваемое равенством $ U f = f \mid_D, $ является линейным гомеоморфизмом.

Доказательство.

Линейность и сюръективность отображения $ U $ очевидны. Проверим его инъективность.
Если для $ f \in \mathcal P_{k,\sigma}^{l,d,m,D} $ значение $ U f =0, $
то ввиду представления (5.1.3) для $ x \in D $ выполняется равенство
$$
\sum_{\nu \in N_{k,\sigma}^{d,m,D}} f_\nu(x) g_{k,\nu}^{m,d}(x) =0,
$$
а, следовательно, в силу того, что $ Q_{k,\nu}^{m,d} \cap
\supp g_{k,\nu^{\prime}}^{m,d} = \emptyset $
для $ \nu, \nu^{\prime} \in (\sigma +(m +1) \Z^d): \nu \ne \nu^{\prime}, $
при $ \nu \in N_{k,\sigma}^{d,m,D} $ для $ x \in D \cap Q_{k,\nu}^{m,d} =
D \cap \inter (\supp g_{k,\nu}^{m,d}) $ соблюдается равенство
$ f_\nu(x) g_{k,\nu}^{m,d}(x) =0, $ и, значит (см. (1.3.1)), $ f_\nu(x) =0. $
Отсюда, учитывая включение $ f_\nu \in \mathcal P^{l,d} $ и открытость непустого
множества $ D \cap Q_{k,\nu}^{m,d}, $ заключаем, что $ f_\nu(x) =0 $ для
$ x \in \R^d, \ \nu \in N_{k,\sigma}^{d,m,D}, $ т.е. $ f =0. $
Принимая во внимание сказанное, видим, что функционал $ \| U f \|_{W_q^{\mathpzc m}(D)} $
является нормой на конечномерном линейном пространстве
$ \mathcal P_{k,\sigma}^{l,d,m,D}. $ А поскольку любые нормы на конечномерном
линейном пространстве эквивалентны, то существует константа
$ c_{15}(l,d,m,D,q,\mathpzc m,k,\sigma) >0 $
такая, что для $ f \in \mathcal P_{k,\sigma}^{l,d,m,D} $ справедливо неравенство
\begin{equation*} \tag{5.1.12}
c_{15} \| f \|_{W_q^{\mathpzc m}(\R^d)} \le \| f \mid_D \|_{W_q^{\mathpzc m}(D)} \le
\| f \|_{W_q^{\mathpzc m}(\R^d)},
\end{equation*}
что завершает доказательство леммы. $ \square $

Предложение 5.1.6

Пусть $ l \in \Z_+, \ d,m \in \N, \ 1 \le q \le \infty, \ D $ -- ограниченная
область в $ \R^d. $ Тогда существует константа $ c_{16}(l,d,m,D,q) >0 $
такая, что при $ \mathpzc m \in \Z_+: \mathpzc m \le m, \ k \in \Z_+, \ \sigma \in \Nu_{0,m}^d $
для $ f \in \mathcal P_{k,\sigma}^{l,d,m,D} $ справедливо неравенство
\begin{equation*} \tag{5.1.13}
\|f\|_{W_q^{\mathpzc m}(\R^d)} \le c_{16} 2^{k \mathpzc m} \|f\|_{L_q(\R^d)}.
\end{equation*}

Доказательство.

Доказательство проведём при $ q \ne \infty. $

В условиях предложения в силу тех же соображений, что при выводе (5.1.7),
принимая во внимание (1.4.33), неравенство Гёльдера, (1.3.10) и
(1.4.1), для $ f  \in \mathcal P_{k,\sigma}^{l,d,m,D} $
при $ \lambda \in \Z_{+ m}^d $ имеем
\begin{multline*} \tag{5.1.14}
\| \D^{\lambda} f\|_{L_q(\R^d)}^q = \int_{\R^d}
| \D^{\lambda} f|^q dx = \sum_{\nu \in N_{k,\sigma}^{d,m,D}}
\int_{Q_{k,\nu}^{m,d}} \biggl| \D^{\lambda} \biggl(\sum_{\nu^{\prime} \in
N_{k,\sigma}^{d,m,D}} f_{\nu^{\prime}} g_{k,\nu^{\prime}}^{m,d} \biggr)\biggr|^q dx \\
= \sum_{\nu \in N_{k,\sigma}^{d,m,D}} \int_{ Q_{k,\nu}^{m,d}}
| \D^{\lambda}(f_{\nu} g_{k,\nu}^{m,d})|^q dx
= \sum_{\nu \in N_{k,\sigma}^{d,m,D}} \biggl\|\sum_{\mu \in \Z_+^d(\lambda)}
C_{\lambda}^{\mu} \D^{\mu} f_{\nu}
\D^{\lambda -\mu} g_{k,\nu}^{m,d}\biggr\|_{L_q( Q_{k,\nu}^{m,d})}^q \\
\le \sum_{ \nu \in N_{k,\sigma}^{d,m,D}} \biggl(\sum_{\mu \in \Z_+^d(\lambda)}
C_{\lambda}^{\mu} \| \D^{\mu} f_{\nu} \D^{\lambda -\mu} g_{k,\nu}^{m,d}
\|_{L_q( Q_{k,\nu}^{m,d})}\biggr)^q \\
\le \sum_{\nu \in N_{k,\sigma}^{d,m,D}} c_{17} \biggl(\sum_{\mu \in
\Z_+^d(\lambda)} \| \D^{\mu} f_{\nu} \D^{\lambda -\mu} g_{k,\nu}^{m,d}
\|_{L_q(Q_{k,\nu}^{m,d})}^q \biggr) \\
\le c_{17} \sum_{\nu \in N_{k,\sigma}^{d,m,D}} \sum_{\mu \in
\Z_+^d(\lambda)} \|\D^{\lambda -\mu} g_{k,\nu}^{m,d}\|_{L_{\infty}(\R^d)}^q
\| \D^{\mu} f_{\nu}\|_{L_q(Q_{k,\nu}^{m,d})}^q \\
\le c_{17} \sum_{\nu \in N_{k,\sigma}^{d,m,D}} \sum_{\mu \in
\Z_+^d(\lambda)} \biggl(c_{18} 2^{k|\lambda -\mu|} \biggr)^q
\left(c_{19} 2^{k|\mu|} \|f_{\nu}\|_{L_q( Q_{k,\nu}^{m,d})} \right)^q \\
\le c_{20} 2^{k |\lambda| q} \sum_{\nu \in N_{k,\sigma}^{d,m,D}}
\| f_{\nu} \|_{L_q( Q_{k,\nu}^{m,d})}^q.
\end{multline*}

Из (5.1.14) и (5.1.6) вытекает (5.1.13). $\square$

Для формулировки леммы 5.1.7 введём следующие обозначения.

Для $ d \in \N, \ l \in \Z_+, \ m \in \N, $ ограниченной области $ D $ в $ \R^d $
при $ k \in \Z_+, $ для которого существует $ \nu \in \Z^d:
\overline Q_{k,\nu}^d \subset D, $ определим непрерывный линейный оператор
$ V_k^{l,d,m,D}: L_1(D) \mapsto
\mathcal P_k^{l,d,m,D} \mid_D \cap W_\infty^m(D), $ полагая для
$ f \in L_1(D) $ значение
$$
V_k^{l,d,m,D} f = (E_k^{l,d,m,D} f) \mid_D \ (\text{см. } (1.4.9),
(1.4.33), (1.4.1), (1.4.3), (1.3.10)).
$$

Лемма 5.1.7

Пусть $ d \in \N, \ \alpha \in \R_+, \ l = l(\alpha), \ m \in \N, \ \mathpzc m \in \Z_+:
\mathpzc m \le m, \ D \subset \R^d $ -- ограниченная область $ \e $-типа.
Пусть ещё $ 1 < p < \infty, \ 1 \le q \le \infty $ и соблюдается условие (3.1.2).
Тогда существует константа $ c_{21}(d,\alpha,m,D,\mathpzc m,p,q) >0 $ такая, что для
любой функции $ f \in (\mathcal H_p^\alpha)^\prime(D) $ при $ k \in \Z_+:
k \ge K^0 $ (см. п. 1.2., определение 1 при $ \alpha = \e $) соблюдается
неравенство
\begin{equation*} \tag{5.1.15}
\| f -V_k^{l -1,d,m,D} f \|_{W_q^{\mathpzc m}(D)} \le
c_{21} 2^{-k(\alpha -\mathpzc m -(d /p -d /q)_+)}.
\end{equation*}

Доказательство.

В условиях леммы, сопоставляя (3.1.2) с (1.4.71) и применяя (1.4.72), для
$ f \in (\mathcal H_p^\alpha)^\prime(D), \ \lambda \in \Z_{+ \mathpzc m}^d $ при
$ k \in \Z_+: k \ge K^0, $ имеем
\begin{multline*}
\| \D^\lambda f -\D^\lambda V_k^{l -1,d,m,D} f \|_{L_q(D)} =
\| \D^\lambda f -\D^\lambda (E_k^{l -1,d,m,D} f) \mid_D \|_{L_q(D)} \le \\
c_{22} 2^{-k(\alpha -|\lambda| -(d /p -d /q)_+)} \le
c_{21} 2^{-k(\alpha -\mathpzc m -(d /p -d /q)_+)},
\end{multline*}
что влечёт (5.1.15). $ \square $

Непосредственным следствием леммы 5.1.7 является одно из утверждений леммы
5.1.8, для формулировки которой потребуются следующие обозначения.

Для $ d \in \N, \ l \in \Z_+, \ m \in \N, $ ограниченной области $ D $ в $ \R^d $
при $ k \in \N, $ для которого существует $ \nu^\prime \in \Z^d:
\overline Q_{k -1,\nu^\prime}^d \subset D, $ с учётом (1.4.40) определим
непрерывный линейный оператор
$ \mathcal V_k^{l,d,m,D}: L_1(D) \mapsto \mathcal P_k^{l,d,m,D} \mid_D \cap
W_\infty^m(D) $ равенством
\begin{equation*}
\mathcal V_k^{l,d,m,D} f = (\mathcal E_k^{l,d,m,D} f) \mid_D, \ f \in L_1(D).
\end{equation*}

При соблюдении указанных выше условий с учётом (1.4.40), (1.3.11), (1.3.16),
а также (5.1.1) имеем
\begin{multline*} \tag{5.1.16}
\mathcal E_k^{l,d,m,D} f =
E_k^{l,d,m,D} f -H_{k -1}^{l,d,m,D} (E_{k -1}^{l,d,m,D} f) = \\
\sum_{\nu \in N_k^{d,m,D}}
(S_{k, \nu_k^D(\nu)}^{l,d} f) g_{k, \nu}^{m,d} -
H_{k -1}^{l,d,m,D} \biggl(\sum_{\nu^\prime \in N_{k -1}^{d,m,D}}
(S_{k -1, \nu_{k -1}^D(\nu^\prime)}^{l,d} f)
g_{k -1, \nu^\prime}^{m,d}\biggr) = \\
\sum_{\nu \in N_k^{d,m,D}}
(S_{k, \nu_k^D(\nu)}^{l,d} f) g_{k, \nu}^{m,d} -
\sum_{\nu \in N_k^{d,m,D}}
\biggl(\sum_{\m \in \M^{m,d}(\nu)} A_{\m}^{m,d}
S_{k -1, \nu_{k -1}^D(\n(\nu,\m))}^{l,d} f\biggr)
g_{k, \nu}^{m,d} = \\
\sum_{\nu \in N_k^{d,m,D}}
\biggl((S_{k, \nu_k^D(\nu)}^{l,d} f) -
\sum_{\m \in \M^{m,d}(\nu)} A_{\m}^{m,d}
S_{k -1, \nu_{k -1}^D(\n(\nu,\m))}^{l,d} f\biggr) g_{k, \nu}^{m,d} = \\
\sum_{\nu \in \cup_{\sigma \in \Nu_{0,m}^d} N_{k,\sigma}^{d,m,D}}
\biggl((S_{k, \nu_k^D(\nu)}^{l,d} f) -
\sum_{\m \in \M^{m,d}(\nu)} A_{\m}^{m,d}
S_{k -1, \nu_{k -1}^D(\n(\nu,\m))}^{l,d} f\biggr) g_{k, \nu}^{m,d} = \\
\sum_{\sigma \in \Nu_{0,m}^d}
\sum_{\nu \in  N_{k,\sigma}^{d,m,D}}
\biggl((S_{k, \nu_k^D(\nu)}^{l,d} f) -
\sum_{\m \in \M^{m,d}(\nu)} A_{\m}^{m,d}
S_{k -1, \nu_{k -1}^D(\n(\nu,\m))}^{l,d} f\biggr) g_{k, \nu}^{m,d}.
\end{multline*}

При тех же условиях на $ d,l,m,D $ и $ k \in \N, \ \sigma \in \Nu_{0,m}^d, $ что и
выше, определим ещё непрерывные линейные операторы
$ \mathcal E_{k,\sigma}^{l,d,m,D}: L_1(D) \mapsto
\mathcal P_{k,\sigma}^{l,d,m,D} \cap W_\infty^m(\R^d) $ и
$ \mathcal V_{k,\sigma}^{l,d,m,D}: L_1(D) \mapsto
\mathcal P_{k,\sigma}^{l,d,m,D} \mid_D \cap W_\infty^m(D) $ равенствами
\begin{equation*}
\mathcal E_{k,\sigma}^{l,d,m,D} f = \sum_{\nu \in N_{k,\sigma}^{d,m,D}}
\biggl((S_{k, \nu_k^D(\nu)}^{l,d} f) -
\sum_{\m \in \M^{m,d}(\nu)} A_{\m}^{m,d}
S_{k -1, \nu_{k -1}^D(\n(\nu,\m))}^{l,d} f\biggr) g_{k, \nu}^{m,d}, \
f \in L_1(D),
\end{equation*}
и
\begin{equation*}
\mathcal V_{k,\sigma}^{l,d,m,D} f =
(\mathcal E_{k,\sigma}^{l,d,m,D} f) \mid_D, \ f \in L_1(D).
\end{equation*}
Из (5.1.16) выводим
\begin{equation*}
\mathcal E_k^{l,d,m,D} f = \sum_{\sigma \in \Nu_{0,m}^d}
\mathcal E_{k,\sigma}^{l,d,m,D} f, \ f \in L_1(D),
\end{equation*}
и
\begin{multline*} \tag{5.1.17}
\mathcal V_k^{l,d,m,D} f = \biggl(\sum_{\sigma \in \Nu_{0,m}^d}
\mathcal E_{k,\sigma}^{l,d,m,D} f \biggr) \mid_D =
\sum_{\sigma \in \Nu_{0,m}^d} (\mathcal E_{k,\sigma}^{l,d,m,D} f) \mid_D = \\
\sum_{\sigma \in \Nu_{0,m}^d} \mathcal V_{k,\sigma}^{l,d,m,D} f, \ f \in L_1(D).
\end{multline*}

Ещё раз обращаясь к (1.4.40), ввиду (1.3.13) получаем
\begin{multline*} \tag{5.1.18}
\mathcal V_k^{l,d,m,D} f =
(E_k^{l,d,m,D} f -H_{k -1}^{l,d,m,D} (E_{k -1}^{l,d,m,D} f)) \mid_D = \\
(E_k^{l,d,m,D} f) \mid_D -(H_{k -1}^{l,d,m,D} (E_{k -1}^{l,d,m,D} f)) \mid_D = \\
(E_k^{l,d,m,D} f) \mid_D -(E_{k -1}^{l,d,m,D} f) \mid_D =
V_k^{l,d,m,D} f -V_{k -1}^{l,d,m,D} f, \ f \in L_1(D).
\end{multline*}

Лемма 5.1.8

В условиях леммы 5.1.7

1) для $ f \in (\mathcal H_p^\alpha)^\prime(D) $ при $ k \in \Z_+: k \ge K^0, $
в $ W_q^{\mathpzc m}(D) $ имеет место равенство
\begin{equation*} \tag{5.1.19}
f = V_k^{l -1,d,m,D} f +\sum_{j = 1}^\infty \mathcal V_{k +j}^{l -1,d,m,D} f,
\end{equation*}
и

2) существует константа $ c_{23}(d,\alpha,m,D,\mathpzc m,p,q) >0 $ такая, что
для $ f \in (\mathcal H_p^\alpha)^\prime(D) $ при $ \lambda \in
\Z_{+ \mathpzc m}^d, \ \sigma \in \Nu_{0,m}^d, \
k \in \Z_+: k > K^0, $ справедливо неравенство
\begin{equation*} \tag{5.1.20}
\| \D^\lambda \mathcal E_{k,\sigma}^{l -1,d,m,D} f \|_{L_q(\R^d)}
\le c_{23} 2^{-k (\alpha -| \lambda| -(d /p -d /q)_+)}.
\end{equation*}

Доказательство.

Начнём с доказательства (5.1.19). При соблюдении условий леммы для
$ f \in (\mathcal H_p^\alpha)^\prime(D), \ k \in \Z_+: k \ge K^0 (K^0 $ см. в п. 1.2.),
в силу (5.1.18), (5.1.15) выполняется соотношение
\begin{multline*}
\| f -(V_k^{l -1,d,m,D} f +\sum_{j =1}^J \mathcal V_{k +j}^{l -1,d,m,D} f)
\|_{W_q^{\mathpzc m}(D)} = \\
\| f -(V_k^{l -1,d,m,D} f +\sum_{j =1}^J (V_{k +j}^{l -1,d,m,D} f -
V_{k +j -1}^{l -1,d,m,D} f))\|_{W_q^{\mathpzc m}(D)} = \\
\| f -V_{k +J}^{l -1,d,m,D} f \|_{W_q^{\mathpzc m}(D)} \le \\
c_{21} 2^{-(k +J)(\alpha -\mathpzc m -(d /p -d /q)_+)} \to 0 \text{ при }
J \to \infty,
\end{multline*}
т.е. справедливо (5.1.19).

Перейдём к проверке справедливости неравенства (5.1.20). В условиях леммы,
повторяя вывод неравенства (1.4.41) с заменой $ N_k^{d,m,D} $ на
$ N_{k,\sigma}^{d,m,D}, $ для $ f \in L_p(D) $ при $ \lambda \in \Z_{+ \mathpzc m}^d,
\ \sigma \in \Nu_{0,m}^d, \ k \in \Z_+: k > K^0, $ получим оценку
\begin{equation*}
\| \D^\lambda \mathcal E_{k,\sigma}^{l -1,d,m,D} f \|_{L_q(\R^d)}
\le c_{24} 2^{k (| \lambda| +(d /p -d /q)_+)}
\Omega^{\prime l}(f, 2^{-k +1})_{L_p(D)}.
\end{equation*}
Отсюда, учитывая, что для $ f \in (\mathcal H_p^\alpha)^\prime(D) $
при $ l = l(\alpha) $ соблюдается неравенство
$$
\Omega^{\prime l}(f, t)_{L_p(D)} \le t^\alpha, \ t \in \R_+,
$$
приходим к (5.1.20). $ \square $

С применением леммы 5.1.7 устанавливается также предложение 5.1.9.

Предложение 5.1.9

В условиях леммы 5.1.7 множество $ B((H_p^\alpha)^\prime(D)) $ --
компактно в $ W_q^{\mathpzc m}(D). $

Доказательство.

Ввиду непрерывности и конечномерности линейного оператора
$ V_k^{l -1,d,m,D}: L_1(D) \mapsto W_{\infty}^{\mathpzc m}(D) $ (см. (1.4.34)),
множество $ V_k^{l -1,d,m,D} (B((H_p^\alpha)^\prime(D))), \ k \in \Z_+: k \ge K^0, $ --
вполне ограничено относительно нормы $ W_{\infty}^{\mathpzc m}(D), $ а , следовательно,
и относительно нормы $ W_q^{\mathpzc m}(D). $

Далее, для $ f \in B((H_p^\alpha)^\prime(D)) $ при $ k \in \Z_+: k \ge K^0, $
выполняется (5.1.15), причём, ввиду (3.1.2)
\linebreak $ 2^{-k(\alpha -\mathpzc m -(d /p -d /q)_+)} \to 0 $ при $ k \to \infty. $
Поэтому для любого $ \epsilon >0 $ существует $ k \in \Z_+ $
такое, что вполне ограниченное множество $ V_k^{l -1,d,m,D} (B((H_p^\alpha)^\prime(D))) $
образует $ \epsilon $-сеть для $ B((H_p^\alpha)^\prime(D)) $ в $ W_q^{\mathpzc m}(D), $
что влечёт полную ограниченность $ B((H_p^\alpha)^\prime(D)) $ в $ W_q^{\mathpzc m}(D). $

В силу теоремы Хаусдорфа, чтобы убедиться в компактности
$ B((H_p^\alpha)^\prime(D)) $ в $ W_q^{\mathpzc m}(D), $ остаётся проверить
замкнутость $ B((H_p^\alpha)^\prime(D)) $ в $ W_q^{\mathpzc m}(D). $ Для этого
рассмотрим функцию $ f \in W_q^{\mathpzc m}(D), $ для которой имеется
последовательность $ \{ f_n \in B((H_p^\alpha)^\prime(D)): \ n \in \N\}, $
сходящаяся к $ f $ в $ W_q^{\mathpzc m}(D), $ а, следовательно, и в смысле сходимости
обобщённых функций.
При $ 1 < p < \infty $ в силу секвенциальной компактности
шара $ B(L_p(D)) $ относительно $ * $-слабой топологии в $ L_p(D) =
(L_{p^{\prime}}(D))^*, $ существует функция $ f_0 \in L_p(D) $ и
подпоследовательность последовательности $ \{f_n, \ n \in \N\}, $ слабо
сходящаяся к $ f_0 $ в $ L_p(D). $ Сопоставляя сказанное, заключаем, что
$ f = f_0, $ и переходя, если нужно, к подпоследовательности, получаем, что
для любой функции $ g \in L_{p^{\prime}}(D) $ справедливо равенство
\begin{equation*} \tag{5.1.21}
\int_{ D } f(x) g(x) dx = \lim_{ n \to \infty}
\int_{ D} f_n(x) g(x) dx.
\end{equation*}
Из равенства (5.1.21), учитывая, что в силу неравенства Гёльдера при
$ n \in \N $ соблюдается оценка
$$
\int_{ D} f_n(x) g(x) dx \le \| f_n \|_{L_p(D)} \|g\|_{L_{p^\prime}(D)} \le
\|g\|_{L_{p^\prime}(D)},
$$
находим, что
\begin{equation*} \tag{5.1.22}
\| f \|_{L_p(D)} = \sup_{ g \in B(L_{p^\prime}(D))} \int_{ D} f(x) g(x) dx \le 1.
\end{equation*}

Далее, заметим, что для любого $ n \in \N $ при $ t \in \R_+ $ выполняется
неравенство
$$
\biggl((2 t)^{-d} \int_{ t B^d} \int_{D_{l \xi}} |\Delta_{\xi}^l f_n(x)|^p
dx d\xi\biggr)^{1 /p} = \Omega^{\prime l}(f_n,t)_{L_p(D)} \le t^\alpha.
$$
При $ t \in \R_+ $ полагая $ \mathfrak D_t = \{(\xi,x): \ \xi \in t B^d, \
x \in D_{l \xi}\}, $ из последнего неравенства на основании теоремы Фубини
выводим
$$
\| \Delta_\xi^l f_n(x) \|_{L_p(\mathfrak D_t)} = \biggl(\int_{ t B^d} \int_{D_{l \xi}}
|\Delta_{\xi}^l f_n(x)|^p dx d\xi\biggr)^{1 /p} \le (2 t)^{d /p} t^\alpha, \
n \in \N.
$$
Отсюда в силу теоремы Фубини и неравенства Гёльдера при $ t \in \R_+, \
n \in \N $ для $ g \in L_{p^\prime}(\mathfrak D_t) $ соблюдается неравенство
\begin{multline*} \tag{5.1.23}
\int_{ t B^d} \int_{D_{l \xi}} \Delta_{\xi}^l f_n(x) g(\xi,x) dx d\xi =
\int_{\mathfrak D_t} \Delta_{\xi}^l f_n(x) g(\xi,x) d\xi dx \le \\
\| \Delta_\xi^l f_n(x) \|_{L_p(\mathfrak D_t)} \|g\|_{L_{p^\prime}(\mathfrak D_t)} \le
(2 t)^{d /p} t^\alpha \|g\|_{L_{p^\prime}(\mathfrak D_t)}.
\end{multline*}

Покажем, что при $ t \in \R_+ $ для любой функции $ g \in L_{p^\prime}(\mathfrak D_t) $
последовательность
$ \{ \int_{ t B^d} \int_{D_{l \xi}} \Delta_{\xi}^l f_n(x) g(\xi,x) dx d\xi, \ n \in \N\} $
сходится к
$$
\int_{ t B^d} \int_{D_{l \xi}} \Delta_{\xi}^l f(x) g(\xi,x) dx d\xi.
$$
В самом деле, поскольку для $ g \in L_{p^\prime}(\mathfrak D_t) $ вследствие
теоремы Фубини почти для всех $ \xi \in t B^d $ функция $ g(\xi, \cdot) \in
L_{p^\prime}(D_{l \xi}), $ то в силу (5.1.21) почти для всех $ \xi \in t B^d $
справедливо равенство
\begin{multline*} \tag{5.1.24}
\lim_{n \to \infty} \int_{D_{l \xi}} \Delta_{\xi}^l f_n(x) g(\xi,x) dx =
\lim_{n \to \infty} \int_{D_{l \xi}}
(\sum_{ k =0}^l C_l^k (-1)^{l -k} f_n(x +k \xi)) g(\xi,x) dx = \\
\lim_{n \to \infty} \sum_{ k =0}^l C_l^k (-1)^{l -k}
\int_{D_{l \xi}} f_n(x +k \xi) g(\xi,x) dx = \\
\lim_{n \to \infty} \sum_{ k =0}^l C_l^k (-1)^{l -k}
\int_{D_{l \xi} +k \xi} f_n(x) g(\xi,x -k \xi) dx = \\
\sum_{ k =0}^l C_l^k (-1)^{l -k} \lim_{n \to \infty}
\int_{D_{l \xi} +k \xi} f_n(x) g(\xi,x -k \xi) dx = \\
\sum_{ k =0}^l C_l^k (-1)^{l -k}
\int_{D_{l \xi} +k \xi} f(x) g(\xi,x -k \xi) dx = \\
\sum_{ k =0}^l C_l^k (-1)^{l -k}
\int_{D_{l \xi} } f(x +k \xi) g(\xi,x) dx = \\
\int_{D_{l \xi} } (\sum_{ k =0}^l C_l^k (-1)^{l -k} f(x +k \xi)) g(\xi,x) dx = \\
\int_{D_{l \xi}} \Delta_{\xi}^l f(x) g(\xi,x) dx,
\end{multline*}
а при $ n \in \N $ благодаря неравенству Гёльдера почти
для всех $ \xi \in t B^d $ имеет место неравенство
\begin{multline*} \tag{5.1.25}
\biggl|\int_{D_{l \xi}} \Delta_{\xi}^l f_n(x) g(\xi,x) dx\biggr| \le
\|\Delta_\xi^l f_n\|_{L_p(D_{l \xi})}
\biggl(\int_{D_{l \xi}} |g(\xi,x)|^{p^\prime} dx\biggr)^{1 /p^\prime} \le \\
2^l \| f_n \|_{L_p(D)} \biggl(\int_{D_{l \xi}} |g(\xi,x)|^{p^\prime} dx\biggr)^{1 /p^\prime} \le \\
2^l \biggl(\int_{D_{l \xi}} |g(\xi,x)|^{p^\prime} dx\biggr)^{1 /p^\prime} \in
L_{p^\prime}(t B^d) \subset L_1(t B^d).
\end{multline*}

Учитывая (5.1.24), (5.1.25), на основании теоремы Лебега о предельном
переходе под знаком интеграла заключаем, что при $ t \in \R_+ $ для любой
функции $ g \in L_{p^\prime}(\mathfrak D_t) $ выполняется равенство
\begin{equation*}
\lim_{n \to \infty} \int_{ t B^d} \int_{D_{l \xi}} \Delta_{\xi}^l f_n(x)
g(\xi,x) dx d\xi =
\int_{ t B^d} \int_{D_{l \xi}} \Delta_{\xi}^l f(x) g(\xi,x) dx d\xi.
\end{equation*}

Сопоставляя это равенство с (5.1.23), приходим к неравенству
\begin{equation*}
\int_{ t B^d} \int_{D_{l \xi}} \Delta_{\xi}^l f(x) g(\xi,x) dx d\xi \le
(2 t)^{d /p} t^\alpha \|g\|_{L_{p^\prime}(\mathfrak D_t)}, \
g \in L_{p^\prime}(\mathfrak D_t), \ t \in \R_+.
\end{equation*}
Поэтому
\begin{equation*}
\biggl(\int_{ t B^d} \int_{D_{l \xi}} |\Delta_{\xi}^l f(x)|^p dx d\xi\biggr)^{1 /p} =
\sup_{g \in B(L_{p^\prime}(\mathfrak D_t))}
\int_{ t B^d} \int_{D_{l \xi}} \Delta_{\xi}^l f(x) g(\xi,x) dx d\xi \le
(2 t)^{d /p} t^\alpha,
\end{equation*}
или
$$
\Omega^{\prime l}(f,t)_{L_p(D)} =
\biggl((2 t)^{-d} \int_{ t B^d} \int_{D_{l \xi}} |\Delta_{\xi}^l f(x)|^p
dx d\xi\biggr)^{1 /p} \le t^\alpha, \ t \in \R_+.
$$
Отсюда и из (5.1.22) видим, что $ f \in  B((H_p^\alpha)^\prime(D)), $ а,
следовательно, $ B((H_p^\alpha)^\prime(D)) $ -- замкнуто в $ W_q^{\mathpzc m}(D). \square $
\bigskip

5.2. В этом пункте даётся описание слабой асимптотики изучаемых в
настоящем параграфе поперечников, но сначала напомним их
определения (см. [13]).

Пусть $ C $ -- подмножество банахова пространства $X$ и
$ n \in \Z_+.$ Тогда $n$-поперечником по Колмогорову множества $C$
в пространстве $X$ называется величина
$$
d_n(C,X) = \inf_{M \in \mathcal M_n(X)} \sup_{x \in C} \inf_{y \in M} \|x -y\|_X,
$$
где $ \mathcal M_n(X) $ -- совокупность всех плоскостей $ M $ в $ X, $ у которых
$ \dim M \le n. $

$ n$-поперечником по Гельфанду множества $ C $ в пространстве $ X $ называется
величина
$$
d^n(C,X) = \inf_{M \in \mathcal M^n(X)} \sup_{x \in C \cap M} \|x\|_X,
$$
где $ \mathcal M^n(X)$ -- множество всех замкнутых линейных подпространств $ M $
в $ X, $ у которых $ \codim M \le n. $

Линейным $ n$-поперечником множества $ C $ в пространстве $ X $
называется величина
$$
\lambda_n(C,X) = \inf_{A \in \mathcal L_n(X)} \sup_{x \in C} \|x -Ax\|_X,
$$
где $ \mathcal L_n(X)$ -- множество всех непрерывных аффинных
отображений $ A: X \mapsto X, $ у которых $ \Rank A \le n. $

$ n$-поперечником по Александрову  множества $ C $ в пространстве $ X $
называется величина
$$
a_n(C,X) = \inf_{\phi \in \mathcal A_n(C,X)} \sup_{x \in C} \|x -\phi(x)\|_X,
$$
где $ \mathcal A_n(C,X)$ -- множество  всех непрерывных отображений
$ \phi: C \mapsto X, $ для которых существует компакт
$ K_{\phi} \subset X $ такой, что $ \Im \phi \subset K_{\phi} $ и
$ \dim K_{\phi} \le n. $

Введём ещё следующую величину
$$
a^n(C,X) = \inf_{\phi \in \mathcal A^n(C,X)} \sup_{y \in \Im \phi}
\diam(\phi^{-1}(y))_X,
$$
где $ \mathcal A^n(C,X)$ -- совокупность всех непрерывных относительно метрики
пространства $ X $ отображений $ \phi:  C \mapsto K_{\phi}, $ для которых
$ K_{\phi}$ -- метрический компакт размерности не больше $ n, $ а
$ \diam(E)_X$ -- диаметр множества $ E $ в метрике пространства $ X. $

Наконец, энтропийным $n$-поперечником множества $ C $ в пространстве $ X $
будем называть величину
$$
\epsilon_n(C,X) = \inf\{\epsilon >0: \exists \ x^1,\ldots,x^{2^n} \in
X: C \subset \cup_{j=1}^{2^n}(x^j +\epsilon B(X))\}.
$$

С помощью утверждений 5.2.1 и 5.2.2, взятых из [14], доказывается лемма 5.2.3.

Предложение 5.2.1

Пусть $ U: X \mapsto Y$ -- непрерывное линейное отображение
банахова пространства $ X $ в банахово пространство $ Y $ и $ L \subset X$ --
замкнутое линейное подпространство, для которого сужение $ U \mid_L $ является
гомеоморфизмом $ L $ на $ U(L), $ и пусть $ C \subset L$ -- некоторое
множество. Тогда для $ p_n(\cdot,\cdot) = d_n(\cdot,\cdot), d^n(\cdot,\cdot),
\lambda_n(\cdot,\cdot),  a^n(\cdot,\cdot), \epsilon_n(\cdot,\cdot) $ при
$ n \in \Z_+ $ имеет место неравенство
\begin{equation*} \tag{5.2.1}
p_n(U(C), Y) \le \| U \|_{\mathcal B(X,Y)} p_n(C,X).
\end{equation*}

Предложение 5.2.2

Пусть $ C $ -- подмножество банахова  пространства $ X, $ для которого
существует последовательность непрерывных линейных отображений
$ V_j: X \mapsto X_j $ пространства $ X $ в замкнутые линейные подпространства
$ X_j, \ j \in \Z_+, $ такая, что для каждого $ x \in C $ имеет место
представление $ x = \sum_{j =0}^{\infty} V_j x. $ Тогда для любых $ n,j_0 \in
\Z_+ $ и любого набора $ \{n_j \in \Z_+, \ j =0,\ldots,j_0\} $ таких, что
$ n \ge \sum_{j =0}^{j_0} n_j, $ для $ p_n(\cdot,\cdot) = d_n(\cdot,\cdot),
d^n(\cdot,\cdot), \lambda_n(\cdot,\cdot), a^n(\cdot,\cdot), \epsilon_n(\cdot,\cdot) $
выполняется неравенство
\begin{equation*} \tag{5.2.2}
p_n(C, X) \le \sum_{j =0}^{j_0} p_{n_j}(V_j(C), X_j)
+2\sum_{j =j_0 +1}^{\infty} d^0(V_j(C), X_j).
\end{equation*}

Лемма 5.2.3

Пусть выполнены условия леммы 5.1.7.  Тогда для
$ C = B((H_p^\alpha)^\prime(D)) $ и $ X = W_q^{\mathpzc m}(D) $
справедливы следующие утверждения:

1) существует константа $ c_1(d,\alpha,m,D,p,q,\mathpzc m) >0 $ такая, что для
любых $ n,j_0,k \in \Z_+: \ k \ge K^0, $ и любого набора чисел
$ \{n_{j,\sigma} \in \Z_+, \ j =1,\ldots,j_0, \ \sigma \in \Nu_{0,m}^d \} $
таких, что
$ n \ge R_k^{l -1,d,m,D} +\sum_{j =1}^{j_0} \sum_{\sigma \in \Nu_{0,m}^d}
n_{j,\sigma}, $ для $ p_n(\cdot,\cdot) = d_n(\cdot,\cdot), d^n(\cdot,\cdot), \lambda_n(\cdot,\cdot),
a^n(\cdot,\cdot) $ имеет место неравенство
\begin{multline*} \tag{5.2.3}
p_n(C, X) \le c_1 \biggl(2^{-k(\alpha -\mathpzc m -d /p +d /q)}
\sum_{j =1}^{j_0} \sum_{\sigma \in \Nu_{0,m}^d} 2^{-j(\alpha -\mathpzc m -d /p +d /q)} \\
p_{n_{j,\sigma}} \bigl(B(l_p^{R_{k +j,\sigma}^{l -1,d,m,D}}),
l_q^{R_{k +j,\sigma}^{l -1,d,m,D}}\bigr) +2^{-(k +j_0)(\alpha -\mathpzc m -(d /p -d /q)_+)} \biggr);
\end{multline*}

2) существует константа $ c_2(d,\alpha,m,D,p,q,\mathpzc m) >0 $ такая,
что для любых $ n,j_0 \in \Z_+ $ и любого набора чисел
$ \{n_{j,\sigma} \in \Z_+, \ j =0,\ldots, j_0, \ \sigma \in \Nu_{0,m}^d \} $
таких, что $ n \ge \sum_{j =0}^{j_0} \sum_{\sigma \in \Nu_{0,m}^d}
n_{j,\sigma}, $ соблюдается
неравенство
\begin{multline*} \tag{5.2.4}
\epsilon_n(C,X) \le c_2 \biggl( \sum_{j =0}^{j_0} \sum_{\sigma \in \Nu_{0,m}^d}
2^{-j(\alpha -\mathpzc m -d /p +d /q)} \\
\epsilon_{n_{j,\sigma}} \bigl(B(l_p^{R_{K^0 +j,\sigma}^{l -1,d,m,D}}),
l_q^{R_{K^0 +j,\sigma}^{l -1,d,m,D}}\bigr) +2^{-j_0(\alpha -\mathpzc m -(d /p -d /q)_+)} \biggr).
\end{multline*}

Доказательство.

Прежде всего, в силу (5.1.19), (5.1.17) при $ k \in \Z_+: k \ge K^0, \
j_0 \in \Z_+ $ для $ f \in C $ в $ X $ имеет место равенство
\begin{equation*}
f = V_k^{l -1,d,m,D} f
+\sum_{j =1}^{j_0} \sum_{\sigma \in \Nu_{0,m}^d} \mathcal V_{k +j,\sigma}^{l -1,d,m,D} f
+\sum_{j =j_0 +1}^\infty \mathcal V_{k +j}^{l -1,d,m,D} f.
\end{equation*}

Поэтому для $ n \in \Z_+ $ и $ \{n_{j,\sigma} \in \Z_+, \ j =1,\ldots,j_0, \
\sigma \in \Nu_{0,m}^d\} $
таких, что $ n \ge R_k^{l -1,d,m,D}
+\sum_{j =1}^{j_0} \sum_{\sigma \in \Nu_{0,m}^d} n_{j,\sigma}, $ ввиду
(5.2.2) справедлива оценка
\begin{multline*} \tag{5.2.5}
p_n(C,X) \le p_{R_k^{l -1,d,m,D}}
\bigl(V_k^{l -1,d,m,D}(C),
(\mathcal P_k^{l -1,d,m,D}) \mid_D \cap X\bigr) \\
+\sum_{j =1}^{j_0} \sum_{\sigma \in \Nu_{0,m}^d} p_{n_{j,\sigma}}
\bigl(\mathcal V_{k +j,\sigma}^{l -1,d,m,D}(C),
(\mathcal P_{k +j,\sigma}^{l -1,d,m,D}) \mid_D \cap W_q^{\mathpzc m}(D)\bigr) \\
+2\sum_{j =j_0 +1}^\infty d^0\bigl(\mathcal V_{k +j}^{l -1,d,m,D}(C),
(\mathcal P_{k +j}^{l -1,d,m,D}) \mid_D \cap W_q6{\mathpzc m}(D)\bigr).
\end{multline*}

Ясно, что
\begin{multline*} \tag{5.2.6}
p_{R_k^{l -1,d,m,D}}
\bigl(V_k^{l -1,d,m,D}(C),
(\mathcal P_k^{l -1,d,m,D}) \mid_D \cap X \bigr) = \\
p_{R_k^{l -1,d,m,D}}\bigl((E_k^{l -1,d,m,D}(C)) \mid_D,
(\mathcal P_k^{l -1,d,m,D}) \mid_D \cap W_q^{\mathpzc m}(D) \bigr) =0.
\end{multline*}

Далее, при $ j =1,\ldots,j_0, \ \sigma \in \Nu_{0,m}^d, $ применяя сначала с
учётом леммы 5.1.5 неравенства (5.2.1) и (5.1.12), а затем -- (5.2.1) и
(5.1.13), а также (5.2.1) и (5.1.11), наконец, используя (5.1.11) и (5.1.20)
(при $ \lambda =0, \ q = p $), получаем
\begin{equation*} \tag{5.2.7}
\begin{split}
p_{n_{j,\sigma}} \bigl(\mathcal V_{k +j,\sigma}^{l -1,d,m,D}(C),
(\mathcal P_{k +j,\sigma}^{l -1,d,m,D}) \mid_D \cap W_q^{\mathpzc m}(D) \bigr) = \\
p_{n_{j,\sigma}} \bigl((\mathcal E_{k +j,\sigma}^{l -1,d,m,D}(C)) \mid_D,
(\mathcal P_{k +j,\sigma}^{l -1,d,m,D}) \mid_D \cap W_q^{\mathpzc m}(D) \bigr) \le \\
p_{n_{j,\sigma}} \bigl(\mathcal E_{k +j,\sigma}^{l -1,d,m,D}(C),
\mathcal P_{k +j,\sigma}^{l -1,d,m,D} \cap W_q^{\mathpzc m}(\R^d) \bigr) \le \\
c_3 2^{(k +j) \mathpzc m} p_{n_{j,\sigma}} \bigl(\mathcal E_{k +j,\sigma}^{l -1,d,m,D}(C),
\mathcal P_{k +j,\sigma}^{l -1,d,m,D} \cap L_q(\R^d) \bigr) \\
\le c_4 2^{-(k +j)(\alpha -\mathpzc m -d /p +d /q)}
p_{n_{j,\sigma}} \bigl(B(l_p^{R_{k +j,\sigma}^{l -1,d,m,D}}),
l_q^{R_{k +j,\sigma}^{l -1,d,m,D}}\bigr).
\end{split}
\end{equation*}

Наконец, учитывая лемму 5.1.5 и соотношения (5.1.12), (5.2.1) и
(1.4.41), имеем
\begin{multline*} \tag{5.2.8}
\sum_{j =j_0 +1}^\infty d^0\bigl(\mathcal V_{k +j}^{l -1,d,m,D}(C),
(\mathcal P_{k +j}^{l -1,d,m,D}) \mid_D \cap W_q^{\mathpzc m}(D) \bigr) = \\
\sum_{j =j_0 +1}^\infty d^0\bigl((\mathcal E_{k +j}^{l -1,d,m,D}(C)) \mid_D,
(\mathcal P_{k +j}^{l -1,d,m,D}) \mid_D \cap W_q^{\mathpzc m}(D) \bigr) \le \\
\sum_{j =J_0 +1}^\infty
d^0 \bigl(\mathcal E_{k +j}^{l -1,d,m,D}(C),
\mathcal P_{k +j}^{l -1,d,m,D} \cap W_q^{\mathpzc m}(\R^d) \bigr) \le \\
c_5 2^{-(k +j_0)(\alpha -\mathpzc m -(d /p -d /q)_+)}.
\end{multline*}

Соединяя (5.2.5) -- (5.2.8), приходим к (5.2.3).

Вывод (5.2.4) проводится аналогично. $ \square $

Используя факты, применявшиеся при выводе результатов в [6], [7], несложно
показать, что имеет место такая лемма.

Лемма 5.2.4

Пусть $ d \in \N, \ \alpha \in \R_+, \ \mathpzc m \in \Z_+, \ 1 \le p,q \le \infty $
удовлетворяют условию (3.1.2) и $ \theta \in \R: \ 1 \le \theta < \infty, \
\l = \l(\alpha). $ Тогда существует
константа $ c_6(d,\alpha,p,\theta,q,\mathpzc m) >0 $ такая, что для $ C =
(B((B_{p,\theta}^\alpha)^{\l}(\R^d)) \cap C_0^\infty(I^d)) \mid_{I^d}, \
X = W_q^{\mathpzc m}(I^d) $ и $ p_n(\cdot,\cdot) = d_n(\cdot,\cdot), d^n(\cdot,\cdot),
\lambda_n(\cdot,\cdot), a^n(\cdot,\cdot), \epsilon_n(\cdot,\cdot) $
для любых натуральных чисел $ n, N: \ n < N, $ справедливо неравенство
\begin{equation*} \tag{5.2.9}
p_n(C,X) \ge c_6 N^{-((\alpha -\mathpzc m) /d -1 /p +1 /q)}
p_n(B(l_p^N), l_q^N).
\end{equation*}

С использованием леммы 5.2.4 доказывается

Лемма 5.2.5

Пусть $ d \in \N, \ \alpha \in \R_+, \ D $ -- ограниченная область $ \e $-типа
в $ \R^d, \ \mathpzc m \in \Z_+, \ 1 < p < \infty, \ 1 \le q \le \infty $
удовлетворяют условию (3.1.2) и $ \theta \in \R: 1 \le \theta < \infty. $
Тогда существует константа $ c_7(d,\alpha,D,p,\theta,q,\mathpzc m) >0 $ такая, что
для $ C = B((B_{p,\theta}^\alpha)^\prime(D)), \ X = W_q^{\mathpzc m}(D) $ и
$ p_n(\cdot,\cdot) = d_n(\cdot,\cdot), d^n(\cdot,\cdot), \lambda_n(\cdot,\cdot),
a^n(\cdot,\cdot), \epsilon_n(\cdot,\cdot) $
для любых натуральных чисел $ n, N : \ n < N, $ выполняется неравенство (5.2.9)
с константой $ c_7 $ вместо $ c_6. $

Доказательство.

Фиксируем точку $ x^0 \in \R^d $ и число $ \delta \in \R_+ $ такие, что
$ Q = (x^0 +\delta I^d) \subset D. $
Сначала заметим, что в условиях леммы при $ n, N \in \N: \ n < N, \ \l = \l(\alpha), $
благодаря (5.2.9) и (1.1.6), имеет место неравенство
\begin{multline*} \tag{5.2.10}
c_6 N^{-((\alpha -\mathpzc m) /d -1 /p +1 /q)}
p_n(B(l_p^N), l_q^N) \le \\
p_n((B((B_{p,\theta}^\alpha)^{\l}(\R^d)) \cap C_0^\infty(I^d)) \mid_{I^d}, W_q^{\mathpzc m}(I^d)) \le \\
c_8 p_n((B((B_{p,\theta}^\alpha)^\prime(\R^d)) \cap C_0^\infty(I^d)) \mid_{I^d}, W_q^{\mathpzc m}(I^d)).
\end{multline*}

Далее, принимая во внимание, что в силу (2.2.6) справедливо включение
\begin{multline*}
(B((B_{p,\theta}^\alpha)^\prime(\R^d)) \cap C_0^\infty(I^d)) \mid_{I^d} =
\{ f \mid_{I^d}: f \in
(B((B_{p,\theta}^\alpha)^\prime(\R^d)) \cap C_0^\infty(I^d)) \} = \\
h_{\delta,x^0} (\{ h_{\delta,x^0}^{-1}(f \mid_{I^d}): f \in
(B((B_{p,\theta}^\alpha)^\prime(\R^d)) \cap C_0^\infty(I^d)) \}) = \\
h_{\delta,x^0} (\{ (h_{\delta,x^0}^{-1} f) \mid_Q:
f \in (B((B_{p,\theta}^\alpha)^\prime(\R^d)) \cap C_0^\infty(I^d)) \}) = \\
h_{\delta,x^0} (\{ F \mid_Q:
F \in h_{\delta,x^0}^{-1}(B((B_{p,\theta}^\alpha)^\prime(\R^d)) \cap C_0^\infty(I^d)) \}) = \\
h_{\delta,x^0} (\{ F \mid_Q:
F \in h_{\delta,x^0}^{-1}(B((B_{p,\theta}^\alpha)^\prime(\R^d))) \cap C_0^\infty(Q) \}) \subset \\
h_{\delta,x^0} (\{ F \mid_Q:
F \in (c_9 B((B_{p,\theta}^\alpha)^\prime(\R^d))) \cap C_0^\infty(Q) \}) = \\
c_9 h_{\delta,x^0} (\{ F \mid_Q:
F \in B((B_{p,\theta}^\alpha)^\prime(\R^d)) \cap C_0^\infty(Q) \}),
\end{multline*}
благодаря (5.2.1), (2.2.2), имеем
\begin{multline*} \tag{5.2.11}
p_n((B((B_{p,\theta}^\alpha)^\prime(\R^d)) \cap C_0^\infty(I^d)) \mid_{I^d}, W_q^{\mathpzc m}(I^d)) \le \\
p_n(c_9 h_{\delta,x^0} (\{ F \mid_Q:
F \in B((B_{p,\theta}^\alpha)^\prime(\R^d)) \cap C_0^\infty(Q) \}), W_q^{\mathpzc m}(I^d)) \le \\
c_{10} p_n(\{ F \mid_Q:
F \in B((B_{p,\theta}^\alpha)^\prime(\R^d)) \cap C_0^\infty(Q) \}, W_q^{\mathpzc m}(Q)).
\end{multline*}

Теперь, определяя в $ W_q^{\mathpzc m}(D) $ замкнутое подпространство
$$
W_{0 \ q}^{\mathpzc m}(Q, D) = \{f \in W_q^{\mathpzc m}(D): \ f = \chi_Q f\},
$$
заметим, что
$$
\{ F \mid_D:
F \in B((B_{p,\theta}^\alpha)^\prime(\R^d)) \cap C_0^\infty(Q) \} \subset
W_{0 \ q}^{\mathpzc m}(Q, D),
$$
а для оператора $ U: \ W_q^{\mathpzc m}(D) \ni f \mapsto U f = f
\mid_Q \in W_q^{\mathpzc m}(Q), $ его сужение $ U \mid_{W_{0 \
q}^{\mathpzc m}(Q, D)} $ является изометрическим изоморфизмом
пространства $ W_{0 \ q}^{\mathpzc m}(Q, D) \cap W_q^{\mathpzc
m}(D) $ на $ W_{0 \ q}^{\mathpzc m}(Q, D) \mid_Q \cap
W_q^{\mathpzc m}(Q). $ С учётом сказанного применяя (5.2.1),
получаем
\begin{multline*} \tag{5.2.12}
p_n(\{ F \mid_Q:
F \in B((B_{p,\theta}^\alpha)^\prime(\R^d)) \cap C_0^\infty(Q) \}, W_q^{\mathpzc m}(Q)) = \\
p_n(\{ (F \mid_D) \mid_Q:
F \in B((B_{p,\theta}^\alpha)^\prime(\R^d)) \cap C_0^\infty(Q) \}, W_q^{\mathpzc m}(Q)) = \\
p_n(U(\{ F \mid_D:
F \in B((B_{p,\theta}^\alpha)^\prime(\R^d)) \cap C_0^\infty(Q) \}), W_q^{\mathpzc m}(Q)) \le \\
p_n(\{ F \mid_D:
F \in B((B_{p,\theta}^\alpha)^\prime(\R^d)) \cap C_0^\infty(Q) \}, W_q^{\mathpzc m}(D)) \le \\
p_n(\{ f: \ f \in B((B_{p,\theta}^\alpha)^\prime(D)), W_q^{\mathpzc m}(D)).
\end{multline*}

Соединяя (5.2.12), (5.2.11), (5.2.10), видим, что для $ C $ и $ X $ из
формулировки леммы соблюдается (5.2.9) с константой $ c_7 $ вместо $ c_6. \square $

Опираясь на соотношения (5.1.8), (5.1.9), (5.1.15), (5.2.3), (5.2.4) и
(5.2.9) с константой $ c_7 $ вместо $ c_6, $ подобно тому, как это сделано
в [7], устанавливается слабая асимптотика колмогоровского,
гельфандовского, линейного, александровского и энтропийного
поперечников классов $ B((H_p^\alpha)^\prime(D)) $ и
$ B((B_{p,\theta}^\alpha)^\prime(D)) $ в пространстве $ W_q^{\mathpzc m}(D), $
заданных в ограниченной области $ D \  \e $-типа.

Теорема 5.2.6

Пусть $ d \in \N, \ \alpha \in \R_+, \ \mathpzc m \in \Z_+, \ 1 < p < \infty, \
1 \le q \le \infty $ удовлетворяют условию (3.1.2) и $ \theta \in \R:
\ 1 \le \theta < \infty, $ а $ D \subset \R^d $ -- ограниченная область
$ \e $-типа. Тогда для $ C = B((H_p^\alpha)^\prime(D)), \
B((B_{p,\theta}^\alpha)^\prime(D)) $ и $ X = W_q^{\mathpzc m}(D) $
справедливо соотношение
\begin{equation*}
d_n(C,X) \asymp
\begin{cases} n^{-(\alpha -\mathpzc m) /d
+(p^{-1} -q^{-1})_+}, \text{ при } q \le p \text{ или } p < q \le 2; \\
n^{-(\alpha -\mathpzc m) /d +(p^{-1} -2^{-1})_+},\parbox[t]{.4\textwidth} { при $ q > \max(2,p),\\
\alpha -\mathpzc m -d /p -(d /2 -d /p)_+ >0.$ }
\end{cases}
\end{equation*}

Теорема 5.2.7

Пусть $ d, \alpha, \mathpzc m,p,q, \theta, D, $ а также $ C $ и $ X $ имеют тот же
смысл, что и в теореме 5.2.6. Тогда имеет место
соотношение
\begin{equation*}
d^n(C,X) \asymp
\begin{cases} n^{-(\alpha -\mathpzc m) /d +(p^{-1} -q^{-1})_+},
\parbox[t]{.4\textwidth} { при $(q \le p$ или $ 2 \le p < q)$
и соблюдении условия (3.1.2);} \\
n^{-(\alpha -\mathpzc m) /d +(2^{-1} -q^{-1})_+},\
\parbox[t]{.4\textwidth} { при $p < \min(2,q),
\alpha -\mathpzc m -d /q^{\prime} -(d /2 -d /q^{\prime})_+ >0. $}
\end{cases}
\end{equation*}

Теорема 5.2.8

В условиях теоремы 5.2.6 верно соотношение
\begin{equation*}
\lambda_n(C,X) \asymp
\begin{cases} n^{-(\alpha -\mathpzc m) /d +(p^{-1} -q^{-1})_+},
\parbox[t]{.4\textwidth} { при $(q  \le p$
или $p < q \le 2$ или
$2 \le p < q)$;} \\
n^{-(\alpha -\mathpzc m) /d +1/2 -1/q +(1/p +1/q -1)_+}, \parbox[t]{.4\textwidth} { при $p < 2 < q,
\alpha -\mathpzc m -d /p -(d /p^{\prime} - d/q)_+ >0.$}
\end{cases}
\end{equation*}

Теорема 5.2.9

В условиях теоремы 5.2.6 соблюдается соотношение
\begin{equation*}
a_n(C,X) \asymp n^{-(\alpha -\mathpzc m) /d}.
\end{equation*}

Теорема 5.2.10

В условиях теоремы 5.2.6 при $ \alpha -\mathpzc m -d >0 $ выполняется соотношение
\begin{equation*}
\epsilon_n(C,X) \asymp n^{-(\alpha -\mathpzc m) /d}.
\end{equation*}

\end{document}